\documentclass{amsart}

\usepackage{color}
\usepackage{setspace}

\usepackage{amssymb}   

\usepackage{amsmath}

\usepackage{mathrsfs}

\usepackage{stmaryrd}

\usepackage{amsthm}
\usepackage{newlfont}
\usepackage{amscd}
\usepackage{mathtools}
\usepackage{bm}
\usepackage{hyperref}
\usepackage{enumitem}

\usepackage[all]{xy}

\usepackage{textcomp}

\usepackage{amsbsy}

\newtheorem{thm}{Theorem}[section]

\newtheorem{thm-defn}[thm]{Theorem/Definition}
\newtheorem{lem}[thm]{Lemma}
\newtheorem{prop}[thm]{Proposition}
\newtheorem{cor}[thm]{Corollary}

\theoremstyle{definition}
\newtheorem{defn}[thm]{Definition}

\newtheorem{eg}[thm]{Example}

\newtheorem{construction}[thm]{Construction}
\newtheorem{set-up}[thm]{Set-up}
\newtheorem{notation}[thm]{Notation}
\newtheorem{assumption}[thm]{Assumption}

\theoremstyle{remark}
\newtheorem{rem}[thm]{Remark}

\numberwithin{equation}{section}

\usepackage{relsize}
\usepackage[bbgreekl]{mathbbol}
\usepackage{amsfonts}

\DeclareSymbolFontAlphabet{\mathbb}{AMSb}
\DeclareSymbolFontAlphabet{\mathbbl}{bbold}
\newcommand{\Prism}{{\mathlarger{\mathbbl{\Delta}}}}

\DeclareMathOperator{\Gal}{Gal}

\DeclareMathOperator{\Spec}{Spec}
\DeclareMathOperator{\Spa}{Spa}
\DeclareMathOperator{\Spf}{Spf}

\DeclareMathOperator{\Frac}{Frac}
\DeclareMathOperator{\Ker}{Ker}
\DeclareMathOperator{\rank}{rank}
\DeclareMathOperator{\Fil}{Fil}

\newcommand{\ur}{\mathrm{ur}}
\newcommand{\cris}{\mathrm{cris}}

\newcommand{\et}{\mathrm{\acute{e}t}}
\newcommand{\proet}{\mathrm{pro\acute{e}t}}

\newcommand{\N}{\mathbf{N}}
\newcommand{\Z}{\mathbf{Z}}
\newcommand{\Q}{\mathbf{Q}}
\newcommand{\C}{\mathbf{C}}
\newcommand{\A}{\mathbf{A}}
\newcommand{\B}{\mathbf{B}}
\newcommand{\OA}{\mathbf{OA}}
\newcommand{\OB}{\mathbf{OB}}
\newcommand{\Ainf}{\mathbf{A}_{\mathrm{inf}}}
\newcommand{\calOA}{\mathcal{O}\mathbb{A}}
\newcommand{\calOB}{\mathcal{O}\mathbb{B}}

\newcommand{\fkm}{\mathfrak{m}}

\newcommand{\fkL}{{\mathfrak L}}
\newcommand{\fkM}{{\mathfrak M}}
\newcommand{\fkN}{{\mathfrak N}}

\newcommand{\fkS}{{\mathfrak S}}

\newcommand{\fkU}{{\mathfrak U}}

\newcommand{\fkX}{{\mathfrak X}}
\newcommand{\fkY}{{\mathfrak Y}}

\newcommand{\calE}{\mathcal{E}}
\newcommand{\calF}{\mathcal{F}}
\newcommand{\calG}{\mathcal{G}}

\newcommand{\calI}{\mathcal{I}}

\newcommand{\calM}{\mathcal{M}}

\newcommand{\calO}{\mathcal{O}}

\newcommand{\calS}{\mathcal{S}}

\newcommand{\calV}{\mathcal{V}}

\newcommand{\bA}{\mathbb{A}}
\newcommand{\bB}{\mathbb{B}}

\newcommand{\bL}{\mathbb{L}}

\begin{document}

\pagenumbering{arabic}

\title{Completed prismatic $F$-crystals and crystalline $\Z_p$-local systems}
\author{Heng Du, Tong Liu, Yong Suk Moon, Koji Shimizu} 
\date{}

\begin{abstract} 
We introduce the notion of completed $F$-crystals on the absolute prismatic site of a smooth $p$-adic formal scheme. 
We define a functor from the category of completed prismatic $F$-crystals to that of crystalline \'etale $\Z_p$-local systems on the generic fiber of the formal scheme and show that it gives an equivalence of categories.
This generalizes the work of Bhatt and Scholze, which treats the case of a mixed characteristic complete discrete valuation ring with perfect residue field.
\end{abstract}

\maketitle

\tableofcontents

\section{Introduction}

Let $p$ be a prime.
In \cite{bhatt-scholze-prismaticcohom}, Bhatt and Scholze introduce the notion of prisms and the relative prismatic ringed site $((\fkX/(A,I))_\Prism,\calO_\Prism)$ for a bounded prism $(A,I)$ and a smooth $p$-adic formal scheme $\fkX$ over $A/I$. Surprisingly, the cohomology $R\Gamma((\fkX/(A,I))_\Prism,\calO_\Prism)$ gives a good integral $p$-adic cohomology of $\fkX$: it recovers the crystalline cohomology of the special fiber as well as the \'etale cohomology of the generic fiber. The prismatic formalism also gives a site-theoretic construction of the $A_{\mathrm{inf}}$-cohomology and the Breuil--Kisin cohomology when $(A,I)=(A_{\mathrm{inf}}(\calO_{\C_p}),\ker\theta)$ and $(\fkS,(E))$, respectively (see \cite[Ex.~1.9]{bhatt-scholze-prismaticcohom} for the details). 

Another advantage of this site-theoretic approach is that it provides a natural framework of the coefficient theory.
In the case of the relative prismatic site, Tian \cite{Tian} studies the cohomology of prismatic crystals when $\fkX$ is proper over $A/I$.

One can study crystals on the absolute prismatic site as well.
Let $\calO_K$ be a complete discrete valuation ring of mixed characteristic $(0,p)$ with perfect residue field $k$, and let $\fkX$ be a smooth $p$-adic formal scheme over $\calO_K$.
In the subsequent paper \cite{bhatt-scholze-prismaticFcrystal}, Bhatt and Scholze study sheaves on the absolute prismatic site $\fkX_\Prism$. Recall that the site $\fkX_\Prism$ has a sheaf $\calO_\Prism$ of rings equipped with a Frobenius $\varphi$ and an ideal sheaf $\calI_\Prism$ (see Definition~\ref{defn:abs-prism-site} for the details). They introduce the category $\mathrm{Vect}^\varphi(\fkX_\Prism)$ of prismatic $F$-crystals of vector bundles on $(\fkX_\Prism,\calO_\Prism)$ as well as the category $\mathrm{Vect}(\fkX_\Prism,\calO_\Prism[1/\calI_\Prism]^\wedge_p)^{\varphi=1}$ of so-called Laurent $F$-crystals on $\fkX$.

The main theorem \cite[Thm.~1.2]{bhatt-scholze-prismaticFcrystal} of Bhatt and Scholze states that $\mathrm{Vect}^\varphi((\calO_K)_\Prism)$ is equivalent to the category of lattices in crystalline representations of $K$. They also show that, for general $\fkX$,  $\mathrm{Vect}(\fkX_\Prism,\calO_\Prism[1/\calI_\Prism]^\wedge_p)^{\varphi=1}$ is equivalent to the category of $\Z_p$-local systems on the generic fiber of $\fkX$.
Part of their work is reproved or generalized by Du--Liu \cite{du-liu-prismaticphiGhatmodule}, Wu \cite{wu-Gal-rep-prism-F-cryst}, and Min--Wang \cite{min-wang-rel-phi-gamma-prism-F-crys}. For other works on the prismatic site, we refer the reader to the recent survey \cite{bhatt-ICM}.

The present article studies the relationship between lattices in crystalline representations and suitable $F$-crystals on the absolute prismatic site in the relative situation. For this, we need to enlarge the category $\mathrm{Vect}^\varphi(\fkX_\Prism)$ of prismatic $F$-crystals on $\fkX$. To explain the enlarged category, we first focus on the small affine case.
More precisely, let $R_0$ be the $p$-adic completion of an integral domain that is \'etale over $W(k)[T_1^{\pm 1},\ldots, T_d^{\pm 1}]$ for some $d\geq 0$ and set $R\coloneqq R_0\otimes_{W(k)}\calO_K$. We also fix a uniformizer $\pi$ of $\calO_K$ with monic minimal polynomial $E(u)$ over $W(k)$.
We consider the following type of sheaves on the absolute prismatic site $R_\Prism$.

\begin{defn}[(Definition~\ref{defn:category-good-prism-completed-F-crystal}, Remark~\ref{rem:definition of completed prismatic F-crystals})]\label{defn:category-good-prism-completed-F-crystal-intro}
A \emph{completed prismatic $F$-crystal} on $R$ is a sheaf $\calF$ of $\calO_\Prism$-modules on $R_\Prism$ together with a $\varphi$-semilinear endomorphism $\varphi_{\calF}\colon \calF\rightarrow \calF$ that satisfies the following properties:
\begin{enumerate}
    \item for each $(A,I)\in R_\Prism$, the $A$-module $\calF_A\coloneqq \calF(A,I)$ is finitely generated and classically $(p,I)$-complete;
    \item for any morphism $(A,IA)\rightarrow (B,IB)$ of bounded prisms over $R$, the map $B\widehat{\otimes}_A\calF_A\rightarrow \calF_B$ is an isomorphism;
    \item for the Breuil--Kisin prism $(\fkS=R_0[\![u]\!],(E(u)))\in R_\Prism$ (see Example~\ref{eg:prism-S2S3}), $\calF_{\fkS}$ is torsion free, $\calF_{\fkS}[E^{-1}]^\wedge_p$ is finite projective over $\fkS[E^{-1}]^\wedge_p$, and $\calF_{\fkS}=\calF_{\fkS}[p^{-1}]\cap \calF_{\fkS}[E^{-1}]^\wedge_p$;
    \item the cokernel of $1\otimes \varphi_{\calF_\fkS}\colon \varphi^\ast\calF_{\fkS}\rightarrow \calF_{\fkS}$ is killed by $E^r$ for a non-negative integer $r$.
\end{enumerate}
We write $\mathrm{CR}^{\wedge,\varphi}(R_\Prism)$ for the category of completed prismatic $F$-crystals on $R$.
\end{defn}

Condition (iii) in the definition is technical but plays a crucial role in our theory. 
The category $\mathrm{CR}^{\wedge,\varphi}(R_\Prism)$ contains, as a full subcategory, the category $\mathrm{Vect}_{\mathrm{eff}}^\varphi(R_\Prism)$ of effective prismatic $F$-crystals of vector bundles on $R_\Prism$ in the sense of \cite[Def.~4.1]{bhatt-scholze-prismaticFcrystal}.
We note that there exists a completed prismatic $F$-crystal which is not a prismatic $F$-crystal of vector bundles (see below and Example~\ref{eg:non-p-divisible crystalline rep}), and thus $\mathrm{CR}^{\wedge,\varphi}(R_\Prism)$ is strictly larger than $\mathrm{Vect}_{\mathrm{eff}}^\varphi(R_\Prism)$ in general.

The main goal of this paper is to describe lattices in crystalline representations of the Galois group $\calG_R$ of $R[p^{-1}]$ in terms of completed prismatic $F$-crystals.
Here is our main result:
\begin{thm}[({Theorem~\ref{thm:main}})]\label{thm:main in intro}
There is a contravariant equivalence of categories 
\[
T\colon \mathrm{CR}^{\wedge,\varphi}(R_\Prism)\xrightarrow{\cong} \mathrm{Rep}_{\Z_p,\geq 0}^{\cris}(\calG_R)
\]
from the category of completed prismatic $F$-crystals on $R$ to the category of crystalline $\Z_p$-representations of $R[p^{-1}]$ with non-negative Hodge--Tate weights, which is functorial in $R$.
\end{thm}

Following \cite{bhatt-scholze-prismaticFcrystal}, we call $T$ the \emph{\'etale realization functor}. Since $T$ is compatible with Breuil--Kisin and Tate twists, one can further enlarge $\mathrm{CR}^{\wedge,\varphi}(R_\Prism)$ and obtain an equivalence of categories with the category of \emph{all} crystalline $\Z_p$-representations of $R[p^{-1}]$ (see Example~\ref{eg:Breuil-Kisin twists} and Remark~\ref{rem:equivalence with all crystalline reps}).

Once Theorem~\ref{thm:main in intro} is obtained, we can globalize our equivalence to the one for a smooth $p$-adic formal scheme over $\calO_K$.

\begin{thm}[(Theorem~\ref{thm:main-global})] \label{thm:main-intro-global}
Let $\fkX$ be a smooth $p$-adic formal scheme over $\calO_K$. Then there is a natural equivalence of categories
\[
T\colon \mathrm{CR}^{\wedge,\varphi}(\fkX_\Prism)\xrightarrow{\cong}\mathrm{Loc}_{\Z_p,\geq 0}^\cris(\fkX_\eta)
\]
between the category of completed prismatic $F$-crystals on $\fkX$ and the category of crystalline $\Z_p$-local systems with non-negative Hodge--Tate weights on the adic generic fiber $\fkX_\eta$ of $\fkX$ (see \S~\ref{sec:globalization} for the precise definitions).
\end{thm}

The main theorems give a prismatic description of crystalline $\Z_p$-representations in the relative case.
Note that when $R=\calO_K$, Kisin \cite{kisin-crystalline} gave a description of lattices in crystalline representations of $K$ in terms of Breuil--Kisin modules. His work was generalized by Brinon--Trihan \cite{brinon-trihan} to the case of complete discrete valuation rings with \emph{imperfect} residue field with a finite $p$-basis. Furthermore, Kim \cite{kim-groupscheme-relative} introduced the notion of Kisin $\fkS$-modules over $R$ as a generalization of Breuil--Kisin modules of $E$-height $\leq 1$ in the relative case. Kim attached to a $p$-divisible group over a general $R$ a Kisin $\fkS$-module and showed that the category of $p$-divisible groups over $R$ is equivalent to the category of Kisin $\fkS$-modules when $p \geq 3$ \cite[Cor.~3]{kim-groupscheme-relative}. 
However, it has not yet been known how to describe crystalline $\Z_p$-representations of $R[p^{-1}]$ with non-negative Hodge--Tate weights in terms of suitable Breuil--Kisin type modules in the relative case. In fact, even a suitable description of rational crystalline representations of $R[p^{-1}]$ has not been given yet in general: while crystalline $\mathbf{Q}_p$-representations of $K$ can be classified by weakly admissible filtered $\varphi$-modules \cite{colmez-fontaine}, the correct weakly admissibility has not been found in the relative case. We hope that the notion of completed prismatic $F$-crystals clarifies these complications.

Examples of crystalline $\Z_p$-representations of $R[p^{-1}]$ with Hodge--Tate weights in $[0,1]$ arise from $p$-divisible groups over $R$.
Ansch\"utz and Le Bras \cite{anschutz-lebras-prismaticdieudonne} developed the prismatic Dieudonn\'e theory. It follows from their work that the category of $p$-divisible groups over $R$ is equivalent to the category of effective prismatic $F$-crystals of vector bundles of $\calI_\Prism$-height $\leq 1$ (see \S~\ref{subsec:height one} for the details). It is easy to see that their formulation is compatible with ours. Following \cite{vasiu-zink-purity}, we will also provide an example of a completed prismatic $F$-crystal over $R$ that does not arise from a $p$-divisible group over $R$ (Example~\ref{eg:non-p-divisible crystalline rep}).
This implies that our category $\mathrm{CR}^{\wedge,\varphi}(R_\Prism)$ is strictly larger than the subcategory of effective prismatic $F$-crystals of vector bundles on $R$ and that the former category  is necessary to describe crystalline $\Z_p$-representations in the relative case. It is an interesting question whether a completed prismatic $F$-crystal on $R$ becomes a prismatic $F$-crystal of vector bundles on $\fkX$ by the pullback along an admissible blow-up $\fkX\rightarrow \Spf R$. A related question is whether a crystalline $\Z_p$-representation of $R[p^{-1}]$ with Hodge--Tate weights in $[0,1]$ comes from a $p$-divisible group on $\fkX$ for some admissible blow-up $\fkX\rightarrow \Spf R$.
We note that admissible blow-ups $\fkX\rightarrow \Spf R$ usually yield \emph{non-smooth} $p$-adic formal schemes $\fkX$ and thus these questions diverge from our current work.

Let us now explain the construction of the \'etale realization functor $T$ in Theorem~\ref{thm:main in intro} (see Proposition~\ref{prop:etale-realization}).
This will be explained best in the following commutative diagram:
\[
\xymatrix{
\mathrm{Vect}_{\mathrm{eff}}^\varphi(R_\Prism)\ar[r]\ar@{^{(}->}[d]
&\mathrm{Vect}(R_\Prism,\calO_\Prism[1/\calI_\Prism]^\wedge_p)^{\varphi=1}\ar[r]^-{\cong}
&\mathrm{Rep}_{\Z_p}^{\mathrm{pr}}(\calG_R)\\
\mathrm{CR}^{\wedge,\varphi}(R_\Prism).\ar@{.>}[ur]^-{\calF\mapsto\calF_{\et}}\ar@{.>}[urr]_{T^{\vee}}
&&
}
\]
Here $\mathrm{Vect}(R_\Prism,\calO_\Prism[1/\calI_\Prism]^\wedge_p)^{\varphi=1}$ denotes the category of Laurent $F$-crystals, namely,  prismatic $F$-crystals of vector bundles of $\calO_\Prism[1/\calI_\Prism]^\wedge_p$-modules on $R$, and $\mathrm{Rep}_{\Z_p}^{\mathrm{pr}}(\calG_R)$ denotes the category of finite free $\Z_p$-representations of the Galois group $\calG_R$ of $R[p^{-1}]$.
The functor $\mathrm{Vect}_{\mathrm{eff}}^\varphi(R_\Prism)\rightarrow\mathrm{Vect}(R_\Prism,\calO_\Prism[1/\calI_\Prism]^\wedge_p)^{\varphi=1}$ is the scalar extension functor. Bhatt--Scholze \cite{bhatt-scholze-prismaticFcrystal} and Min--Wang \cite{min-wang-rel-phi-gamma-prism-F-crys} show the (covariant) equivalence of categories $\mathrm{Vect}(R_\Prism,\calO_\Prism[1/\calI_\Prism]^\wedge_p)^{\varphi=1}\cong \mathrm{Rep}_{\Z_p}^{\mathrm{pr}}(\calG_R)$. Hence, to define the contravariant functor $T\colon \mathrm{CR}^{\wedge,\varphi}(R_\Prism)\rightarrow \mathrm{Rep}_{\Z_p}^{\mathrm{pr}}(\calG_R)$ or its dual $T^\vee$, it suffices to show that the functor 
$\mathrm{Vect}_{\mathrm{eff}}^\varphi(R_\Prism)\rightarrow\mathrm{Vect}(R_\Prism,\calO_\Prism[1/\calI_\Prism]^\wedge_p)^{\varphi=1}$ extends to a functor 
\[
\mathrm{CR}^{\wedge,\varphi}(R_\Prism)\rightarrow\mathrm{Vect}(R_\Prism,\calO_\Prism[1/\calI_\Prism]^\wedge_p)^{\varphi=1},\quad\calF\mapsto \calF_{\et}.
\]
The construction of the latter functor uses the following fact on the Breuil--Kisin prism $(\fkS,(E(u)))$: it covers the final object of $\mathrm{Shv}(R_\Prism)$, and thus a sheaf on $R_\Prism$ is described by a descent datum involving the self-product $(\fkS^{(1)},(E(u)))$ and the self-triple-product $(\fkS^{(2)},(E(u)))$ of the Breuil--Kisin prism. 
In particular, completed prismatic $F$-crystals are described by the following data:

\begin{prop}[(Proposition~\ref{prop:equivalence-to-descent-datum})]\label{prop:equivalence-to-descent-datum in intro}
The association $\calF\mapsto \calF_\fkS$ gives rise to an equivalence of categories $\mathrm{CR}^{\wedge,\varphi}(R_{\Prism})\xrightarrow{\cong}\mathrm{DD}_\fkS$. Here $\mathrm{DD}_\fkS$ consists of triples $(\mathfrak{M}, \varphi_\fkM, f)$ where
\begin{enumerate}
\item $\fkM$ is a finite $\fkS$-module satisfying condition (iii) of Definition~\ref{defn:category-good-prism-completed-F-crystal-intro} in place of $\calF_\fkS$;
\item $\varphi_\fkM\colon \mathfrak{M} \rightarrow \mathfrak{M}$ is a $\varphi$-semi-linear endomorphism such that the cokernel of $1\otimes\varphi_\fkM\colon \varphi^*\mathfrak{M} \rightarrow \mathfrak{M}$ is killed by $E^r$ for a non-negative integer $r$;
\item $f\colon \mathfrak{S}^{(1)}\otimes_{p_1,\mathfrak{S}}\mathfrak{M} \xrightarrow{\cong}\mathfrak{S}^{(1)}\otimes_{p_2,\mathfrak{S}}\mathfrak{M}$ is an isomorphism of $\fkS^{(1)}$-modules that is compatible with Frobenii and satisfies the cocycle condition over $\mathfrak{S}^{(2)}$.
\end{enumerate}
\end{prop}
Since $\mathrm{Vect}(R_\Prism,\calO_\Prism[1/\calI_\Prism]^\wedge_p)^{\varphi=1}$ has a similar description in terms of descent data involving $\fkS^{(1)}[E^{-1}]^\wedge_p$ and $\fkS^{(2)}[E^{-1}]^\wedge_p$, the base change along the map $\fkS^{(1)}\rightarrow \fkS^{(1)}[E^{-1}]^\wedge_p$ yields the desired functor $\mathrm{CR}^{\wedge,\varphi}(R_\Prism)\rightarrow\mathrm{Vect}(R_\Prism,\calO_\Prism[1/\calI_\Prism]^\wedge_p)^{\varphi=1}$.
To put things together, the contravariant functor $T$ is explicitly given by
\[
T(\mathcal{F})\coloneqq 
\bigl((\mathcal{F}_{\et}(\mathbf{A}_{\mathrm{inf}}(\overline{R}), (\xi)))^{\varphi_{\calF_{\et}} = 1}\bigr)^\vee.
\]
See Example~\ref{eg:prism-Ainf} for the definition of the $\A_\mathrm{inf}$-prism $(\mathbf{A}_{\mathrm{inf}}(\overline{R}), (\xi))$.
Once $T$ is defined, it is not difficult to see that $T$ is fully faithful and that $T(\calF)[p^{-1}]$ is a crystalline $\Q_p$-representation of $\calG_R$ with non-negative Hodge--Tate weights.

The hardest part of the proof of Theorem~\ref{thm:main in intro} concerns the essential surjectivity of $T$, and \S~\ref{sec:quasi-kisin-mods-cryst-loc-systs} is devoted to proving it.
For this, take $T_0\in \mathrm{Rep}_{\Z_p,\geq 0}^{\cris}(\calG_R)$. We will attach to $T_0$ an object $(\fkM,\varphi_\fkM,f)\in \mathrm{DD}_\fkS$.

To explain the outline of the construction, let us introduce several rings.
Let $\calO_\calE$ denote the $p$-adic completion $\fkS[E^{-1}]^\wedge_p$ of $\fkS[E^{-1}]$. We write $\calO_{L_0}$ (resp.~$\calO_L$)  for the $p$-adic completion of the localization of $R_0$ at the prime $(p)$ (resp.~the localization of $R$ at the prime $(\pi)$): $\calO_{L_0}$ is an absolutely unramified complete discrete valuation ring with imperfect residue field having a finite $p$-basis, and $\calO_L=\calO_{L_0}\otimes_{W(k)}\calO_K$. We set $\fkS_L\coloneqq \calO_{L_0}[\![u]\!]$ and $\calO_{\calE_L}\coloneqq \fkS_L[E^{-1}]^{\wedge}_{p}$.

On the one hand, the theory of \'etale $\varphi$-modules attaches to $T_0$ a finite free $\calO_\calE$-module $\calM$ together with a  Frobenius and a descent datum. 
On the other hand, Brinon--Trihan's theory \cite{brinon-trihan} of Breuil--Kisin modules associates with $T_0|_{\Gal(\overline{L}/L)}$ a finite free $\fkS_L$-module $\fkM_L$ with a Frobenius.
We set $\fkM\coloneqq \calM\cap \fkM_L$ inside $\calO_{\calE_L}\otimes_{\calO_\calE}\calM=\calO_{\calE_L}\otimes_{\fkS_L}\fkM_L$. Naturally, $\fkM$ is equipped with a Frobenius $\varphi_\fkM$.
With careful study of the structures, we are able to show that 
the pair $(\fkM,\varphi_\fkM)$ satisfies conditions (i) and (ii) of Proposition~\ref{prop:equivalence-to-descent-datum in intro}.
Finally, the connection on $D_{\cris}(T_0[p^{-1}])$ equips $\fkM[p^{-1}]$ with a descent datum. Combined with the descent datum on $\calM$, it yields a descent datum $f$ on $\fkM$ and thus an object $(\fkM,\varphi_\fkM,f)\in \mathrm{DD}_\fkS$.
The associated completed prismatic $F$-crystal $\calF$ satisfies $T(\calF)\cong T_0$.

\begin{rem}
After we posted our paper on arXiv, we learned that Guo--Reinecke independently proved Theorem~\ref{thm:main-intro-global} \cite[Thm. A, Rem.~1.8]{guo-reinecke-prism-F-crys}. Their proof generalizes the method in \cite{bhatt-scholze-prismaticFcrystal} and is different from ours. Our category $\mathrm{CR}^{\wedge,\varphi}(\fkX_\Prism)$ corresponds to the category of \textit{effective analytic prismatic} $F$-\textit{crystals} on $\fkX$ in their terminology \cite[Def. 3.2]{guo-reinecke-prism-F-crys}. Note that one can deduce from Theorem~\ref{thm:main-intro-global} together with the compatibility of $T$ with Breuil--Kisin and Tate twists (see Example~\ref{eg:Breuil-Kisin twists} and Remark~\ref{rem:equivalence with all crystalline reps}) that the category of analytic prismatic $F$-crystals on $\fkX$ is equivalent to the category of all crystalline $\Z_p$-local systems on $\fkX_\eta$ as in \cite[Thm. A]{guo-reinecke-prism-F-crys}.
\end{rem}

\medskip
\noindent
\textbf{Organization of the paper.}
Section~\ref{sec:review of crystallline rep and etale phi modules} reviews basic concepts in relative $p$-adic Hodge theory.
In \S~\ref{sec-basering}, we explain the assumptions on our base ring $R$ and objects attached to $R$, which we will use throughout this article.  We review crystalline representations developed by Brinon \cite{brinon-relative} in \S~\ref{sec:crystalline representations}, and \'etale $\varphi$-modules in \S~\ref{sec:etale phi-module}. The topics in the latter two subsections are standard, and the reader may skip them. 

Section~\ref{sec:prism-cryst-loc-syst} introduces the notion of completed prismatic $F$-crystals and states the main theorems. In \S~\ref{sec:prism-site}, we recall the definition of the absolute prismatic site of a $p$-adic formal scheme and explain key examples of prisms in the small affine case. In \S~\ref{sec:completed-prism-F-crystals}, we define finitely generated completed prismatic crystals and completed prismatic $F$-crystals in the small affine case. Then we describe the category of completed prismatic $F$-crystals in terms of descent data in \S~\ref{sec:descent data}. Section~\ref{sec:etale-realization-main-thm} introduces the \'etale realization functor $T\colon \mathrm{CR}^{\wedge,\varphi}(R_\Prism)\rightarrow \mathrm{Rep}_{\Z_p,\geq 0}^{\cris}(\calG_R)$ and states the main theorem in the small affine case (Theorem~\ref{thm:main}). We also prove part of the main theorem that $T$ is fully faithfully and $T(\calF)$ is crystalline for $\calF\in\mathrm{CR}^{\wedge,\varphi}(R_\Prism)$ in this subsection.  In \S~\ref{subsec:height one}, we consider the height one case and compare the \'etale realization functor with prismatic Dieudonn\'e theory by Ansch\"utz--Le Bras \cite{anschutz-lebras-prismaticdieudonne}. We also present an example of a completed prismatic $F$-crystal that is not an effective prismatic $F$-crystal of vector bundles (Example~\ref{eg:non-p-divisible crystalline rep}). In \S~\ref{sec:globalization}, we define the notion of completed prismatic $F$-crystals on a general smooth $p$-adic formal scheme and the \'etale realization functor. We end the subsection with the main theorem in this general case (Theorem~\ref{thm:main-global}), which is a direct consequence of Theorem~\ref{thm:main}.

Section~\ref{sec:quasi-kisin-mods-cryst-loc-systs} is devoted to the proof of the remaining part of the main theorem (the essential surjectivity of Theorem~\ref{thm:main}). In \S~\ref{sec:quasi-kisin-mod-rational-descent-data}, we define the notion of quasi-Kisin modules and show that such an object yields a rational Kisin descent datum.  Section~\ref{sec:quasi-kisin-mod-projectivity} proves the general fact that for a finite torsion free $\varphi$-module $(\fkM,\varphi_\fkM)$ of finite $E$-height over $\fkS$, $\fkM[p^{-1}]$ is projective over $\fkS[p^{-1}]$ (Proposition~\ref{prop:rational-projectivity-etale-over-torus-case}). In \S~\ref{sec:cryst-rep-cdvr}, we consider the CDVR case. With these preparations, we attach a quasi-Kisin module to a lattice in a crystalline representation in \S~\ref{sec:quasi-kisin-mod-construction} and \ref{sec:quasi-kisin-mod-connection}, and  complete the proof of Theorem~\ref{thm:main} in \S~\ref{sec:equivalence-categories}.

Appendix~\ref{sec:crystalline local systems} follows the work of Tan and Tong \cite{Tan-Tong} and defines the notion of crystalline local systems on the generic fiber of a smooth $p$-adic formal scheme.

\medskip
\noindent
\textbf{Notation and conventions}.
Let $p$ be a prime and let $k$ be a perfect field of characteristic $p$. Write $W=W(k)$ and let $K$ be a finite totally ramified extension of $K_0\coloneqq W[p^{-1}]$. 
Fix a uniformizer $\pi$ of $K$ and let $E(u)\in W[u]$ denote the monic minimal polynomial of $\pi$.

For derived completions and relevant concepts, we refer the reader to \cite[\S 1.2]{bhatt-scholze-prismaticcohom}. In this article, most rings are classically $p$-complete, and we also call them \emph{$p$-adically complete}. Similarly, a \emph{$p$-adically completed \'etale map} from a $p$-adically complete ring $A$ refers to the (classical) $p$-adic completion of an \'etale map from $A$.

We also follow \cite{bhatt-scholze-prismaticcohom} for the definitions of $\delta$-rings and prisms. However, to avoid confusion, we say that a map of prisms $(A,I)\rightarrow (B,J)$ is \emph{$(p,I)$-completely (faithfully) flat} if the map $A\rightarrow B$ is $(p,I)$-completely faithfully flat (compare \cite[Def.~3.2]{bhatt-scholze-prismaticcohom}). 

We write $W\langle T_1^{\pm 1},\ldots,T_d^{\pm 1}\rangle$ for the $p$-adic completion of the Laurent polynomial ring $W[T_1^{\pm 1},\ldots,T_d^{\pm}]$.
 In this article, the braces $\{\cdots\}$ denote the $p$-adically completed divided power polynomials, and for a fixed prism $(A,I)$, the notation $\{\cdots\}^\wedge_\delta$ stands for adjoining elements in the category of derived $(p,I)$-complete simplicial $\delta$-$A$-algebras.

For an element $a$ of a $\Q$-algebra $A$ and $n\geq 0$, write $\gamma_n(a)$ for the element $\frac{a^n}{n!}\in A$.

Our convention is that the cyclotomic character $\Z_p(1)\coloneqq T_p(\mu_{p^\infty})$ has Hodge--Tate weight one.

\medskip
\noindent
\textbf{Acknowledgements}.	
We thank H\'el\`ene Esnault, Mark Kisin, and Peter Scholze for their valuable comments on earlier versions of this article. We are also grateful to Abhinandan for pointing out an error in the previous version of Lemma~\ref{lem:intersection-witt-rings}. Lastly, we thank the anonymous referee for carefully reading the paper and making many valuable suggestions. The third and fourth authors are partially supported by AMS--Simons Travel Grant.

\section{Review of crystalline representations and \'etale \texorpdfstring{$\varphi$}{phi}-modules}\label{sec:review of crystallline rep and etale phi modules}

 \subsection{Base ring}\label{sec-basering}
In this subsection, we will introduce our base ring $R$.

\begin{defn} \label{defn:small-over-OK}
A $p$-adically complete $\calO_K$-algebra is called \emph{small and smooth} (or \emph{small} for short) if it is $p$-adically completed \'etale over $\calO_K\langle T_1^{\pm 1},\ldots,T_d^{\pm 1}\rangle$ for some $d\geq 0$.
\end{defn}

\begin{rem} \label{rem:p-complete-etale}
Let $R$ be small over $\calO_K$. 
Since $R/\pi R$ is \'etale over $k[T_1^{\pm 1}\mspace{-5mu},\ldots,T_d^{\pm 1}]$, there exists a subalgebra $R_0\subset R$ such that $R_0$ is $p$-adically completed \'etale over $W\langle T_1^{\pm 1},\ldots,T_d^{\pm 1}\rangle$ and $R=R_0\otimes_W\calO_K$.

Let $R'$ be $p$-adically completed \'etale over $R$.
If one fixes a subring $R_0\subset R$ as above, then the \'etale map $R_0/pR_0=R/\pi R\rightarrow R'/\pi R'$ lifts uniquely to a $p$-adically completed \'etale map $R_0\rightarrow R_0'$. Moreover, $R_0'\otimes_W\calO_K$ is isomorphic to $R'$ as $R$-algebras.
\end{rem}

In this paper, we use the crystalline period rings developed in \cite{brinon-relative}.
For this, we consider the following class of $p$-adic rings that contains connected small $\calO_K$-algebras:

\begin{set-up}\label{set-up:base ring}
A connected $p$-adically complete $\calO_K$-algebra $R$ is said to be a \emph{base ring} if it is of the form $R\coloneqq R_0\otimes_W\calO_K$, where $R_0$ is an integral domain obtained from $W\langle T_1^{\pm 1},\ldots,T_d^{\pm 1}\rangle$ by a finite number of iterations of the following operations:
\begin{itemize}
 \item $p$-adic completion of an \'etale extension;
 \item $p$-adic completion of a localization;
 \item completion with respect to an ideal containing $p$.
\end{itemize}
\end{set-up}

\begin{rem}
To apply Faltings' almost purity theorem, Brinon \cite[p.~7]{brinon-relative} further assumes that $W[T_1^{\pm 1},\ldots,T_d^{\pm 1}]\rightarrow R_0$ has geometrically regular fibers and that $k\rightarrow R\otimes_{\calO_K}k$ is geometrically integral. By \cite[Prop.~5.12]{andreatta} and the fact that any ideal-adic completion of an excellent ring is excellent \cite[Main Thm.~2]{Kurano-Shimomoto}, we see that $W[T_1^{\pm 1},\ldots,T_d^{\pm 1}]\rightarrow R_0$ has geometrically regular fibers for a ring $R_0$ as in Set-up~\ref{set-up:base ring}.
We also note that the latter assumption can be dropped. Indeed, if we let $k'$ be the integral closure of $k$ inside $\Frac (R\otimes_{\calO_K}k)$,  then $R\otimes_{\calO_K}k$ is geometrically connected (and thus geometrically integral) over $k'$, and $R$ is an $\calO_K\otimes_WW(k')$-algebra.
The claim now follows since Brinon's period rings for $R$ are defined without any reference to $\calO_K$. Finally, we note that if $R$ is a base ring, then it satisfies Brinon's good reduction condition (BR) in \cite[p.~9]{brinon-relative}.
\end{rem}

Let $R$ be a base ring as defined in Set-up~\ref{set-up:base ring}. In the rest of this subsection, we introduce basic objects attached to $R$ that we will use throughout this article.

Let $\varphi\colon R_0\rightarrow R_0$ denote the lift of the Frobenius on $R_0/pR_0$ with $\varphi(T_i)=T_i^p$; this uniquely determines $\varphi$.  Let $\widehat{\Omega}_{R_0}$ denote the module of continuous K\"ahler differentials $\varprojlim_n\Omega_{(R_0/p^nR_0)/(W/p^nW)}$. By \cite[Prop.~2.0.2]{brinon-relative}, we have $\widehat{\Omega}_{R_0} = \bigoplus_{i=1}^d R_0 \cdot d\log{T_i}$.

Let $\overline{R}$ denote the union of finite $R$-subalgebras $R'$ of a fixed algebraic closure of $\Frac{R}$ such that $R'[p^{-1}]$ is \'etale over $R[p^{-1}]$.  Set 
\[
\calG_{R}\coloneqq \Gal(\overline{R}[p^{-1}]/R[p^{-1}]).
\]
Let $\mathrm{Rep}_{\Q_p}(\calG_R)$ denote the category of finite-dimensional $\Q_p$-vector spaces with continuous $\calG_R$-action. We call its objects \emph{$\Q_p$-representations of $\calG_R$} for short. Similarly, let $\mathrm{Rep}_{\mathbf{Z}_p}(\calG_R)$ denote the category of finite $\mathbf{Z}_p$-modules equipped with continuous $\calG_R$-action and let $\mathrm{Rep}_{\mathbf{Z}_p}^{\mathrm{pr}}(\calG_R)$ denote the full subcategory consisting of finite free objects.

\begin{rem}\label{rem:Galois rep and local system}
Assume that $R$ is of topologically finite type over $\calO_K$ (for example, $R$ is small over $\calO_K$). If we equip $R$ with the $p$-adic topology, then $\mathrm{Rep}_{\Z_p}^\mathrm{pr}(\calG_R)$ (resp.~$\mathrm{Rep}_{\Q_p}(\calG_R)$) is equivalent to the category of $\Z_p$-local systems (resp.~isogeny $\Z_p$-local systems) on the \'etale site of the adic space $\Spa(R[p^{-1}],R)$ by \cite[Ex.~1.6.6 ii)]{Huber-etale} and \cite[Rem.~1.4.4]{kedlaya-liu-relative-padichodge}. Note also that $\Spa(R[p^{-1}],R)$ is the adic generic fiber of $\Spf R$.
\end{rem}

Let $\overline{R}^\wedge$ be the $p$-adic completion of $\overline{R}$ and let $\overline{R}^{\flat}$ be its tilt $\varprojlim_\varphi \overline{R}/p\overline{R}$. 
Set $\Ainf(\overline{R}) \coloneqq W(\overline{R}^{\flat})$.
The first projection $\overline{R}^{\flat}\rightarrow \overline{R}/p\overline{R}$ lifts uniquely to a surjective $W$-algebra homomorphism $\theta\colon \Ainf(\overline{R})\rightarrow \overline{R}^\wedge$. 

\begin{notation}\label{notation:Sigma}
Let $\mathfrak{S}=\mathfrak{S}_R \coloneqq R_0[\![u]\!]$ equipped with the Frobenius given by $\varphi(u) = u^p$. Let $\mathcal{O}_{\mathcal{E}}$ be the $p$-adic completion of $\mathfrak{S}[u^{-1}]$, equipped with the Frobenius $\varphi$ extending that on $\mathfrak{S}$. Note that $E$ is invertible in $\calO_\calE$ and the map $\fkS[E^{-1}]\rightarrow \calO_\calE$ induces an isomorphism $\fkS[E^{-1}]^\wedge_p\xrightarrow{\cong} \calO_\calE$.
\end{notation}

We recall a result about the Frobenius on $\mathfrak{S}$.

\begin{lem} \label{lem:Frobenius-faith-flat}
The map $\varphi\colon \mathfrak{S} \rightarrow \mathfrak{S}$ is classically faithfully flat. Moreover, $\mathfrak{S}$ as a module over itself via $\varphi$ is finite free.	
\end{lem}

\begin{proof}
By \cite[Lem.~7.1.5]{brinon-relative}, $\varphi \colon \mathfrak{S} \rightarrow \mathfrak{S}$ is classically flat. Let $\mathfrak{q}\subset \mathfrak{S}$ be any maximal ideal. Since $\mathfrak{S}$ is $p$-adically complete, we have $p \in \mathfrak{q}$. Thus, $\varphi(\mathfrak{q}) \subset \mathfrak{q}$, which implies $\varphi \colon \mathfrak{S} \rightarrow \mathfrak{S}$ is classically faithfully flat.

For the second part, consider $\mathfrak{S}$ as a module over itself via $\varphi$. Note that $\mathfrak{S}/(p)$ has a finite $p$-basis. By Nakayama's lemma, a lift of a $p$-basis to $\mathfrak{S}$ generates $\mathfrak{S}$. There cannot be any non-trivial relation among such a lift, since $\mathfrak{S}$ is $p$-torsion free.   	
\end{proof}

\begin{notation}\label{notation:L}
Note that $(p)$ (resp.~$(\pi)$) is a prime ideal of $R_0$ (resp.~$R$).
We let $\calO_{L_0}$ (resp.~$\calO_L$) denote the $p$-adic completion of the localization $(R_0)_{(p)}$ (resp.~$R_{(\pi)}$).
Then $\calO_{L_0}$ and $\calO_L$ are complete discrete valuation rings (CDVR's for short) with the same residue field $\Frac (R_0/(p))=\Frac (R/(\pi))$.
Set $L_0\coloneqq \calO_{L_0}[p^{-1}]$ and $L\coloneqq\calO_L[p^{-1}]$.
The Frobenius $\varphi$ on $R_0$ extends to $\varphi\colon \calO_{L_0}\rightarrow \calO_{L_0}$. Note that $\calO_L$ is also a base ring. When we work on $\calO_L$ for a fixed base ring $R$, we simply write $R=\calO_L$ by abuse of notation.

Define $\calO_{K_{0,g}}$ to be the $p$-adic completion of $\varinjlim_{\varphi}\calO_{L_0}$ and let $\calO_{K_g}\coloneqq \calO_{K_{0,g}}\otimes_{W}\calO_K$. Set $K_{0,g}\coloneqq \calO_{K_{0,g}}[p^{-1}]$ and $K_g\coloneqq\calO_{K_g}[p^{-1}]$.
Note that there is a unique $\varphi$-compatible isomorphism $\calO_{K_{0,g}}\cong W(k_g)$ that reduces to the identity modulo $p$, where $k_g$ denotes $\varinjlim_{\varphi}\Frac(R_0/(p))$. 
Hence $\calO_{K_{0,g}}$ and $\calO_{K_g}$ are CDVR's with the same \emph{perfect} residue field $k_g$.
Note that the structure map $R_0\rightarrow \calO_{K_{0,g}}$ factors through $\calO_{L_0}\rightarrow \calO_{K_{0,g}}=W(k_g)$. 

We will often deduce our statements over $\mathfrak{S}$ from those over $\mathcal{O}_{\mathcal{E}}$ and $\mathcal{O}_L$ (or $\mathcal{O}_{K_g}$) by taking certain intersections of modules (e.g. Construction~\ref{construction:quasi-Kisin module} and the proof of Theorem \ref{thm:main} (i)). For proofs, we need the following localization method (cf.~\cite[\S~3.3]{brinon-relative}): fix an algebraic closure $\overline{K_g}$ and let $\mathcal{O}_{\overline{K_g}}$ denote its ring of integers. Let $\mathcal{P}$ be the set of minimal prime ideals of $\overline{R}$ containing $p$. For each $\mathfrak{p}\in \mathcal{P}$, fix a continuous ring homomorphism $(\overline{R}_{\mathfrak{p}})^{\wedge} \rightarrow (\calO_{\overline{K_g}})^{\wedge}$ extending $R_{(\pi)}\rightarrow \calO_{\overline{K_g}}$, where $(\cdots)^{\wedge}$ denotes the $p$-adic completion.
Taking the product over the $\mathfrak{p}$'s induces injective maps
\[
\overline{R}^{\wedge} \rightarrow \prod_{\mathfrak{p} \in \mathcal{P}} (\overline{R}_{\mathfrak{p}})^{\wedge} \rightarrow \prod_{\mathfrak{p} \in \mathcal{P}} (\calO_{\overline{K_g}})^{\wedge}
\quad\text{and}\quad
\overline{R}^{\flat} \rightarrow \prod_{\mathfrak{p} \in \mathcal{P}} \calO_{\overline{K_g}}^{\flat}.
\]

\end{notation}

In \S~\ref{sec:prism-cryst-loc-syst}, we will consider the absolute prismatic site on a $p$-adic formal scheme. In the affine case, we usually make the following additional assumption:

\begin{assumption} \label{assumption:base-ring-sec-3.4}
The base ring $R$ is small over $\mathcal{O}_K$ or $R = \mathcal{O}_L$.  We equip $R$ with $p$-adic topology. In particular, $\Spf R$ is a smooth $p$-adic formal scheme over $\calO_K$ (or $\calO_L$ in the second case).
Note that $\mathrm{Rep}^\mathrm{pr}_{\Z_p}(\calG_R)$ is equivalent to the category of \'etale $\Z_p$-local systems on the adic generic fiber $\Spa(R[p^{-1}],R)$ of $\Spf R$ (cf.~Remark~\ref{rem:Galois rep and local system}).
\end{assumption}

\subsection{Crystalline representations}\label{sec:crystalline representations}
Let $R$ be a base ring.
In this subsection, we will review the crystalline period ring $\OB_{\cris}(\overline{R})$ and the notion of crystalline representations of the Galois group $\calG_R$ of $R[p^{-1}]$ developed in \cite[Chap.~6]{brinon-relative}.

Recall the surjective $W$-algebra homomorphism $\theta\colon \Ainf(\overline{R})\coloneqq W(\overline{R}^{\flat})\rightarrow \overline{R}^\wedge$. 
Define $\A_{\cris}(\overline{R})$ to be the $p$-adic completion of the divided power envelope of $\Ainf(\overline{R})$ with respect to $\Ker\theta$. 
Choose a non-trivial compatible system of $p$-power roots of unity: $\varepsilon_n\in \overline{R}$ with $\varepsilon_0=1, \varepsilon_1\neq 1$, and $\varepsilon_n=\varepsilon_{n+1}^p$. Set $\varepsilon = (\varepsilon_n)_n\in \overline{R}^{\flat}$ and $t\coloneqq \log [\varepsilon]\in \A_{\cris}(\overline{R})$. 
Define $\B_{\cris}(\overline{R})\coloneqq \A_{\cris}(\overline{R})[p^{-1},t^{-1}]$.

Extend the map $\theta$ to $\theta_{R_0}\colon R_0\otimes_W \Ainf(\overline{R})\rightarrow \overline{R}^{\wedge}$. 
Define $\OA_{\cris}(\overline{R})$ to be the $p$-adic completion of the divided power envelope of $R_0\otimes_W \Ainf(\overline{R})$ with respect to $\Ker \theta_{R_0}$. Define $\OB_{\cris}(\overline{R})\coloneqq \OA_{\cris}(\overline{R})[p^{-1},t^{-1}]$.

\begin{rem}\hfill
\begin{enumerate}
    \item Our period rings $\B_{\cris}(\overline{R})$ and $\OB_{\cris}(\overline{R})$ are written as $\mathrm{B}_{\cris}^\nabla(R)$ and $\mathrm{B}_{\cris}(R)$ respectively in \cite{brinon-relative}. 
    \item When $K$ is absolutely unramified and $R$ is of topologically finite type over $\calO_K=W$, Tan and Tong define the crystalline period sheaves $\mathbb{B}_{\cris}$ and $\calO\mathbb{B}_{\cris}$ on the pro-\'etale site of $\Spa(R[p^{-1}],R)$ \cite[Def.~2.4, 2.9]{Tan-Tong}. In this case, $U\coloneqq \Spa(\overline{R}^\wedge[p^{-1}],\overline{R}^\wedge)$ is an affinoid perfectoid object of the pro-\'etale site. We then have $\B_{\cris}(\overline{R})=\mathbb{B}_{\cris}(U)$ and $\OB_{\cris}(\overline{R})=\calO\mathbb{B}_{\cris}(U)$.
    See Proposition~\ref{prop:Tan-Tong Cor.2.19}.
    \item The ring $\calO_L=(R_{(\pi)})^\wedge$ is also a  base ring. In this case, $\B_{\cris}(\overline{\calO_L})$ and $\OB_{\cris}(\overline{\calO_L})$ are studied in \cite{Brinon-cris} and written as $\mathrm{B}_{\cris}^\nabla$ and $\mathrm{B}_{\cris}$, respectively. The notation $\mathrm{B}_{\cris}$ is also used in \cite{brinon-trihan}.
\end{enumerate}
\end{rem}

The crystalline period ring $\OB_{\cris}(\overline{R})$ has a natural $\calG_R$-action and a Frobenius endomorphism $\varphi$ extending those on $R_0\otimes_W \Ainf(\overline{R})$, and there is a natural $\B_{\cris}(\overline{R})$-linear integrable connection $\nabla\colon \OB_{\cris}(\overline{R})\rightarrow \OB_{\cris}(\overline{R})\otimes_{R_0}\widehat{\Omega}_{R_0}$.
Moreover, $R\otimes_{R_0}\OB_{\cris}(\overline{R})$ is equipped with a filtration by $R[p^{-1}]$-modules, which is compatible with the natural PD-filtration on $\A_{\cris}(\overline{R})$. See \cite[Chap.~6]{brinon-relative} for the detail of these structures.

The following result on the crystalline period ring will be used later:

\begin{lem}[({\cite[Prop.~6.1.5]{brinon-relative}})] \label{lem:cryst-period-ring}
Choose a compatible system $(T_{i,n})$ of $p$-power roots of $T_i$ in $\overline{R}$ with $T_{i,0}=T_i$, and let $T_i^\flat\in \overline{R}^\flat$ denote the corresponding element. 
The map $X_i \mapsto T_i\otimes1-1\otimes[T_i^\flat]$ induces an $\mathbf{A}_{\mathrm{cris}}(\overline{R})$-linear isomorphism
\[
\mathbf{A}_{\mathrm{cris}}(\overline{R})\{X_1, \ldots, X_d\} \cong \mathbf{OA}_{\mathrm{cris}}(\overline{R}),
\]	
where the former ring denotes the $p$-adically completed divided power polynomial with variables $X_i$ and coefficients in $\mathbf{A}_{\mathrm{cris}}(\overline{R})$. 
\end{lem}

Let us recall the  definition of crystalline representations.
\begin{defn}\label{defn:crystalline representations}
 For $V\in\mathrm{Rep}_{\Q_p}(\calG_R)$, set
\[
D_{\cris}(V)\coloneqq (\OB_{\cris}(\overline{R})\otimes_{\Q_p}V)^{\calG_R}\quad \text{and}\quad 
D_{\cris}^\vee(V)\coloneqq \mathrm{Hom}_{\calG_R}(V,\OB_{\cris}(\overline{R})).
\]
Then $D_{\cris}(V)$ is a finite projective $R_0[p^{-1}]$-module of rank at most $\dim_{\Q_p}V$ equipped with a natural $\varphi$ and $\nabla$ structure induced from $\OB_{\cris}(\overline{R})$, and $R\otimes_{R_0}D_{\cris}(V)$ has a filtration induced from $R\otimes_{R_0}\OB_{\cris}(\overline{R})$. The natural map
\[
\alpha_{\cris}(V)\colon \OB_{\cris}(\overline{R})\otimes_{R_0[p^{-1}]}D_{\cris}(V) \rightarrow \OB_{\cris}(\overline{R})\otimes_{\Q_p}V
\]
is injective by \cite[Prop.~8.2.6]{brinon-relative}.
We say that $V$ is \emph{$R_0$-crystalline} if $\alpha_{\cris}(V)$ is an isomorphism.
By \cite[Prop.~8.3.5]{brinon-relative}, this notion depends only on $R$, not on $R_0$. Hence we simply say that $V$ is crystalline from now on.

By \cite[Thm.~8.4.2]{brinon-relative}, $V$ is crystalline if and only if $V^{\vee}$ is crystalline. Note also that $D_{\cris}^\vee(V)=D_{\cris}(V^\vee)=\mathrm{Hom}_{R_0[p^{-1}]}(D_{\cris}(V),R_0[p^{-1}])$. We will mainly use $D_{\cris}^\vee(V)$ in this paper.

A finite free $\Z_p$-representation $T\in \mathrm{Rep}_{\Z_p}^{\mathrm{pr}}(\calG_R)$ is called \emph{crystalline} if the associated $\Q_p$-representation $T\otimes_{\Z_p}\Q_p$ is crystalline.
\end{defn}

Finally, let us explain the functoriality.
Let $R'=R_0'\otimes_{W}\calO_K$ be another base ring and assume that there exists a $\varphi$-equivariant ring homomorphism $g\colon R_0\rightarrow R_0'$ that extends to $g\colon \overline{R}\rightarrow \overline{R'}$. By the change of paths for \'etale fundamental groups, $g$ induces a continuous group homomorphism $\calG_{R'}\rightarrow \calG_R$ and thus a natural $\otimes$-functor $\mathrm{Rep}_{\Q_p}(\calG_R)\rightarrow \mathrm{Rep}_{\Q_p}(\calG_{R'})$.
The map $g$ also induces a ring homomorphism $\OB_{\cris}(\overline{R})\rightarrow \OB_{\cris}(\overline{R'})$, and the latter is compatible with Frobenii and Galois actions.
 
\begin{lem}\label{lem:base change map for Dcris}
With the notation as above, for $V\in \mathrm{Rep}_{\Q_p}(\calG_R)$, the map $R_0'[p^{-1}]\otimes_{R_0[p^{-1}]}\OB_\cris(\overline{R})\otimes_{\Q_p}V\rightarrow \OB_\cris(\overline{R'})\otimes_{\Q_p}V$ induces a $\varphi$-equivariant morphism of $R_0'[p^{-1}]$-modules
\begin{equation}\label{eq:base change map for Dcris}
R_0'[p^{-1}]\otimes_{R_0[p^{-1}]}D_{\cris}(V)\rightarrow D_{\cris}(V|_{\calG_{R'}}).    
\end{equation}
Moreover, if $V$ is crystalline, then $V|_{\calG_{R'}}$ is crystalline and the above map is an isomorphism.
\end{lem}

\begin{proof}
The first assertion is obvious. Now assume that $V$ is crystalline. Consider the composite of $\OB_\cris(\overline{R'})$-linear maps
\begin{align*}
\alpha\colon \OB_\cris(\overline{R'})\otimes_{R_0'[p^{-1}]}(R_0'[p^{-1}] \otimes_{R_0[p^{-1}]}D_{\cris}(V))&\rightarrow \OB_\cris(\overline{R'})\otimes_{R_0'[p^{-1}]}D_{\cris}(V|_{\calG_{R'}})\\
&\xhookrightarrow{\alpha_\cris(V|_{\calG_{R'}})}\OB_\cris(\overline{R'})\otimes_{\Q_p}V.    
\end{align*}
Observe that $\alpha$ is the base change of $\alpha_\cris(V)$ along the map $\OB_\cris(\overline{R})\rightarrow\OB_\cris(\overline{R'})$. Since $V$ is crystalline, $\alpha$ is an isomorphism. Moreover, the second map $\alpha_\cris(V|_{\calG_{R'}})$ in $\alpha$ is injective. Hence $\alpha_\cris(V|_{\calG_{R'}})$ is an isomorphism and thus $V|_{\calG_{R'}}$ is crystalline. We also see that the first map in $\alpha$ is an isomorphism. Since the map $R_0'[p^{-1}]\rightarrow \OB_\cris(\overline{R'})$ is faithfully flat by \cite[Thm.~6.3.8]{brinon-relative}, the morphism \eqref{eq:base change map for Dcris} is an isomorphism.
\end{proof}

\subsection{\'Etale \texorpdfstring{$\varphi$}{phi}-modules} \label{sec:etale phi-module}

The classical theory of \'etale $\varphi$-modules and Galois representations is generalized to our relative setting in \cite{kim-groupscheme-relative}. We briefly review some necessary facts discussed in \cite{kim-groupscheme-relative} and \cite{liu-moon-rel-crys-rep-p-div-gps-small-ramification}. Recall Notation~\ref{notation:Sigma}: $\mathfrak{S}=\mathfrak{S}_R \coloneqq R_0[\![u]\!]$ equipped with the Frobenius given by $\varphi(u) = u^p$; $\mathcal{O}_{\mathcal{E}}$ is the $p$-adic completion of $\mathfrak{S}[u^{-1}]$, equipped with the Frobenius $\varphi$ extending that on $\mathfrak{S}$.

\begin{defn} \label{defn:etale phi-module}
An \textit{\'{e}tale} $(\varphi, \mathcal{O}_\mathcal{E})$-\textit{module} is a pair $(\mathcal{M}, \varphi_{\mathcal{M}})$ where $\mathcal{M}$ is a finitely generated $\mathcal{O}_\mathcal{E}$-module and $\varphi_{\mathcal{M}}\colon \mathcal{M} \rightarrow \mathcal{M}$ is a $\varphi$-semi-linear endomorphism such that $1\otimes\varphi_{\mathcal{M}}\colon \varphi^*\mathcal{M} \rightarrow \mathcal{M}$ is an isomorphism. We say that an \'{e}tale $(\varphi, \mathcal{O}_{\mathcal{E}})$-module is \textit{projective} (resp. \textit{torsion}) if the underlying $\mathcal{O}_{\mathcal{E}}$-module $\mathcal{M}$ is projective (resp. $p$-power torsion).	

Let $\mathrm{Mod}_{\mathcal{O}_{\mathcal{E}}}$ denote the category of \'{e}tale $(\varphi, \mathcal{O}_{\mathcal{E}})$-modules whose morphisms are $\mathcal{O}_\mathcal{E}$-linear maps compatible with Frobenii. Let  $\mathrm{Mod}_{\mathcal{O}_{\mathcal{E}}}^{\mathrm{pr}}$ and  $\mathrm{Mod}_{\mathcal{O}_{\mathcal{E}}}^{\mathrm{tor}}$, respectively, denote the full subcategories of projective and torsion objects. Note that we have a natural notion of tensor products for \'{e}tale $(\varphi, \mathcal{O}_\mathcal{E})$-modules, and duals are defined for projective and torsion objects.
\end{defn}

We use \'{e}tale $(\varphi, \mathcal{O}_{\mathcal{E}})$-modules to study certain Galois representations as follows. We refer the reader to \cite{scholze-perfectoid} for definitions and facts on perfectoid algebras. Recall that $\pi$ denotes a uniformizer in  $\mathcal{O}_K$. For integers $n \geq 0$, compatibly choose $\pi_n \in \overline{K}$ such that $\pi_0 = \pi$ and $\pi_{n+1}^p = \pi_n$, and let $K_{\infty}$ be the $p$-adic completion of $\bigcup_{n \geq 0} K(\pi_n)$. Then $K_{\infty}$ is a perfectoid field, and $(\overline{R}^\wedge[p^{-1}], \overline{R}^\wedge)$ is a perfectoid affinoid $K_{\infty}$-algebra. Let $K_{\infty}^\flat$ denote the tilt of $K_{\infty}$, and set $\pi^\flat \coloneqq (\pi_n) \in K_{\infty}^\flat$. 
 
Let $E_{R_\infty}^+ = \mathfrak{S}/p\mathfrak{S}$, and let $\tilde{E}_{R_\infty}^+$ be the $u$-adic completion of $\varinjlim_{\varphi}E_{R_\infty}^+$. Let $E_{R_\infty} = E_{R_\infty}^+[u^{-1}]$ and $\tilde{E}_{R_\infty} = \tilde{E}_{R_\infty}^+[u^{-1}]$. By \cite[Prop.~5.9]{scholze-perfectoid}, $(\tilde{E}_{R_\infty}, \tilde{E}_{R_\infty}^+)$ is a perfectoid affinoid $K_{\infty}^\flat$-algebra, and we have a natural injective map $(\tilde{E}_{R_\infty}, \tilde{E}_{R_\infty}^+) \hookrightarrow (\overline{R}^\flat[(\pi^\flat)^{-1}], \overline{R}^\flat)$ given by $u \mapsto \pi^\flat$. 

Consider 
\begin{equation} \label{eq:R-infinity}
\tilde{R}_{\infty} \coloneqq W(\tilde{E}_{R_\infty}^+)\otimes_{W(K_{\infty}^{\flat \circ}), ~\theta} \mathcal{O}_{K_{\infty}}.	
\end{equation}
By \cite[Rem.~5.19]{scholze-perfectoid}, $(\tilde{R}_\infty[p^{-1}], \tilde{R}_\infty)$ is a perfectoid affinoid $K_{\infty}$-algebra whose tilt is $(\tilde{E}_{R_\infty}, \tilde{E}_{R_\infty}^+)$. Furthermore, we have a natural injective map $(\tilde{R}_\infty[p^{-1}], \tilde{R}_\infty) \hookrightarrow (\overline{R}^\wedge[p^{-1}], \overline{R}^\wedge)$ whose tilt is $(\tilde{E}_{R_\infty}, \tilde{E}_{R_\infty}^+) \hookrightarrow (\overline{R}^\flat[(\pi^{\flat})^{-1}], \overline{R}^\flat)$. If we write $\mathcal{G}_{\tilde{R}_\infty}$ for $\pi_1^{\text{\'{e}t}}(\operatorname{Spec}\tilde{R}_\infty[p^{-1}])$, we then have a continuous map of Galois groups $\mathcal{G}_{\tilde{R}_\infty} \rightarrow \mathcal{G}_R$, which is a closed embedding by \cite[Prop.~5.4.54]{gabber-almost}. By \cite[Thm.~7.12]{scholze-perfectoid}, we can canonically identify $\overline{R}^\flat[(\pi^{\flat})^{-1}]$ with the $\pi^\flat$-adic completion of the affine ring of a universal pro-\'etale covering of $\operatorname{Spec}\tilde{E}_{R_\infty}$. 
Let $\mathcal{G}_{\tilde{E}_{R_\infty}}$ be the Galois group corresponding to the universal pro-\'etale covering. Then we have a canonical isomorphism $\mathcal{G}_{\tilde{E}_{R_\infty}} \cong \mathcal{G}_{\tilde{R}_\infty}$. 

There exists a unique $W(k)$-linear map $R_0 \rightarrow W(\overline{R}^{\flat})$ which maps $T_i$ to $[T_i^\flat]$ and is compatible with Frobenii (see Lemma~\ref{lem:cryst-period-ring} for the definition of $[T_i^\flat]$). This induces a $\varphi$-equivariant embedding $\mathfrak{S} \rightarrow W(\overline{R}^{\flat})$ given by $u \mapsto [\pi^\flat]$, which further extends to an embedding $\mathcal{O}_{\mathcal{E}} \rightarrow W(\overline{R}^\flat[(\pi^{\flat})^{-1}])$. Let $\mathcal{O}_\mathcal{E}^{\mathrm{ur}}$ be the union of finite \'etale $\mathcal{O}_\mathcal{E}$-subalgebras of $W(\overline{R}^\flat[(\pi^{\flat})^{-1}])$, and let $\widehat{\mathcal{O}}_{\mathcal{E}}^{\mathrm{ur}}$ be its $p$-adic completion. We also define $\widehat{\mathfrak{S}}^{\mathrm{ur}} \coloneqq \widehat{\mathcal{O}}_{\mathcal{E}}^{\mathrm{ur}} \cap W(\overline{R}^\flat) \subset W(\overline{R}^\flat[\pi^\flat]^{-1})$. We note that the definitions of these rings in \cite[p.~8201]{kim-groupscheme-relative} are incorrect but that the results concerning these rings in \textit{ibid.} hold with the correct definitions: since $(\mathcal{O}_{\mathcal{E}}, (p))$ is a henselian pair, $\mathcal{O}_\mathcal{E}^{\mathrm{ur}}/(p)=\widehat{\mathcal{O}}_{\mathcal{E}}^{\mathrm{ur}}/(p)$ is the union of finite \'etale $E_{R_\infty}$-subalgebras of $\overline{R}^\flat[(\pi^{\flat})^{-1}]$. In particular, we have $\mathrm{Aut}_{\mathcal{O}_\mathcal{E}}(\mathcal{O}_\mathcal{E}^{\mathrm{ur}}) \cong \mathcal{G}_{E_{R_\infty}} \coloneqq \pi_1^{\text{\'et}}(\operatorname{Spec}E_{R_\infty})$. By \cite[Prop.~5.4.54]{gabber-almost} and \cite[Lem.~7.5]{scholze-perfectoid}, we have $\mathcal{G}_{E_{R_\infty}} \cong \mathcal{G}_{\tilde{E}_{R_\infty}} \cong \mathcal{G}_{\tilde{R}_\infty}$. This induces $\mathcal{G}_{\tilde{R}_\infty}$-action on $\widehat{\mathcal{O}}_{\mathcal{E}}^{\mathrm{ur}}$. The following is proved in \cite{kim-groupscheme-relative}.

\begin{lem}[(cf.~{\cite[Lem.~7.5 and 7.6]{kim-groupscheme-relative}})] \label{lem:etale-ring-frob-galois} 
We have $(\widehat{\mathcal{O}}_{\mathcal{E}}^{\mathrm{ur}})^{\mathcal{G}_{\tilde{R}_\infty}} = \mathcal{O}_\mathcal{E}$ and the same holds modulo $p^n$. Furthermore, there exists a unique $\mathcal{G}_{\tilde{R}_\infty}$-equivariant ring endomorphism $\varphi$ on $\widehat{\mathcal{O}}_{\mathcal{E}}^{\mathrm{ur}}$ lifting the $p$-th power Frobenius on $\widehat{\mathcal{O}}_{\mathcal{E}}^{\mathrm{ur}}/(p)$ and extending $\varphi$ on $\mathcal{O}_\mathcal{E}$. The inclusion $\widehat{\mathcal{O}}_{\mathcal{E}}^{\mathrm{ur}} \hookrightarrow W(\overline{R}^\flat[(\pi^{\flat})^{-1}])$ is $\varphi$-equivariant where the latter ring is given the Witt vector Frobenius.
\end{lem}

Let $\mathrm{Rep}_{\mathbf{Z}_p}(\mathcal{G}_{\tilde{R}_\infty})$ denote the category of finite $\mathbf{Z}_p$-modules equipped with continuous $\mathcal{G}_{\tilde{R}_\infty}$-action, and let $\mathrm{Rep}_{\mathbf{Z}_p}^{\mathrm{pr}}(\mathcal{G}_{\tilde{R}_\infty})$ and $\mathrm{Rep}_{\mathbf{Z}_p}^{\mathrm{tor}}(\mathcal{G}_{\tilde{R}_\infty})$, respectively, denote the full subcategories of free and torsion objects. For $\mathcal{M} \in \mathrm{Mod}_{\mathcal{O}_\mathcal{E}}$ and $T \in \mathrm{Rep}_{\mathbf{Z}_p}(\mathcal{G}_{\tilde{R}_\infty})$, define 
\[
T(\mathcal{M}) \coloneqq (\widehat{\mathcal{O}}_{\mathcal{E}}^{\mathrm{ur}}\otimes_{\mathcal{O}_\mathcal{E}}\mathcal{M})^{\varphi = 1} \quad\text{and}\quad
\mathcal{M}(T) \coloneqq (\widehat{\mathcal{O}}_{\mathcal{E}}^{\mathrm{ur}}\otimes_{\mathbf{Z}_p}T)^{\mathcal{G}_{\tilde{R}_\infty}}.
\]
For a torsion \'etale $\varphi$-module $\mathcal{M} \in \mathrm{Mod}_{\mathcal{O}_{\mathcal{E}}}^{\mathrm{tor}}$, we define its \textit{length} to be the length of $(\mathcal{O}_{\mathcal{E}})_{(p)}\otimes_{\mathcal{O}_{\mathcal{E}}}\mathcal{M}$ as an $(\mathcal{O}_{\mathcal{E}})_{(p)}$-module. The following equivalence is proved in \cite{kim-groupscheme-relative} (see also \cite[Prop.~2.5]{liu-moon-rel-crys-rep-p-div-gps-small-ramification}).

\begin{prop}[(cf.~{\cite[Prop.~7.7]{kim-groupscheme-relative}, \cite[Prop.~2.5]{liu-moon-rel-crys-rep-p-div-gps-small-ramification}})]  \label{prop:etale-gal-equiv}
The assignments $T(\cdot)$ and $\mathcal{M}(\cdot)$ are exact equivalences \emph{(}quasi-inverse of each other\emph{)} of $\otimes$-categories between $\mathrm{Mod}_{\mathcal{O}_{\mathcal{E}}}$ and $\mathrm{Rep}_{\mathbf{Z}_p}(\mathcal{G}_{\tilde{R}_\infty})$. Moreover, $T(\cdot)$ and $\mathcal{M}(\cdot)$ restrict to rank-preserving equivalence of categories between $\mathrm{Mod}_{\mathcal{O}_{\mathcal{E}}}^{\mathrm{pr}}$ and $\mathrm{Rep}_{\mathbf{Z}_p}^{\mathrm{pr}}(\mathcal{G}_{\tilde{R}_\infty})$ and length-preserving equivalence of categories between $\mathrm{Mod}_{\mathcal{O}_{\mathcal{E}}}^{\mathrm{tor}}$ and $\mathrm{Rep}_{\mathbf{Z}_p}^{\mathrm{tor}}({\mathcal{G}}_{\tilde{R}_\infty})$. In both cases, $T(\cdot)$ and $\mathcal{M}(\cdot)$ commute with taking duals. For $T\in \mathrm{Rep}_{\mathbf{Z}_p}(\mathcal{G}_{\tilde{R}_\infty})$, the natural map $\widehat{\mathcal{O}}_{\mathcal{E}}^{\mathrm{ur}}\otimes_{\mathcal{O}_\mathcal{E}}\mathcal{M}(T) \rightarrow \widehat{\mathcal{O}}_{\mathcal{E}}^{\mathrm{ur}}\otimes_{\mathbf{Z}_p}T$ is an $\widehat{\mathcal{O}}_{\mathcal{E}}^{\mathrm{ur}}$-linear isomorphism compatible with $\varphi$ and $\mathcal{G}_{\tilde{R}_\infty}$-actions, and a similar statement holds for $\mathrm{Mod}_{\mathcal{O}_{\mathcal{E}}}$.  
\end{prop}

\begin{proof}
All statements except the last one are in \cite[Prop.~7.7]{kim-groupscheme-relative}, whose proof is based on \cite[Prop.~4.1.1]{katz}. The last statement also follows from the proof of \cite[Prop.~4.1.1]{katz} by the standard d\'evissage.
\end{proof}

For $\mathcal{M} \in \mathrm{Mod}_{\mathcal{O}_\mathcal{E}}^{\mathrm{pr}}$ and $T \in \mathrm{Rep}_{\mathbf{Z}_p}^{\mathrm{pr}}(\mathcal{G}_{\tilde{R}_\infty})$, we can consider the contravariant functors
\[
T^{\vee}(\mathcal{M}) \coloneqq \mathrm{Hom}_{\mathcal{O}_\mathcal{E}, \varphi}(\mathcal{M}, \widehat{\mathcal{O}}_{\mathcal{E}}^{\mathrm{ur}})
\quad\text{and}\quad
\mathcal{M}^{\vee}(T) \coloneqq \mathrm{Hom}_{\mathcal{G}_{\tilde{R}_\infty}}(T, \widehat{\mathcal{O}}_{\mathcal{E}}^{\mathrm{ur}}).
\]
We have natural isomorphisms
\[
T^{\vee}(\mathcal{M}) \cong T(\mathcal{M}^{\vee})
\quad\text{and}\quad
\mathcal{M}^{\vee}(T) \cong \mathcal{M}(T^{\vee}),
\]
and these contravariant functors give equivalences of categories between $\mathrm{Mod}_{\mathcal{O}_{\mathcal{E}}}^{\mathrm{pr}}$ and $\mathrm{Rep}_{\mathbf{Z}_p}^{\mathrm{pr}}(\mathcal{G}_{\tilde{R}_\infty})$ by Proposition~\ref{prop:etale-gal-equiv}.

We now explain certain functoriality of above constructions. Let $R_0'$ be another base ring over $W(k)\langle T_1^{\pm 1}, \ldots, T_d^{\pm 1}\rangle$ as in \S~2.1 equipped with Frobenius, and suppose a $\varphi$-equivariant map $R_0 \to  R_0'$ of $W(k)\langle T_1^{\pm 1}, \ldots, T_d^{\pm 1}\rangle$-algebras is given. 
Consider the induced $\mathcal{O}_K$-linear extension $R = R_0\otimes_{W(k)}\mathcal{O}_K \rightarrow R' \coloneqq R_0'\otimes_{W(k)}\mathcal{O}_K$. By fixing an algebraic closure of $\mathrm{Frac}(R')$, we have a map $\overline{R} \rightarrow \overline{R'}$, and this induces $\tilde{R}_{\infty} \rightarrow \tilde{R'}_{\infty}$ by the constructions given in \eqref{eq:R-infinity}. Hence, the corresponding map of Galois groups $\mathcal{G}_{R'} \rightarrow \mathcal{G}_R$ restricts to $\mathcal{G}_{\tilde{R'}_\infty} \rightarrow \mathcal{G}_{\tilde{R}_\infty}$. Let $\mathfrak{S}_{R'} = R_0'[\![u]\!]$ and let $\mathcal{O}_{\mathcal{E}, R'}$ be the $p$-adic completion of $\mathfrak{S}_{R'}[u^{-1}]$. Let $\mathcal{M}_{R'}(\cdot)$ be the functor for the base ring $R'$ constructed similarly as above. If $T \in \mathrm{Rep}_{\mathbf{Z}_p}^{\mathrm{pr}}(\mathcal{G}_{\tilde{R}_\infty})$, then $T$ can be also considered as a $\mathcal{G}_{\tilde{R'}_\infty}$-representation via the map $\mathcal{G}_{\tilde{R'}_\infty} \rightarrow \mathcal{G}_{\tilde{R}_\infty}$. We claim that there is a natural isomorphism $\mathcal{O}_{\mathcal{E}, R'}\otimes_{\mathcal{O}_\mathcal{E}}\mathcal{M}(T) \xrightarrow{\cong} \mathcal{M}_{R'}(T)$ of \'etale $(\varphi, \mathcal{O}_{\mathcal{E}, R'})$-modules. Indeed, the $W(k)\langle T_1^{\pm 1}, \ldots, T_d^{\pm 1}\rangle$-algebra homomorphism $\mathcal{O}_\mathcal{E}\rightarrow \mathcal{O}_{\mathcal{E}, R'}$ extends to a map $\widehat{\mathcal{O}}_{\mathcal{E}}^{\mathrm{ur}}\rightarrow \widehat{\mathcal{O}}_{\mathcal{E},R'}^{\mathrm{ur}}$, which defines the desired map $\mathcal{O}_{\mathcal{E}, R'}\otimes_{\mathcal{O}_\mathcal{E}}\mathcal{M}(T) \rightarrow \mathcal{M}_{R'}(T)$. To see that this is an isomorphism, observe 
\[
\widehat{\mathcal{O}}_{\mathcal{E}, R'}^{\mathrm{ur}} \otimes_{ \calO_{\calE , R'}} (\calO_{\calE, R'}\otimes_{\mathcal{O}_{\mathcal{E}}}\mathcal{M}(T)) = \widehat{\mathcal{O}}_{\mathcal{E}, R'}^{\mathrm{ur}} \otimes_{ \widehat \calO^\mathrm{ur} _{\calE}} (\widehat \calO^\mathrm{ur}_{\calE}\otimes_{\mathcal{O}_{\mathcal{E}}}\mathcal{M}(T) ) \xrightarrow{\cong} \widehat{\mathcal{O}}_{\mathcal{E}, R'}^{\mathrm{ur}} \otimes_{ \widehat \calO^\mathrm{ur} _{\calE}} (\widehat \calO^\mathrm{ur}_{\calE}\otimes_{\mathbf{Z}_p} T)= \widehat{\mathcal{O}}_{\mathcal{E}, R'}^{\mathrm{ur}}\otimes_{\mathbf{Z}_p} T. 
\]
Hence we conclude that $\mathcal{O}_{\mathcal{E}, R'}\otimes_{\mathcal{O}_\mathcal{E}}\mathcal{M}(T) \rightarrow \mathcal{M}_{R'}(T)$ is an isomorphism since it is so after the base change along the faithfully flat map $\calO_{\calE, R'}\rightarrow \widehat{\mathcal{O}}_{\mathcal{E}, R'}^{\mathrm{ur}}$.

We will use this functoriality for the maps of base rings $R \rightarrow \mathcal{O}_L$  and $R \rightarrow \mathcal{O}_{K_g}$ as in Notation \ref{notation:L} in later sections.
For $\mathcal{O}_L$ and $\mathcal{O}_{K_g}$, the relevant rings will be denoted by $\widehat{\mathcal{O}}_{\mathcal{E}, L}^{\mathrm{ur}}$, $\widehat{\mathfrak{S}}_L^{\mathrm{ur}}$, $\widehat{\mathcal{O}}_{\mathcal{E}, g}^{\mathrm{ur}}$, and $\widehat{\mathfrak{S}}_g^{\mathrm{ur}}$.

\section{Completed prismatic \texorpdfstring{$F$}{F}-crystals and crystalline representations} \label{sec:prism-cryst-loc-syst}

This section introduces the notion of completed prismatic $F$-crystals on the absolute prismatic site of $R$ and formulates the main theorem. In \S~\ref{sec:prism-site}, we recall the definition of absolute prismatic site and consider some important examples of prisms. In \S~\ref{sec:completed-prism-F-crystals}-\ref{sec:descent data}, we define completed prismatic $F$-crystals and study their basic properties in the small affine case. In \S~\ref{sec:etale-realization-main-thm}, we study the \'etale realization and formulate our main theorem. In \S~\ref{subsec:height one}, we consider the special case where crystalline representations have Hodge--Tate weights in $[0, 1]$ and study the relation to $p$-divisible groups. Finally, we globalize the \'etale realization functor and the main theorem in \S~\ref{sec:globalization}.

We will frequently use the following lemma.

\begin{lem} \label{lem:intersection-modules-flat-base-change}
Let $A$ be a ring. 
\begin{enumerate}
    \item Let $M$ be a flat $A$-module, and $N_1$, $N_2$ submodules of an $A$-module $N$. Then as submodules of $M\otimes_A N$, we have
\[
M\otimes_A (N_1 \cap N_2) = (M\otimes_A N_1) \cap (M\otimes_A N_2).
\]

 \item 
 Let $M$ be a finite projective $A$-module, and $N$ an $A$-module. Let $\mathcal{I}$ be a (possibly infinite) index set. Suppose for each $i \in \mathcal{I}$, we are given an $A$-submodule $N_i$ of $N$. Then as submodules $M\otimes_A N$, we have
\[
M\otimes_A \bigl(\,\bigcap_{i \in \mathcal{I}} N_i\bigr) = \bigcap_{i \in \mathcal{I}} (M\otimes_A N_i).
\]
\end{enumerate}
\end{lem}

\begin{proof}
(i) is well known. For (ii), it suffices to show that the natural injective map 
\[
f\colon M\otimes_A (\bigcap_{i \in \mathcal{I}} N_i) \rightarrow \bigcap_{i \in \mathcal{I}} (M\otimes_A N_i)
\]
is also surjective. Let $M'$ be an $A$-module such that $M\oplus M'$ is finite free over $A$. Then the map
\[
(M\otimes_A (\bigcap_{i \in \mathcal{I}} N_i))\oplus (M'\otimes_A (\bigcap_{i \in \mathcal{I}} N_i))=(M\oplus M')\otimes_A (\bigcap_{i \in \mathcal{I}} N_i)  \rightarrow \bigcap_{i \in \mathcal{I}} (M\oplus M')\otimes_A N_i = \bigcap_{i \in \mathcal{I}} ((M\otimes_A N_i) \oplus (M'\otimes_A N_i)) 
\]
is an isomorphism. This implies that the above map $f$ is also surjective. 
\end{proof}

\subsection{The absolute prismatic site} \label{sec:prism-site}

We first recall the definition of the absolute prismatic site from \cite{bhatt-scholze-prismaticcohom} and \cite{bhatt-scholze-prismaticFcrystal}. 
Let $\fkX$ be a smooth $p$-adic formal scheme over $\calO_K$ (or a CDVR of mixed characteristic $(0,p)$ such as $\calO_L$ in Notation~\ref{notation:L}).

\begin{defn}(\cite[Def.~2.3]{bhatt-scholze-prismaticFcrystal}) \label{defn:abs-prism-site}
The \emph{absolute prismatic site} $\fkX_\Prism$ of $\fkX$ consists of the pairs $((A,I),\Spf A/I\rightarrow \fkX)$, where $(A,I)$ is a bounded prism and  $\Spf A/I\rightarrow \fkX$ is a morphism of $p$-adic formal schemes. For simplicity, we often omit the structure map $\Spf A/I\rightarrow \fkX$ and simply write $(A,I)$ for an object of $\fkX_\Prism$. 
The morphisms are the opposite of morphisms of bounded prisms over $\fkX$, i.e., the ones compatible with the structure morphisms to $\fkX$. We equip $\fkX_\Prism$ with the topology given by $(p,I)$-completely faithfully flat maps of prisms $(A,I) \rightarrow (B,J)$ over $\fkX$. 
If $\fkX=\Spf R$ is affine, then we also write $R_\Prism$ for $\fkX_\Prism$. Note that the associated topos is replete by \cite[Rem.~2.4]{bhatt-scholze-prismaticFcrystal}.

The prismatic site $\fkX_\Prism$ has a sheaf $\mathcal{O}_{\Prism}$ of rings defined by $\mathcal{O}_{\Prism}(A, I) = A$ and an ideal sheaf $\mathcal{I}_{\Prism} \subset \mathcal{O}_\Prism$ given by $\mathcal{I}_{\Prism}(A, I) = I$ (cf.~\cite[Cor.~3.12]{bhatt-scholze-prismaticcohom}).
A similar argument shows that for each $n\geq 1$, the association $(A,I)\mapsto A/(p,I)^n$ defines a sheaf $\calO_{\Prism,n}$ on $\fkX_\Prism$. Moreover, we have $\calO_\Prism\xrightarrow{\cong}\varprojlim_n\calO_{\Prism,n}\cong \operatorname{Rlim}\calO_{\Prism,n}$ (see Lemma \ref{lem:calF vs (calF_n) for completed crystals} below). Finally, the $\delta$-structure on each $(A, I) \in \fkX_{\Prism}$ induces a ring endomorphism $\varphi\colon \mathcal{O}_{\Prism} \rightarrow \mathcal{O}_{\Prism}$.
\end{defn}

Let us explain the functoriality of the prismatic topoi. Let $f\colon\fkY\rightarrow \fkX$ be a morphism of smooth $p$-adic formal schemes over $\calO_K$.
Then $f$ induces a cocontinuous functor
\[
\fkY_\Prism \rightarrow \fkX_\Prism, ~~((A,I),\iota\colon \Spf A/I\rightarrow \fkY)\mapsto ((A,I),f\circ\iota\colon \Spf A/I\rightarrow \fkX).
\]
Hence we have a morphism of topoi
\[
f_\Prism=(f_\Prism^{-1},f_{\Prism,\ast})\colon \mathrm{Shv}(\fkY_\Prism)\rightarrow \mathrm{Shv}(\fkX_\Prism).
\]
Observe that if $f$ is an open immersion, then $\mathrm{Shv}(\fkY_\Prism)$ is an open subtopos of $\mathrm{Shv}(\fkX_\Prism)$ by $f_\Prism$.

\begin{lem}\label{lem:pushout for prisms along fflat map}
Let $(B,IB)\xleftarrow{b} (A,I)\xrightarrow{c}(C,IC)$ be a diagram of maps of bounded prisms over $\fkX$ with $b$ being $(p,I)$-completely faithfully flat, and let $(B\otimes_AC)^\wedge_{(p,I)}$ denote the classical $(p,I)$-completion of $B\otimes_AC$. Then the pushout of $b$ along $c$ is represented by the map $(C,IC)\rightarrow \bigl((B\otimes_AC)^\wedge_{(p,I)}, I(B\otimes_AC)^\wedge_{(p,I)}\bigr)$ and is $(p,IC)$-completely faithfully flat.
\end{lem}

\begin{proof}
By the proof of \cite[Cor.~3.12]{bhatt-scholze-prismaticcohom}, the pushout of the diagram is represented by the derived $(p, I)$-completion $\widehat{B\otimes^{\mathbb{L}}_A C}$ of $B\otimes^{\mathbb{L}}_A C$. Moreover, it is discrete and classically $(p, I)$-complete, and the map from $C$ is $(p,IC)$-completely faithfully flat. By the proof of \cite[Prop.~3.2]{wu-Gal-rep-prism-F-cryst}, we also have
$H^0(\widehat{B\otimes^{\mathbb{L}}_A C}) = (B\otimes_A C)^{\wedge}_{(p, I)}$, i.e., $\widehat{B\otimes^{\mathbb{L}}_A C}$ is nothing but the classical $(p, I)$-completion of  $B\otimes_A C$.
\end{proof}

Let $R$ be small over $\calO_K$ or $R=\calO_L$ (Assumption~\ref{assumption:base-ring-sec-3.4}).
We now explain several objects of $R_\Prism$ that we will use later.

\begin{eg}[The Breuil--Kisin prism and its self-products] \label{eg:prism-S2S3}
Consider the pair $(\mathfrak{S},(E))$ where $\mathfrak{S} = R_0[\![u]\!]$ and $E = E(u)$ is the Eisenstein polynomial for $\pi \in \mathcal{O}_K$ over $W$ as before. Equip $\mathfrak{S}$ with the $\delta$-structure defined by extending the fixed Frobenius $\varphi$ on $R_0$ to $\mathfrak{S}$ via $\varphi(u)=u^p$. Then $(\mathfrak{S},(E)) \in R_\Prism$, where the structure map $R \rightarrow \mathfrak{S}/(E)$ is given by the natural isomorphism $R \cong \mathfrak{S}/(E)$. We call $(\fkS,(E))$ the \emph{Breuil--Kisin prism} attached to $\pi$ and $R_0$.

The self-product of $(\mathfrak{S},(E))$ exists in $R_{\Prism}$ as follows. Consider the $p$-adically complete tensor-product $\mathfrak{S}\widehat{\otimes}_{\mathbf{Z}_p} \mathfrak{S}$ equipped with the induced $\otimes$-product Frobenius. We have a projection $d\colon \mathfrak{S}\widehat{\otimes}_{\mathbf{Z}_p} \mathfrak{S} \rightarrow R$ given by the composite of the multiplication $\mathfrak{S}\widehat{\otimes}_{\mathbf{Z}_p} \mathfrak{S} \rightarrow \mathfrak{S}$ and the natural projection $\mathfrak{S} \rightarrow \mathfrak{S}/(E)\cong R$. Let $J$ be the kernel of $d$, and let
\[
\mathfrak{S}^{(1)} \coloneqq (\mathfrak{S}\widehat{\otimes}_{\mathbf{Z}_p} \mathfrak{S})\biggl \{\frac{J}{E}\biggr\}_\delta^{\wedge}.
\]
Here $\mathfrak{S}\widehat{\otimes}_{\mathbf{Z}_p} \mathfrak{S}$ is regarded as an $\mathfrak{S}$-algebra via $a \mapsto a\otimes 1$, and $\{\cdot\}_\delta^{\wedge}$ means adjoining elements in the category of derived $(p, E)$-complete simplicial $\delta$-$\mathfrak{S}$-algebras. Note that the $E$ in $\bigl\{\frac{J}{E}\bigr\}_\delta^{\wedge}$ denotes $E\otimes 1$ but using $1\otimes E$ instead also gives the same $\mathfrak{S}^{(1)}$ (see \cite[Construction~7.13]{bhatt-scholze-prismaticFcrystal}). By \cite[Cor.~3.14]{bhatt-scholze-prismaticcohom-v3}, $(\mathfrak{S}^{(1)}, (E))$ is a $(p,E)$-completely flat prism over $(\mathfrak{S},(E))$. Furthermore, $(\mathfrak{S}^{(1)},(E))$ is bounded by \cite[Lem.~3.7 (2)]{bhatt-scholze-prismaticcohom}, so $(\mathfrak{S}^{(1)},(E)) \in R_{\Prism}$. Let $(B, I) \in R_{\Prism}$. If we are given maps $f_1, f_2\colon (\mathfrak{S}, (E)) \rightarrow (B, I)$ such that two maps $R \cong \mathfrak{S}/E \rightarrow B/I$ induced by $f_1$ and $f_2$ agree, then we have a natural induced map $f_1\otimes f_2\colon \mathfrak{S}\widehat{\otimes}_{\mathbf{Z}_p} \mathfrak{S} \rightarrow B$ of $\delta$-rings, and $(f_1\otimes f_2)(J) \subset I$. Thus, by the universal property of prismatic envelope (\cite[Cor.~3.14]{bhatt-scholze-prismaticcohom-v3}), we obtain a map $(\mathfrak{S}^{(1)},(E)) \rightarrow (B, I)$ in $R_{\Prism}$ uniquely determined by $f_1$, $f_2$. So $(\mathfrak{S}^{(1)},(E))$ is the self-product of $(\mathfrak{S},(E))$ in $R_{\Prism}$. Similarly, the self-triple-product $(\mathfrak{S}^{(2)},(E))$ of $(\mathfrak{S},(E))$ exists in $R_{\Prism}$. Write $p_1$, $p_2$ (resp.~$q_1$, $q_2$, $q_3$) for the maps from $(\mathfrak{S},(E))$ to $(\mathfrak{S}^{(1)},(E))$ (resp. to $(\mathfrak{S}^{(2)},(E))$).

A little more explicit description is given in \cite[\S~4.1]{du-liu-prismaticphiGhatmodule} as follows. Recall that $R_0$ is a $W\langle T_1^{\pm1},\ldots,T_d^{\pm1}\rangle$-algebra. Let $B^{\widehat{\otimes}[1]}$ denote the completion of $\fkS\otimes_{\mathbf{Z}_p}\fkS$ with respect to the ideal $\Ker(\fkS\otimes_{\mathbf{Z}_p}\fkS\rightarrow \fkS)$.
We have two natural maps $p_1,p_2\colon \fkS\rightarrow B^{\widehat{\otimes}[1]}$ and regard $B^{\widehat{\otimes}[1]}$ as an $\fkS$-algebra via $p_1$. If we set $s_j\coloneqq p_2(T_j)$ and $y\coloneqq p_2(u)$, then  
\[
\fkS^{\widehat{\otimes}[1]} \coloneqq \mathfrak{S} [\![y-u, s_1 - T_1, \dots , s_d - T_d]\!]
\]
can be naturally considered as an $\mathfrak{S}$-subalgebra of $B^{\widehat{\otimes}[1]}$. Similarly, let $B^{\widehat{\otimes}[2]}$ denote the completion of $\fkS\otimes_{\mathbf{Z}_p}\fkS\otimes_{\mathbf{Z}_p}\fkS$ with respect to the ideal $\Ker(\fkS\otimes_{\mathbf{Z}_p}\fkS\otimes_{\mathbf{Z}_p}\fkS\rightarrow \fkS)$.
We have maps $q_1,q_2,q_3\colon \fkS\rightarrow B^{\widehat{\otimes}[2]}$. Via $q_1$, we can naturally consider 
\[
\fkS^{\widehat{\otimes}[2]} \coloneqq \mathfrak{S} [\![y-u, w-u,\{s_j - T_j,r_j-T_j\}_{j=1,\ldots,d}]\!]
\]
as an $\fkS$-subalgebra of $B^{\widehat{\otimes}[2]}$, where $s_j\coloneqq q_2(T_i), r_j\coloneqq q_3(T_j), y\coloneqq q_2(u)$ and $w\coloneqq q_3(u)$.
Let 
\begin{align*}
J^{(1)}&= (E, y-u , \{s_j- T_j \}_{j = 1, \dots,  d} )\subset \mathfrak{S}^{\widehat\otimes [1]}\quad \text{and} \\
J^{(2)}&= (E, y-u , w-u, \{s_j-T_j, r_j - T_j\}_{j = 1, \dots, d}) \subset \mathfrak{S}^{\widehat\otimes [2]}.
\end{align*}
For $i=1,2$, $\fkS^{\widehat{\otimes}[i]}$ is naturally a $(p,E)$-completed $\delta$-$\fkS$-algebra and $\mathfrak{S}^{(i)}\cong \mathfrak{S}^{\widehat\otimes [i]}\bigl \{\frac{J ^{(i)}}{E}\bigr\}_\delta^{\wedge}$. 
\end{eg}

\begin{lem}\label{lem:AtoA2A3-faithful-flat}
The maps $p_i\colon \mathfrak{S} \rightarrow \mathfrak{S}^{(1)}$ \emph{(}resp. $q_i\colon \mathfrak{S} \rightarrow \mathfrak{S}^{(2)}$\emph{)} for $i = 1, 2$ \emph{(}resp. $i=1, 2, 3$\emph{)} are classically faithfully flat. 	
\end{lem}

\begin{proof}
We only consider $p_i\colon \mathfrak{S} \rightarrow \mathfrak{S}^{(1)}$. The proof for $q_i\colon \mathfrak{S} \rightarrow \mathfrak{S}^{(2)}$ is similar. Note that $\mathfrak{S}^{(1)}$ is classically $(p, E)$-complete by \cite[Lem.~3.7 (1)]{bhatt-scholze-prismaticcohom}, and $p_i\colon \mathfrak{S} \rightarrow \mathfrak{S}^{(1)}$ is $(p, E)$-completely flat. In particular, the induced map $\mathfrak{S}/(p, E)^n \rightarrow \mathfrak{S}^{(1)}/(p, E)^n$ is flat for each $n \geq 1$. Since $\mathfrak{S}$ is noetherian, $p_i\colon \mathfrak{S} \rightarrow \mathfrak{S}^{(1)}$ is classically flat by \cite[Tag 0912]{stacks-project}. Note that $p_i$ is a section of the  diagonal map $\mathfrak{S}^{(1)} \twoheadrightarrow \mathfrak{S}$. So if $N$ is any non-zero $\mathfrak{S}$-module, then $\mathfrak{S}^{(1)}\otimes_{p_i,\mathfrak{S}}N \neq 0$. Thus, $p_i$ is classically faithfully flat.	
\end{proof}

\begin{cor} \label{cor:A2-properties}
$\mathfrak{S}^{(1)}$ is $p$-torsion free and $E$-torsion free. Furthermore, 
\[
\mathfrak{S}^{(1)}[p^{-1}] \cap \mathfrak{S}^{(1)}[E^{-1}] = \mathfrak{S}^{(1)},
\]
and $\mathfrak{S}^{(1)}[E^{-1}]$ is $p$-adically separated.
\end{cor}

\begin{proof}
Since $\mathfrak{S}$ is torsion free and $p_1\colon \mathfrak{S} \rightarrow \mathfrak{S}^{(1)}$ is classically flat, $\mathfrak{S}^{(1)}$ is $p$-torsion free and $E$-torsion free. We deduce by Lemma~\ref{lem:intersection-modules-flat-base-change} (i) and $\mathfrak{S}[p^{-1}] \cap \mathfrak{S}[E^{-1}] = \mathfrak{S}$ that
\[
\mathfrak{S}^{(1)}[p^{-1}] \cap \mathfrak{S}^{(1)}[E^{-1}] = (\mathfrak{S}^{(1)}\otimes_{p_1, \mathfrak{S}} \mathfrak{S}[p^{-1}]) \cap (\mathfrak{S}^{(1)}\otimes_{p_1, \mathfrak{S}}\mathfrak{S}[E^{-1}]) = \mathfrak{S}^{(1)}.
\]
Since $\mathfrak{S}^{(1)}$ is $p$-adically complete, this also implies that $\mathfrak{S}^{(1)}[E^{-1}]$ is $p$-adically separated.
\end{proof}

\begin{eg}[The $\A_{\mathrm{inf}}$-prism] \label{eg:prism-Ainf}
Let $(\xi)$ be the kernel of $\theta\colon \A_{\mathrm{inf}}(\overline{R}) \rightarrow \overline{R}^{\wedge}$. Then $(\A_{\mathrm{inf}}(\overline{R}), (\xi)) \in R_\Prism$, with the structure map $R \rightarrow \A_{\mathrm{inf}}(\overline{R}) / (\xi)$ given by the natural inclusion $R \rightarrow \overline{R}^{\wedge}$. Note that the map $f_{\pi^\flat,T_i^\flat}\colon\mathfrak{S} \rightarrow \A_{\mathrm{inf}}(\overline{R})$ given by $u \mapsto [\pi^\flat]$ and $T_i \mapsto [T_i^\flat]$ induces a map of prisms $(\mathfrak{S}, (E)) \rightarrow (\A_{\mathrm{inf}}(\overline{R}), (\xi))$ over $R$.
Moreover, each $\sigma\in\calG_R$ induces a map of prisms $(\A_{\mathrm{inf}}(\overline{R}), (\xi))\rightarrow (\A_{\mathrm{inf}}(\overline{R}), (\xi))$ satisfying $\sigma\circ f_{\pi^\flat,T_i^\flat}=f_{\sigma(\pi^\flat),\sigma(T_i^\flat)}$.
\end{eg}

\begin{eg}[The $\OA_\cris$-prism and its Frobenius twists] \label{eg:prism-OAcris}
Consider the surjective map $\theta_{R_0}\colon \OA_{\cris}(\overline{R}) \rightarrow \overline{R}^\wedge$. The map 
\[
\varphi\colon \OA_{\cris}(\overline{R})/(p) \rightarrow \OA_{\cris}(\overline{R})/(p)
\]
factors through 
\[
\OA_{\mathrm{cris}}(\overline{R})/(p) \rightarrow \OA_{\mathrm{cris}}(\overline{R})/((p)+\mathrm{ker}(\theta_{R_0})) \cong \overline{R}/(p) \xrightarrow{h} \OA_{\mathrm{cris}}(\overline{R})/(p).
\]  
The pair $(\OA_{\mathrm{cris}}(\overline{R}), (p))$ defines a prism in $R_\Prism$, where the structure map $R \rightarrow \OA_{\mathrm{cris}}(\overline{R})/(p)$ is given by the composite of $R\rightarrow \overline{R}/(p)$ and $h$ (defined in above factorization). Consider the composite $\Ainf(\overline{R}) \xrightarrow{\varphi} \Ainf(\overline{R}) \rightarrow \OA_{\mathrm{cris}}(\overline{R})$, where the second map is the natural inclusion. This induces a map of prisms $(\Ainf(\overline{R}), (\xi)) \xrightarrow{\varphi} (\OA_{\mathrm{cris}}(\overline{R}), (p))$ over $R$, which is compatible with Frobenii and $\mathcal{G}_R$-actions.

On the other hand, consider the prism $(R_0, (p))$ in $R_{\Prism}$, where the structure map $R \rightarrow R_0/(p)$ is given by the natural projection $R \rightarrow R/(\pi) \cong R_0/(p)$. For any integer $j \geq 1$, write $(\phi_j \OA_{\mathrm{cris}}(\overline{R}), (p))$ for the prism in $R_{\Prism}$ whose underlying $\delta$-pair is $(\OA_{\mathrm{cris}}(\overline{R}),(p))$ and the structure map $R \rightarrow \OA_{\mathrm{cris}}(\overline{R})/(p)$ is given the structure map of $(\OA_{\mathrm{cris}}(\overline{R}),(p))$ composed with $\varphi^j\colon \OA_{\mathrm{cris}}(\overline{R})\rightarrow\OA_{\mathrm{cris}}(\overline{R})$. For a sufficiently large $j$, there exists a map of prisms $(R_0, (p)) \rightarrow (\phi_j\OA_{\mathrm{cris}}(\overline{R}), (p))$ given by Dwork's trick. Indeed, let $e = [K: K_0]$ be the ramification index, and choose an integer $l$ such that $p^l \geq e$. Consider $\varphi^{l+1}\colon \OA_{\mathrm{cris}}(\overline{R}) \rightarrow \OA_{\mathrm{cris}}(\overline{R})$. Taking modulo the ideal $(p)\subset \OA_{\mathrm{cris}}(\overline{R})$, this induces a map $\overline{R}/(p) \rightarrow \OA_{\mathrm{cris}}(\overline{R})/(p)$ as above, which further factors through $\overline{R}/(\pi)$. Thus, the ring map $R_0 \xrightarrow{\varphi^{l+1}} \OA_{\mathrm{cris}}(\overline{R})$ induces a map $(R_0, (p)) \rightarrow (\phi_{l}\OA_{\mathrm{cris}}(\overline{R}), (p))$ of prisms over $R$.
\end{eg}

\begin{eg}[The Breuil prism] \label{eg:prism-PD-S}
Let $S$ denote the $p$-adically completed PD-envelope of $\mathfrak{S}$ with respect to $(E)$, equipped with the Frobenius extending $\varphi$ on $\mathfrak{S}$. Note that $c \coloneqq \frac{\varphi(E)}{p}$ is a unit in $S$. So $(S, (p)) \in R_{\Prism}$ with the structure map given by $R \cong \mathfrak{S}/(E) \xrightarrow{\varphi} S/(p)$, and we have a map of prisms $\varphi\colon (\mathfrak{S}, (E)) \rightarrow (S, (p))$.
\end{eg}

For $i = 1, 2$, let $S^{(i)} \coloneqq D_{\mathfrak{S}^{\widehat{\otimes}[i]}} (J^{(i)})^{\wedge}$ be the $p$-adically completed PD-envelope of $\mathfrak{S}^{\widehat{\otimes}[i]}$ with respect to the ideal $J^{(i)}$. We set
\[
z_0 \coloneqq \frac{y-u}{E}\quad\text{and}\quad
z_j \coloneqq \frac{s_j-T_j}{E} \quad(j = 1, \ldots, d).
\]
Let $A_{\mathrm{max}}^{(1)}$ be the $p$-adic completion of the $\mathfrak{S}$-subalgebra of $(\mathfrak{S}[p^{-1}])[z_0, z_1, \ldots, z_d]$ generated by $\frac{E}{p}$ and $\{\gamma_n(z_j)\}_{n \geq 1, 0\leq j\leq d}$. By \cite[\S 2.2]{du-liu-prismaticphiGhatmodule}\footnote{We warn the reader that our ring $A_{\mathrm{max}}^{(1)}$ is denoted by $A_{\mathrm{max}}^{(2)}$ in \textit{ibid.}. }, we have a ring endomorphism $\varphi\colon A_{\mathrm{max}}^{(1)} \rightarrow A_{\mathrm{max}}^{(1)}$ extending $\varphi\colon \mathfrak{S} \rightarrow \mathfrak{S}$ and satisfying
\[
\varphi(z_0) = \frac{y^p-u^p}{\varphi(E)} \quad\text{and}\quad
\varphi(z_j) = \frac{s_j^p-T_j^p}{\varphi(E)}
\quad(j = 1, \ldots, d).
\]
Since $\varphi\colon \mathfrak{S} \rightarrow \mathfrak{S}$ is injective, $\varphi\colon A_{\mathrm{max}}^{(1)} \rightarrow A_{\mathrm{max}}^{(1)}$ is injective. By \textit{ibid.}, we have a natural ring map $\mathfrak{S}^{(1)} \rightarrow A_{\mathrm{max}}^{(1)}$ which is injective and compatible with $\varphi$. In fact, $A_{\mathrm{max}}^{(1)}$ is isomorphic to $\mathfrak{S}^{(1)}\langle \frac{E}{p} \rangle$, the $p$-adic completion of $\mathfrak{S}^{(1)}[ \frac{E}{p} ]$ by \cite[Rem.~2.2.11]{du-liu-prismaticphiGhatmodule}.

Let $S_1 \coloneqq \mathfrak{S}^{\widehat{\otimes}[1]} [\gamma_n(E), \gamma_n(y-u), \{\gamma_n (s_j-T_j)\}_{n \geq 1, ~j = 0, \ldots, d}] \subset \mathfrak{S}^{\widehat{\otimes}[1]} [p^{-1}]$. Note that $S_1 \subset A_{\mathrm{max}}^{(1)}$ since $\gamma_n(y-u) =\gamma_n(z_0)E^n \in A_{\mathrm{max}}^{(1)}$ and similarly for $\gamma_n(s_i-T_i)$. Since $E$, $y-u$, and $\{s_j-T_j\}_{j = 1, \ldots, d}$ form a regular sequence in $\mathfrak{S}^{\widehat{\otimes}[1]}$, $S_1$ is the PD-envelope of $\mathfrak{S}^{\widehat{\otimes}[1]}$ for $J^{(1)}$ by \cite[Cor.~2.39]{bhatt-scholze-prismaticcohom}. Then $S^{(1)}$ is the $p$-adic completion of $S_1$. As a subring of $(R_0[p^{-1}])[\![u, y-u, s_1-T_1, \ldots, s_d-T_d]\!]$, we have
\begin{align*}
S^{(1)} =  \bigl\{ \sum a_{i_0, \ldots, i_{d+1}} &\gamma_{i_0}(E) \gamma_{i_1}(y-u)\gamma_{i_2}(s_1-T_1)\cdots \gamma_{i_{d+1}} (s_d-T_d) ~|\\
&a_{i_0, \dots , i_{d+1}} \in \mathfrak{S}^{\widehat{\otimes}[1]}, ~a_{i_0, \dots , i_{d+1}} \rightarrow 0 \quad(\text{as } i_0+\cdots+i_{d+1}\to\infty)
\bigr\}
\end{align*}
where the sum goes over the multi-index $(i_0, \ldots, i_{d+1})$ of non-negative integers and $a_{i_0, \dots , i_{d+1}} \rightarrow 0$ means in the $p$-adic topology. Note that $S^{(1)}$ is a $\delta$-ring by \cite[Cor.~2.39]{bhatt-scholze-prismaticcohom}. We similarly construct a $\delta$-ring $S^{(2)}$. 

\begin{lem} \label{lem:A2-to-S2-embedding}
For $i = 1, 2$, we have an embedding $\mathfrak{S}^{(i)} \xhookrightarrow{\varphi} S^{(i)}$.
\end{lem}

\begin{proof}
We prove the statement for $\varphi\colon \mathfrak{S}^{(1)} \rightarrow S^{(1)}$, and the proof for $\varphi\colon \mathfrak{S}^{(2)} \rightarrow S^{(2)}$ is analogous. It suffices to show $\varphi(\delta^n(z_j)) \in S^{(1)}$ for $j = 0, \ldots, d$. Since $S^{(1)}$ is a $\delta$-ring, we have
\[
\varphi(z_0) = c^{-1}\frac{\varphi(y-u)}{p} = c^{-1}(\frac{(y-u)^p}{p}+\delta(y-u)) \in S^{(1)},
\]
and similarly $\varphi(z_j) \in S^{(1)}$ for $j = 1, \ldots, d$. Again since $S^{(1)}$ is a $\delta$-ring, we have $\varphi(\delta^n(z_j)) = \delta^n(\varphi(z_j)) \in S^{(1)}$ for any $n \geq 0$. 
\end{proof}

\subsection{Completed prismatic \texorpdfstring{$F$}{F}-crystals in the small affine case} \label{sec:completed-prism-F-crystals}

In this subsection, we introduce completed prismatic crystals and completed prismatic $F$-crystals on the absolute prismatic site.

We first introduce the notion of finitely generated completed prismatic crystals. Let $\fkX$ be a smooth $p$-adic formal scheme over $\calO_K$ (or a CDVR of mixed characteristic $(0,p)$).

\begin{defn}\label{defn:completed-crystal:global case}
A \textit{finitely generated completed crystal of $\mathcal{O}_\Prism$-modules} on $\fkX_\Prism$ is a sheaf $\mathcal{F}$ of $\mathcal{O}_\Prism$-modules on $\fkX_\Prism$
such that
\begin{enumerate}
    \item for each $(A,I)\in \fkX_{\Prism}$, the evaluation $\mathcal{F}_A\coloneqq \mathcal{F}(A,I)$ of $\mathcal{F}$ on $(A,I)$ is a finitely generated and classically $(p,I)$-complete $A$-module;
    \item for any morphism $(A,I)\rightarrow (B,IB)$ of bounded prisms over $\fkX$, the canonical linearized transition map
    \[
    B\widehat{\otimes}_A \mathcal{F}_A\rightarrow \mathcal{F}_B
    \]
    is an isomorphism, where $B\widehat{\otimes}_A \mathcal{F}_A$ denotes the completed tensor product $\varprojlim_n (B\otimes_A\calF_A)/(p,I)^n(B\otimes_A\calF_A)$.  
%the $(p,I)$-adically completed tensor product 
\end{enumerate}
We also call such a sheaf a \textit{finitely generated completed prismatic crystal} on $\fkX$, or a \emph{completed prismatic crystal} on $\fkX$ for short.

Similarly, a \textit{finitely generated crystal of $\mathcal{O}_{\Prism,n}$-modules} on $\fkX_\Prism$ is a sheaf $\mathcal{F}_n$ of $\mathcal{O}_{\Prism,n}$-modules on $\fkX_\Prism$
such that
\begin{enumerate}
    \item for each $(A,I)\in \fkX_{\Prism}$, the evaluation $\mathcal{F}_{n,A}\coloneqq \mathcal{F}_n(A,I)$ of $\mathcal{F}_n$ on $(A,I)$ is a finitely generated $A/(p,I)^n$-module;
    \item for any morphism $(A,I)\rightarrow (B,IB)$ of bounded prisms over $\fkX$, the canonical linearized transition map $B\otimes_A\mathcal{F}_{n,A}\rightarrow \mathcal{F}_{n,B}$ is an isomorphism.
\end{enumerate}
\end{defn}

\begin{rem}\label{rem:conditions for uncompleted base change}
Let $\calF$ be a finitely generated completed prismatic crystal on $\fkX$ and let $(A,I)\rightarrow (B,IB)$ be a map of bounded prisms over $\fkX$.
Since $\calF_A$ is a finitely generated $A$-module and $B$ is classically $(p,IB)$-complete, the natural map
\begin{equation}\label{eq:rem:conditions for uncompleted base change}
B\otimes_A\calF_A\rightarrow B\widehat{\otimes}_A\calF_A\xrightarrow{\cong}\calF_B    
\end{equation}
is surjective. Since $(p,IB)$ is a finitely generated ideal of $B$, the map \eqref{eq:rem:conditions for uncompleted base change} induces an isomorphism $B/(p,IB)^n\otimes_A\calF_A\xrightarrow{\cong}\calF_B/(p,IB)^n\calF_B$ by \cite[Thm.~1.2 (2)]{yekutieli-flatnesscompletion}.
Moreover, the map \eqref{eq:rem:conditions for uncompleted base change} is an isomorphism if $A$ and $B$ are both noetherian or if $A$ is noetherian and the map $A\rightarrow B$ is classically flat. The latter case follows from \cite[Tag 0912]{stacks-project}.
It is also an isomorphism if $\calF_A$ is a finite projective $A$-module.
\end{rem}

\begin{lem}\label{lem:calF vs (calF_n) for completed crystals}
Let $\calF$ be a finitely generated completed crystal of $\calO_\Prism$-modules on $\fkX_\Prism$. Then for each $n\geq 1$, the association $(A,I)\mapsto \calF_A/(p,I)^n\calF_A$ represents the quotient sheaf $\calF/(p,\calI_\Prism)^n\calF$ and defines a finitely generated crystal $\calF_n$ of $\calO_{\Prism,n}$-modules on $\fkX_\Prism$. Moreover, we have isomorphisms of $\calO_\Prism$-modules $\calO_{\Prism,n}\otimes_{\calO_{\Prism,n+1}}\calF_{n+1}\xrightarrow{\cong}\calF_n$ and $\calF\cong \varprojlim_n \calF_n\cong \operatorname{Rlim}\calF_n$.

Conversely, let $(\calF_n)_n$ be an inverse system of sheaves of $\calO_\Prism$-modules such that $\calF_n$ is a finitely generated crystal of $\calO_{\Prism,n}$-modules and such that the projection $\calF_{n+1}\rightarrow\calF_n$ induces an isomorphism $\calO_{\Prism,n}\otimes_{\calO_{\Prism,n+1}}\calF_{n+1}\xrightarrow{\cong}\calF_n$ for each $n$. Then $\calF\coloneqq \varprojlim_n\calF_n$ is a finitely generated completed crystal of $\calO_\Prism$-modules on $\fkX_\Prism$, and we have isomorphisms of $\calO_\Prism$-modules $\calO_{\Prism,n}\otimes_{\calO_\Prism}\calF \cong \calF_n$ and $\calF\cong \operatorname{Rlim}\calF_n$.
\end{lem}

\begin{proof}
Let $\calF$ be a finitely generated completed crystal of $\calO_\Prism$-modules on $\fkX_\Prism$.
Let $(A,I)\rightarrow (B,IB)$ be a $(p,I)$-completely faithfully flat map of bounded prisms over $\fkX$.
Set $B'=(B\otimes_AB)^\wedge_{(p,I)}$. By Lemma~\ref{lem:pushout for prisms along fflat map}, $B'$ is $(p,I)$-completely faithfully flat over $A$ and
$(B',IB')\in \fkX_\Prism$ is the self-fiber product of $(B,IB)$ over $(A,I)$.
Let $p_1,p_2\colon (B,IB)\rightarrow (B',IB')$ be the two maps of bounded prisms over $\fkX$.
Since $B/(p,I)^nB$ is classically faithfully flat over $A/(p,I)^n$ and since  $B/(p,I)^nB\otimes_{A/(p,I)^n}B/(p,I)^nB\cong B'/(p,I)^nB'$, we have an exact sequence
\[
0\rightarrow \calF_{n,A}\rightarrow B/(p,I)^nB\otimes_{A/(p,I)^n}\calF_{n,A} \xrightarrow{p_1\otimes 1-p_2\otimes 1} B'/(p,I)^nB'\otimes_{A/(p,I)^n}\calF_{n,A}.
\]
On the other hand, since $\calF$ is a completed prismatic crystal,  we have the isomorphisms
$B\widehat{\otimes}_A\calF_A\cong \calF_B$ and $B'\widehat{\otimes}_A\calF_A\cong \calF_{B'}$. It follows that the above exact sequence is identified with
\[
0\rightarrow \calF_{n,A}\rightarrow \calF_{n,B}\xrightarrow{p_1^\ast-p_2^\ast} \calF_{n,B'}.
\]
This implies that $\calF_n$ is a sheaf on $\fkX_\Prism$, representing the quotient sheaf $\calF/(p,\calI_\Prism)^n\calF$.
Since $\calF_A$ is a finitely generated classically $(p,I)$-complete $A$-module, $\calF_n$ is a finitely generated crystal of $\calO_{\Prism,n}$-modules. Moreover, we have isomorphisms of $\calO_\Prism$-modules $\calO_{\Prism,n}\otimes_{\calO_{\Prism,n+1}}\calF_{n+1}\xrightarrow{\cong}\calF_n$ and $\calF\xrightarrow{\cong} \varprojlim_n \calF_n$. Finally, since the absolute prismatic topos is replete and $\calF_{n+1}\rightarrow\calF_n$ is surjective for every $n$, we obtain $\varprojlim_n\calF_n\cong \operatorname{Rlim}\calF_n$ by \cite[Prop.~3.1.10]{bhatt-scholze-proetale}.

Conversely, let $(\calF_n)_n$ be an inverse system of sheaves of $\calO_\Prism$-modules satisfying the properties as in the lemma. An argument similar to the previous paragraph shows that the association $(A,I)\mapsto A/(p,I)^n\otimes_{A/(p,I)^{n+1}}\calF_{n+1,A}$ represents the sheaf $\calO_{\Prism,n}\otimes_{\calO_{\Prism,n+1}}\calF_{n+1}$. Set $\calF\coloneqq \varprojlim_n\calF_n$ and take any $(A,I)\in \fkX_\Prism$. Then we have $\calF(A,I)=(\varprojlim_n\calF_n)(A,I)=\varprojlim_n \calF_{n,A}$ and $A/(p,I)^n\otimes_{A/(p,I)^{n+1}}\calF_{n+1,A}\cong \calF_{n,A}$.  
It follows from \cite[Thm~2.8]{yekutieli-flatnesscompletion} that $\calF(A,I)$ is a finitely generated and classically $(p,I)$-complete $A$-module with $\calF(A,I)/(p,I)^n\calF(A,I)\cong \calF_{A,n}$.
Moreover, for a morphism $(A,I)\rightarrow (B,IB)$ of bounded prisms over $\fkX$, we have $B\otimes_A\calF_{n,A}\cong \calF_{n,B}$. It follows that the canonical map $B\widehat{\otimes}_A\calF(A,I)\rightarrow \calF(B,IB)$ is an isomorphism. Hence $\calF$ is a finitely generated completed crystal of $\calO_\Prism$-modules. Now the remaining assertions follow easily.
\end{proof}

To define completed prismatic $F$-crystals in the affine case, let us first introduce the following terminologies.

\begin{defn} \label{defn:proj-away-(p, E)-saturated}
Let $R$ be a base ring and keep the notation as in \S~\ref{sec-basering}.
\begin{enumerate}
\item We say that a finite $\mathfrak{S}$-module $N$ is \textit{projective away from} $(p, E)$ if $N$ is torsion free, $N[p^{-1}]$ is projective over $\mathfrak{S}[p^{-1}]$, and $N[E^{-1}]^{\wedge}_p$ is projective over $\mathfrak{S}[E^{-1}]^{\wedge}_p = \mathcal{O}_{\mathcal{E}}$.

\item We say that a finite $\mathfrak{S}$-module $N$ is \textit{saturated} if $N$ is torsion free and 
\[
N = N[p^{-1}] \cap N[E^{-1}].
\]
\item Let $r$ be a non-negative integer and let $N$ be an $\fkS$-module equipped with a $\varphi$-semi-linear endomorphism $\varphi_{N}\colon N \rightarrow N$. We say that the pair $(N,\varphi_N)$ has \emph{$E$-height $\leq r$} if 
\[
1\otimes\varphi_{N} \colon \mathfrak{S}\otimes_{\varphi, \mathfrak{S}} N \rightarrow N
\] 	
is injective and its cokernel is killed by $E(u)^r$.
We say that $(N,\varphi_N)$ has finite $E$-height if it has $E$-height $\leq r$ for some $r$.
\end{enumerate}
We will also use these terminologies for a finite module over an $\fkS$-algebra.
\end{defn}

\begin{rem} \label{rem:saturated-regular-seq}
\hfill
\begin{enumerate}
 \item By the Beauville--Laszlo theorem, any finitely generated $\mathfrak{S}$-module which is saturated and projective away from $(p, E)$ is the pushforward of a vector bundle on $\operatorname{Spec}\mathfrak{S}\setminus V(p, E)$ to $\operatorname{Spec}\mathfrak{S}$. In fact, let $N$ be a torsion free finite $\mathfrak{S}$-module.  Then $N$ is saturated if and only if the natural map
\[
N/pN \rightarrow  N[E^{-1}]/pN[E^{-1}]
\]
is injective. Since $N[E^{-1}]/pN[E^{-1}]\cong N[E^{-1}]^{\wedge}_p/pN[E^{-1}]^{\wedge}_p \cong N[u^{-1}]/pN[u^{-1}]$, we deduce that $N$ is saturated if and only if $N = N[p^{-1}] \cap N[E^{-1}]^\wedge_p$, or equivalently, $N=N[p^{-1}]\cap N[u^{-1}]$. Moreover, if $N$ is projective away from $(p, E)$, then $N[E^{-1}]$ is finite projective over $\mathfrak{S}[E^{-1}]$.
 \item Assume that either $R$ is small over $\calO_K$ or $R=\calO_L$.
We will show that if $(N,\varphi_N)$ is a torsion free finite $\mathfrak{S}$-module with Frobenius of finite $E$-height, then $N[p^{-1}]$ is projective over $\fkS[p^{-1}]$ (Proposition~\ref{prop:rational-projectivity-etale-over-torus-case}).

\end{enumerate}
\end{rem}

 We now introduce the notion of completed prismatic $F$-crystals on $R$, which will be our main object of study. We make Assumption~\ref{assumption:base-ring-sec-3.4}: $R$ is small over $\calO_K$ or $R=\calO_L$.

\begin{defn} \label{defn:category-good-prism-completed-F-crystal}
A \emph{completed $F$-crystal of $\calO_\Prism$-modules} on $R_\Prism$ is a pair
 $(\mathcal{F}, \varphi_{\mathcal{F}})$, where $\mathcal{F}$ is a finitely generated completed crystal of $\mathcal{O}_{\Prism}$-modules on $R_\Prism$ and 
\[
\varphi_{\mathcal{F}}\colon \mathcal{F} \rightarrow \mathcal{F}
\]
is a $\varphi$-semilinear morphism of $\mathcal{O}_{\Prism}$-modules such that 
\begin{enumerate}
\item $\mathcal{F}_{\mathfrak{S}} \coloneqq \mathcal{F}(\mathfrak{S}, E)$ is projective away from $(p, E)$ and saturated;
\item the pair $(\calF_\fkS,\varphi_{\calF_\fkS})$ has finite $E$-height.
\end{enumerate}
We also call such an object a \textit{completed prismatic $F$-crystal} on $R$.
The morphisms between completed $F$-crystals of $\calO_\Prism$-modules are $\mathcal{O}_{\Prism}$-module maps compatible with Frobenii. 

We write $\mathrm{CR}^{\wedge,\varphi}(R_{\Prism})$ for the category of completed $F$-crystals of $\calO_\Prism$-modules on $R_\Prism$.
Let $\mathrm{Vect}^{\varphi}_{\mathrm{eff}}(R_{\Prism})$ denote the full subcategory of $\mathrm{CR}^{\wedge,\varphi}(R_{\Prism})$ consisting of objects $(\calF,\varphi_{\calF})$ where $\calF$ is a locally free $\calO_\Prism$-module. 
For a fixed non-negative integer $r$, we let $\mathrm{CR}^{\wedge,\varphi}_{[0,r]}(R_{\Prism})$ and $\mathrm{Vect}^{\varphi}_{[0,r]}(R_{\Prism})$ denote the full subcategories consisting of objects for which $(\calF_\fkS,\varphi_{\calF_\fkS})$ has $E$-height $\leq r$.
\end{defn}

\begin{rem}\label{rem:definition of completed prismatic F-crystals}
When $R$ is small over $\calO_K$, the above definition agrees with Definition~\ref{defn:category-good-prism-completed-F-crystal-intro} by Remark~\ref{rem:saturated-regular-seq} (i) and (ii). In \S~\ref{sec:globalization}, we will define completed prismatic $F$-crystals on a smooth $p$-adic formal scheme by gluing.
\end{rem}

\begin{rem} \label{rem:vector-bundle-CDVF-case}
When $R = \mathcal{O}_L\coloneqq R_{(\pi)}^\wedge$ (i.e., a CDVR with residue field having a finite $p$-basis and a uniformizer finite over $W(k)$), any finite $\mathfrak{S}_L$-module which is projective away from $(p, E)$ and saturated is free over $\mathfrak{S}_L$, since $\mathfrak{S}_L$ is a regular local ring of dimension 2 (e.g. \cite[Cor.~4.1.1]{horrocks-vectorbundles}). Thus, by Proposition~\ref{prop:equivalence-to-descent-datum} below, the category $\mathrm{CR}^{\wedge,\varphi}((\mathcal{O}_L)_{\Prism})$ is equal to the category $\mathrm{Vect}^{\varphi}_{\mathrm{eff}}((\mathcal{O}_L)_{\Prism})$. Furthermore, when $R = \mathcal{O}_K$ (i.e., a CDVR with perfect residue field), our category $\mathrm{Vect}^{\varphi}_{\mathrm{eff}}((\mathcal{O}_K)_{\Prism})$ coincides with the full subcategory of $\mathrm{Vect}^{\varphi}(\mathrm{Spf}(\mathcal{O}_K)_{\Prism}, \mathcal{O}_{\Prism})$ defined in \cite[Def.~4.1]{bhatt-scholze-prismaticFcrystal} consisting of \emph{effective} prismatic $F$-crystals of vector bundles. 
\end{rem}
Let us explain that the definition of completed prismatic $F$-crystals is independent of the choice of a Breuil--Kisin prism, namely, a uniformizer $\pi\in\calO_K$ and a $W$-subalgebra $R_0\subset R$ (Corollary~\ref{cor:independence of Breuil--Kisin condition}). For this, we need the following two lemmas.

\begin{lem} \label{lem:independence of Breuil--Kisin condition}
Let $\mathcal{F}$ be a finitely generated completed prismatic crystal on $R$ equipped with a morphism $1\otimes\varphi_{\mathcal{F}}\colon \varphi^*\mathcal{F} \rightarrow \mathcal{F}$ of $\mathcal{O}_{\Prism}$-modules. 
Fix a uniformizer $\pi\in\calO_K$ with minimal polynomial $E(u)$ and associated Breuil--Kisin prism $(\fkS,(E))$, and let $(\fkS,(E))\rightarrow (B,EB)$ be a classically flat map of bounded prisms over $R$.
Then the following properties hold:
\begin{enumerate}
    \item if $\mathcal{F}_{\mathfrak{S}}$ is  projective away from $(p, E)$ and saturated as an $\mathfrak{S}$-module, then $\mathcal{F}_B$ is projective away from $(p, E)$ and saturated as a $B$-module;
    \item for a non-negative integer $r$, if the pair $(\calF_\fkS,\varphi_{\calF_\fkS})$ has $E$-height $\leq r$, then $(\calF_B,\varphi_{\calF_B})$ has $E$-height $\leq r$.
\end{enumerate}
Moreover, the converse also holds if $\fkS\rightarrow B$ is classically faithfully flat.
\end{lem}

\begin{proof}
Note $B\otimes_\fkS \calF_\fkS\cong \calF_B$ by Remark~\ref{rem:conditions for uncompleted base change}.

(i) Suppose $\mathcal{F}_{\mathfrak{S}}$ is projective away from $(p, E)$ and saturated as an $\mathfrak{S}$-module. 
Then $\calF_B$ is $p$-torsion free, and $\calF_B[p^{-1}]$ is projective over $B[p^{-1}]$. It follows that $\calF_B\subset \calF_B[p^{-1}]$ is torsion free. Since $\fkS$ is noetherian and $\calF_\fkS$ is finitely generated, the induced map $\mathfrak{S}[E^{-1}]^{\wedge}_p \rightarrow B[E^{-1}]^{\wedge}_p$ is classically flat and $B[E^{-1}]^{\wedge}_p\otimes_{\mathfrak{S}[E^{-1}]^{\wedge}_p}\mathcal{F}_\fkS[E^{-1}]^{\wedge}_p\cong\mathcal{F}_B[E^{-1}]^{\wedge}_p$ by \cite[Tag 0912]{stacks-project}. We deduce that $\mathcal{F}_B[E^{-1}]^{\wedge}_p$ is projective over $B[E^{-1}]^{\wedge}_p$. 
Thus $\mathcal{F}_B$ is projective away from $(p, E)$.
Since $\mathcal{F}_{\mathfrak{S}}$ is saturated, Lemma~\ref{lem:intersection-modules-flat-base-change} (i) implies that
\[
\calF_B=B\otimes_\fkS \calF_\fkS=B\otimes_\fkS (\calF_\fkS[p^{-1}]\cap \calF_\fkS[E^{-1}])=\calF_B[p^{-1}]\cap\calF_B[E^{-1}].
\]
This means that $\calF_B$ is saturated.

(ii) The assertion follows from $\operatorname{Coker}(1\otimes \varphi_{\calF_\fkS})\otimes_\fkS B=\operatorname{Coker}(1\otimes \varphi_{\calF_B})$.

Finally, if the map $\fkS\rightarrow B$ is classically faithfully flat, then so is $\mathfrak{S}[E^{-1}]^{\wedge}_p \rightarrow B[E^{-1}]^{\wedge}_p$. Hence the converse direction follows similarly.
\end{proof}

Suppose $R$ is small over $\mathcal{O}_K$. Let $\pi' \in \mathcal{O}_K$ be another uniformizer of $\calO_K$, $E'(y) \in W[y]$ the Eisenstein polynomial for $\pi'$, and $R_0'$  a $W\langle (T_1')^{\pm 1}, \ldots, (T_d')^{\pm 1}\rangle$-algebra with $R_0'\otimes_W\calO_K=R$ as in Remark \ref{rem:p-complete-etale}. Set $\fkS'\coloneqq R_0'[\![y]\!]$ equipped with Frobenius given by $\varphi(T_i') = (T_i')^p$ and $\varphi(y) = y^p$. Then we have a Breuil--Kisin prism $(\fkS',(E'))\in R_\Prism$ with the structure map $R\xrightarrow{\cong} \fkS'/(E')$.

\begin{lem}\label{lem:productoffkSandfkSprime}
\hfill
\begin{enumerate}
    \item The absolute product of $(\mathfrak{S}, (E))$ and $(\mathfrak{S}', (E'))$ exists in $R_{\Prism}$. Write $(\mathfrak{S}_{\pi, \pi'}^{(1)}, I)$ for the absolute product. We also have $I=E\mathfrak{S}_{\pi, \pi'}^{(1)}=E'\mathfrak{S}_{\pi, \pi'}^{(1)}$.
    \item The maps $\mathfrak{S} \rightarrow \mathfrak{S}_{\pi, \pi'}^{(1)}$ and $\mathfrak{S}' \rightarrow \mathfrak{S}_{\pi, \pi'}^{(1)}$ are classically faithfully flat.
\end{enumerate}
\end{lem}

\begin{proof}
(i) Consider the $p$-adically complete tensor-product $\mathfrak{S}\widehat{\otimes}_{\mathbf{Z}_p} \mathfrak{S}'$, and let
\[
d\colon \mathfrak{S}\widehat{\otimes}_{\mathbf{Z}_p} \mathfrak{S}' \rightarrow R
\]
be the composite of the natural projection $\mathfrak{S}\widehat{\otimes}_{\mathbf{Z}_p} \mathfrak{S}' \rightarrow \mathfrak{S}/(E)\widehat{\otimes}_{\mathbf{Z}_p} \mathfrak{S}'/(E') \cong R\widehat{\otimes}_{\mathbf{Z}_p} R$ and the multiplication $R\widehat{\otimes}_{\mathbf{Z}_p} R \rightarrow R$. Let $J$ be the kernel of $d$. We claim that the absolute product of $(\mathfrak{S}, (E))$ and $(\mathfrak{S}', (E'))$ in $R_{\Prism}$ is given by 
\[
\mathfrak{S}_{\pi, \pi'}^{(1)} = \mathfrak{S}\widehat{\otimes}_{\mathbf{Z}_p} \mathfrak{S}'\biggl \{\frac{J}{E}\biggr\}_\delta^{\wedge},
\]
where $\{\cdot\}_\delta^{\wedge}$ means adjoining elements in the category of derived $(p, E)$-complete simplicial $\delta$-$\mathfrak{S}$-algebras. Indeed, by \cite[Cor.~3.14]{bhatt-scholze-prismaticcohom-v3}, $(\mathfrak{S}_{\pi, \pi'}^{(1)}, (E))$ is a $(p,E)$-completely flat prism over $(\mathfrak{S}, (E))$. We have a natural map of prisms $(\mathfrak{S}', (E')) \rightarrow (\mathfrak{S}_{\pi, \pi'}^{(1)}, (E))$, and thus $(E')\mathfrak{S}_{\pi, \pi'}^{(1)} = (E)\mathfrak{S}_{\pi, \pi'}^{(1)}$ by \cite[Lem.~3.5]{bhatt-scholze-prismaticcohom}. So the construction is symmetric, and $(\mathfrak{S}_{\pi, \pi'}^{(1)}, (E))$ is also a $(p,E')$-completely flat prism over $(\mathfrak{S}', (E'))$. By \cite[Lem.~3.7 (2)]{bhatt-scholze-prismaticcohom}, $(\mathfrak{S}_{\pi, \pi'}^{(1)}, (E))$ is bounded. The universal property can be checked similarly as in Example~\ref{eg:prism-S2S3}.

(ii) The classical flatness follows by a similar argument as in the proof of Lemma~\ref{lem:AtoA2A3-faithful-flat}: note that $\mathfrak{S}_{\pi, \pi'}^{(1)}$ is classically $(p, E)$-complete by \cite[Lem.~3.7 (1)]{bhatt-scholze-prismaticcohom}. Since $\mathfrak{S} \rightarrow \mathfrak{S}_{\pi, \pi'}^{(1)}$ is $(p, E)$-completely flat and $\mathfrak{S}$ is noetherian, $\mathfrak{S} \rightarrow \mathfrak{S}_{\pi, \pi'}^{(1)}$ is classically flat by \cite[Tag 0912]{stacks-project}. 

Consider the composite $R_0 \rightarrow \mathfrak{S} \rightarrow \mathfrak{S}_{\pi, \pi'}^{(1)}$, where the first map is given by the natural inclusion $R_0 \hookrightarrow R_0[\![u]\!]=\fkS$ (which is classically faithfully flat). Since $\mathfrak{S}_{\pi, \pi'}^{(1)}$ is classically $(p, E(u))$-complete, it is $u$-complete and $u$ lies in the radical of $\mathfrak{S}_{\pi, \pi'}^{(1)}$. Thus, to prove $\mathfrak{S} \rightarrow \mathfrak{S}_{\pi, \pi'}^{(1)}$ is classically faithfully flat, it suffices to show that $R_0 \rightarrow \mathfrak{S}_{\pi, \pi'}^{(1)}$ is classically faithfully flat by \cite[Tag 00HQ]{stacks-project}. Let $\mathfrak{P} \subset R$ be a maximal ideal, and let $\mathfrak{m} = R_0 \cap \mathfrak{P}$ and $\mathfrak{m}' = R_0' \cap \mathfrak{P}$ be the corresponding maximal ideals of $R_0$ and $R_0'$ respectively. Let $(R_0)_{\mathfrak{m}}^{\wedge}$ denote the $\mathfrak{m}$-adic completion of the localization $(R_0)_{\mathfrak{m}}$. It is shown in the proof of Proposition \ref{prop:rational-projectivity-etale-over-torus-case} below that $(R_0)_{\mathfrak{m}}^{\wedge}$ is equipped with the Frobenius induced from $R_0$, and that $(R_0)_{\mathfrak{m}}^{\wedge} \cong W(k_1)[\![t_1, \ldots, t_d]\!]$ where $k_1 \coloneqq R/\mathfrak{P}$ is a finite extension of $k$. Similarly, we have $(R_0')_{\mathfrak{m}'}^{\wedge} \cong W(k_1)[\![t_1', \ldots, t_d']\!]$.

Let $A_{\mathfrak{P}}^{(1)}$ be the absolute product of $((R_0)_{\mathfrak{m}}^{\wedge}[\![u]\!], (E))$ and $((R_0')_{\mathfrak{m}'}^{\wedge}[\![y]\!], (E'))$ constructed as in (i) with $R_0$ (resp. $R_0'$) replaced by $(R_0)_{\mathfrak{m}}^{\wedge}$ (resp. $(R_0')_{\mathfrak{m}'}^{\wedge}$). Note that the map
\begin{equation} \label{eq:map-on-localized-rings}
f_{\mathfrak{m}}\colon W(k_1)[\![t_1, \ldots, t_d]\!] \cong (R_0)_{\mathfrak{m}}^{\wedge} \rightarrow A_{\mathfrak{P}}^{(1)}    
\end{equation}
is classically flat similarly as above. Consider the induced map 
\[
W(k_1)[\![t_1, \ldots, t_d]\!] / (t_1, \ldots, t_d) \cong W(k_1) \rightarrow A_{\mathfrak{P}}^{(1)}/(t_1, \ldots, t_d)A_{\mathfrak{P}}^{(1)}.
\]
From the explicit construction of the absolute product $A_{\mathfrak{P}}^{(1)}$ in (i), we deduce that $1 \notin (t_1, \ldots, t_d)A_{\mathfrak{P}}^{(1)}$, and so $A_{\mathfrak{P}}^{(1)}/(t_1, \ldots, t_d)A_{\mathfrak{P}}^{(1)}$ is not the zero ring. Furthermore, since $A_{\mathfrak{P}}^{(1)}$ is classically $p$-complete, $p$ lies in the radical of $A_{\mathfrak{P}}^{(1)}/(t_1, \ldots, t_d)A_{\mathfrak{P}}^{(1)}$. Thus, $A_{\mathfrak{P}}^{(1)}$ has a maximal ideal which lies over the maximal ideal $(p, t_1, \ldots, t_d)$ of $(R_0)_{\mathfrak{m}}^{\wedge}$, and the map $f_{\mathfrak{m}}$ in (\ref{eq:map-on-localized-rings}) is classically faithfully flat. 

Now, consider the map $(R_0)_{\mathfrak{m}}^{\wedge} \rightarrow \mathfrak{S}_{\pi, \pi'}^{(1)}\otimes_{R_0} (R_0)_{\mathfrak{m}}^{\wedge}$ induced from $R_0 \rightarrow \mathfrak{S}_{\pi, \pi'}^{(1)}$. We claim that $(R_0)_{\mathfrak{m}}^{\wedge} \rightarrow \mathfrak{S}_{\pi, \pi'}^{(1)}\otimes_{R_0} (R_0)_{\mathfrak{m}}^{\wedge}$ is classically faithfully flat. The classical flatness is clear. Note that by \cite[Cor. 3.14]{bhatt-scholze-prismaticcohom-v3}, the construction of the absolute product in (i) commutes with $(p, E)$-completely flat base change. Thus, the map $f_{\mathfrak{m}}\colon (R_0)_{\mathfrak{m}}^{\wedge} \rightarrow A_{\mathfrak{P}}^{(1)}$ in  (\ref{eq:map-on-localized-rings}) naturally factors through $(R_0)_{\mathfrak{m}}^{\wedge} \rightarrow \mathfrak{S}_{\pi, \pi'}^{(1)}\otimes_{R_0} (R_0)_{\mathfrak{m}}^{\wedge}$. Since $f_\fkm$ is classically faithfully flat, so is the flat map $(R_0)_{\mathfrak{m}}^{\wedge} \rightarrow \mathfrak{S}_{\pi, \pi'}^{(1)}\otimes_{R_0} (R_0)_{\mathfrak{m}}^{\wedge}$. Now since the claim holds for any maximal ideal $\mathfrak{m} \subset R_0$, $R_0 \rightarrow \mathfrak{S}_{\pi, \pi'}^{(1)}$ is classically faithfully flat.

By symmetry, $\mathfrak{S}' \rightarrow \mathfrak{S}_{\pi, \pi'}^{(1)}$ is also classically faithfully flat.
\end{proof}

\begin{rem}
When $R=\mathbf{Z}_p$, the above lemma follows from \cite[Prop.~2.4.5 and 2.4.9]{Bhatt-Lurie-AbsolutePC} and \cite[Tag 0912]{stacks-project}: for any $(A,I)$ and $(B,J)$ in $R_\Prism$ with $(A,I)$ nonzero and transversal in the sense of \cite[Def.~2.1.3]{Bhatt-Lurie-AbsolutePC}, the product of $(A,I)$ and $(B,J)$ exists in $R_\Prism$, and it covers $(B,J)$. 
\end{rem}

\begin{cor}\label{cor:independence of Breuil--Kisin condition}
Definition~\ref{defn:category-good-prism-completed-F-crystal} of completed prismatic $F$-crystals is independent of the choice of a uniformizer $\pi\in\calO_K$ and a $W$-subalgebra $R_0$ of $R$.
\end{cor}

\begin{proof}
This follows from Lemmas~\ref{lem:independence of Breuil--Kisin condition} and \ref{lem:productoffkSandfkSprime}.
\end{proof}

\begin{rem}[Restriction of completed prismatic $F$-crystals] \label{rem:restriction-of-completed-F-crystals} 
Suppose $R$ is small over $\calO_K$, and let $R\rightarrow R'$ be a $p$-adically completed \'etale map. Let $R_0' \subset R'$ such that $R_0'\otimes_{W} \mathcal{O}_K \cong R'$ with a $p$-adically completed \'etale map $R_0 \rightarrow R_0'$ as in Remark~\ref{rem:p-complete-etale}. Note that the Frobenius on $R_0$ extends uniquely to a Frobenius on $R_0'$. For $\mathcal{F} \in  \mathrm{CR}^{\wedge,\varphi}(R_{\Prism})$, consider its restriction $\mathcal{F} |_{(R')_\Prism}$ to $(R')_\Prism$. Since $\fkS=R_0[\![u]\!]\rightarrow R_0'[\![u]\!]$ is classically flat, we deduce from Remark~\ref{rem:conditions for uncompleted base change}  and Lemma~\ref{lem:independence of Breuil--Kisin condition} that $\calF|_{(R')_\Prism}(R_0'[\![u]\!],(E))=R_0'[\![u]\!]\otimes_{\fkS}\calF_\fkS$ and $\mathcal{F} |_{(R')_{\Prism}} \in \mathrm{CR}^{\wedge,\varphi}((R')_{\Prism})$. We similarly have the restriction of completed prismatic $F$-crystals for the maps $R \rightarrow \mathcal{O}_L$ and $R \rightarrow \mathcal{O}_{K_g}$ as in Notation~\ref{notation:L}. 
\end{rem}

We now study some properties of completed prismatic $F$-crystals on $R$. Let $\mathfrak{S}_L \coloneqq \mathcal{O}_{L_0}[\![u]\!]$ equipped with Frobenius given by $\varphi(u) = u^p$. Note that $(\mathfrak{S}_L, (E)) \in R_{\Prism}$ with $R \rightarrow \mathfrak{S}_L/(E)=\calO_L=R_{(\pi)}^\wedge$ and that the natural map $\mathfrak{S} \rightarrow \mathfrak{S}_L$ induces a map of prisms $(\mathfrak{S}, (E)) \rightarrow (\mathfrak{S}_L, (E))$ over $R$. Let $\mathcal{O}_{\mathcal{E}, L}$ denote the $p$-adic completion of $\mathfrak{S}_L[u^{-1}]$.

\begin{lem} \label{lem:completed-crystals-basic-properties}
Let $\mathcal{F} \in \mathrm{CR}^{\wedge,\varphi}(R_{\Prism})$. 
Then the following properties hold.
\begin{enumerate}
\item We have $\mathcal{F}_{\mathfrak{S}_L} \cong \mathfrak{S}_L\otimes_{\mathfrak{S}}\mathcal{F}_{\mathfrak{S}}$. Furthermore, $\mathcal{F}_{\mathfrak{S}_L}$ is finite free over $\mathfrak{S}_L$.

\item We have $\mathcal{F}_{\mathfrak{S}} = \mathcal{F}_{\mathfrak{S}_L} \cap \mathcal{F}_{\mathfrak{S}}[E^{-1}]^{\wedge}_p$ as submodules of $\mathcal{O}_{\mathcal{E}, L}\otimes_{\mathfrak{S}}\mathcal{F}_{\mathfrak{S}}$.

\item The natural map
\[
\mathfrak{S}^{(1)}\otimes_{p_i,\mathfrak{S}}\mathcal{F}_{\mathfrak{S}} \rightarrow  \mathfrak{S}^{(1)}[E^{-1}]^{\wedge}_p\otimes_{p_i,\mathfrak{S}}\mathcal{F}_{\mathfrak{S}}
\]
is injective for $i = 1, 2$.
\item For any map of bounded prisms $(\mathfrak{S}, (E)) \rightarrow (A, EA)$ over $R$, the natural map 
\[
A[p^{-1}]\otimes_{\mathfrak{S}}\mathcal{F}_{\mathfrak{S}} \rightarrow 
\mathcal{F}_A [p^{-1}]
\]
is a $\varphi$-compatible isomorphism of $A[p^{-1}]$-modules. Similarly, the natural map 
\[
A[E^{-1}]\otimes_{\mathfrak{S}}\mathcal{F}_{\mathfrak{S}} \rightarrow 
\mathcal{F}_A [E^{-1}]
\]
is an isomorphism of $A[E^{-1}]$-modules. Furthermore, the classical $p$-adic completions $(A[E^{-1}]\otimes_{\mathfrak{S}}\mathcal{F}_{\mathfrak{S}})^{\wedge}_p$ and $(\mathcal{F}_A [E^{-1}])^{\wedge}_p$ have naturally induced Frobenii, and the induced isomorphism $(A[E^{-1}]\otimes_{\mathfrak{S}}\mathcal{F}_{\mathfrak{S}})^{\wedge}_p \stackrel{\cong}{\rightarrow} (\mathcal{F}_A [E^{-1}])^{\wedge}_p$ is $\varphi$-compatible.
\end{enumerate}
\end{lem}

\begin{proof}
(i) Since $\mathfrak{S}_L$ is noetherian and $\mathfrak{S} \rightarrow \mathfrak{S}_L$ is classically flat, we deduce by a similar argument as in Remark~\ref{rem:restriction-of-completed-F-crystals} that $\mathcal{F}_{\mathfrak{S}_L} \cong \mathfrak{S}_L\otimes_{\mathfrak{S}} \mathcal{F}_{\mathfrak{S}}$, $\mathcal{F}_{\mathfrak{S}_L}$ is torsion free, and
\[
\mathcal{F}_{\mathfrak{S}_L}[p^{-1}] \cap \mathcal{F}_{\mathfrak{S}_L}[E^{-1}] = \mathcal{F}_{\mathfrak{S}_L}.
\]
So by Remark~\ref{rem:vector-bundle-CDVF-case}, $\mathcal{F}_{\mathfrak{S}_L}$ is finite free over $\mathfrak{S}_L$.

(ii) It suffices to show $\mathcal{F}_{\mathfrak{S}_L} \cap \mathcal{F}_{\mathfrak{S}}[E^{-1}]^{\wedge}_p \subset \mathcal{F}_{\mathfrak{S}}$. Since $\mathcal{F}_{\mathfrak{S}}[p^{-1}]$ is projective over $\mathfrak{S}[p^{-1}]$ and $\mathfrak{S}_L \cap \mathcal{O}_{\mathcal{E}} = \mathfrak{S}$, we have by Lemma~\ref{lem:intersection-modules-flat-base-change} (i) that
\begin{align*}
\mathcal{F}_{\mathfrak{S}_L}[p^{-1}] \cap \mathcal{F}_{\mathfrak{S}}[E^{-1}]^{\wedge}_p[p^{-1}] &= (\mathfrak{S}_L[p^{-1}]\otimes_{\mathfrak{S}[p^{-1}]}\mathcal{F}_{\mathfrak{S}}[p^{-1}]) \cap (\mathcal{O}_{\mathcal{E}}[p^{-1}]\otimes_{\mathfrak{S}[p^{-1}]}\mathcal{F}_{\mathfrak{S}}[p^{-1}]) \\
    &= \mathcal{F}_{\mathfrak{S}}[p^{-1}].
\end{align*}
Thus,
\[
\mathcal{F}_{\mathfrak{S}_L} \cap \mathcal{F}_{\mathfrak{S}}[E^{-1}]^{\wedge}_p \subset \mathcal{F}_{\mathfrak{S}}[p^{-1}] \cap \mathcal{F}_{\mathfrak{S}}[E^{-1}]^{\wedge}_p = \mathcal{F}_{\mathfrak{S}}.
\]

(iii) The natural map
\[
\mathfrak{S}^{(1)}\otimes_{p_i,\mathfrak{S}}\mathcal{F}_{\mathfrak{S}} \rightarrow  \mathfrak{S}^{(1)}[p^{-1}]\otimes_{p_i,\mathfrak{S}}\mathcal{F}_{\mathfrak{S}}
\]
is injective since $\mathcal{F}_{\mathfrak{S}} \rightarrow \mathcal{F}_{\mathfrak{S}}[p^{-1}]$ is injective and $p_i\colon \mathfrak{S} \rightarrow \mathfrak{S}^{(1)}$ is classically flat by Lemma~\ref{lem:AtoA2A3-faithful-flat}. Furthermore, since $\mathcal{F}_{\mathfrak{S}}[p^{-1}]$ is projective over $\mathfrak{S}[p^{-1}]$ and $\mathfrak{S}^{(1)}[p^{-1}] \rightarrow \mathfrak{S}^{(1)}[E^{-1}]^{\wedge}_p[p^{-1}]$ is injective by Corollary~\ref{cor:A2-properties}, the natural map
\[
\mathfrak{S}^{(1)}[p^{-1}]\otimes_{p_i,\mathfrak{S}[p^{-1}]}\mathcal{F}_{\mathfrak{S}}[p^{-1}] \rightarrow  \mathfrak{S}^{(1)}[E^{-1}]^{\wedge}_p[p^{-1}]\otimes_{p_i,\mathfrak{S}[p^{-1}]}\mathcal{F}_{\mathfrak{S}}[p^{-1}]
\]
is injective. So the composite map
\[
\mathfrak{S}^{(1)}\otimes_{p_i,\mathfrak{S}}\mathcal{F}_{\mathfrak{S}} \rightarrow  \mathfrak{S}^{(1)}[p^{-1}]\otimes_{p_i,\mathfrak{S}}\mathcal{F}_{\mathfrak{S}} \rightarrow  \mathfrak{S}^{(1)}[E^{-1}]^{\wedge}_p[p^{-1}]\otimes_{p_i,\mathfrak{S}}\mathcal{F}_{\mathfrak{S}}
\]
is injective. The composite factors through the map
\[
\mathfrak{S}^{(1)}\otimes_{p_i,\mathfrak{S}}\mathcal{F}_{\mathfrak{S}} \rightarrow \mathfrak{S}^{(1)}[E^{-1}]^{\wedge}_p\otimes_{p_i,\mathfrak{S}}\mathcal{F}_{\mathfrak{S}},
\]
which is therefore injective.

(iv)  We first show that the first map is an isomorphism of $A[p^{-1}]$-modules.
By the definition of a finitely generated completed prismatic crystal, 
the map $A\widehat{\otimes}_{\mathfrak{S}}\mathcal{F}_{\fkS}\rightarrow \calF_A$ is an isomorphism. Since $\mathcal{F}_{\mathfrak{S}}$ is finitely generated over $\mathfrak{S}$ and $A$ is classically $(p, E)$-complete, the natural map
\[
A\otimes_{\mathfrak{S}}\mathcal{F}_{\mathfrak{S}} \rightarrow A\widehat\otimes_{\mathfrak{S}}\mathcal{F}_{\mathfrak{S}}
\]
is surjective. So it suffices to show that the induced surjective map $(A\otimes_{\mathfrak{S}}\mathcal{F}_{\mathfrak{S}})[p^{-1}] \rightarrow (A\widehat{\otimes}_{\mathfrak{S}} \mathcal{F}_{\mathfrak{S}})[p^{-1}]$ is also injective.

Since $\mathcal{F}_{\mathfrak{S}}[p^{-1}]$ is finite projective over $\mathfrak{S}[p^{-1}]$, there exists an $\mathfrak{S}[p^{-1}]$-module $Q$ such that $\mathcal{F}_{\mathfrak{S}}[p^{-1}]\oplus Q$ is finite free over $\mathfrak{S}[p^{-1}]$. We have an $\mathfrak{S}$-submodule $N \subset \mathcal{F}_{\mathfrak{S}}[p^{-1}]\oplus Q$ with $N[p^{-1}] = \mathcal{F}_{\mathfrak{S}}[p^{-1}]\oplus Q$ such that $N$ is free over $\mathfrak{S}$ and that the inclusion $\mathcal{F}_{\mathfrak{S}} \hookrightarrow N[p^{-1}]$ factors through $\mathcal{F}_{\mathfrak{S}} \hookrightarrow N \subset N[p^{-1}]$.

Consider the induced map $A\otimes_{\mathfrak{S}}\mathcal{F}_{\mathfrak{S}} \rightarrow A\otimes_{\mathfrak{S}}N$. Note that $A\otimes_{\mathfrak{S}}N$ is $(p, E)$-complete since $N$ is finite free over $\mathfrak{S}$. Thus, this map factors through
\[
A\otimes_{\mathfrak{S}}\mathcal{F}_{\mathfrak{S}} \rightarrow A\widehat{\otimes}_{\mathfrak{S}}\mathcal{F}_{\mathfrak{S}}  \rightarrow A\otimes_{\mathfrak{S}}N.
\] 
On the other hand, since $\mathcal{F}_{\mathfrak{S}}[p^{-1}]$ is a direct summand of $N[p^{-1}]$, the map $(A\otimes_{\mathfrak{S}} \mathcal{F}_{\mathfrak{S}})[p^{-1}] \rightarrow (A\otimes_{\mathfrak{S}} N)[p^{-1}]$ is injective. Since it factors through
\[
(A\otimes_{\mathfrak{S}}\mathcal{F}_{\mathfrak{S}})[p^{-1}] \rightarrow (A\widehat{\otimes}_{\mathfrak{S}}\mathcal{F}_{\mathfrak{S}})[p^{-1}]  \rightarrow (A\otimes_{\mathfrak{S}}N)[p^{-1}],
\] 
the map $(A\otimes_{\mathfrak{S}}\mathcal{F}_{\mathfrak{S}})[p^{-1}] \rightarrow (A\widehat{\otimes}_{\mathfrak{S}} \mathcal{F}_{\mathfrak{S}})[p^{-1}]$ in question is also injective.

Similarly, the second map is an isomorphism of $A[E^{-1}]$-modules since $\mathcal{F}_{\mathfrak{S}}[E^{-1}]$ is finite projective over $\mathfrak{S}[E^{-1}]$ by Remark~\ref{rem:saturated-regular-seq} (i). Hence it remains to show the statements for $\varphi$-compatibility. Note that $\varphi((p, E)^m) \subset (p, E)^m$ for each $m \geq 1$. It follows that  $A\widehat{\otimes}_{\mathfrak{S}}\mathcal{F}_{\mathfrak{S}}$ admits a Frobenius endomorphism induced from that on $A\otimes_{\mathfrak{S}}\mathcal{F}_{\mathfrak{S}}$. Thus, the natural map $A\otimes_{\mathfrak{S}}\mathcal{F}_{\mathfrak{S}} \rightarrow A\widehat\otimes_{\mathfrak{S}}\mathcal{F}_{\mathfrak{S}}$ is $\varphi$-compatible, and so are
\[
(A\otimes_{\mathfrak{S}}\mathcal{F}_{\mathfrak{S}})[p^{-1}] \to (A\widehat{\otimes}_{\mathfrak{S}}\mathcal{F}_{\mathfrak{S}})[p^{-1}]
\quad\text{and}\quad
((A\otimes_{\mathfrak{S}}\mathcal{F}_{\mathfrak{S}})[E^{-1}])/p^n \to ((A\widehat{\otimes}_{\mathfrak{S}}\mathcal{F}_{\mathfrak{S}})[E^{-1}])/p^n
\quad\text{for each $n \geq 1$.}
\]
\end{proof}

\subsection{Completed prismatic \texorpdfstring{$F$}{F}-crystals in terms of descent data}\label{sec:descent data}
Keep Assumption~\ref{assumption:base-ring-sec-3.4}.
We can explicitly describe the category $\mathrm{CR}^{\wedge,\varphi}(R_{\Prism})$ in terms of certain descent data as follows. 

\begin{defn}
Let $\mathrm{DD}_\fkS$ denote the category consisting of 
triples $(\mathfrak{M}, \varphi_\fkM, f)$ where
\begin{enumerate}
\item $\mathfrak{M}$ is a finite $\mathfrak{S}$-module that is projective away from $(p, E)$ and saturated;
\item $\varphi_\fkM\colon \mathfrak{M} \rightarrow \mathfrak{M}$ is a $\varphi$-semi-linear endomorphism such that $(\fkM,\varphi_\fkM)$ has finite $E$-height;
\item $f\colon \mathfrak{S}^{(1)} \otimes_{p_1,\mathfrak{S}}\mathfrak{M} \xrightarrow{\cong} \mathfrak{S}^{(1)} \otimes_{p_2,\mathfrak{S}}\mathfrak{M}$ is an isomorphism of $\fkS^{(1)}$-modules that is compatible with Frobenii and satisfies the cocycle condition over $\mathfrak{S}^{(2)}$.
\end{enumerate}
The morphisms of $\mathrm{DD}_\fkS$ are $\mathfrak{S}$-linear maps compatible with all structures.

For a fixed non-negative integer $r$, let $\mathrm{DD}_{\fkS,[0, r]}$ denote the full subcategory consisting of objects for which $(\fkM,\varphi_\fkM)$ has $E$-height $\leq r$.

We call an object of $\mathrm{DD}_\fkS$  an \emph{integral Kisin descent datum}.
One can also consider a triple $(\fkM,\varphi_\fkM,f)$ where $(\fkM,\varphi_\fkM)$ is as above and $f\colon \mathfrak{S}^{(1)}[p^{-1}] \otimes_{p_1,\mathfrak{S}}\mathfrak{M} \xrightarrow{\cong} \mathfrak{S}^{(1)}[p^{-1}] \otimes_{p_2,\mathfrak{S}}\mathfrak{M}$ is an isomorphism of $\fkS^{(1)}[p^{-1}]$-modules that is compatible with Frobenii and satisfies the cocycle condition over $\mathfrak{S}^{(2)}[p^{-1}]$. Such an object is called a \emph{rational Kisin descent datum}.
\end{defn}

\begin{prop} \label{prop:equivalence-to-descent-datum}
The association $\calF\mapsto \calF_\fkS = \calF(\fkS,(E))$ gives rise to a functor $\mathrm{CR}^{\wedge,\varphi}(R_{\Prism})\rightarrow\mathrm{DD}_\fkS$ and induces equivalences of categories 
\[
 \mathrm{CR}^{\wedge,\varphi}(R_{\Prism})\cong\mathrm{DD}_{\fkS}\quad\text{and}\quad
 \mathrm{CR}^{\wedge,\varphi}_{[0,r]}(R_{\Prism})\cong\mathrm{DD}_{\fkS,[0,r]}.
\]
Furthermore, under this equivalence, $(\mathfrak{M}, \varphi, f)$ corresponds to an object in $\mathrm{Vect}^{\varphi}_{\mathrm{eff}}(R_{\Prism})$ if and only if $\mathfrak{M}$ is finite projective over $\mathfrak{S}$. 
\end{prop}

\begin{proof}
Let $\mathcal{F} \in \mathrm{CR}^{\wedge,\varphi}(R_{\Prism})$. By Lemma~\ref{lem:AtoA2A3-faithful-flat} and Remark~\ref{rem:conditions for uncompleted base change}, we have an isomorphism of $\mathfrak{S}^{(1)}$-modules
\[
f\colon  \mathfrak{S}^{(1)}\otimes_{p_1,\mathfrak{S}}\mathcal{F}_{\mathfrak{S}} \xrightarrow{\cong}\mathcal{F}_{\mathfrak{S}^{(1)}}\stackrel{\cong}{\leftarrow}\mathfrak{S}^{(1)}\otimes_{p_2,\mathfrak{S}} \mathcal{F}_{\mathfrak{S}}
\]
satisfying the cocycle condition over $\mathfrak{S}^{(2)}$. Thus, any completed crystal in $\mathrm{CR}^{\wedge,\varphi}(R_{\Prism})$ naturally gives an object in $\mathrm{DD}_\fkS$ via $\mathcal{F} \mapsto \mathcal{F}_{\mathfrak{S}}$ , which gives a functor from $\mathrm{CR}^{\wedge,\varphi}(R_{\Prism})$ to $\mathrm{DD}_\fkS$. 

Conversely, let $(\mathfrak{M}, \varphi_\fkM, f) \in \mathrm{DD}_\fkS$. Take any prism $(A, I) \in R_{\Prism}$. By \cite[Lem.~4.1.8]{du-liu-prismaticphiGhatmodule}, there exists a prism $(B, IB)\in R_\Prism$ which covers $(A, I)$ and admits a map $(\mathfrak{S}, (E))\rightarrow (B,IB)$ over $R$. 
By Lemma~\ref{lem:pushout for prisms along fflat map}, the pushout of the diagram $(B,IB)\leftarrow (A,I)\rightarrow (B,IB)$ of maps of bounded prism over $R$ is represented by 
$(B\otimes_A B)^{\wedge}_{(p, I)}$, and $(B\otimes_A B\otimes_AB)^{\wedge}_{(p, I)}$ satisfies a similar property for the self-triple cofiber product. By the universal property of $\fkS^{(1)}$ and $\fkS^{(2)}$, we have maps $\fkS^{(1)}\rightarrow (B\otimes_A B)^{\wedge}_{(p, I)}$ and $\fkS^{(2)}\rightarrow (B\otimes_A B\otimes_AB)^{\wedge}_{(p, I)}$.

Consider the $B$-module $B\otimes_\fkS \fkM$.
The base change of the descent datum
$f\colon \mathfrak{S}^{(1)}\otimes_{p_1,\mathfrak{S}}\mathfrak{M} \xrightarrow{\cong} \mathfrak{S}^{(1)}\otimes_{p_2,\mathfrak{S}}\mathfrak{M}$ along $\fkS^{(1)}\rightarrow (B\otimes_A B)^{\wedge}_{(p, I)}$ gives a descent datum of $B\otimes_\fkS \fkM$, namely, a $(B\otimes_A B)^{\wedge}_{(p, I)}$-linear isomorphism
\[
f_B\colon (B\otimes_A B)^{\wedge}_{(p, I)}\otimes_{p_1, B} (B\otimes_{\mathfrak{S}}\mathfrak{M}) \xrightarrow{\cong}  (B\otimes_A B)^{\wedge}_{(p, I)}\otimes_{p_2,B} (B\otimes_{\mathfrak{S}}\mathfrak{M})
\]
satisfying the cocycle condition over $(B\otimes_A B\otimes_AB)^{\wedge}_{(p, I)}$.
By reducing modulo $(p,I)^n$, $f_B$ induces a compatible system of isomorphisms
\[
f_{B,n}\colon (B\otimes_A B)/(p, I)^n\otimes_{p_1,B}(B\otimes_{\mathfrak{S}}\mathfrak{M}) \xrightarrow{\cong} (B\otimes_A B)/(p, I)^n\otimes_{p_2,B}(B\otimes_{\mathfrak{S}}\mathfrak{M}) 
\]
satisfying the cocycle condition over $(B\otimes_A B\otimes_A B)/(p, I)^n$ for each $n \geq 1$. 

Since $A \rightarrow B$ is $(p, I)$-completely faithfully flat, each $f_{B,n}$ defines a finitely generated $A/(p,I)^n$-module $\calF_{n,A}$ by the usual faithfully flat descent.
We claim that $\calF_{n,A}$ is independent of the choice of the cover $(A,I)\rightarrow (B,IB)$ and that the association $(A,I)\mapsto \calF_{n,A}$ defines a sheaf $\calF_n$ of $\calO_\Prism$-modules on $R_\Prism$.
To see the former, take another prism $(B',IB')\in R_\Prism$ which covers $(A, I)$ and admits a map $(\mathfrak{S}, (E))\rightarrow (B',IB')$ over $R$. Let $\calF_{n,A}'$ denote the finitely generated $A/(p,I)^n$-module given by the descent of $\bigl((B'\otimes_\fkS\fkM)/(p,I)^n,f_{B',n}\bigr)$.
By Lemma~\ref{lem:pushout for prisms along fflat map}, the pushout of the diagram $(B,IB)\leftarrow (A,I)\rightarrow (B',IB')$ of maps of bounded prism is represented by  $(B\otimes_AB')^\wedge_{(p,I)}$. Since the maps $B\rightarrow (B\otimes_AB')^\wedge_{(p,I)}$ and $B'\rightarrow (B\otimes_AB')^\wedge_{(p,I)}$ are $(p,I)$-completely faithfully flat, we can canonically identify both $\calF_{n,A}$ and $\calF_{n,A}'$ with the descent of $\bigl(( (B\otimes_AB')^\wedge_{(p,I)}\otimes_\fkS\fkM)/(p,I)^n,f_{(B\otimes_AB')^\wedge_{(p,I)},n}\bigr)$. 
To see that the association $(A,I)\mapsto \calF_{n,A}$ defines a sheaf $\calF_n$ of $\calO_\Prism$-modules on $R_\Prism$, take a $(p,I)$-completely faithfully flat map of prisms $(A,I)\rightarrow (A',IA')$ over $R$. Then the pushout of the diagram $(A',IA')\leftarrow (A,I)\rightarrow (A',IA')$ is represented by  $(A'\otimes_AA')^\wedge_{(p,I)}$.
Hence we need to show the exactness of the sequence
\begin{equation}\label{eq:exact seq in the proof of descent datum over fkS}
0\rightarrow \calF_{n,A}\rightarrow \calF_{n,A'}\rightarrow \calF_{n,(A'\otimes_AA')^\wedge_{(p,I)}}.    
\end{equation}
On the other hand, we see that $(A'\otimes_AB)^\wedge_{(p,I)}$ (resp.~ $(A'\otimes_AA'\otimes_AB)^\wedge_{(p,I)}$) together with the ideal generated by $I$ gives a bounded prism over $R$ that admits a map from $(\fkS,(E))$ over $R$ and covers $(A',IA')$ (resp.~$((A'\otimes_AA')^\wedge_{(p,I)}, I(A'\otimes_AA')^\wedge_{(p,I)})$). By construction, we have a left exact sequence 
\[
0\rightarrow (B\otimes_\fkS \fkM)/(p,I)^n\rightarrow (A'\otimes B\otimes_\fkS \fkM)/(p,I)^n\rightarrow (A'\otimes_AA'\otimes_AB\otimes_\fkS \fkM )/(p,I)^n.
\]
Since this left exact sequence is the base change of the sequence \eqref{eq:exact seq in the proof of descent datum over fkS}
along the classically faithfully flat map $A/(p,I)^n\rightarrow B/(p,I)^n$, we conclude that the sequence \eqref{eq:exact seq in the proof of descent datum over fkS} is left exact. This completes the verification of the claim.

The sheaf $\mathcal{F}_n$ is equipped with an induced Frobenius, since $\varphi((p, I)^n) \subset (p, I)^n$. Furthermore, $\mathcal{F}_n$ is a finitely generated crystal of $\calO_{\Prism,n}$-modules. This follows from a similar argument as in the above paragraph and the verification is left to the reader.
We also remark that $\{\mathcal{F}_n\}_{n \geq 1}$ forms an inverse system of sheaves of $\mathcal{O}_{\Prism}$-modules such that $\calO_{\Prism,n+1}\otimes_{\calO_{\Prism,n}}\calF_{n+1}\cong \calF_n$. Hence $\mathcal{F} \coloneqq \varprojlim_n \mathcal{F}_n$ is a completed prismatic crystal on $R$ equipped with Frobenius by Lemma~\ref{lem:calF vs (calF_n) for completed crystals}.
By construction, we see $\calF(\fkS,(E))=\fkM$.
As a result, $\mathcal{F} \in \mathrm{CR}^{\wedge,\varphi}(R_{\Prism})$. This proves the essential surjectivity. 

The fully faithfulness also follows directly from a similar argument as above (alternatively, one can check that the above two functors are quasi-inverse to each other). Obviously, this equivalence also induces $\mathrm{CR}^{\wedge,\varphi}_{[0,r]}(R_{\Prism})\cong\mathrm{DD}_{\fkS,[0,r]}$. The last assertion follows from \cite[Tag 0D4B]{stacks-project}.
\end{proof}

\subsection{\'Etale realization and the main theorem in the small affine case} \label{sec:etale-realization-main-thm}

We now formulate our main theorem. 
For this, we first attach to a completed prismatic $F$-crystal $\calF$ on $R$ a finite free $\Z_p$-representation $T(\calF)$ of $\calG_R$. This will be based on the results in \cite[\S 3]{bhatt-scholze-prismaticFcrystal} (see also \cite{min-wang-rel-phi-gamma-prism-F-crys}). Keep Assumption~\ref{assumption:base-ring-sec-3.4}: $R$ is small over $\calO_K$ or $R=\calO_L$.

Recall that $\mathrm{Vect}(R_{\Prism}, \mathcal{O}_{\Prism}[1/\mathcal{I_{\Prism}}]^{\wedge}_p)^{\varphi = 1}$ denotes the category of Laurent $F$-crystals, i.e., crystals of vector bundles $\calV$ on $(R_\Prism,\calO_\Prism[1/\calI_\Prism]^\wedge_p)$
together with isomorphisms $\varphi_{\calV}\colon \varphi^\ast\calV\cong \calV$ \cite[Def.~3.2]{bhatt-scholze-prismaticFcrystal}. There is an equivalence of categories
\[
\mathrm{Vect}(R_{\Prism}, \mathcal{O}_{\Prism}[1/\mathcal{I_{\Prism}}]^{\wedge}_p)^{\varphi = 1}\cong \mathrm{Rep}^\mathrm{pr}_{\Z_p}(\calG_R)
\]
given by $(\calV,\varphi_{\calV})\mapsto \calV(\A_{\mathrm{inf}}(\overline{R}),(\xi))^{\varphi_{\calV}=1}$ (see \cite[Cor.~3.8]{bhatt-scholze-prismaticFcrystal},  \cite[Thm.~ 3.2]{min-wang-rel-phi-gamma-prism-F-crys}), which is functorial in $R$.

\begin{prop} \label{prop:etale-realization} \hfill
\begin{enumerate}
\item 
The assignment $\calF\mapsto \calF_\et \coloneqq \varprojlim_n\calO_\Prism[1/\calI_\Prism]/p^n\otimes_{\calO_\Prism}\mathcal{F}$ 
defines a faithful functor
\[
\mathrm{CR}^{\wedge,\varphi}(R_\Prism)\rightarrow \mathrm{Vect}(R_{\Prism}, \mathcal{O}_{\Prism}[1/\mathcal{I_{\Prism}}]^{\wedge}_p)^{\varphi = 1}.
\]
Moreover, if $\mathcal{F}\in\mathrm{Vect}^\varphi_{\mathrm{eff}}(R_\Prism)$, then
the canonical morphism $\calO_\Prism[1/\calI_\Prism]^\wedge_p\otimes_{\calO_\Prism}\calF\rightarrow\mathcal{F}_{\mathrm{\acute{e}t}}$ is an isomorphism.
\item Define a contravariant functor $T\colon \mathrm{CR}^{\wedge,\varphi}(R_\Prism)\rightarrow \mathrm{Rep}_{\Z_p}^\mathrm{pr}(\calG_R)$
by
\[
T(\mathcal{F})\coloneqq 
\bigl((\mathcal{F}_{\et}(\mathbf{A}_{\mathrm{inf}}(\overline{R}), (\xi)))^{\varphi_{\calF_{\et}} = 1}\bigr)^\vee.
\]
Then it satisfies the following properties:
\begin{enumerate}
\item there is a $\mathcal{G}_R$-equivariant identification
\[
T(\mathcal{F})[p^{-1}]^{\vee} = (W(\overline{R}^{\flat}[(\pi^{\flat})^{-1}])[p^{-1}] \otimes_{\mathbf{A}_{\mathrm{inf}}(\overline{R})} \mathcal{F}_{\mathbf{A}_{\mathrm{inf}}(\overline{R})})^{\varphi = 1},
\]
where the $\mathcal{G}_R$-action on the right hand side is the tensor product of those on $W(\overline{R}^{\flat}[(\pi^{\flat})^{-1}])$ and on $\mathcal{F}_{\mathbf{A}_{\mathrm{inf}}(\overline{R})} = \mathcal{F}(\mathbf{A}_{\mathrm{inf}}(\overline{R}), (\xi))$;
\item $\mathcal{F}_{\mathfrak{S}}[E^{-1}]^{\wedge}_p$ is the \'etale $\varphi$-module associated with $T(\mathcal{F})^{\vee}|_{\mathcal{G}_{\tilde{R}_{\infty}}}$ via \emph{Proposition~\ref{prop:etale-gal-equiv}}. 
\item If $R$ is small over $\calO_K$ and if $R\rightarrow R'$ is a $p$-adically completed \'etale map together with a compatible $W$-map $R_0\rightarrow R_0'$, then $T$ is compatible with the restrictions $\mathrm{CR}^{\wedge,\varphi}(R_\Prism)\rightarrow \mathrm{CR}^{\wedge,\varphi}(R'_\Prism)$ (see Remark~\ref{rem:restriction-of-completed-F-crystals}) and $\mathrm{Rep}_{\Z_p}^\mathrm{pr}(\calG_R)\rightarrow \mathrm{Rep}_{\Z_p}^\mathrm{pr}(\calG_{R'})$. We also have the analogous compatibility for the base changes along $R \rightarrow \mathcal{O}_L$ and $R \rightarrow \mathcal{O}_{K_g}$.
\end{enumerate}
\end{enumerate}
\end{prop}

We call the functor $T$ the \emph{\'etale realization functor}.
We remark that our functor is contravariant and it is the dual of the covariant \'etale realization functor in \cite{bhatt-scholze-prismaticFcrystal}. Our contravariant convention agrees with the one in the theory of Breuil--Kisin modules \cite{kisin-crystalline, brinon-trihan}, which is heavily used in this paper.

\begin{proof}
(i) 
We start with a description of $\mathrm{Vect}(R_{\Prism}, \mathcal{O}_{\Prism}[1/\mathcal{I_{\Prism}}]^{\wedge}_p)^{\varphi = 1}$ in terms of the category of certain descent data: let $\mathrm{DD}_{\mathcal{O}_{\mathcal{E}}}$ denote the category of triples $(\mathcal{M},\varphi_{\mathcal{M}},g)$ where $(\mathcal{M},\varphi_{\mathcal{M}})$ is a finite projective \'etale $\varphi$-module over $\mathcal{O}_{\mathcal{E}}$ as in \S~\ref{sec:etale phi-module}, and $g$ is an isomorphism of $\mathfrak{S}^{(1)}[E^{-1}]^\wedge_p$-modules 
\[
g\colon \mathfrak{S}^{(1)}[E^{-1}]^\wedge_p \otimes_{p_1,\mathcal{O}_{\mathcal{E}}}\mathcal{M} \xrightarrow{\cong} \mathfrak{S}^{(1)}[E^{-1}]^\wedge_p \otimes_{p_2,\mathcal{O}_{\mathcal{E}}}\mathcal{M}
\]
that is compatible with Frobenii and satisfies the cocycle condition over $\mathfrak{S}^{(2)}[E^{-1}]^\wedge_p$. We claim that evaluating on the diagram $\mathfrak{S}\xrightarrow{p_1} \mathfrak{S}^{(1)}\xleftarrow{p_2}\mathfrak{S}$ gives an equivalence of categories from $\mathrm{Vect}(R_{\Prism}, \mathcal{O}_{\Prism}[1/\mathcal{I_{\Prism}}]^{\wedge}_p)^{\varphi = 1}$ to $\mathrm{DD}_{\mathcal{O}_{\mathcal{E}}}$;
for any prism $(A, I) \in R_{\Prism}$, take a prism $(B, IB)\in R_\Prism$ which covers $(A, I)$ and admits a map $(\mathfrak{S}, (E))\rightarrow (B,IB)$ over $R$. Then $(B\otimes_A B)^{\wedge}_{(p, I)}$ and $(B\otimes_A B\otimes_AB)^{\wedge}_{(p, I)}$ represent the self cofiber product and the self-triple cofiber product of the map $(A,I)\rightarrow (B,IB)$ of bounded prisms over $R$. Moreover, for each $n$, the map $A/p^n\rightarrow B/p^n$ is $I$-completely faithfully flat, and the self cofiber product and the self-triple cofiber product for $I$-completely flat topology are given by $(B/p^n\otimes_{A/p^n}B/p^n)^\wedge_I$ and $(B/p^n\otimes_{A/p^n}B/p^n\otimes_{A/p^n}B/p^n)^\wedge_I$, respectively. Now the claim follows as in the proof of Proposition~\ref{prop:equivalence-to-descent-datum} with faithfully flat descent replaced by \cite[Thm.~7.8]{Mathew-flatdescent} (see also \cite[\S 3]{wu-Gal-rep-prism-F-cryst} when $R=\mathcal{O}_K$).
Moreover, the proof shows the following: if $\mathcal{G}$ is the Laurent $F$-crystal associated to $(\mathcal{M},\varphi_{\mathcal{M}},g)$, then $\mathcal{G}(B,J)=\varprojlim_n B/p^n\otimes_{\mathfrak{S}}\mathcal{M}$ for any prism $(B,J)$ that admits a map from $(\mathfrak{S},(E))$.

Let $\mathcal{F} \in \mathrm{CR}^{\wedge,\varphi}(R_{\Prism})$. By Proposition~\ref{prop:equivalence-to-descent-datum}, we have an isomorphism of $\fkS^{(1)}$-modules
\[
f\colon  \mathfrak{S}^{(1)}\otimes_{p_1, \mathfrak{S}}\mathcal{F}_{\mathfrak{S}} \xrightarrow{\cong} \mathfrak{S}^{(1)} \otimes_{p_2,\mathfrak{S}}\mathcal{F}_{\mathfrak{S}}
\]
satisfying the cocycle condition over $\mathfrak{S}^{(2)}$. Let $\mathcal{M} = \mathcal{F}_{\mathfrak{S}}[E^{-1}]^{\wedge}_p$, which is a finite projective \'etale $\varphi$-module over $\mathcal{O}_{\mathcal{E}}$. By extending scalars, $f$ induces a descent datum $(\mathcal{M},\varphi_{\mathcal{M}},g)$ with
\[
g\colon \mathfrak{S}^{(1)}[E^{-1}]^{\wedge}_p\otimes_{p_1,\mathcal{O}_{\mathcal{E}}}\mathcal{M} \xrightarrow{\cong}  \mathfrak{S}^{(1)}[E^{-1}]^{\wedge}_p\otimes_{p_2,\mathcal{O}_{\mathcal{E}}}\mathcal{M}.
\]
Via the equivalence $\mathrm{Vect}(R_{\Prism}, \mathcal{O}_{\Prism}[1/\mathcal{I_{\Prism}}]^{\wedge}_p)^{\varphi = 1}\cong \mathrm{DD}_{\mathcal{O}_{\mathcal{E}}}$, it defines a Laurent $F$-crystal $\mathcal{F}_{\et}'$ such that $\mathcal{F}_\et'(B,J)=\varprojlim_nB/p^n\otimes_{\mathfrak{S}}\mathcal{M}$ for any prism $(B, J) \in R_{\Prism}$ with a map $(\mathfrak{S}, (E))\rightarrow (B,J)$. For such a prism, Lemma~\ref{lem:completed-crystals-basic-properties} (iv) gives the natural isomorphism $B[E^{-1}]/p^n\otimes_{\mathfrak{S}}\mathcal{M}/p^n\xrightarrow{\cong} \mathcal{F}_{B}[E^{-1}]/p^n$ and thus yields an identification $\mathcal{F}_\et'(B,J)=\mathcal{F}_B[E^{-1}]^\wedge_p$. Since every prism over $R$ admits a map from $(\mathfrak{S},(E))$ locally in the $(p,I)$-completely faithfully flat topology, one can check that $\mathcal{F}_\et'$ coincides with $\mathcal{F}_{\et}\coloneqq \varprojlim_n\calO_\Prism[1/\calI_\Prism]/p^n\otimes_{\calO_\Prism}\mathcal{F}$, which is obtained as the sheafification of the presheaf $(A,I)\mapsto \varprojlim_n A[1/I]/p^n\otimes_A\mathcal{F}_A$.

Let us verify the remaining assertions: the faithfulness of $\calF\mapsto \calF_\et\cong \calF_\et'$ follows from the construction of $\mathcal{F}_{\et}'$ and Lemma~\ref{lem:completed-crystals-basic-properties} (iii). 
Now assume $\mathcal{F}\in\mathrm{Vect}^\varphi_{\mathrm{eff}}(R_\Prism)$ and consider 
the canonical morphism $\calO_\Prism[1/\calI_\Prism]^\wedge_p\otimes_{\calO_\Prism}\calF\rightarrow\mathcal{F}_{\mathrm{\acute{e}t}}$. We need to show that it is an isomorphism, which can be checked locally on $R_\Prism$. Since $\mathcal{F}$ is locally a direct summand of $\mathcal{O}_\Prism^{\oplus m}$ for some $m$, the statement follows from the case $\mathcal{F}=\mathcal{O}_\Prism$.

(ii) 
By the paragraph before the proposition,
$\bigl((\mathcal{F}_{\et}(\mathbf{A}_{\mathrm{inf}}(\overline{R}), (\xi)))^{\varphi_{\calF_{\et}} = 1}\bigr)^\vee$ is a finite free $\Z_p$-representation of $\calG_R$ for $\calF\in\mathrm{CR}^{\wedge,\varphi}(R_\Prism)$,
and $T$ is well-defined.

First we verify (a). 
By construction in (i), we have
\[
\mathcal{F}_{\et}(\mathbf{A}_{\mathrm{inf}}(\overline{R}), (\xi)) \cong  \mathbf{A}_{\mathrm{inf}}(\overline{R})[E^{-1}]^{\wedge}_p\otimes_{\mathcal{O}_{\mathcal{E}}}\mathcal{M} \cong  W(\overline{R}^{\flat}[(\pi^{\flat})^{-1}])\otimes_{\mathfrak{S}}\mathcal{F}_{\mathfrak{S}}.
\]
On the other hand, it follows from Lemma~\ref{lem:completed-crystals-basic-properties} (iv) and Example~\ref{eg:prism-Ainf} that
\[
\mathcal{F}(\mathbf{A}_{\mathrm{inf}}(\overline{R}), (\xi))[p^{-1}] \cong  \mathbf{A}_{\mathrm{inf}}(\overline{R})[p^{-1}]\otimes_{\mathfrak{S}}\mathcal{F}_{\mathfrak{S}}.
\]
Thus, we deduce
\begin{align*}
T(\mathcal{F})[p^{-1}]^{\vee} 
&=(\mathcal{F}_{\et}(\mathbf{A}_{\mathrm{inf}}(\overline{R}), (\xi)))^{\varphi_{\calF_{\et}} = 1}[p^{-1}]
=(\mathcal{F}_{\et}(\mathbf{A}_{\mathrm{inf}}(\overline{R}), (\xi))[p^{-1}])^{\varphi_{\calF_{\et}} = 1}\\
&\cong (W(\overline{R}^{\flat}[(\pi^{\flat})^{-1}])[p^{-1}]\otimes_{\mathfrak{S}} \mathcal{F}_{\mathfrak{S}})^{\varphi = 1}\\
&\cong (W(\overline{R}^{\flat}[(\pi^{\flat})^{-1}])[p^{-1}]\otimes_{\mathbf{A}_{\mathrm{inf}}(\overline{R})}\mathcal{F}_{\mathbf{A}_{\mathrm{inf}}(\overline{R})})^{\varphi = 1}.
\end{align*}

Since $\mathcal{F}_{\et}$ is a crystal, the $\mathcal{G}_R$-action on the prism $(\mathbf{A}_{\mathrm{inf}}(\overline{R}), (\xi))$ induces the $\mathcal{G}_R$-action on the last term $( W(\overline{R}^{\flat}[(\pi^{\flat})^{-1}])[p^{-1}]\otimes_{\mathbf{A}_{\mathrm{inf}}(\overline{R})} \mathcal{F}_{\mathbf{A}_{\mathrm{inf}}(\overline{R})})^{\varphi = 1}$, for which $T(\mathcal{F})[p^{-1}]^{\vee} \cong (W(\overline{R}^{\flat}[(\pi^{\flat})^{-1}])[p^{-1}]\otimes_{\mathbf{A}_{\mathrm{inf}}(\overline{R})}\mathcal{F}_{\mathbf{A}_{\mathrm{inf}}(\overline{R})})^{\varphi = 1}$ in the above paragraph becomes $\calG_R$-equivariant.

Next we prove (b). Since $T(\mathcal{F})^{\vee} \cong (W(\overline{R}^{\flat}[(\pi^{\flat})^{-1}])\otimes_{\mathcal{O}_{\mathcal{E}}}\mathcal{M})^{\varphi = 1}$, it suffices to show that the natural injective map
\[
(\widehat{\mathcal{O}}_{\mathcal{E}}^{\mathrm{ur}}\otimes_{\mathcal{O}_\mathcal{E}}\mathcal{M})^{\varphi = 1} \rightarrow (W(\overline{R}^{\flat}[(\pi^{\flat})^{-1}])\otimes_{\mathcal{O}_{\mathcal{E}}}\mathcal{M} )^{\varphi = 1}
\]
is bijective. Indeed, this holds for any \'etale $(\varphi,\calO_\calE)$-module; as in the proof of \cite[Lem.~2.1.4]{gao-liu-loose-cryst-lifts}, one can reduce it to the $p$-torsion case, where the \'etale $(\varphi,\calO_\calE)$-module is finite projective over $\calO_\calE/(p)$ by \cite[p.~8200]{kim-groupscheme-relative} and thus the result follows from $\mathbf{F}_p=(\tilde{E}_{R_\infty})^{\varphi=1}=(\overline{R}^{\flat}[(\pi^{\flat})^{-1}])^{\varphi=1}$.

We now prove (c). Let $R \rightarrow R'$ be a $p$-complete \'etale map. From the above construction, we have an induced map of $\mathbf{Z}_p$-modules 
\[
T(\mathcal{F})^{\vee} \rightarrow T(\mathcal{F} |_{R'})^{\vee}
\]
which is compatible with $\mathcal{G}_{R'}$-actions. By part (b) and the functoriality of \'etale $\varphi$-modules as in the end of \S~\ref{sec:etale phi-module}, this map $T(\mathcal{F})^{\vee} \rightarrow T(\mathcal{F} |_{R'})^{\vee}$ is an isomorphism. The statements for $R \rightarrow \mathcal{O}_L$ and $R \rightarrow \mathcal{O}_{K_g}$ follow from a similar argument.
\end{proof}

\begin{eg}\label{eg:Breuil-Kisin twists}
Recall the Breuil--Kisin twist $\calO_\Prism\{1\}\in \mathrm{Vect}^\varphi(R_\Prism)$ from \cite[Ex.~4.5]{bhatt-scholze-prismaticFcrystal}. It is an invertible $\calO_\Prism$-module with $\varphi^\ast\calO_\Prism\{1\}\cong\calI_\Prism^{-1}\calO_\Prism\{1\}$ and is given informally by $\calO_\Prism\{1\}\coloneqq \bigotimes_{i\geq 0}(\varphi^i)^\ast\calI_\Prism$. For $n \in \mathbb{Z}$, set $\calO_\Prism\{n\}\coloneqq  \bigl(\calO_\Prism\{1\}\bigr)^{\otimes n}$. This is an invertible $\calO_\Prism$-module such that $\varphi^\ast\calO_\Prism\{n\}\cong \calI_\Prism^{-n}\calO_\Prism\{n\}$. By \cite[Ex. 4.9]{bhatt-scholze-prismaticFcrystal}, we have $T(\calO_\Prism\{n\})=\Z_p(-n)$.\footnote{Recall that our \'etale realization functor $T$ is the dual of that of \cite{bhatt-scholze-prismaticFcrystal}.}

For $\mathcal{F} \in \mathrm{CR}^{\wedge,\varphi}(R_\Prism)$, consider a sheaf of $\mathcal{O}_{\Prism}$-modules $\calF\{n\}\coloneqq \calF\otimes_{\calO_\Prism}\calO_\Prism\{n\}$. Suppose that the image of the induced map $\varphi^*(\calF\{n\})_{\mathfrak{S}} \rightarrow (\calF\{n\})_{\mathfrak{S}}[E^{-1}]$ lies in $(\calF\{n\})_{\mathfrak{S}}$. It follows directly from the definition that $\calF\{n\} \in \mathrm{CR}^{\wedge,\varphi}(R_\Prism)$. We claim that $T(\calF\{n\}) \cong T(\calF)\otimes_{\Z_p}\mathbf{Z}_p(-n)$. To see this, note that we have a natural $\mathcal{G}_R$-equivariant map $T(\mathcal{F})^{\vee}\otimes_{\mathbf{Z}_p}\mathbf{Z}_p(n)=T(\mathcal{F})^{\vee}\otimes_{\mathbf{Z}_p}T(\calO_\Prism\{n\})^\vee \rightarrow T(\calF\{n\})^{\vee}$. Since the equivalence in Proposition~\ref{prop:etale-gal-equiv} is compatible with tensor products and duals, it follows from the proof of Proposition~\ref{prop:etale-realization} (iib) that this map is bijective, which shows the claim.
\end{eg}

Now we state our main theorem.
Let $\mathrm{Rep}_{\mathbf{Z}_p, \geq 0}^{\mathrm{cris}}(\mathcal{G}_R)$ (resp.~
$\mathrm{Rep}_{\mathbf{Z}_p, [0, r]}^{\mathrm{cris}}(\mathcal{G}_R)$) denote the category of crystalline $\mathbf{Z}_p$-representations of $\calG_R$ with non-negative Hodge--Tate weights (resp.~ with Hodge--Tate weights in $[0, r]$). Note that $\Z_p(1)$ has Hodge--Tate weight one by our convention.\\ 

\begin{thm} \label{thm:main}
We keep Assumption~\ref{assumption:base-ring-sec-3.4}.
\begin{enumerate}
\item The \'etale realization $T$ as in Proposition~\ref{prop:etale-realization} gives a fully faithful functor from $\mathrm{CR}^{\wedge,\varphi}(R_{\Prism})$ to $\mathrm{Rep}_{\mathbf{Z}_p, \geq 0}^{\mathrm{cris}}(\mathcal{G}_R)$.
Moreover, $T$ restricts to $\mathrm{CR}^{\wedge,\varphi}_{[0,r]}(R_{\Prism})\rightarrow\mathrm{Rep}_{\mathbf{Z}_p, [0,r]}^{\mathrm{cris}}(\mathcal{G}_R)$.
\item  The functor $T$ gives an equivalence $\mathrm{CR}^{\wedge,\varphi}_{[0,r]}(R_{\Prism})\cong\mathrm{Rep}_{\mathbf{Z}_p, [0,r]}^{\mathrm{cris}}(\mathcal{G}_R)$, which is functorial in $R$.
\end{enumerate}
\end{thm}

\begin{rem}\label{rem:equivalence with all crystalline reps}
Note that for every crystalline $\Z_p$-representation $T_0$ of $\calG_R$, there exists $n\in\Z$ such that $T_0\otimes_{\Z_p}\Z_p(n)\in \mathrm{Rep}_{\Z_p,\geq 0}^\cris(\calG_R)$, and that the \'etale realization functor $T$ is compatible with Breuil--Kisin twists by Example~\ref{eg:Breuil-Kisin twists}. Hence as in \cite[\S 1.2]{Kisin-Sh}, one can extend the definition of completed prismatic $F$-crystals in a way that the resulting category is equivalent to $\mathrm{Rep}_{\Z_p}^\cris(\calG_R)$, the category of $\Z_p$-crystalline representations of $\calG_R$. We leave it to the reader to make a precise formulation.
\end{rem}

The functor $T$ is functorial in $R$ since so are the \'etale realization for Laurent $F$-crystals and the functor $\mathcal{F}\mapsto \mathcal{F}_{\mathrm{\acute{e}t}}$. We prove the first part here. The essential surjectivity in the second part will be proved in the next section. 

\begin{proof}[Proof of Theorem~\ref{thm:main} (i)]
Let $\mathcal{F} \in \mathrm{CR}^{\wedge,\varphi}_{[0, r]}(R_{\Prism})$. Consider the map $R \rightarrow \mathcal{O}_{K_g}$ as in Notation~\ref{notation:L}. By Remarks~\ref{rem:restriction-of-completed-F-crystals} and \ref{rem:vector-bundle-CDVF-case}, we have $\mathcal{F} |_{(\mathcal{O}_{K_g})_{\Prism}} \in \mathrm{CR}^{\wedge,\varphi}_{[0, r]}((\mathcal{O}_{K_g})_{\Prism}) = \mathrm{Vect}^{\varphi}_{[0, r]}((\mathcal{O}_{K_g})_{\Prism})$. Thus, by \cite[Prop.~5.3]{bhatt-scholze-prismaticFcrystal} (see also \cite[Thm.~4.1.10]{du-liu-prismaticphiGhatmodule}), we have $T(\mathcal{F}|_{(\mathcal{O}_{K_g})_{\Prism}}) \in \mathrm{Rep}_{\mathbf{Z}_p, [0, r]}^{\mathrm{cris}}(G_{K_g})$ where $G_{K_g} \coloneqq \mathcal{G}_{\mathcal{O}_{K_g}}$. Note that by Proposition~\ref{prop:etale-realization} (iic), $T(\mathcal{F}|_{(\mathcal{O}_{K_g})_{\Prism}})$ is equal to $T(\mathcal{F}) |_{G_{K_g}}$.

We first show that the essential image of $T$ is contained in $\mathrm{Rep}_{\mathbf{Z}_p, \geq 0}^{\mathrm{cris}}(\mathcal{G}_R)$. Let $V(\mathcal{F}) = T(\mathcal{F})[p^{-1}]$ 
denote the corresponding $\mathbf{Q}_p$-representation of $\mathcal{G}_R$. By Proposition~\ref{prop:etale-realization} (iia), we see
\[
V(\mathcal{F})^{\vee} \cong ( W(\overline{R}^{\flat}[(\pi^{\flat})^{-1}])[p^{-1}]\otimes_{\mathbf{A}_{\mathrm{inf}}(\overline{R})}\mathcal{F}_{\mathbf{A}_{\mathrm{inf}}(\overline{R})} )^{\varphi = 1}.
\]
By Lemma~\ref{lem:completed-crystals-basic-properties} (iv), we have $\mathcal{F}_{\mathbf{A}_{\mathrm{inf}}(\overline{R})}[p^{-1}] \cong \mathbf{A}_{\mathrm{inf}}(\overline{R})[p^{-1}]\otimes_{\mathfrak{S}[p^{-1}]}\mathcal{F}_{\mathfrak{S}}[p^{-1}]$,  which is finite projective over $\mathbf{A}_{\mathrm{inf}}(\overline{R})[p^{-1}]$. Since $\mathcal{F}_{\et}(\mathbf{A}_{\mathrm{inf}}(\overline{R}), (\xi))$ is an \'etale $\varphi$-module finite projective over $W(\overline{R}^{\flat}[(\pi^{\flat})^{-1}])$, we obtain  
\begin{equation} \label{eq:BKF-mod}
 W(\overline{R}^{\flat}[(\pi^{\flat})^{-1}])[p^{-1}]\otimes_{\mathbf{Q}_p}V(\mathcal{F})^{\vee} \cong  W(\overline{R}^{\flat}[(\pi^{\flat})^{-1}])[p^{-1}]\otimes_{\mathbf{A}_{\mathrm{inf}}(\overline{R})[p^{-1}]}\mathcal{F}_{\mathbf{A}_{\mathrm{inf}}(\overline{R})}[p^{-1}].
\end{equation}

Consider each $\mathfrak{p} \in \mathcal{P}$ and $\overline{R}^{\wedge}\rightarrow (\overline{R}_{\mathfrak{p}})^{\wedge} \rightarrow (\calO_{\overline{K_g}})^{\wedge}$ as in Notation \ref{notation:L}. Equation~\eqref{eq:BKF-mod} induces 
\[
 W(\mathcal{O}_{\overline{K_g}}^{\flat}[(\pi^{\flat})^{-1}])[p^{-1}]\otimes_{\mathbf{Q}_p}V(\mathcal{F})^{\vee} = W(\mathcal{O}_{\overline{K_g}}^{\flat}[(\pi^{\flat})^{-1}])[p^{-1}]\otimes_{\mathbf{A}_{\mathrm{inf}}(\overline{R})[p^{-1}]} \mathcal{F}_{\mathbf{A}_{\mathrm{inf}}(\overline{R})}[p^{-1}]
\]
by the base change along $W(\overline{R}^{\flat}[(\pi^{\flat})^{-1}]) \rightarrow W(\mathcal{O}_{\overline{K_g}}^{\flat}[(\pi^{\flat})^{-1}])$. Set $T_g\coloneqq T(\mathcal{F} |_{(\mathcal{O}_{K_g})_{\Prism}})$. Since $T_g \in \mathrm{Rep}_{\mathbf{Z}_p, [0, r]}^{\mathrm{cris}}(G_{K_g})$, the proof of \cite[Lem.~4.26]{bhatt-morrow-scholze-integralpadic} shows $ \mathcal{F}_{\Ainf(\mathcal{O}_{\overline{K_g}})} \subset T_g^{\vee}\otimes_{\Z_p} \Ainf (\calO_{\overline{K _g}})$ and $T_g^\vee \subset \mathcal{F}_{\Ainf(\mathcal{O}_{\overline{K_g}})}\otimes_{\Ainf(\mathcal{O}_{\overline{K_g}})} \Ainf(\mathcal{O}_{\overline{K_g}})\{r\}$. Here $\Ainf(\mathcal{O}_{\overline{K_g}})\{r\}$ denotes the base change of the $r$-th Tate twist $\mu^{-r}\Ainf(\mathcal{O}_{\overline{K_g}})\otimes_{\mathbf{Z}_p} \mathbf{Z}_p(r)$, defined in \cite[Ex.~4.24]{bhatt-morrow-scholze-integralpadic} with $\mu \coloneqq [\varepsilon] - 1$, along the Frobenius inverse $\varphi^{-1}\colon \Ainf(\mathcal{O}_{\overline{K_g}})\rightarrow \Ainf(\mathcal{O}_{\overline{K_g}})$; note that our functor $T$ is contravariant and the cokernel of $1\otimes \varphi_{\mathcal{F}_{\Ainf(\mathcal{O}_{\overline{K_g}})}}$ is supported on $(E)$, while \cite[Def.~4.22]{bhatt-morrow-scholze-integralpadic} uses $(\varphi(E))$. Since $\varphi^{-1}(\mu)$ divides $\mu$, we obtain inclusions of $W(\mathcal{O}_{\overline{K_g}}^{\flat})[p^{-1}]$-modules
\begin{equation} \label{eq:BKF-mod-DVR}
W(\mathcal{O}_{\overline{K_g}}^{\flat})[p^{-1}]\otimes_{\mathbf{A}_{\mathrm{inf}}(\overline{R})[p^{-1}]} \mathcal{F}_{\mathbf{A}_{\mathrm{inf}}(\overline{R})}[p^{-1}] \subset  W(\mathcal{O}_{\overline{K_g}}^{\flat})[p^{-1}]\otimes_{\mathbf{Q}_p}V(\mathcal{F})^{\vee} \subset \frac{1}{\mu ^r}W(\mathcal{O}_{\overline{K_g}}^{\flat})[p^{-1}]\otimes_{\mathbf{A}_{\mathrm{inf}}(\overline{R})[p^{-1}]} \mathcal{F}_{\mathbf{A}_{\mathrm{inf}}(\overline{R})}[p^{-1}].
\end{equation}
Since the $\mathbf{A}_{\mathrm{inf}}(\overline{R})[p^{-1}]$-module $\mathcal{F}_{\mathbf{A}_{\mathrm{inf}}(\overline{R})}[p^{-1}]$ is finite projective, \eqref{eq:BKF-mod} and \eqref{eq:BKF-mod-DVR} together with Lemma~\ref{lem:intersection-modules-flat-base-change} (ii) and Lemma~\ref{lem:intersection-witt-rings} below yield
\begin{equation} \label{eq:gal-rep-BKF-mod}
 \mathbf{A}_{\mathrm{inf}}(\overline{R})[p^{-1}][\mu^{-1}]\otimes_{\mathbf{Q}_p}V(\mathcal{F})^{\vee} = \mathcal{F}_{\mathbf{A}_{\mathrm{inf}}(\overline{R})}[p^{-1}][\mu^{-1}].
\end{equation}

Consider the map of prisms $(\mathfrak{S}, (E)) \rightarrow (R_0, (p))$ over $R$ given by $u \mapsto 0$. Let $D(\mathcal{F}) \coloneqq \mathcal{F}_{R_0}[p^{-1}]$. We have $D(\mathcal{F}) \cong R_0[p^{-1}]\otimes_{\mathfrak{S}} \mathcal{F}_{\mathfrak{S}}$ by Lemma~\ref{lem:completed-crystals-basic-properties} (iv), and $1\otimes \varphi\colon \varphi^*D(\mathcal{F}) \rightarrow D(\mathcal{F})$ is an isomorphism. Choose a positive integer $l$ with $p^l \geq e$ as in Example~\ref{eg:prism-OAcris} so that we have the map of prisms $(R_0, (p)) \xrightarrow{\varphi_{l+1}} (\phi_{l}\OA_{\mathrm{cris}}(\overline{R}), (p))$ over $R$. Thus,
\begin{align*}
\OA_{\mathrm{cris}}(\overline{R})[p^{-1}]\otimes_{\varphi^{l+1},R_0}D(\mathcal{F})  &\cong \mathcal{F}(\phi_{l}\OA_{\mathrm{cris}}(\overline{R}), (p))[p^{-1}] \\
    &\cong \OA_{\mathrm{cris}}(\overline{R})[p^{-1}] \otimes_{\varphi^l,\OA_{\mathrm{cris}}(\overline{R})} \mathcal{F}(\OA_{\mathrm{cris}}(\overline{R}), (p))
\end{align*}
by Lemma~\ref{lem:completed-crystals-basic-properties} (iv), and we obtain
\begin{equation} \label{eq:comparison-D}
\OB_{\mathrm{cris}}(\overline{R})\otimes_{\varphi^{l+1},R_0}D(\mathcal{F})  \cong  \OB_{\mathrm{cris}}(\overline{R})\otimes_{\varphi^l,\OA_{\mathrm{cris}}(\overline{R})} \mathcal{F}(\OA_{\mathrm{cris}}(\overline{R}), (p)).
\end{equation}

On the other hand, again by Lemma~\ref{lem:completed-crystals-basic-properties} (iv),
\[
\mathcal{F}(\phi_{l}\OA_{\mathrm{cris}}(\overline{R}), (p))[p^{-1}] \cong  \OA_{\mathrm{cris}}(\overline{R})[p^{-1}]\otimes_{\varphi^{l+1},\mathbf{A}_{\mathrm{inf}}(\overline{R})}\mathcal{F}_{\mathbf{A}_{\mathrm{inf}}(\overline{R})}[p^{-1}].  
\]
So Equation \eqref{eq:gal-rep-BKF-mod} gives
\begin{equation} \label{eq:comparison-V}
\OB_{\mathrm{cris}}(\overline{R})\otimes_{\mathbf{Q}_p}V(\mathcal{F})^{\vee} \cong  \OB_{\mathrm{cris}}(\overline{R})\otimes_{\varphi^l,\OA_{\mathrm{cris}}(\overline{R})} \mathcal{F}(\OA_{\mathrm{cris}}(\overline{R}), (p)).     
\end{equation}
By Equations \eqref{eq:comparison-D}, \eqref{eq:comparison-V}, and  $\OB_{\mathrm{cris}}(\overline{R})\otimes_{ \varphi^{l+1},R_0} D(\mathcal{F}) \cong \OB_{\mathrm{cris}}(\overline{R})\otimes_{R_0} D(\mathcal{F})$ obtained by $(l+1)$-times iterations of the isomorphism $1\otimes \varphi$, we deduce the isomorphism
\begin{equation}\label{eq:crystalline comparison in proof of main thm}
\OB_{\mathrm{cris}}(\overline{R})\otimes_{R_0} D(\mathcal{F}) \cong  \OB_{\mathrm{cris}}(\overline{R})\otimes_{\mathbf{Q}_p}V(\mathcal{F})^{\vee}
\end{equation}
that is compatible with $\mathcal{G}_R$-actions and $\varphi$. 
Since $(\OB_\cris(\overline{R}))^{\calG_R}=R_0[p^{-1}]$, we deduce from this isomorphism that $ \OB_{\mathrm{cris}}(\overline{R})\otimes_{\mathbf{Q}_p}V(\mathcal{F})^{\vee}$ is spanned by its $\calG_R$-invariants as an $\OB_{\cris}(\overline{R})$-module.
It follows that $\alpha_{\mathrm{cris}}(V(\mathcal{F})^{\vee})$ is surjective and thus $V(\mathcal{F})^{\vee}$ is crystalline. So $V(\mathcal{F})$ is crystalline. Note that $V(\mathcal{F})$ has Hodge--Tate weights in $[0, r]$, since it has Hodge--Tate weights in $[0, r]$ considered as a representation of $G_{K_g}$. 

The faithfulness follows from the construction of the \'etale realization $T$ in Proposition~\ref{prop:etale-realization}. For the fullness, let $\mathcal{F}_1, \mathcal{F}_2 \in \mathrm{CR}^{\wedge,\varphi}(R_{\Prism})$, and suppose we have a map $h\colon T(\mathcal{F}_1) \rightarrow T(\mathcal{F}_2)$ of representations of $\mathcal{G}_R$. By Proposition~\ref{prop:etale-realization} (iib), $h|_{\mathcal{G}_{\tilde{R}_{\infty}}}$ induces a map 
\[
(\mathcal{F}_2)_{\mathfrak{S}}[E^{-1}]^{\wedge}_p \rightarrow (\mathcal{F}_1)_{\mathfrak{S}}[E^{-1}]^{\wedge}_p
\]
of \'etale $\varphi$-modules over $\mathcal{O}_{\mathcal{E}}$. On the other hand, by Lemma~\ref{lem:completed-crystals-basic-properties} (i) and \cite[Prop.~4.2.7]{gao-integral-padic-hodge-imperfect}, $h|_{\mathcal{G}_{\widetilde{\mathcal{O}}_{L, \infty}}}$ induces a $\varphi$-equivariant map
\[
(\mathcal{F}_{2})_{\mathfrak{S}_L} \rightarrow (\mathcal{F}_{1})_{\mathfrak{S}_L}
\]
of $\mathfrak{S}_L$-modules. These two maps are compatible after the base changes to $\mathcal{O}_{\mathcal{E}, L}$, and thus we obtain an induced $\varphi$-equivariant map of $\mathfrak{S}$-modules
\[
(\mathcal{F}_2)_{\mathfrak{S}}[E^{-1}]^{\wedge}_p \cap (\mathcal{F}_{2})_{\mathfrak{S}_L} \rightarrow (\mathcal{F}_1)_{\mathfrak{S}}[E^{-1}]^{\wedge}_p \cap (\mathcal{F}_{1})_{\mathfrak{S}_L},
\]
i.e., a map $f\colon (\mathcal{F}_2)_{\mathfrak{S}} \rightarrow (\mathcal{F}_1)_{\mathfrak{S}}$ by Lemma~\ref{lem:completed-crystals-basic-properties} (ii). 

By the construction of the \'etale realization, $f$ is compatible with the descent data
\[
\mathfrak{S}^{(1)}[E^{-1}]^{\wedge}_p\otimes_{p_1,\mathcal{O}_{\mathcal{E}}} (\mathcal{F}_i)_{\mathfrak{S}}[E^{-1}]^{\wedge}_p \xrightarrow{\cong} \mathfrak{S}^{(1)}[E^{-1}]^{\wedge}_p\otimes_{p_2,\mathcal{O}_{\mathcal{E}}}(\mathcal{F}_i)_{\mathfrak{S}}[E^{-1}]^{\wedge}_p
\]
for $i = 1, 2$. So by Lemma~\ref{lem:completed-crystals-basic-properties} (iii), the map $f\colon (\mathcal{F}_2)_{\mathfrak{S}} \rightarrow (\mathcal{F}_1)_{\mathfrak{S}}$ is compatible with the descent data
\[
\mathfrak{S}^{(1)}\otimes_{p_1,\mathfrak{S}}(\mathcal{F}_i)_{\mathfrak{S}}  \xrightarrow{\cong} \mathfrak{S}^{(1)}\otimes_{p_2,\mathfrak{S}}(\mathcal{F}_i)_{\mathfrak{S}}
\]
for $i = 1, 2$. Thus, the fullness follows from Proposition~\ref{prop:equivalence-to-descent-datum}. 
\end{proof}

\begin{rem}\label{rem:how to recover Dcris-1}
For $\calF\in\mathrm{CR}^{\wedge,\varphi}(R_\Prism)$, 
the isomorphism \eqref{eq:crystalline comparison in proof of main thm} in the above proof shows that there is an isomorphism $\calF(R_0,(p))[p^{-1}]\cong D_\cris^\vee(T(\calF)[p^{-1}])$ as $\varphi$-modules over $R_0[p^{-1}]$.
Since $\varphi$ is an isomorphism on $D_\cris^\vee(T(\calF)[p^{-1}])$, Lemma~\ref{lem:completed-crystals-basic-properties} (iv) for the map of prisms $(\fkS,E)\xrightarrow{u\mapsto 0}(R_0,(p)) \xrightarrow{\varphi}(R_0,(p)) $ gives
a $\varphi$-equivalent $R_0[p^{-1}]$-linear isomorphism
\[
(R_0\otimes_{\varphi,R_0}\calF_\fkS/u\calF_\fkS)[p^{-1}]\cong D_\cris^\vee(T(\calF)[p^{-1}]).
\]
In Remark~\ref{rem:how to recover Dcris-2}, we will explain how to obtain the connection on $D_\cris^\vee(T(\calF)[p^{-1}])$ and the filtration on $R\otimes_{R_0}D_\cris^\vee(T(\calF)[p^{-1}])$ from $\calF$ under the above isomorphism. 
\end{rem}

We used the following lemma in the proof of Theorem~\ref{thm:main} (i).

\begin{lem} \label{lem:intersection-witt-rings}
As subrings of $\prod_{\mathfrak{p} \in \mathcal{P}} W(\overline{K_g}^{\flat})$, 
\[
W(\overline{R}^\flat[(\pi^{\flat})^{-1}]) \cap \prod_{\mathfrak{p} \in \mathcal{P}} W(\mathcal{O}_{\overline{K_g}}^{\flat}) = \A_\mathrm{inf}(\overline{R}).
\]
Furthermore, we have
\[
W(\overline{R}^\flat[(\pi^{\flat})^{-1}])[p^{-1}] \cap \prod_{\mathfrak{p} \in \mathcal{P}} W(\mathcal{O}_{\overline{K_g}}^{\flat})[p^{-1}] = \A_\mathrm{inf}(\overline{R})[p^{-1}].
\]
\end{lem}

\begin{proof}
Recall that $\A_\mathrm{inf}(\overline{R})$ denotes $W(\overline{R}^\flat)$. For any $x \in \overline{R}^\flat[(\pi^{\flat})^{-1}]$, if its $\pi^{\flat}$-adic valuation as element in $\overline{K_g}^{\flat}$ is nonnegative for all $\mathfrak{p} \in \mathcal{P}$, then $x \in \overline{R}^\flat$. So $\overline{R}^\flat[(\pi^{\flat})^{-1}] \cap \prod_{\mathfrak{p} \in \mathcal{P}}\mathcal{O}_{\overline{K_g}}^{\flat} = \overline{R}^\flat$ as subrings of $(\prod_{\mathfrak{p} \in \mathcal{P}}\mathcal{O}_{\overline{K_g}}^{\flat})[(\pi^{\flat})^{-1}]$. By considering the Teichm\"{u}ller expansion of $p$-typical Witt vectors, we deduce both statements.
\end{proof}

\subsection{Height one case} \label{subsec:height one}

This subsection discusses the case where crystalline representations have Hodge--Tate weights in $[0, 1]$, and studies the relation to $p$-divisible groups. We keep Assumption~\ref{assumption:base-ring-sec-3.4}: $R$ is small over $\mathcal{O}_K$ or $R = \mathcal{O}_L$. We first recall a main result in \cite{anschutz-lebras-prismaticdieudonne} on classifying $p$-divisible groups over $R$ via prismatic $F$-crystals on $R$. Let $\mathrm{BT}(R)$ denote the category of $p$-divisible groups over $R$

For a $p$-complete $R$-algebra with bounded $p^{\infty}$-torsion, let
$R_{\mathrm{QSYN}}$ denote the big quasi-syntomic site of $R$ (cf.~\cite[\S 3.3]{anschutz-lebras-prismaticdieudonne}, \cite[\S 4]{bhatt-morrow-scholze-tophochschildhom}). By \cite[Cor.~3.24]{anschutz-lebras-prismaticdieudonne}, the functor $R_{\Prism} \rightarrow R_{\mathrm{QSYN}}$ sending $(A, I)$ to $R \rightarrow A/I$ is cocontinuous, so it defines a morphism of topoi 
\[
u\colon \mathrm{Shv}(R_{\Prism}) \rightarrow \mathrm{Shv}(R_{\mathrm{QSYN}}).
\]
For a $p$-divisible group $H$ over $R$, we consider the sheaf $\mathfrak{M}_{\Prism}(H) \coloneqq \mathcal{E}\mathit{xt}^1_{R_{\Prism}}(u^{-1}H, \mathcal{O}_{\Prism})$ on $R_{\Prism}$. Let $\varphi_{\mathfrak{M}_{\Prism}(H)}$ be the endomorphism of $\mathfrak{M}_{\Prism}(H)$ induced from $\varphi$ on $\mathcal{O}_{\Prism}$. The following is proved in \cite{anschutz-lebras-prismaticdieudonne}.

\begin{thm}[(Ansch\"utz--Le Bras)] \label{thm:p-div-gp-classification}
The assignment 
\[
H \mapsto (\mathfrak{M}_{\Prism}(H), \varphi_{\mathfrak{M}_{\Prism}(H)})
\]	
gives an equivalence of categories from $\mathrm{BT}(R)$ to $\mathrm{Vect}^{\varphi}_{[0,1]}(R_{\Prism})$.
\end{thm}

\begin{proof}
Let $H \in \mathrm{BT}(R)$. By \cite[Thm.~4.71, Lem.~4.38]{anschutz-lebras-prismaticdieudonne}, we see that $(\mathfrak{M}_{\Prism}(H), \varphi_{\mathfrak{M}_{\Prism}(H)})$ is an object in $\mathrm{Vect}^{\varphi}_{[0,1]}(R_{\Prism})$. So the assignment defines a functor from $\mathrm{BT}(R)$ to $\mathrm{Vect}^{\varphi}_{[0,1]}(R_{\Prism})$. From \cite[Thm.~4.74, Prop.~5.10]{anschutz-lebras-prismaticdieudonne}, we deduce that this is an equivalence: note that the proof of \cite[Prop.~5.10]{anschutz-lebras-prismaticdieudonne} uses the existence of a quasi-syntomic cover $R_{\infty}$ of $R$, which is constructed in the proof of \cite[Prop.~5.8]{anschutz-lebras-prismaticdieudonne}. In our case, $R = R_0[\![u]\!]/(E(u))$ with $R_0$ unramified, so we can do a similar construction by extracting $p$-power roots of $u$. 
\end{proof}

\begin{rem} \label{rem:equiv-holds-quasi-syntomic-regular-ring}
In \cite{anschutz-lebras-prismaticdieudonne}, an equivalence between $\mathrm{BT}(R)$ and the category of admissible prismatic Dieudonn\'{e} crystals over $R$ is proved for any quasi-syntomic ring $R$.
\end{rem}

For $H \in \mathrm{BT}(R)$, we write $T_p(H)$ for its Tate module. Note that we have a natural $\mathcal{G}_R$-equivariant isomorphism
\[
T_p(H) \cong \mathrm{Hom}_{\mathrm{BT}(\overline{R}^{\wedge})} ((\mathbf{Q}_p/\mathbf{Z}_p)_{\overline{R}^{\wedge}}, H_{\overline{R}^{\wedge}}).
\]  

\begin{prop} \label{prop:Tate-mod}
There exists a natural $\mathcal{G}_R$-equivariant isomorphism 
\[
T_p(H) \cong T(\mathfrak{M}_{\Prism}(H)),
\]	
where $T(\mathfrak{M}_{\Prism}(H))$ is the \'etale realization of $\mathfrak{M}_{\Prism}(H) \in \mathrm{Vect}^{\varphi}_{[0,1]}(R_{\Prism}) \subset \mathrm{CR}^{\wedge,\varphi}(R_{\Prism})$ as in Proposition~\ref{prop:etale-realization}. 
\end{prop}

\begin{proof}
Since $\overline{R}^{\wedge}$ is an integral perfectoid ring, the prism $(\mathbf{A}_{\mathrm{inf}}(\overline{R}), (\xi))$ is the final object of  $(\overline{R}^{\wedge})_{\Prism}$. Let $M = \mathfrak{M}_{\Prism}(H)(\mathbf{A}_{\mathrm{inf}}(\overline{R}), (\xi))$. By \cite[Prop.~1.39]{morrowtsujigeneralised} (where the covariant version of $\mathfrak{M}_{\Prism}(\cdot)$ is used), we have a $\mathcal{G}_R$-equivariant isomorphism
\[
T_p(H) \cong (M^{\vee})^{\varphi = 1}.
\]

We claim that the natural injective map
\[
(M^{\vee})^{\varphi = 1} \rightarrow (M^{\vee}\otimes_{\mathbf{A}_{\mathrm{inf}}(\overline{R})} W(\overline{R}^{\flat}[(\pi^{\flat})^{-1}]))^{\varphi = 1} 
\]
is also surjective. Note that the natural map $(M^{\vee})^{\varphi=1}/p \rightarrow (M^{\vee}/p)^{\varphi=1}$ is injective, and $(M^{\vee}/p)^{\varphi=1} \subset ((M^{\vee}/p)[(\pi^{\flat})^{-1}])^{\varphi=1}$. On the other hand, by Proposition~\ref{prop:etale-gal-equiv}, we have
\[
\mathrm{dim}_{\mathbf{F}_p} ((M^{\vee}/p)[(\pi^{\flat})^{-1}])^{\varphi=1} = \mathrm{rank}_{\mathbf{A}_{\mathrm{inf}}(\overline{R})} M = \mathrm{rank}_{\mathbf{Z}_p}T_p(H).
\]
Thus, we deduce that the map $(M^{\vee})^{\varphi=1}/p \rightarrow ((M^{\vee}/p)[(\pi^{\flat})^{-1}])^{\varphi=1}$ is bijective, and the claim follows. Since $(M^{\vee}\otimes_{\mathbf{A}_{\mathrm{inf}}(\overline{R})} W(\overline{R}^{\flat}[(\pi^{\flat})^{-1}]))^{\varphi = 1} \cong ((M\otimes_{\mathbf{A}_{\mathrm{inf}}(\overline{R})} W(\overline{R}^{\flat}[(\pi^{\flat})^{-1}]))^{\varphi = 1})^{\vee}$, it follows from the definition of the \'etale realization functor that $T_p(H) \cong T(\mathfrak{M}_{\Prism}(H))$.
\end{proof}

Based on an example in \cite[\S~5.4]{vasiu-zink-purity}, we now present an example of a crystalline representation with Hodge--Tate weights in $[0, 1]$ that does not come from a $p$-divisible group. By Theorems~\ref{thm:main}, \ref{thm:p-div-gp-classification}, and Proposition~\ref{prop:Tate-mod}, such an example implies that the inclusion $\mathrm{Vect}^\varphi_{\mathrm{eff}}(R_\Prism) \subset  \mathrm{CR}^{\wedge,\varphi}(R_{\Prism})$ is strict in general.

\begin{eg}\label{eg:non-p-divisible crystalline rep}
Let $R_0= W(k)\langle T^{\pm 1}\rangle$ and $R = R_0\otimes_{W(k)}\mathcal{O}_K$. Suppose $p \geq 3$ and the ramification index $[K : K_0]$ is $p$. Let $\mathfrak{M}_1$ be a free $\mathfrak{S}/p\mathfrak{S}$-module with a basis $\{e_1, e_2, e_3\}$, equipped with Frobenius given by
\[
\varphi = \begin{pmatrix} (1+T)^{p-1}-u & 0 & u\\ 0 & u^{p-1} & (1+T)^{p}u^{p-1}-(1+T)\\ 1 & 0 & 0 \end{pmatrix}
\]
and with the trivial connection $\nabla(e_i) = 0$ for $i = 1, 2, 3$. Note that for 
\[
\psi = \begin{pmatrix} 0 & 0 & u^p\\ (1+T)-(1+T)^pu^{p-1} & u & (u-(1+T)^{p-1})((1+T)-(1+T)^pu^{p-1}) \\ u^{p-1} & 0 & u^{p-1}(u-(1+T)^{p-1})\end{pmatrix},
\]
we have $\varphi\psi = \psi\varphi = u^p I_3$. Thus $\mathfrak{M}_1 \in (\mathrm{Mod  ~FI})_{\mathfrak{S}}^{\mathrm{Ki}}(\varphi, \nabla^0)$ in the sense of \cite[Def.~9.2]{kim-groupscheme-relative}. By \cite[Thm.~9.8]{kim-groupscheme-relative}, $\mathfrak{M}_1$ is associated with a finite flat group scheme $H_1$ over $R$.

Let $\mathfrak{M}_2$ be a free $\mathfrak{S}/p\mathfrak{S}$-module with a basis $\{f\}$, equipped with Frobenius given by $\varphi(f) = f$ and with the trivial connection $\nabla(f) = 0$. Then $\mathfrak{M}_2 \in (\mathrm{Mod  ~FI})_{\mathfrak{S}}^{\mathrm{Ki}}(\varphi, \nabla^0)$ and it is associated with a finite flat group scheme $H_2$ over $R$.

Let $h\colon \mathfrak{M}_1 \rightarrow \mathfrak{M}_2$ be a map of torsion Kisin modules given by $\begin{pmatrix} 1+T & u & (1+T)u\end{pmatrix}$. Since $h$ is not surjective, the associated map $H_2 \rightarrow H_1$ of finite flat group schemes is not a monomorphism. On the other hand, the induced maps $H_2[p^{-1}] \rightarrow H_1[p^{-1}]$ and $H_2\times_R \mathcal{O}_L \rightarrow H_1\times_R \mathcal{O}_L$ are monomorphisms of finite flat group schemes over $R[p^{-1}]$ and $\mathcal{O}_L$, respectively.

By \cite[Thm.~3.1.1]{berthelot-breen-messing}, there exists $a \in R$ with $a \notin (\pi, 1+T) \subset R$ such that $(H_1)_{R'}$ can be embedded into some $p$-divisible group $H$ over $R'$, the $p$-adic completion of $R[a^{-1}]$. Consider the $p$-divisible group $H'$ over $R'[p^{-1}]$ given by
\[
H' \coloneqq H_{R'[p^{-1}]} / (H_2)_{R'[p^{-1}]},
\] 
and let $V$ be the associated representation of $\mathcal{G}_{R'}$. Since $V \cong T_p(H)[p^{-1}]$ by construction, $V$ is a crystalline representation of $\mathcal{G}_{R'}$ with Hodge--Tate weights in $[0, 1]$.

On the other hand, we claim that $H'$ cannot be extended to a $p$-divisible group over $R'$. Suppose otherwise, i.e., suppose that $H'$ extends to a $p$-divisible group $H'_{R'}$ over $R'$. Let $R_1$ be the $(\pi, 1+T)$-adic completion of the localization $R'_{(\pi, 1+T)}$. Let $\mathfrak{m} \in \Spec R_1$ be the closed point, and let $U \coloneqq \Spec R_1 - \mathfrak{m}$ be the open subscheme of $\Spec R_1$. Note that by the construction of $h\colon \mathfrak{M}_1 \rightarrow \mathfrak{M}_2$ above, the induced map $(H_2)_{R_1} \rightarrow (H_1)_{R_1}$ is not a monomorphism whereas the restriction $(H_2)_U \rightarrow (H_1)_U$ to $U$ is a monomorphism. In particular, by \cite[Thm.~ 4]{tate}, we have $H'_{R'}\times_{R'} U = (H\times_{R'} U) / (H_2)_U$ as $p$-divisible groups over $U$. The isogeny $H\times_{R'} U \rightarrow H'_{R'}\times_{R'} U$ extends to an isogeny $i\colon H\times_{R'} R_1 \rightarrow H'_{R'}\times_{R'} R_1$. The kernel of $i$ is a finite flat group scheme over $R_1$ whose restriction to $U$ is $(H_2)_U$. Since $R_1$ is a regular local ring of dimension $2$, the kernel of $i$ is then equal to $(H_2)_{R_1}$, and $(H_2)_{R_1}$ embeds into $H\times_{R'} R_1$. This contradicts that $(H_2)_{R_1} \rightarrow (H_1)_{R_1}$ is not a monomorphism, so $H$ cannot be extended to a $p$-divisible group over $R'$.
\end{eg}

In the above example, the ramification index $e \coloneqq [K : K_0]$ needs to be large. In fact, when $e<p-1$, such an example does not exist.

\begin{thm}[({\cite[Thm.~1.2]{liu-moon-rel-crys-rep-p-div-gps-small-ramification}})] \label{thm:liu-moon-small-ramification}
Suppose $e < p-1$ (so that $p \geq 3$). Assume moreover that $R_0/pR_0$ is a UFD, and that $R_0$ is complete with respect to some ideal $J \subset R_0$ containing $p$ such that $R_0/J$ is finitely generated over some field. Then for any $T \in \mathrm{Rep}_{\mathbf{Z}_p, [0, 1]}^{\mathrm{cris}}(\mathcal{G}_R)$, there exists a $p$-divisible group $G$ over $R$ such that $T_p(G) \cong T$ as representations of $\mathcal{G}_R$.
\end{thm}

In particular, we have $\mathrm{Vect}^\varphi_{[0, 1]}(R_\Prism) = \mathrm{CR}^{\wedge,\varphi}_{[0, 1]}(R_{\Prism})$ under the assumptions of the above theorem.

\begin{rem}
 In fact, we have a little stronger result: $\mathrm{Vect}^\varphi_{[0, 1]}(R_\Prism) = \mathrm{CR}^{\wedge,\varphi}_{[0, 1]}(R_{\Prism})$ if $e < p-1$ and $R_0$ is small over $\calO_K$.  To see this, we use arguments in our proof of Theorem~\ref{thm:main} (ii) (the essential surjectivity) in \S~\ref{sec:quasi-kisin-mods-cryst-loc-systs}. More precisely, for $T \in \mathrm{Rep}_{\mathbf{Z}_p, [0, 1]}^{\mathrm{cris}}(\mathcal{G}_R)$, the associated completed prismatic $F$-crystal $\mathcal{F} \in \mathrm{CR}^{\wedge,\varphi}_{[0, 1]}(R_{\Prism})$ satisfies $\mathcal{F}_{\mathfrak{S}} = \mathfrak{M}$ where $\mathfrak{M}$ is given by the construction in \S~\ref{sec:quasi-kisin-mod-construction}. By Remark~\ref{rem:small-ramification-case}, $\mathcal{F}_{\mathfrak{S}}$ is projective over $\mathfrak{S}$ when $e < p-1$. Hence Proposition~\ref{prop:equivalence-to-descent-datum} implies $\mathrm{Vect}^\varphi_{[0, 1]}(R_\Prism) = \mathrm{CR}^{\wedge,\varphi}_{[0, 1]}(R_{\Prism})$ when $e < p-1$ (even without the additional assumptions on $R_0$ in Theorem~\ref{thm:liu-moon-small-ramification}). By \cite[Rem.~4.6]{liu-moon-rel-crys-rep-p-div-gps-small-ramification}, we can similarly deduce $\mathrm{Vect}^\varphi_{[0, r]}(R_\Prism) = \mathrm{CR}^{\wedge,\varphi}_{[0, r]}(R_{\Prism})$ when $er < p-1$.  In particular, when $r = 0$, we have $\mathrm{Vect}^\varphi_{[0, 0]}(R_\Prism) = \mathrm{CR}^{\wedge,\varphi}_{[0, 0]}(R_{\Prism})$ for any $e$.  \end{rem}

\subsection{Completed prismatic \texorpdfstring{$F$}{F}-crystals on a smooth \texorpdfstring{$p$}{p}-adic formal scheme}\label{sec:globalization}

This subsection globalizes the construction and the main theorem in \S~\ref{sec:etale-realization-main-thm}. Let $\fkX$ be a smooth $p$-adic formal scheme over $\calO_K$.
To define the category $\mathrm{CR}^{\wedge,\varphi}(\fkX_\Prism)$ by gluing, we need to show the descent property of completed prismatic $F$-crystals with respect to Zariski open coverings.

\begin{lem}\label{lem:Zariski descent of completed prismatic crystals}
Let $\fkX=\bigcup_{\lambda\in \Lambda} \Spf R_\lambda$ be an affine open covering of $\fkX$ .
For a sheaf $\calF$ of $\calO_\Prism$-modules on $\fkX_\Prism$, it is a finitely generated completed prismatic crystal on $\fkX$ if and only if $\calF|_{\Spf R_\lambda}$ is a finitely generated completed prismatic crystal on $\Spf R_\lambda$ for every $\lambda$.
\end{lem}

\begin{proof}
The necessity is obvious. 
To show the sufficiency, assume that $\calF|_{\Spf R_\lambda}$ is a finitely generated completed prismatic crystal on $\Spf R_\lambda$ for every $\lambda$. 
Consider the quotient sheaf $\calF_n\coloneqq \calF/(p,\calI_\Prism)^n\calF=\calO_{\Prism,n}\otimes_{\calO_\Prism}\calF$ on $\fkX_\Prism$ for each $n\in\N$. Then $(\calF_n)_n$ forms an inverse system of $\calO_\Prism$-modules and the natural morphism $\calF\rightarrow \varprojlim_n\calF_n$ is an isomorphism since it is so on $(\Spf R_\lambda)_\Prism$ for each $\lambda$.
By Lemma~\ref{lem:calF vs (calF_n) for completed crystals}, it is enough to show that $\calF_n$ is a finitely generated crystal of $\calO_{\Prism,n}$-modules for every $n\in\N$.

Take any $(A,I)\in\fkX_\Prism$.
Then there exist a finite affine open covering $\Spf A/I=\bigcup_{j=1}^l \Spf R_j$ and an element $\lambda_j\in\Lambda$ for each $j=1,\ldots,l$ such that the map $\Spf R_j\rightarrow \Spf A/I\rightarrow \fkX$ factors through $\Spf R_{\lambda_j}\subset \fkX$.
Since $A/I\rightarrow R_j$ is $p$-completely \'etale map, it lifts uniquely to a $(p,I)$-completely \'etale map $A\rightarrow A_j$ of $\delta$-rings (cf.~\cite[Construction~4.4]{bhatt-scholze-prismaticcohom}) and defines $(A_j,IA_j)\in (\Spf R_{\lambda_i})_\Prism \subset\fkX_\Prism$.
Set $B\coloneqq \prod_{j=1}^l A_j$. Then $B$ admits a natural $\delta$-structure and  $(B,IB)\in \fkX_\Prism$.
Moreover, $(A,I)\rightarrow (B,IB)$ is $(p,I)$-completely faithfully flat.
Let $(B',IB')$ be the object of $\fkX_\Prism$ corresponding to the pushout of the diagram $(B,IB)\leftarrow (A,IA)\rightarrow (B,IB)$ of maps of bounded prisms over $\fkX$. Note $B'/(p, I)^nB'=B/(p, I)^nB\otimes_{A/(p, I)^n}B/(p, I)^nB$ by Lemma~\ref{lem:pushout for prisms along fflat map}. Let $p_1$ and $p_2$ denote the two structure maps $B\rightarrow B'$.

Since $\calF_n$ is a sheaf on $\fkX_\Prism$, we have an exact sequence
\[
0\rightarrow \calF_n(A,I)\rightarrow \calF_n(B,IB)\xrightarrow{p_1^\ast-p_2^\ast}\calF_n(B',IB').
\]
By definition of $B$, we also have an identification $\calF_n(B,IB)=\prod_{j=1}^l\calF_n(A_j,IA_j)$.
By assumption, $\calF_n|_{(\Spf R_\lambda)_\Prism}$ is a finitely generated crystal of $\calO_{\Prism,n}$-modules. Hence we have a $B'/(p,I)^nB'$-linear isomorphism
\[
\eta\colon B'\otimes_{p_1,B}\calF_n(B,IB)\cong \calF_n(B',IB')\cong B'\otimes_{p_2,B}\calF_n(B,IB)
\]
satisfying the cocycle condition over $B/(p,I)^nB\otimes_{A}B/(p,I)^nB\otimes_AB/(p,I)^nB$.
Since $A/(p,I)^n\rightarrow B/(p,I)^nB$ is classically faithfully flat, it follows from faithfully flat descent that $\calF_n(A,I)$ is a finitely generated $A/(p,I)^n$-module and $B\otimes_A\calF_n(A,I)\cong \calF_n(B,IB)$. 

Let $(A,I)\rightarrow (\tilde{A},I\tilde{A})$ be a map of bounded prisms over $\fkX$. Set $\tilde{B}\coloneqq \tilde{A}\,\widehat{\otimes}_AB$. Then $\tilde{B}$ admits a natural $\delta$-structure and  $(\tilde{B},I\tilde{B})\in \fkX_\Prism$.
Moreover, $(\tilde{A},I\tilde{A})\rightarrow (\tilde{B},I\tilde{B})$ is $(p,I)$-completely faithfully flat. By the same argument as above, $\calF_n(\tilde{A},I\tilde{A})$ is an $\tilde{A}/(p,I)^n\tilde{A}$-module with $\tilde{B}\otimes_{\tilde{A}}\calF_n(\tilde{A},I\tilde{A})\cong \calF_n(\tilde{B},I\tilde{B})$.
Since $\calF_n|_{(\Spf R_\lambda)_\Prism}$ is a finitely generated crystal of $\calO_{\Prism,n}$-modules, we also have $\tilde{B}\otimes_B\calF_n(B,IB)\cong \calF_n(\tilde{B},I\tilde{B})$.
Hence the natural map $\tilde{A}\otimes_A\calF_n(A,I)\rightarrow \calF_n(\tilde{A},I\tilde{A})$ is an isomorphism since it is so after tensored with $\tilde{B}/(p,I)\tilde{B}$ over $\tilde{A}/(p,I)^n\tilde{A}$.
Therefore $\calF_n$ is a finitely generated crystal of $\calO_{\Prism,n}$-modules on $\fkX_\Prism$.
\end{proof}

\begin{rem}
 An analogue of Lemma~\ref{lem:Zariski descent of completed prismatic crystals} holds for an \'etale covering of $\fkX$ in place of an affine open covering. The verification is left to the reader.
\end{rem}

Recall that for an integral domain $R$ that is small over $\calO_K$, we defined the category $\mathrm{CR}^{\wedge,\varphi}(R_\Prism)$ of completed prismatic $F$-crystals on $R$ in Definition~\ref{defn:category-good-prism-completed-F-crystal}.

\begin{lem}\label{lem:zariskidescentoveraffinebases}
Assume that $\fkX=\Spf R$ is an affine formal scheme that is connected and small over $\calO_K$ and let $\calF$ be a sheaf of $\calO_\Prism$-modules on $\fkX_\Prism$ together with $1\otimes \varphi_\calF\colon \varphi^\ast \calF\rightarrow \calF$.
Then $(\calF,\varphi_{\calF})\in \mathrm{CR}^{\wedge,\varphi}(R_\Prism)$ if and only if
there exists an affine open covering $\fkX=\bigcup_{\lambda\in \Lambda} \Spf R_\lambda$ such that for each $\lambda$,  $R_\lambda$ is connected and small over $\calO_K$, and $(\calF|_{(\Spf R_\lambda)_\Prism},\varphi_{\mathcal{F}}|_{(\Spf R_\lambda)_\Prism})\in \mathrm{CR}^{\wedge,\varphi}((R_\lambda)_{\Prism})$.
\end{lem}

\begin{proof}
The necessity is straightforward. For the sufficiency, choose a $p$-complete \'etale map $R_0\rightarrow R_{\lambda,0}$ that induces $R\rightarrow R_\lambda$ after the base change along $W\rightarrow \calO_K$.
Set $R_0'\coloneqq \prod_{\lambda\in \Lambda}R_{\lambda,0}$ and extend the Frobenius on $R_0$ to $R_0'$.
Let $\fkS'\coloneqq R_0'[\![u]\!]$ and equip it with Frobenius $\varphi$ extending the one on $R_0'$ by $\varphi(u)=u^p$.
Via $\fkS'/(E)\cong \prod_{\lambda\in \Lambda}R_{\lambda}\leftarrow R$, we regard $(\fkS',(E))$ as an object of $R_\Prism$. Since $(\fkS,(E))\rightarrow (\fkS',(E))$ is a classically faithfully flat map of bounded prisms over $R$, the sufficiency follows from Lemmas~\ref{lem:independence of Breuil--Kisin condition} and \ref{lem:Zariski descent of completed prismatic crystals}. 
\end{proof}

\begin{defn} \label{defn:category-good-prism-completed-F-crystal:global case}
Let $\fkX$ be a smooth $p$-adic formal scheme over $\calO_K$.
A \emph{completed $F$-crystal of $\calO_\Prism$-modules} on $\fkX_\Prism$ is a pair
 $(\mathcal{F}, \varphi_{\mathcal{F}})$, where $\mathcal{F}$ is a finitely generated completed crystal of $\mathcal{O}_{\Prism}$-modules on $\fkX_\Prism$ and 
\[
\varphi_{\mathcal{F}}\colon \mathcal{F} \rightarrow \mathcal{F}
\]
is a $\varphi$-semilinear morphism of $\mathcal{O}_{\Prism}$-modules satisfying the following property: there exists an affine open covering $\fkX=\bigcup_{\lambda\in \Lambda} \Spf R_\lambda$  such that each $R_\lambda$ is connected and small over $\calO_K$ in the sense of Definition~\ref{defn:small-over-OK} and such that $(\calF|_{(\Spf R_\lambda)_\Prism},\varphi_\calF|_{(\Spf R_\lambda)_\Prism})\in \mathrm{CR}^{\wedge,\varphi}((R_\lambda)_{\Prism})$.
When $\fkX=\Spf R$ is affine that is connected and small over $\calO_K$, this definition coincides with Definition~\ref{defn:category-good-prism-completed-F-crystal} by Lemma~\ref{lem:zariskidescentoveraffinebases}.

We also call such an object a \textit{completed prismatic $F$-crystal} on $\fkX$.
The morphisms between completed $F$-crystals of $\calO_\Prism$-modules are $\mathcal{O}_{\Prism}$-module maps compatible with Frobenii $\varphi_{\calF}$.

We write $\mathrm{CR}^{\wedge,\varphi}(\fkX_{\Prism})$ for the category of completed $F$-crystals of $\calO_\Prism$-modules on $\fkX_\Prism$.
Let $\mathrm{Vect}^{\varphi}_{\mathrm{eff}}(\fkX_{\Prism})$ denote the full subcategory of $\mathrm{CR}^{\wedge,\varphi}(\fkX_{\Prism})$ consisting of objects $(\calF,\varphi_{\calF})$ where $\calF$ is a locally free $\calO_\Prism$-module. 
For a fixed non-negative integer $r$, we let $\mathrm{CR}^{\wedge,\varphi}_{[0,r]}(\fkX_{\Prism})$ and $\mathrm{Vect}^{\varphi}_{[0,r]}(\fkX_{\Prism})$ denote the full subcategories consisting of objects which locally lie in $\mathrm{CR}_{[0, r]}^{\wedge, \varphi}((R_{\lambda})_{\Prism})$.
\end{defn}

Let $\fkX_\eta$ denote the adic generic fiber of $\fkX$.
Recall that $\mathrm{Vect}(\fkX_{\Prism}, \mathcal{O}_{\Prism}[1/\mathcal{I_{\Prism}}]^{\wedge}_p)^{\varphi = 1}$ denotes the category of crystals of vector bundles $\calV$ on $(\fkX_\Prism,\calO_\Prism[1/\calI_\Prism]^\wedge_p)$
together with isomorphisms $\varphi_{\calV}\colon \varphi^\ast\calV\cong \calV$ \cite[Def.~3.2]{bhatt-scholze-prismaticFcrystal} and that 
there is a natural equivalence of categories
\[
\mathrm{Vect}(\fkX_{\Prism}, \mathcal{O}_{\Prism}[1/\mathcal{I_{\Prism}}]^{\wedge}_p)^{\varphi = 1}\cong \mathrm{Loc}_{\Z_p}(\fkX_\eta),
\]
where $\mathrm{Loc}_{\Z_p}(\fkX_\eta)$ denotes the category of \'etale $\Z_p$-local systems on $\fkX_\eta$ (see \cite[Cor.~ 3.8]{bhatt-scholze-prismaticFcrystal}).

\begin{prop}\label{prop:global-etale-realization}
The assignment $\calF\mapsto \calF_\et \coloneqq \varprojlim_n\calO_\Prism[1/\calI_\Prism]/p^n\otimes_{\calO_\Prism}\mathcal{F}$ 
defines a faithful functor
\[
\mathrm{CR}^{\wedge,\varphi}(\mathfrak{X}_\Prism)\rightarrow \mathrm{Vect}(\mathfrak{X}_{\Prism}, \mathcal{O}_{\Prism}[1/\mathcal{I_{\Prism}}]^{\wedge}_p)^{\varphi = 1}.
\]
\end{prop}

\begin{proof}
Take an affine open covering $\fkX=\bigcup_\lambda \Spf R_\lambda$ such that  $R_\lambda$ is connected and small over $\calO_K$ for each $\lambda$.
Then for each $\lambda$, $\calF|_{(\Spf R_\lambda)_\Prism}$ is naturally an object of $\mathrm{CR}^{\wedge,\varphi}((R_{\lambda})_\Prism)$.
Hence by Proposition~\ref{prop:etale-realization} (i), we obtain an object $(\calF|_{(\Spf R_\lambda)_\Prism})_\et$ of $\mathrm{Vect}((\Spf R_\lambda)_\Prism,\calO_\Prism[1/\calI_\Prism]^\wedge_p)^{\varphi=1}$ together with an identification
\[
(\calF|_{(\Spf R_\lambda)_\Prism})_\et|_{(\Spf R_\lambda\times_{\fkX}\Spf R_{\lambda'})_\Prism}\cong (\calF|_{(\Spf R_{\lambda'})_\Prism})_\et|_{(\Spf R_\lambda\times_{\fkX}\Spf R_{\lambda'})_\Prism}
\]
satisfying the cocycle condition over $(\Spf R_{\lambda}\times_\fkX \Spf R_{\lambda'}\times_\fkX\Spf R_{\lambda''})_\Prism$.
Hence they glue to an object $\calF_\et$ of $\mathrm{Vect}(\fkX_{\Prism}, \mathcal{O}_{\Prism}[1/\mathcal{I_{\Prism}}]^{\wedge}_p)^{\varphi = 1}$.
It is immediate to see that $\calF_\et$ is independent of the choice of the affine open covering and this gives the desired faithful functor $\calF\mapsto \calF_\et$.
\end{proof}

\begin{defn}
Define a contravariant functor $T\colon \mathrm{CR}^{\wedge,\varphi}(\fkX_\Prism)\rightarrow \mathrm{Loc}_{\Z_p}(\fkX_\eta)$
to be the composite
\[
\mathrm{CR}^{\wedge,\varphi}(\fkX_\Prism)\rightarrow \mathrm{Vect}(\fkX_{\Prism}, \mathcal{O}_{\Prism}[1/\mathcal{I_{\Prism}}]^{\wedge}_p)^{\varphi = 1}\cong \mathrm{Loc}_{\Z_p}(\fkX_\eta)\xrightarrow{(\,)^\vee}\mathrm{Loc}_{\Z_p}(\fkX_\eta),
\]
where the last functor sends $\bL$ to its dual $\Z_p$-local system $\bL^\vee$.
We call the functor $T$ the \emph{\'etale realization functor}.
Note that we use the contravariant convention as opposed to the covariant convention in \cite{bhatt-scholze-prismaticFcrystal}.
\end{defn}

\begin{notation}
Let $\mathrm{Loc}_{\Z_p,\geq 0}^\cris(\fkX_\eta)$ denote the full subcategory of $\mathrm{Loc}_{\Z_p}(\fkX_\eta)$ consisting of $\Z_p$-local systems $\bL$ on $\fkX_\eta$ such that $\bL\otimes_{\Z_p}\Q_p$ is a crystalline local system on $\fkX_\eta$ with non-negative Hodge--Tate weights. See Appendix~\ref{sec:crystalline local systems} for the definition of crystalline local systems on $\fkX_\eta$.
\end{notation}

\begin{thm}\label{thm:main-global}
Let $\fkX$ be a smooth $p$-adic formal scheme over $\calO_K$ and let $\fkX_\eta$ denote its adic generic fiber.
The \'etale realization functor $T$ induces the equivalence of categories
\[
T\colon \mathrm{CR}^{\wedge,\varphi}(\fkX_\Prism)
\xrightarrow{\cong}\mathrm{Loc}_{\Z_p,\geq 0}^\cris(\fkX_\eta).
\]
Moreover, $T$ is functorial in $\mathfrak{X}$.
\end{thm}

\begin{proof}
By Theorem~\ref{thm:main} (i), we see that $T$ factors through $\mathrm{Loc}_{\Z_p,\geq 0}^\cris(\fkX_\eta)\subset \mathrm{Loc}_{\Z_p}(\fkX_\eta)$ and $T$ is fully faithful since both properties are of local nature.
Once the full faithfulness is established, it follows from Proposition~\ref{prop:etale-realization} (iic), Theorem~\ref{thm:main} (ii), and gluing that $T\colon \mathrm{CR}^{\wedge,\varphi}(\fkX_\Prism)
\rightarrow\mathrm{Loc}_{\Z_p,\geq 0}^\cris(\fkX_\eta)$ is also essentially surjective. The functoriality follows from the construction of $T$.
\end{proof}

\section{Quasi-Kisin modules associated with crystalline representations} \label{sec:quasi-kisin-mods-cryst-loc-systs}

The goal of this section is to prove the second part of Theorem~\ref{thm:main} (the essential surjectivity of the \'etale realization functor). Recall $\fkS=R_0[\![u]\!]$ as in Notation~\ref{notation:Sigma}.
Given a $\mathbf{Z}_p$-lattice $T$ of a crystalline representation of $\calG_R$, we will construct a certain $\mathfrak{S}$-module equipped with a Frobenius and a connection, which we call a \textit{quasi-Kisin module} associated with $T$. 

In \S~\ref{sec:quasi-kisin-mod-rational-descent-data}, we introduce quasi-Kisin modules (Definition~\ref{defn:quasi-kisin-mod}) and attach a rational Kisin descent datum to a quasi-Kisin module (Construction~\ref{construction:rational S-descent datum} and Propositions~\ref{prop:rational-descent-datum-over-A(2)}, \ref{prop:p=2-rational-descent-datum-over-A(2)}). The proof crucially uses explicit computations of elements in $A_{\mathrm{max}}^{(1)}$ (Lemmas~\ref{lem:h0-filtration} and \ref{lem:p=2-h0-filtration}). Section~\ref{sec:quasi-kisin-mod-projectivity} shows, under Assumption~\ref{assumption:base-ring-sec-3.4}, that if $(\fkM,\varphi_\fkM)$ is a finitely generated torsion free $\varphi$-module  of finite $E$-height over $\fkS$, then $\fkM[p^{-1}]$ is projective over $\fkS[p^{-1}]$ (Proposition~\ref{prop:rational-projectivity-etale-over-torus-case}).
In \S~\ref{sec:cryst-rep-cdvr}, we consider the special case where $R = \mathcal{O}_L$ and establish some preliminary results. In \S~\ref{sec:quasi-kisin-mod-construction}-\ref{sec:quasi-kisin-mod-connection}, we construct a quasi-Kisin module associated with $T \in \mathrm{Rep}_{\mathbf{Z}_p, \geq 0}^{\mathrm{cris}}(\mathcal{G}_R)$. Finally,  \S~\ref{sec:equivalence-categories} completes the proof of Theorem~\ref{thm:main} by spreading the rational Kisin descent datum to an integral Kisin descent datum via the theory of \'etale $\varphi$-modules.

Since some of the arguments work for a general base ring $R$, which may be of some interest, we let $R$ be a base ring over $\calO_K$ as in Set-up~\ref{set-up:base ring} unless otherwise noted.

\subsection{Quasi-Kisin modules and associated rational Kisin descent data} \label{sec:quasi-kisin-mod-rational-descent-data}

Recall that $S$ denotes the $p$-adically completed divided power envelope of $\mathfrak{S}$ with respect to $(E(u))$, equipped with the Frobenius extending that on $\mathfrak{S}$. Let $\mathrm{Fil}^i S$ be the PD-filtration of $S$. Namely, $\Fil^i S$ is the $p$-adically completed ideal of $S$ generated by the divided powers $\gamma_j(E(u))$ ($j \geq i$), where $\gamma_j(x) \coloneqq \frac{x^j}{j!}$.
Let $N_u\colon S \rightarrow S$ be the $R_0$-linear derivation given by $N_u(u) = -u$, and let $\partial_u \colon S \rightarrow S$ be the $R_0$-linear derivation given by $\partial_u (u) = 1$. Note that $-u\partial_u = N_u$. We also have a natural integrable connection $\nabla=\nabla_S\colon S \rightarrow S\otimes_{R_0} \widehat{\Omega}_{R_0}$ given by the universal derivation on $R_0$, which commutes with $N_u$. 

\begin{defn} \label{defn:quasi-kisin-mod}
Let $r$ be a non-negative integer.
A \textit{quasi-Kisin module} over $\mathfrak{S}$ of $E$-height $\leq r$ is a triple $(\mathfrak{M}, \varphi_{\mathfrak{M}}, \nabla_{\mathfrak{M}})$ where
\begin{enumerate}
\item $\mathfrak{M}$ is a finitely generated $\mathfrak{S}$-module that is projective away from $(p,E)$ and saturated;
\item $\varphi_{\mathfrak{M}} \colon \mathfrak{M} \rightarrow \mathfrak{M}$ is a $\varphi$-semi-linear endomorphism such that $(\fkM,\varphi_{\fkM})$ has $E$-height $\leq r$;
\item if we set $M \coloneqq R_0\otimes_{\varphi,R_0} \mathfrak{M}/u\mathfrak{M}$ equipped with the induced tensor-product Frobenius, then 
\[
\nabla_{\mathfrak{M}} \colon M[p^{-1}] \rightarrow M[p^{-1}]\otimes_{R_0} \widehat{\Omega}_{R_0}
\]  
is a topologically quasi-nilpotent integrable connection commuting with Frobenius and satisfies the $S$-Griffiths transversality (see below). 
\end{enumerate}

Let us explain the definition of the $S$-Griffiths transversality.
Set $\mathscr{M}\coloneqq S\otimes_{\varphi,\fkS}\fkM$ and
define a decreasing filtration $F^i\mathscr{M}[p^{-1}]$ by
\[
 F^i\mathscr{M}[p^{-1}]\coloneqq\{x \in \mathscr{M}[p^{-1}] \mid (1\otimes\varphi_{\mathfrak{M}})(x) \in (\mathrm{Fil}^i S[p^{-1}])\otimes_{\mathfrak{S}}\mathfrak{M}\}.
\]
By Lemma~\ref{lem:frobenius-compatible-section} below, we have $\mathscr{M}[p^{-1}]\cong S[p^{-1}]\otimes_{R_0[p^{-1}]}M[p^{-1}]$, which admits a connection 
\[
\nabla_{\mathscr{M}[p^{-1}]}\colon \mathscr{M}[p^{-1}]\rightarrow\mathscr{M}[p^{-1}]\otimes_{R_0}\widehat{\Omega}_{R_0}
\]
given by $\nabla_{\mathscr{M}[p^{-1}]}=\nabla_{S[p^{-1}]}\otimes 1+1\otimes \nabla_{\fkM}$ so that $\varphi$ is horizontal. Let $\partial_u\colon \mathscr{M}[p^{-1}] \rightarrow \mathscr{M}[p^{-1}]$ be the derivation given by $\partial_{u, S[p^{-1}]}\otimes 1$.
We say that the connection $\nabla_\fkM$ or $\nabla_{\mathscr{M}[p^{-1}]}$ satisfies the \emph{$S$-Griffiths transversality} if, for every $i$,
\[
\partial_u(F^{i+1}\mathscr{M}[p^{-1}]) \subset F^i\mathscr{M}[p^{-1}] \text{~~and~~} \nabla_{\mathscr{M}[p^{-1}]}(F^{i+1}\mathscr{M}[p^{-1}])\subset (F^i\mathscr{M}[p^{-1}])\otimes_{R_0}\widehat{\Omega}_{R_0}.
\]
\end{defn}

\begin{lem} \label{lem:frobenius-compatible-section}
Let $(\fkM,\varphi_{\fkM})$ be a $\varphi$-module finite torsion free over $\fkS$ of $E$-height $\leq r$ such that $\fkM[p^{-1}]$ is projective over $\fkS[p^{-1}]$. Let $M \coloneqq R_0\otimes_{\varphi,R_0}\mathfrak{M}/u\mathfrak{M}$ and $\mathscr{M}\coloneqq S\otimes_{\varphi,\mathfrak{S}}\mathfrak{M}$ equipped with the induced Frobenii.
Consider the projection $q\colon \mathscr{M} \twoheadrightarrow M$ induced by the $\varphi$-compatible projection $S \twoheadrightarrow R_0$, $u \mapsto 0$. 
Then $q$ admits a unique $\varphi$-compatible section $s\colon M[p^{-1}] \rightarrow \mathscr{M}[p^{-1}]$. Furthermore, $1\otimes s\colon S[p^{-1}]\otimes_{R_0[p^{-1}]}M[p^{-1}] \rightarrow \mathscr{M}[p^{-1}]$ is an isomorphism.
\end{lem}

\begin{proof}
Since $\fkM$ has $E$-height $\leq r$, the map 
\[
(1\otimes \varphi)[p^{-1}]\colon \varphi^\ast M[p^{-1}] \coloneqq (R_0\otimes_{\varphi, R_0}M)[p^{-1}] \rightarrow M[p^{-1}]
\] 
is an isomorphism, and the preimage of $M$ is contained in $p^{-r}(\varphi^* M)$. It then follows from the standard argument as in the proof of \cite[Lem.~3.14]{kim-groupscheme-relative} that there exists a unique $\varphi$-compatible section $s\colon M[p^{-1}] \rightarrow \mathscr{M}[p^{-1}]$. Furthermore, the map $1\otimes s\colon S[p^{-1}]\otimes_{R_0[p^{-1}]}M[p^{-1}] \rightarrow \mathscr{M}[p^{-1}]$ is a map of projective $S[p^{-1}]$-modules of the same rank. So by a similar argument as in the proof of \cite[Lem.~4.17]{moon-strly-div-latt-cryst-cohom-CDVF}, $1\otimes s$ is an isomorphism.
\end{proof}

Let $(\fkM,\varphi_\fkM,\nabla_\fkM)$ be a quasi-Kisin module of $E$-height $\leq r$. We associate with $(\fkM,\varphi_\fkM,\nabla_\fkM)$ a rational Kisin descent datum, namely, an isomorphism of $\fkS^{(1)}[p^{-1}]$-modules
\[
f\colon \mathfrak{S}^{(1)}[p^{-1}]\otimes_{p_1,\mathfrak{S}}\fkM  \xrightarrow{\cong}\mathfrak{S}^{(1)}[p^{-1}]\otimes_{p_2, \mathfrak{S}}\fkM 
\]
satisfying the cocycle condition over $\mathfrak{S}^{(2)}[p^{-1}]$ and compatible with Frobenius.

First we will construct an isomorphism of $S^{(1)}[p^{-1}]$-modules
\[
f_S\colon S^{(1)}[p^{-1}]\otimes_{p_1, S}\mathscr{M} \xrightarrow{\cong} S^{(1)}[p^{-1}]\otimes_{p_2, S}\mathscr{M}
\]
satisfying the cocycle condition over $S^{(2)}[p^{-1}]$ and compatible with Frobenius and filtration. For each $i = 1, \ldots, d$, let $\partial_{T_i, M}\colon M[p^{-1}] \rightarrow M[p^{-1}]$ be the derivation given by $\nabla_{\mathfrak{M}}\colon M[p^{-1}] \rightarrow M[p^{-1}]\otimes_{R_0} \widehat{\Omega}_{R_0} \cong \bigoplus_{i=1}^d M[p^{-1}]\cdot dT_i$ composed with the projection to the $i$-th factor.

\begin{construction}\label{construction:rational S-descent datum}
Let $(\fkM,\varphi_\fkM,\nabla_\fkM)$ be a quasi-Kisin module of $E$-height $\leq r$.
Identify $\mathscr{M}[p^{-1}]$ with $\mathscr{D} \coloneqq S[p^{-1}]\otimes_{R_0}M$ as in Lemma~\ref{lem:frobenius-compatible-section}. Let $\partial_u \colon \mathscr{D} \rightarrow \mathscr{D}$ be the derivation given by $\partial_{u, S}\otimes1$, and for $i = 1, \ldots, d$, let $\partial_{T_i}\colon \mathscr{D} \rightarrow \mathscr{D}$ be the derivation given by $\partial_{T_i,S}\otimes 1+1\otimes \partial_{T_{i}, M}$. We define $f_S\colon S^{(1)}\otimes_{p_1,S}\mathscr{D} \rightarrow S^{(1)}\otimes_{p_2,S}\mathscr{D}$ by 
\[
f_S(x) = \sum \partial_u^{j_0}\partial_{T_1}^{j_1}\cdots \partial_{T_d}^{j_d} (x) \cdot \gamma_{j_0}(p_2(u)-p_1(u))\prod_{i=1}^d \gamma_{j_i}(p_2(T_i)-p_1(T_i)),
\]
where the sum goes over the multi-index $(j_0, \ldots, j_d)$ of non-negative integers. Note that $\partial_u$ and $\partial_{T_i}$'s are topologically quasi-nilpotent, so the above sum converges. It follows from a standard computation that this defines a $\varphi$-compatible isomorphism of $S^{(1)}[p^{-1}]$-modules $f_S\colon S^{(1)}\otimes_{p_1,S}\mathscr{D} \xrightarrow{\cong} S^{(1)}\otimes_{p_2,S}\mathscr{D}$ satisfying the cocycle condition over $S^{(2)}[p^{-1}]$. 

By the identification $\mathscr{M}[p^{-1}]=\mathscr{D}$, we obtain a descent datum
$f_S\colon S^{(1)}[p^{-1}]\otimes_{p_1,S}\mathscr{M} \xrightarrow{\cong} S^{(1)}[p^{-1}]\otimes_{p_2, S}\mathscr{M}$.
Since $\nabla_{\mathscr{M}[p^{-1}]}$ satisfies the $S$-Griffiths transversality, we see that $f_S$ is compatible with filtrations (see below for the filtration on $S^{(1)}$). 
\end{construction}

To further construct a rational Kisin descent datum, we need to discuss filtrations on subrings of $A_{\mathrm{max}}^{(1)}[p^{-1}]$ such as $\mathfrak{S}^{(1)}$ and $S^{(1)}$. Recall that $A_{\mathrm{max}}^{(1)}$ is defined after Example~\ref{eg:prism-PD-S}. Our argument below can be regarded as a counterpart of the argument about the ``boundedness of descent data at the boundary'' in \cite[\S 6.3]{bhatt-scholze-prismaticFcrystal} via the isomorphism  $A_{\mathrm{max}}^{(1)}\cong\mathfrak{S}^{(1)}\langle \frac{E}{p} \rangle$ (see also \cite[Rem.~2.2.14]{du-liu-prismaticphiGhatmodule}). 
For any subring $B \subset A_{\mathrm{max}}^{(1)}[p^{-1}]$ that is stable under $\varphi_{A_{\mathrm{max}}^{(1)}[p^{-1}]}$, define 
\[\mathrm{Fil}^m B \coloneqq B \cap E^m A_{\mathrm{max}}^{(1)}[p^{-1}].
\]
In particular, we have $\mathrm{Fil}^m \mathfrak{S} = E^m \mathfrak{S}$ and $\mathrm{Fil}^m \mathfrak{S}^{(1)} = E^m \mathfrak{S}^{(1)}$ by \cite[Cor.~2.2.9]{du-liu-prismaticphiGhatmodule}. Note that $\mathrm{Fil}^m S^{(1)}$ is compatible with the PD-filtration on $S^{(1)}$, i.e.,
\begin{align*}
\mathrm{Fil}^m S^{(1)} =
\bigl\{ \sum_{i_0+\cdots+i_{d+1} \geq m} &a_{i_0, \ldots, i_{d+1}} \gamma_{i_0}(E) \gamma_{i_1}(y-u)\gamma_{i_2}(s_1-T_1)\cdots \gamma_{i_{d+1}} (s_d-T_d) ~|\\ 
&a_{i_0, \dots , i_{d+1}} \in \mathfrak{S}^{\widehat{\otimes}[1]}, ~a_{i_0, \dots , i_{d+1}} \rightarrow 0\quad(\text{as } i_0+\cdots+i_{d+1}\to\infty)\bigr\}.
\end{align*}

\begin{lem} \label{lem:h0-filtration} 
Assume $p\geq 3$ and let $r$ be a fixed non-negative integer. There exists an integer $h_0 > r$ such that if $m \geq  h_0$ and $x \in S^{(1)}[E^{-1}]$ with $E^r x \in \Fil^{m} S^{(1)}$, then $\varphi(x) = a + b$ for some $a \in \mathfrak{S}^{(1)}$ and $b \in \Fil^{m+1} S^{(1)}$ (as elements in $A_{\mathrm{max}}^{(1)}$). 
\end{lem}

\begin{proof}
By the explicit description of $\Fil^m S^{(1)}$, since $y-u = Ez_0$, $s_j-T_j = Ez_j$ and $z_j \in \mathfrak{S}^{(1)}$, we can write $E^r x = \sum_{i \geq m} c_i \gamma_i(E) $ for some $c_i \in \mathfrak{S}^{(1)}$ with $c_i \rightarrow 0$ $p$-adically as $i \rightarrow \infty$. So 
\[
\varphi(x)= \sum_{i \geq m } \varphi(c_i) \varphi\Bigl(\frac{\gamma_i(E)}{E^r}\Bigr).
\] 
It suffices to show that there exists $h _0 > r$ such that if $m \geq h_0$ then $\varphi (\frac{E^{m-r}}{m!}) = a_m + b_m$ for some $a_m \in \mathfrak{S}$ and $b_m \in \Fil^{m +1}S$. For this, note that $\varphi(E)= E^p + pt$ for some $t \in \mathfrak{S}$. So 
\[
\varphi(E)^{m-r}= (E^p + pt)^{m-r}= \sum_{i=0}^{m-r} \binom{m-r}{i} E^{p(m -r -i)} (pt)^i.
\] 
Let $v_p(\cdot )$ be the $p$-adic valuation with $v_p(p) = 1$. Since $v_p (m!) < \frac{m}{p-1}$, we have 
\[
a_m \coloneqq \frac{1}{m!} \sum_{i\geq \frac{m}{p-1}}^{m-r} \binom{m-r}{i} E^{p(m -r -i )} (pt)^i\in \mathfrak{S}.
\]
Consider $b_m\coloneqq \frac{1}{m!} \sum_{0\leq i< \frac{m}{p-1}}\binom{m-r}{i} E^{p(m -r -i )} (pt)^i$.
If $p (m-r-\frac {m}{p-1}) \geq m+1$, i.e., 
$m \geq  \frac{(p-1)(pr+1)}{p^2-3p +1}$ (since $p>2$), then $b_m\in \Fil^{m+1}S$. Hence, we can set $h_0 = \bigl\lceil\frac{(p-1)(pr+1)}{p^2-3p +1}\bigr\rceil $.	
\end{proof}

We now return to the discussion on the quasi-Kisin module $(\fkM,\varphi_\fkM,\nabla_\fkM)$. Set $\mathfrak{M}^* \coloneqq \mathfrak{S}\otimes_{\varphi, \mathfrak{S}} \mathfrak{M}$. For $j = 1, 2$, let $\fkM_j^{(1)}\coloneqq \mathfrak{S}^{(1)}\otimes_{p_j, \mathfrak{S}} \mathfrak{M}$,  $\fkM_j^{\ast,(1)} \coloneqq \mathfrak{S}^{(1)}\otimes_{p_j, \mathfrak{S}} \mathfrak{M}^*$, and $\fkM_{\mathrm{max},j}^{\ast,(1)} \coloneqq A_{\mathrm{max}}^{(1)}[p^{-1}]\otimes_{p_j, \mathfrak{S}} \mathfrak{M}^*$. 
If $B$ is a subring of $A_{\mathrm{max}}^{(1)}[p^{-1}]$ stable under $\varphi_{A_{\mathrm{max}}^{(1)}[p^{-1}]}$ and if $p_j\colon \mathfrak{S} \rightarrow A_{\mathrm{max}}^{(1)}[p^{-1}]$ factors through $B$, then define
\[
\mathrm{Fil}^i (B\otimes_{p_j, \mathfrak{S}} \mathfrak{M}^*) \coloneqq \{x \in B\otimes_{p_j, \mathfrak{S}} \mathfrak{M}^* \mid (1\otimes\varphi_{\mathfrak{M}})(x) \in \mathrm{Fil}^i B\otimes_{p_j, \mathfrak{S}} \mathfrak{M}\}.
\] 
Note that
\begin{align*}
\mathrm{Fil}^i \mathfrak{M}^* &= \{x \in \mathfrak{M}^* \mid (1\otimes\varphi)(x) \in E^i \mathfrak{M} \}
\quad\text{and}\\
\mathrm{Fil}^i \fkM_j^{\ast,(1)}  &= \{x \in \fkM_j^{\ast,(1)} \mid (1\otimes\varphi)(x) \in E^i \fkM_j^{(1)}\}. 
\end{align*}

Since $\fkM$ has $E$-height $\leq r$, 
$1\otimes\varphi\colon \mathrm{Fil}^r \mathfrak{M}^* \rightarrow E^r \mathfrak{M}$ is an isomorphism. Let $\varphi_r\colon \mathrm{Fil}^r \mathfrak{M}^* \rightarrow \mathfrak{M}^*$ be the $\varphi$-semi-linear map given by the composite
\[
\varphi_r\colon \mathrm{Fil}^r \mathfrak{M}^* \xrightarrow{1\otimes\varphi}E^r\mathfrak{M} \cong E^r\mathfrak{S}\otimes_{\mathfrak{S}} \mathfrak{M} \xrightarrow{\frac{\varphi\otimes 1}{\varphi(E^r)}} \mathfrak{S}\otimes_{\varphi, \mathfrak{S}} \mathfrak{M} = \mathfrak{M}^*
\]
Note that $\varphi_r(\mathrm{Fil}^r \mathfrak{M}^*)$ generates $\mathfrak{M}^*$ as an $\mathfrak{S}$-module. Similarly, we define the $\varphi$-semi-linear map $\varphi_r\colon \mathrm{Fil}^r \fkM_j^{\ast,(1)} \rightarrow \fkM_j^{\ast,(1)}$.

\begin{lem} \label{lem:filtration-intersection} 
We have $(\mathrm{Fil}^i \fkM_{\mathrm{max}, j}^{\ast,(1)}) \cap \fkM_j^{\ast,(1)} = \mathrm{Fil}^i \fkM_j^{\ast,(1)}$.	
\end{lem}

\begin{proof}
By assumption, $(p, u)$ forms a regular sequence for $\mathfrak{M}$ as an $\mathfrak{S}$-module. So $(p, E)$ is a regular sequence for $\mathfrak{M}$, and $\mathfrak{M}/E\mathfrak{M}$ is $p$-torsion free. Since $p_j\colon \mathfrak{S} \rightarrow \mathfrak{S}^{(1)}$ is classically flat by Lemma \ref{lem:AtoA2A3-faithful-flat}, $\fkM_j^{(1)}/E\fkM_j^{(1)}$ is $p$-torsion free. 
In particular, $E^i \fkM_j^{(1)}[p^{-1}] \cap \fkM_j^{(1)}  = E^i \fkM_j^{(1)}$.

It suffices to show 
\[
(E^i A_{\mathrm{max}}^{(1)}[p^{-1}]\otimes_{p_j, \mathfrak{S}} \mathfrak{M})  \cap \fkM_j^{(1)} = E^i \fkM_j^{(1)}
\] 	
as submodules of $A_{\mathrm{max}}^{(1)}[p^{-1}]\otimes_{p_j, \mathfrak{S}} \mathfrak{M}$, which makes sense since $\mathfrak{M}[p^{-1}]$ is projective over $\mathfrak{S}[p^{-1}]$ by assumption.
Since $E^i A_{\mathrm{max}}^{(1)}[p^{-1}] \cap \mathfrak{S}^{(1)}[p^{-1}]=E^i \mathfrak{S}^{(1)}[p^{-1}]$, Lemma~\ref{lem:intersection-modules-flat-base-change} (i) implies that
\begin{align*}
 (E^i A_{\mathrm{max}}^{(1)}[p^{-1}]\otimes_{p_j, \mathfrak{S}} \mathfrak{M}[p^{-1}])  \cap \fkM_j^{(1)}[p^{-1}] &= E^i \mathfrak{S}^{(1)}[p^{-1}]\otimes_{p_j, \mathfrak{S}[p^{-1}]} \mathfrak{M}[p^{-1}]\\
	&= E^i \fkM_j^{(1)}[p^{-1}].
\end{align*}
Since $E^i \fkM_j^{(1)}[p^{-1}] \cap \fkM_j^{(1)}  = E^i \fkM_j^{(1)}$ by above, the assertion follows.
\end{proof}

We can now show that $f_S$ defines a rational Kisin descent datum when $p \geq 3$. The same result also holds for $p=2$ (Proposition~\ref{prop:p=2-rational-descent-datum-over-A(2)}) with a similar but longer proof, and we postpone the latter case.

\begin{prop} \label{prop:rational-descent-datum-over-A(2)}
Assume $p \geq 3$. Let $(\fkM,\varphi_\fkM,\nabla_\fkM)$ be a quasi-Kisin module of $E$-height $\leq r$. 
There exists a unique rational Kisin descent datum
\[
f\colon \mathfrak{S}^{(1)}[p^{-1}]\otimes_{p_1,\mathfrak{S}}\mathfrak{M}\xrightarrow{\cong} \mathfrak{S}^{(1)}[p^{-1}]\otimes_{p_2,\mathfrak{S}}\mathfrak{M}
\]
such that $\mathrm{id}_{S^{(1)}}\otimes_{\varphi,\mathfrak{S}^{(1)}}f = f_S$, where $f_S$ is defined as in Construction~\ref{construction:rational S-descent datum}.
\end{prop}

\begin{proof}
For $j = 1, 2$, we write $\mathscr{M}_j^{(1)}$ for the image of $S^{(1)}\otimes_{p_j,S}\mathscr{M}$ in $S^{(1)}[p^{-1}]\otimes_{p_j,S}\mathscr{M}$ under the natural map. We will show that there exists a unique $\fkS^{(1)}[p^{-1}]$-linear map
\[
f\colon \mathfrak{S}^{(1)}[p^{-1}]\otimes_{p_1,\mathfrak{S}}\mathfrak{M}  \rightarrow \mathfrak{S}^{(1)}[p^{-1}]\otimes_{p_2,\mathfrak{S}}\mathfrak{M}
\]
such that $\mathrm{id}_{S^{(1)}}\otimes_{\varphi,\mathfrak{S}^{(1)}}f= f_S$.
Let $h_0 > r$ be a constant given as in Lemma~\ref{lem:h0-filtration}. Note that by the explicit description of $\mathrm{Fil}^m S^{(1)}$, for any $x \in S^{(1)}$, we have $p^{h_0}x = y+z$ for some $y \in \mathfrak{S}^{\widehat{\otimes}[1]}$ and $z \in \mathrm{Fil}^{h_0} S^{(1)}$. Thus, we can take a sufficiently large integer $n \geq 0$ such that $f'_S\coloneqq p^n f_S$ satisfies $f'_S(\mathscr{M}_1^{(1)}) \subset \mathscr{M}_2^{(1)}$ and
\[
f'_S(\mathfrak{M}^*) \subset \fkM_2^{\ast,(1)} + \mathrm{Fil}^{h_0} S^{(1)}\cdot  \mathscr{M}_2^{(1)}
\]	
as submodules of $A_{\mathrm{max}}^{(1)}[p^{-1}]  \otimes_{p_2, \mathfrak{S}}\mathfrak{M}^*$.  We claim that 
\[
f'_S(\mathfrak{M}^*) \subset \fkM_2^{\ast,(1)} + \mathrm{Fil}^{m} S^{(1)}\cdot  \mathscr{M}_2^{(1)}
\]	  
for any $m \geq h_0$. We induct on $m$. 
Suppose that the claim holds for $m ~(\geq h_0)$. Let $w \in \mathrm{Fil}^r \mathfrak{M}^*$ viewed as an element in $\mathscr{M}_1^{(1)}$ via $p_1$. Then we can write
\[
f'_S(w) = z+\sum_i a_i w_i
\]
for some $z \in \fkM_2^{\ast,(1)}$, $a_i \in \mathrm{Fil}^{m} S^{(1)}$, $w_i \in \mathfrak{M}^\ast$ that are viewed as elements in $\mathscr{M}_2^{(1)}$ via $p_2$ with a finite index set for $i$. Note that $z \in (\mathrm{Fil}^r \fkM_{\mathrm{max}, 2}^{\ast,(1)}) \cap \fkM_2^{\ast,(1)} = \mathrm{Fil}^r \fkM_2^{\ast,(1)}$ by Lemma~\ref{lem:filtration-intersection}. Let $a_i' = \frac{a_i}{E^r} \in S^{(1)}[E^{-1}]$. Then $f'_S(w) = z+\sum_i a_i' \cdot E^r w_i$ with $E^r w_i \in \mathrm{Fil}^r \mathfrak{M}^*$.

We have
\[
f'_S(\varphi_r(w)) = \varphi_r(z)+\sum_i \varphi(a_i')\varphi_r(E^r w_i).
\]
Since $E^r a_i' \in \mathrm{Fil}^{m} S^{(1)}$, we have $\varphi(a_i') = b_i+c_i$ for some $b_i \in \mathfrak{S}^{(1)}$ and $c_i \in \mathrm{Fil}^{m+1} S^{(1)}$ by Lemma~\ref{lem:h0-filtration}. Thus, $f'_S(\varphi_r(w)) \in \fkM_2^{\ast,(1)} + \mathrm{Fil}^{m+1} S^{(1)}\cdot  \mathscr{M}_2^{(1)}$. Since $\varphi_r(\mathrm{Fil}^r \mathfrak{M}^*)$ generates $\mathfrak{M}^*$ as $\mathfrak{S}$-modules, the claim follows.

Since $\mathfrak{M}[p^{-1}]$ is finite projective over $\mathfrak{S}[p^{-1}]$ by assumption and the filtration $\{\Fil^m S^{(1)}[p^{-1}]\}$ is separated, we deduce that $f'_S(\mathfrak{M}^*) \subset \fkM_2^{\ast,(1)}[p^{-1}]$. By increasing $n$ if necessary, we may further assume $f'_S(\mathfrak{M}^*) \subset \fkM_2^{\ast,(1)}$. Then $f'_S(\fkM_1^{\ast,(1)}) \subset \fkM_2^{\ast,(1)}$, and
\[
f'_S(\mathrm{Fil}^r \fkM_1^{\ast,(1)}) \subset (\mathrm{Fil}^r \fkM_{\mathrm{max}, 2}^{\ast,(1)}) \cap \fkM_2^{\ast,(1)} = \mathrm{Fil}^r \fkM_2^{\ast,(1)} 
\] 
by Lemma~\ref{lem:filtration-intersection}. Consider the composite of the isomorphisms
\[
\mathrm{Fil}^r \mathfrak{M}^* \overset{1\otimes\varphi}{\cong} E(u)^r\mathfrak{M} \cong \mathfrak{M}.
\]
Since $p_j\colon \mathfrak{S} \rightarrow \mathfrak{S}^{(1)}$ is classically faithfully flat by Lemma~\ref{lem:AtoA2A3-faithful-flat}, we obtain the isomorphism $\mathrm{Fil}^r \fkM_j^{\ast,(1)} \cong \fkM_j^{(1)}$ of $\mathfrak{S}^{(1)}$-modules for $j = 1, 2$. Via these isomorphisms, $f'_S\colon\mathrm{Fil}^r \fkM_1^{\ast,(1)} \rightarrow \mathrm{Fil}^r \fkM_2^{\ast,(1)}$ induces $f'\colon \mathfrak{S}^{(1)}[p^{-1}]\otimes_{p_1,\mathfrak{S}}\mathfrak{M}  \rightarrow \mathfrak{S}^{(1)}[p^{-1}]\otimes_{p_2,\mathfrak{S}} \mathfrak{M}$. If we set $f\coloneqq p^{-n}f'$, then we have $\mathrm{id}_{S^{(1)}}\otimes_{\varphi,\mathfrak{S}^{(1)}}f = f_S$.
The uniqueness is obvious.

By applying the same argument to $f_S^{-1}$, we conclude that $f$ is an isomorphism. Hence $f$ is a rational Kisin descent datum.
\end{proof}

We now explain how to obtain a rational Kisin descent datum from $f_S$ when $p = 2$. 
 We consider two auxiliary subrings $\widetilde{S}$, $\widehat{S}$ of $A_{\mathrm{max}}^{(1)}$, defined by
\begin{align*}
\widetilde{S} &\coloneqq \mathfrak{S}^{(1)}\Bigl[\!\!\Bigl[\frac{E^2}{2}\Bigr]\!\!\Bigr] = \Bigl\{\sum_{i \geq 0} a_i\Bigl(\frac{E^2}{2}\Bigr)^i \Bigm| a_i \in \mathfrak{S}^{(1)} \Bigr\}\quad\text{and}\\    
\widehat{S} &\coloneqq \mathfrak{S}^{(1)}\Bigl[\!\!\Bigl[\frac{E^4}{2}\Bigr]\!\!\Bigr] = \Bigl\{\sum_{i \geq 0} a_i\Bigl(\frac{E^4}{2}\Bigr)^i \Bigm| a_i \in \mathfrak{S}^{(1)} \Bigr\}.
\end{align*}
Since $\varphi(E) = E^2+2\delta(E)$ and $\mathfrak{S}^{(1)}$ is $2$-adically complete, both $\widetilde{S}$ and $\widehat{S}$ are stable under the ring endomorphism $\varphi$ on $ A_{\mathrm{max}}^{(1)}$. The following is shown in \cite{du-liu-prismaticphiGhatmodule}.

\begin{lem}[{(cf. {\cite[Lem. 2.2.12]{du-liu-prismaticphiGhatmodule}})}]  \label{lem:tildeShatS-basic-properties}
Suppose $p = 2$. The following properties hold.
\begin{enumerate}
\item $\varphi\bigl(A_{\mathrm{max}}^{(1)}\bigr) \subset \widetilde{S}$ and $\varphi\bigl(\widetilde{S}\bigr) \subset \widehat{S}$.

\item For every positive integer $h$, we have
\[
\mathrm{Fil}^h \widetilde{S} = \Bigl\{ \sum_{i \geq h} a_i \frac{E^i}{2^{\lfloor \frac{i}{2} \rfloor}} \Bigm| a_i \in \mathfrak{S}^{(1)} \Bigr\}\quad\text{and}\quad
\mathrm{Fil}^h \widehat{S} = \Bigl\{ \sum_{i \geq h} a_i \frac{E^i}{2^{\lfloor \frac{i}{4} \rfloor}} \Bigm| a_i \in \mathfrak{S}^{(1)} \Bigr\}.
\]
\end{enumerate}
\end{lem}

\begin{lem} \label{lem:p=2-h0-filtration}
Assume $p = 2$, and let $r$ be a fixed non-negative integer. There exists an integer $h_0 > r$ such that if $m \geq  h_0$ and $x \in \widehat{S}[E^{-1}]$ with $E^r x \in \Fil^{m} \widehat{S}$, then $\varphi(x) = a + b$ for some $a \in \mathfrak{S}^{(1)}$ and $b \in \Fil^{m+1} \widehat{S}$ (as elements in $A_{\mathrm{max}}^{(1)}$). 
\end{lem}

\begin{proof}
By Lemma~\ref{lem:tildeShatS-basic-properties} (ii), we can write
\[
E^r x = \sum_{i \geq m} c_i \frac{E^i}{2^{\lfloor \frac{i}{4} \rfloor}}
\]
for some $c_i \in \mathfrak{S}^{(1)}$. So 
\[
\varphi(x) = \sum_{i \geq m} \varphi(c_i) \frac{\varphi(E^{i-r})}{2^{\lfloor \frac{i}{4} \rfloor}}.
\]

It suffices to show that there exists $h _0 > r$ such that if $m \geq h_0$ then $\displaystyle \frac{\varphi(E^{m-r})}{2^{\lfloor \frac{m}{4} \rfloor}} = a_m + b_m$ for some $a_m \in (2, u)^{m-r-\lfloor \frac{m}{4}\rfloor} \mathfrak{S}$ and $b_m \in \Fil^{m +1} \widehat{S}$. For this, note that 
\[
\varphi(E^{m-r}) = (E^2 + 2\delta(E))^{m-r}= \sum_{i=0}^{m-r} \binom{m-r}{i} E^{2(m -r -i)} (2\delta(E))^i.
\] 
We have 
\[
a_m \coloneqq \frac{1}{2^{\lfloor \frac{m}{4} \rfloor}} \sum_{i\geq \lfloor \frac{m}{4} \rfloor}^{m-r}\binom{m-r}{i} E^{2(m -r -i )} (2\delta(E))^i \in (2, u)^{m-r-\lfloor \frac{m}{4}\rfloor} \mathfrak{S}.
\]
Set $\displaystyle b_m\coloneqq \frac{1}{2^{\lfloor \frac{m}{4} \rfloor}} \sum_{0\leq i \leq \lfloor \frac{m}{4} \rfloor-1}\binom{m-r}{i} E^{2(m -r -i )} (2\delta(E))^i$. If $2(m-r-\lfloor \frac{m}{4} \rfloor+1) \geq m+1$, then $b_m\in \Fil^{m+1}\widehat{S}$. Since 
\[
2\Bigl(m-r-\Bigl\lfloor \frac{m}{4} \Bigr\rfloor+1\Bigr) \geq 2\Bigl(m-r-\frac{m}{4}+1\Bigr) = \frac{3}{2}m-2r+2,
\]
we can set $h_0 = 4r+1$. 
\end{proof}

Using Lemmas~\ref{lem:tildeShatS-basic-properties} and \ref{lem:p=2-h0-filtration}, we now construct a rational Kisin datum when $p=2$.

\begin{prop} \label{prop:p=2-rational-descent-datum-over-A(2)}
Assume $p = 2$. Let $(\fkM,\varphi_\fkM,\nabla_\fkM)$ be a quasi-Kisin module of $E$-height $\leq r$. 
There exists a unique rational Kisin descent datum
\[
f\colon \mathfrak{S}^{(1)}[p^{-1}]\otimes_{p_1,\mathfrak{S}}\mathfrak{M}  \xrightarrow{\cong} \mathfrak{S}^{(1)}[p^{-1}]\otimes_{p_2,\mathfrak{S}}\mathfrak{M}
\]
such that $\mathrm{id}_{S^{(1)}}\otimes_{\varphi,\mathfrak{S}^{(1)}}f = f_S$, where $f_S$ is defined as in Construction~\ref{construction:rational S-descent datum}.
\end{prop}

\begin{proof}
For $j = 1, 2$, write $\mathscr{M}_j^{(1)}$ for the image of $S^{(1)}\otimes_{p_j,S}\mathscr{M}$ in $S^{(1)}[p^{-1}]\otimes_{p_j,S}\mathscr{M}$ under the natural map. We first claim that $f_S(\mathfrak{M}^*) \subset \widetilde{S}[p^{-1}]\otimes_{p_2, \mathfrak{S}} \mathfrak{M}$. For this, take a sufficiently large integer $n \geq 0$ such that $f'_S\coloneqq p^n f_S$ satisfies $f'_S(\mathscr{M}_1^{(1)}) \subset \mathscr{M}_2^{(1)}$ and
\[
f'_S(\mathfrak{M}^*) \subset \fkM_2^{\ast,(1)} + \mathrm{Fil}^{r} S^{(1)}\cdot  \mathscr{M}_2^{(1)}
\]	
as submodules of $A_{\mathrm{max}}^{(1)}[p^{-1}] \otimes_{p_2, \mathfrak{S}}\mathfrak{M}^*$. Let $w \in \mathrm{Fil}^r \mathfrak{M}^*$. We can write
\[
f'_S(w) = z+\sum_i a_i w_i
\]
for some $z \in \fkM_2^{\ast,(1)}$, $a_i \in \mathrm{Fil}^{r} S^{(1)}$, $w_i \in \mathfrak{M}^\ast$ (with finitely many indices $i$). Note that $z \in (\mathrm{Fil}^r \fkM_{\mathrm{max}, 2}^{\ast,(1)} ) \cap \fkM_2^{\ast,(1)} = \mathrm{Fil}^r \fkM_2^{\ast,(1)}$ by Lemma~\ref{lem:filtration-intersection}. Since $a_i \in \mathrm{Fil}^{r} S^{(1)}$, it follows from the explicit description of $\mathrm{Fil}^{r} S^{(1)}$ that $a_i' \coloneqq p^r\frac{a_i}{E^r}$ lies in $A_{\mathrm{max}}^{(1)}$. We have $f'_S(p^rw) = p^r z+\sum_i a_i' \cdot E^r w_i$ with $E^r w_i \in \mathrm{Fil}^r \mathfrak{M}^*$ as elements in $A_{\mathrm{max}}^{(1)}[p^{-1}] \otimes_{p_2, \mathfrak{S}}\mathfrak{M}^*$, and so 
\[
f'_S(\varphi_r(p^r w)) = \varphi_r(p^r z)+\sum_i \varphi(a_i')\varphi_r(E^r w_i).
\]
Note that $\varphi(a_i') \in \widetilde{S}$ by Lemma~\ref{lem:tildeShatS-basic-properties} (i). Thus, we deduce $f_S(\varphi_r(w)) \in \widetilde{S}[p^{-1}]\otimes_{p_2, \mathfrak{S}} \mathfrak{M}$. Since $\varphi_r(\mathrm{Fil}^r \mathfrak{M}^*)$ generates $\mathfrak{M}^*$ as $\mathfrak{S}$-modules, the claim follows.

Let $a \in \mathrm{Fil}^r \widetilde{S}$. Since $\lfloor \frac{i-r}{2} \rfloor - (\lfloor \frac{i}{2} \rfloor -r) \geq 0$, it follows from Lemma~\ref{lem:tildeShatS-basic-properties} (ii) that $p^r\frac{a}{E^r} \in \widetilde{S}$. Furthermore, $\varphi(\widetilde{S}) \subset \widehat{S}$ by Lemma~\ref{lem:tildeShatS-basic-properties} (i). Thus, starting with $f_S(\mathfrak{M}^*) \subset \widetilde{S}[p^{-1}]\otimes_{p_2, \mathfrak{S}} \mathfrak{M}$, we can repeat a similar argument to further obtain
\[
f_S(\mathfrak{M}^*) \subset \widehat{S}[p^{-1}]\otimes_{p_2, \mathfrak{S}} \mathfrak{M}.
\]

As in the proof of Proposition~\ref{prop:rational-descent-datum-over-A(2)} with Lemma~\ref{lem:p=2-h0-filtration} in place of Lemma~\ref{lem:h0-filtration}, we deduce $f_S(\mathfrak{M}^*) \subset \fkM_2^{\ast,(1)}[p^{-1}]$. The rest of the proof proceeds exactly as in the proof of Proposition~\ref{prop:rational-descent-datum-over-A(2)}.
\end{proof}

We end this subsection with a simple lemma.

\begin{lem}\label{lem:intersectin of two base chage over S^(2)}
Let $\mathfrak{M}$ be a finitely generated $\mathfrak{S}$-module which is projective away from $(p,E)$ and saturated.
Then the natural map
\[
\mathfrak{S}^{(1)}\otimes_{p_j,\mathfrak{S}}\mathfrak{M} \rightarrow  (\mathfrak{S}^{(1)}[p^{-1}]\otimes_{p_j,\mathfrak{S}}\mathfrak{M}) \cap (\mathfrak{S}^{(1)}[E^{-1}]^{\wedge}_p\otimes_{p_j,\mathfrak{S}}\mathfrak{M} )
\]
is an isomorphism.
\end{lem}

\begin{proof}
Note first that the maps
\[
\mathfrak{S}^{(1)}\otimes_{p_j,\mathfrak{S}}\mathfrak{M} \rightarrow \mathfrak{S}^{(1)}[p^{-1}] \otimes_{p_j,\mathfrak{S}}\mathfrak{M} 
\quad\text{and}\quad
\mathfrak{S}^{(1)}\otimes_{p_j,\mathfrak{S}}\mathfrak{M} \rightarrow \mathfrak{S}^{(1)}[E^{-1}]^{\wedge}_p\otimes_{p_j,\mathfrak{S}} \mathfrak{M}
\]
are injective by the same argument as in the proof of Lemma~\ref{lem:completed-crystals-basic-properties} (iii).

We need to show that the injective map
\[
\mathfrak{S}^{(1)}\otimes_{p_j,\mathfrak{S}}\mathfrak{M} \hookrightarrow  (\mathfrak{S}^{(1)}[p^{-1}]\otimes_{p_j,\mathfrak{S}}\mathfrak{M}) \cap (\mathfrak{S}^{(1)}[E^{-1}]^{\wedge}_p\otimes_{p_j,\mathfrak{S}}\mathfrak{M} )
\]
is also surjective. Suppose not. 
Set $\mathfrak{L} \coloneqq (\mathfrak{S}^{(1)}[p^{-1}]\otimes_{p_j,\mathfrak{S}}\mathfrak{M}) \cap (\mathfrak{S}^{(1)}[E^{-1}]^{\wedge}_p\otimes_{p_j,\mathfrak{S}}\mathfrak{M})$ for simplicity.
For any $\mathbf{Z}_p$-module $Q$, write $Q/p$ for $Q/pQ$. Then the induced map
\[
\mathfrak{S}^{(1)}/p\otimes_{p_j,\mathfrak{S}}\mathfrak{M} \rightarrow \mathfrak{L} / p
\]
is not injective since $\fkS^{(1)}[p^{-1}] \otimes_{p_j, \fkS} \fkM= \fkL [p^{-1}]$. On the other hand, by the saturation assumption, we have $\mathfrak{M}[p^{-1}] \cap \mathfrak{M}[E^{-1}] = \mathfrak{M}$. So by Lemmas~\ref{lem:intersection-modules-flat-base-change} (i) and \ref{lem:AtoA2A3-faithful-flat}, 
\[
(\mathfrak{S}^{(1)}\otimes_{p_j,\mathfrak{S}}\mathfrak{M}[p^{-1}]) \cap (\mathfrak{S}^{(1)}\otimes_{p_j,\mathfrak{S}}\mathfrak{M}[E^{-1}]) = \mathfrak{S}^{(1)}\otimes_{p_j,\mathfrak{S}}\mathfrak{M}.
\]
This implies that the map
\[
\mathfrak{S}^{(1)} / p\otimes_{p_j,\mathfrak{S}}\mathfrak{M} \rightarrow \mathfrak{S}^{(1)}[E^{-1}]^{\wedge}_p / p\otimes_{p_j,\mathfrak{S}}\mathfrak{M} = \mathfrak{S}^{(1)}[E^{-1}] / p\otimes_{p_j,\mathfrak{S}}\mathfrak{M} 
\]
is injective. This factors through the map $\mathfrak{S}^{(1)} / p\otimes_{p_j,\mathfrak{S}}\mathfrak{M} \rightarrow \mathfrak{L} / p$, which therefore is injective. This gives a contradiction, and the surjectivity follows.
\end{proof}

\subsection{Projectivity of \texorpdfstring{$\fkM[p^{-1}]$}{M[1/p]} under Assumption~\ref{assumption:base-ring-sec-3.4}} \label{sec:quasi-kisin-mod-projectivity}
In this subsection, assume that either $R$ is small over $\calO_K$ or $R=\calO_L$ (Assumption~\ref{assumption:base-ring-sec-3.4}).
We will show that if $(\fkM,\varphi_\fkM)$ is a finitely generated torsion free $\varphi$-module  of finite $E$-height over $\fkS$, then $\fkM[p^{-1}]$ is projective over $\fkS[p^{-1}]$ (Proposition~\ref{prop:rational-projectivity-etale-over-torus-case}). For this, we need two preliminary results.

\begin{lem}\label{lem-recenter}
Let $k_1$ be a perfect field of characteristic $p$, and let $A$ be a power-series ring $W(k_1)[\![s_1, \dots , s_a]\!]$. Suppose that $A$ is equipped with a Frobenius endomorphism $\varphi$ extending the Witt vector Frobenius on $W(k_1)$. Then there exist $t_1,\ldots,t_a \in A$ such that $A =  W(k_1)[\![t_1, \dots , t_a]\!]$ and $\varphi (t_i)$ has zero constant term for each $i$. 
\end{lem}

\begin{proof}
Write $\varphi (s_i) = s_i^p + p (f_i(s_1, \dots , s_a)) + pb _i $ where $f_i(s_1 , \dots , s_a) \in A$ satisfying $f_i(0, \dots , 0)= 0$ and $b_i \in W(k_1)$. Write $v_p(\cdot)$ for the $p$-adic valuation on $W(k_1)$ with $v_p(p) = 1$. Suppose $b_i \neq 0$ for some $i$, and define $I = \{j \mid v_p(b_j) = \min_{1\leq i \leq a} \{ v_p(b_i)\} \}$. 

Let $i_0 \in I$, and let $c_{i_0} \in W(k_1)$ such that $\varphi(c_{i_0}) = b_{i_0}$. We claim that if we replace $s_{i_0}$ by $s_{i_0}- pc_{i_0}$, then $\varphi (s_i) = s_i^p + p (f_i'(s_1, \dots, s_a)) + pb_i'$ satisfying $v_p (b_i')\geq \min \{v_p (b _i), v_p (b_{i_0})+1\}$ for each $i = 1, \ldots, a$, and $v_p (b_{i_0}') \geq v_p (b_{i_0})+1$ and $v_p(b_i') = v_p(b_i)$ if $i_0 \neq i \in I$. Here, $f_i'$ and $b_i'$ denote the corresponding power series and the constant replacing $f_i$ and $b_i$, respectively. To check the claim, note that
\begin{align*}
\varphi (s_{i_0} -p c_{i _0}) &= s_{i_0}^p + p (f_{i_0}(s_1, \dots , s_a)) + pb _{i_0} - \varphi (p c_{i_0})\\
    &= (s_{i_0}-pc_{i_0} + p c_{i _0})^p + p f_{i_0} (s_1, \dots , s_{i_0}-p c_{i_0} + pc_{i_0}, \dots , s_a).   
\end{align*}
Since $v_p(c_{i_0}) = v_p(b_{i_0})$, we have $v_p (b_{i_0}') \geq  v_p (b_{i_0})+1$. For $i \not = i_0$, we have 
\[
\varphi (s_i)=   s_i^p + p (f_i(s_1, \dots , s_{i_0} - p c_{i _0} + p c_{i_0}, \dots,  s_a)) + pb_i.
\]
So $v_p (b_i')\geq \min\{v_p (c_{i_0}) +1, v_p (b_i)\}$. Furthermore, if $i \in I$ (with $i \neq i_0$), then $v_p(b_i') = v_p(b_i)$. This proves the claim.

Thus, if $\# I \geq 2$, then after replacing $s_{i_0}$ by $s_{i_0}- pc_{i_0}$, $\# I$ decreases by $1$. If $\# I = 1$, then after replacing $s_{i_0}$ by $s_{i_0}- pc_{i_0}$, we have 
\[
\min_{1\leq i \leq a} \{ v_p(b_i')\} \geq 1+\min_{1\leq i \leq a} \{ v_p(b_i)\}.
\] 
By repeating the above process, we deduce that there exist $c_1,\ldots,c_a \in W(k_1)$ such that for $t_i = s_i -pc_i$, $\varphi(t_i)$ has zero constant term for each $i$. It is clear that $A = W(k_1)[\![t_1, \dots , t_a]\!]$.
\end{proof}

\begin{lem} \label{lem:rational-projectivity-power-series-case} 
Let $k_1$ be a perfect field of characteristic $p$ over $k$, and let $A$ be a power-series ring $W(k_1)[\![t_1, \ldots, t_a]\!]$. Suppose that $A$ is equipped with a Frobenius endomorphism $\varphi_A$ extending the Witt vector Frobenius on $W(k_1)$ such that $\varphi_A(t_i) \in A$ has zero constant term for each $i$. Let $\mathfrak{S}_A \coloneqq A[\![u]\!]$ equipped with Frobenius extending that on $A$ by $\varphi(u) = u^p$. Let $\mathfrak{N}$ be a finite $\mathfrak{S}_A$-module equipped with a $\varphi$-semi-linear endomorphism $\varphi_{\mathfrak{N}}\colon \mathfrak{N} \rightarrow \mathfrak{N}$ such that the induced map $1\otimes\varphi_{\mathfrak{N}}\colon (\mathfrak{S}_A\otimes_{\varphi, \mathfrak{S}_A} \mathfrak{N})[E(u)^{-1}] \rightarrow \mathfrak{N}[E(u)^{-1}]$ is an isomorphism. Then $\mathfrak{N}[p^{-1}]$ is projective over $\mathfrak{S}_A[p^{-1}]$.
\end{lem}

\begin{proof}
When $a = 0$ (i.e., $A = W(k_1)$), the statement is proved in \cite[Prop.~4.3]{bhatt-morrow-scholze-integralpadic}. For the general case, we prove by reducing to the case $a = 0$ as follows. Suppose $a \geq 1$. Let $J$ be the non-zero Fitting ideal of $\mathfrak{N}$ over $\mathfrak{S}_A$ with the smallest index. It suffices to show that $J\mathfrak{S}_A[p^{-1}] = \mathfrak{S}_A[p^{-1}]$. Assume the contrary. Since Fitting ideals are compatible under base change, we have
\begin{equation} \label{eq:fitting-ideal-1}
J\mathfrak{S}_A[E(u)^{-1}] = \varphi_{\mathfrak{S}_A}(J)\mathfrak{S}_A[E(u)^{-1}]	
\end{equation}
as ideals of $\mathfrak{S}_A[E(u)^{-1}]$, and so
\begin{equation} \label{eq:fitting-ideal}
(\mathfrak{S}_A/J)[E(u)^{-1}] = (\mathfrak{S}_A/\varphi_{\mathfrak{S}_A}(J))[E(u)^{-1}].	
\end{equation}

Write $K_1 = W(k_1)[p^{-1}]$. Let $B$ be the rigid analytic open unit ball in coordinates $(t_1, \ldots, t_a, u)$. Hence the set of $\overline{K_1}$-valued points of $B$ is given by 
\[
\{(t_1, \ldots, t_a, u) \in \overline{K_1}^{a+1} \mid 0 \leq |t_i|, |u| < 1\},
\]
where we use the $p$-adic norm such that $|p| = p^{-1}$. We have a natural map $\mathfrak{S}_A[p^{-1}] \rightarrow \mathcal{O}_B(B)$ whose image is dense. Note that by \cite[Lem.~7.1.9]{deJong-dieudonnemodule}, we have a functorial bijection between the set of maximal ideals of $\mathfrak{S}_A[p^{-1}]$ and the points of $B$. Moreover, the Frobenius $\varphi_{\mathfrak{S}_A}$ on $\fkS_A$ induces an endomorphism on $B$.

For any real number $c$ with $0 < c < 1$, set
\begin{align*}
M_{c} &\coloneqq \{(x_1, \ldots, x_{a+1}) \in \mathbf{R}^{a+1} \mid 0 \leq x_i \leq c\}\quad\text{and}\\
V_c &\coloneqq \{(x_1, \ldots, x_{a+1}) \in \mathbf{R}^{a+1} \mid 0 \leq x_i < 1 \text{ for } 1\leq i \leq a, ~x_{a+1} = c\}.
\end{align*}

Consider the $\overline{K_1}$-valued points of $\operatorname{Spec}(\mathfrak{S}_A[p^{-1}]/J)$, and let $Z = \{(|t_1|, \ldots, |t_a|, |u|)\}$ be the set of corresponding $(a+1)$-tuple norms. Define
\[
Z' = \{(|t_1|, \ldots, |t_a|, |u|) \in \mathbf{R}^{a+1} \mid 0 \leq |t_i|, |u| < 1,  ~(|\varphi_{A}(t_i)|, |u|^p) \in Z\}. 
\]
By Equation \eqref{eq:fitting-ideal}, we have $Z-V_{|\pi|} = Z'-V_{|\pi|}$. For $i = 1, \ldots, a$, let $y_i$, $t_i \in \overline{K_1}$ with $0 \leq |y_i| < 1$ and $0 \leq |t_i| < 1$ such that $\varphi_A(t_i) = y_i$. Note that by the assumption on $\varphi_A(t_i)$'s,
\[
|y_i| \leq \max\{|t_i|^p, p^{-1}|t_1|, \ldots, p^{-1}|t_a|\}
\]
for each $i$. So we have 
\begin{equation} \label{eq:norm-change-under-phi}
\max_{1 \leq i \leq a}\{|y_i|\} \leq \max_{1 \leq i \leq a}\{|t_i|^p, p^{-1}|t_i|\}.
\end{equation}
 
First we show that $Z$ contains a point with $|u| < |\pi|$. Suppose otherwise. Recall that we assume $J\mathfrak{S}_A[p^{-1}] \neq \mathfrak{S}_A[p^{-1}]$ so that $Z \neq \emptyset$. Since $Z-V_{|\pi|} = Z'-V_{|\pi|}$, we deduce that if $Z\cap V_c \neq \emptyset$, then $c = |\pi|^{p^{-n}}$ for some integer $n \geq 0$. Moreover, the rigid analytic $K_1$-space $(\operatorname{Spf} (\mathfrak{S}_A/J))^{\mathrm{rig}}$ has finitely many connected components since they correspond to the idempotents of the noetherian ring $(\mathfrak{S}_A/J)[p^{-1}]$ (see \cite[para. before Lem.~4.13]{kappen-martin-analytic-functions}). So there exists a finite set of non-negative integers $\{n_1, \ldots, n_m\}$ such that $Z\cap V_c \neq \emptyset$ if and only if $c = |\pi|^{p^{-n_i}}$ for some $i$. Without loss of generality, let $n_1$ be maximal among $\{n_1, \ldots, n_m\}$. Since $Z-V_{|\pi|} = Z'-V_{|\pi|}$, we have $Z \cap V_{|\pi|^{p^{-(n_1+1)}}} \neq \emptyset$, which is a contradiction. Thus, $Z$ contains a point with $|u| < |\pi|$.

Next we show $(0, \ldots, 0) \in Z$, i.e., $J\mathfrak{S}_A[p^{-1}] \subset (t_1, \ldots, t_a, u)\mathfrak{S}_A[p^{-1}]$. Suppose otherwise. Then there exists $f(t_1, \ldots, t_a, u) \in J$ whose constant term is non-zero, and let $b$ be the norm of the constant term. Since $Z$ contains a point with $|u| < |\pi|$, we deduce from $Z-V_{|\pi|} = Z'-V_{|\pi|}$ and the inequality \eqref{eq:norm-change-under-phi} that $Z\cap M_{\epsilon} \neq \emptyset$ for any sufficiently small $\epsilon > 0$. But $|f(t_1, \ldots, t_a, u)| = b > 0$ if $(|t_1|, \ldots, |t_a|, |u|) \in M_{\epsilon}$ for any sufficiently small $\epsilon > 0$, which is a contradiction. Thus, $(0, \ldots, 0) \in Z$.   

On the other hand, we claim $J\mathfrak{S}_A[p^{-1}] \not \subset I\mathfrak{S}_A[p^{-1}]$ where $I = (t_1, \ldots, t_a) \subset \mathfrak{S}_A$. Suppose otherwise. 
Take $n \geq 0$ such that $J' \coloneqq p^nJ$ satisfies $J' \subset I\mathfrak{S}_A$. We show by induction that $J' \subset (p, I)^m \cap I$ (as ideals of $\mathfrak{S}_A$) for each $m \geq 0$. The base case $m = 0$ is clear. Suppose $J' \subset (p, I)^m \cap I$. By Equation \eqref{eq:fitting-ideal-1} and the assumption on $\varphi_A(t_i)$'s, we have
\[
E(u)^s J' \subset \varphi((p , I)^{m}\cap I) \subset (p, I)^{m+1} \cap I
\]
for some integer $s \geq 0$. So it suffices to show that if $f \in \mathfrak{S}_A$ satisfies $E(u)f \in (p, I)^{m+1} \cap I$, then $f \in (p, I)^{m+1} \cap I$. For this, choose a set of generators $g_1, \ldots, g_b \in A=W(k_1)[\![t _1 , \dots, t_a]\!]$ of $(p, I)^{m+1} \cap I$. We have $E(u)f = \sum_{i=1}^b g_i h_i$ for some $h_i \in \mathfrak{S}_A$. Note that we can write
\[
h_i = \sum_{j=0}^{e-1}c_{ij}u^j + E(u)h_i' 
\]
for some $c_{ij} \in  W(k_1)[\![t _1 , \dots, t_a]\!]$ and $h_i' \in \mathfrak{S}_A$. So
\[
E(u)f = \sum_{j=0}^{e-1}(\sum_{i=1}^b c_{ij}g_i)u^j + E(u)\sum_{i=1}^b g_i h_i'.
\]
Setting $u = \pi$ in the above equation, we get $ \sum_{j=0}^{e-1}(\sum_{i=1}^b c_{ij}g_i)\pi^j = 0$ as an element in $\mathcal{O}_{K'}[\![t _1 , \dots, t_a]\!]$ where $K' \coloneqq W(k_1)\otimes_{W(k)} K$. This implies $\sum_{i=1}^b c_{ij}g_i = 0$ for each $j = 0, \ldots, e-1$, and thus $f = \sum_{i=1}^b g_i h_i' \in (p, I)^{m+1} \cap I$. Hence $J' \subset (p, I)^m \cap I$ for each $m \geq 0$. Since $\mathfrak{S}_A$ is $(p, I)$-adically separated, we have $J' = 0$ and thus $J=0$, which is a contradiction. This proves the claim.  

Finally, consider the $\varphi$-equivariant projection $\mathfrak{S}_A \rightarrow \mathfrak{S}_{A_0} \coloneqq \mathfrak{S}_A/I\mathfrak{S}_A \cong W(k_1)[\![u]\!]$. Let $J_0 \subset \mathfrak{S}_{A_0}$ be the image of $J$. Since $J\mathfrak{S}_A[p^{-1}] \not\subset I\mathfrak{S}[p^{-1}]$, we have $J_0 \neq (0)$. Moreover, $J_0\mathfrak{S}_{A_0}[p^{-1}] \neq \mathfrak{S}_{A_0}[p^{-1}]$ since $(0, \ldots, 0) \in Z$. On the other hand, Equation \eqref{eq:fitting-ideal} gives via $\mathfrak{S}_A \rightarrow \mathfrak{S}_{A_0}$ 
\[
(\mathfrak{S}_{A_0}/J_0)[E(u)^{-1}] = (\mathfrak{S}_{A_0}/\varphi_{\mathfrak{S}_{A_0}}(J_0))[E(u)^{-1}],
\]
which gives a contradiction by the case $a = 0$. Hence, $J\mathfrak{S}_A[p^{-1}] = \mathfrak{S}_A[p^{-1}]$.
\end{proof}

Let us return to the discussion on the projectivity of $\fkM[p^{-1}]$.

\begin{prop} \label{prop:rational-projectivity-etale-over-torus-case}
Suppose that $R$ satisfies Assumption~\ref{assumption:base-ring-sec-3.4}: $R$ is small over $\calO_K$ or $R=\calO_L$.
If $(\fkM,\varphi_\fkM)$ is a finitely generated $\varphi$-module over $\fkS$ of finite $E$-height, then $\fkM[p^{-1}]$ is projective over $\fkS[p^{-1}]$.
\end{prop}

\begin{proof}
The case where $R=\calO_L$ follows from \cite[Prop.~4.3]{bhatt-morrow-scholze-integralpadic} for $\calO_{K_g}$ (cf. Notation \ref{notation:L}) and the classically faithful flatness of  $\mathfrak{S}_L[p^{-1}] \rightarrow \mathfrak{S}_g[p^{-1}]\coloneqq \calO_{K_{0,g}}[\![u]\!][p^{-1}]$.
Consider the case where $R_0$ is the $p$-adic completion of an \'etale extension of $W(k)\langle T_1^{\pm 1}, \ldots, T_d^{\pm 1}\rangle$. Note that the Krull dimension of $R_0$ is the same as that of $W(k)\langle T_1^{\pm 1}, \ldots, T_d^{\pm 1}\rangle$. Let $\mathfrak{m} \subset R_0$ be any maximal ideal, and let $(R_0)_{\mathfrak{m}}^{\wedge}$ denote the $\mathfrak{m}$-adic completion of the localization $(R_0)_{\mathfrak{m}}$. Since $\mathfrak{m} \cap W(k)\langle T_1^{\pm 1}, \ldots, T_d^{\pm 1}\rangle$ is a maximal ideal of $W(k)\langle T_1^{\pm 1}, \ldots, T_d^{\pm 1}\rangle$, the residue field $k_{\mathfrak{m}} \coloneqq (R_0)_{\mathfrak{m}}^{\wedge} / \mathfrak{m}(R_0)_{\mathfrak{m}}^{\wedge}$ is a finite extension of $k$. Note that since $p \in \mathfrak{m}$ and $\varphi(\mathfrak{m}) \subset \mathfrak{m}$, $(R_0)_{\mathfrak{m}}^{\wedge}$ is equipped with the Frobenius induced from $R_0$. Let $f\colon W(k) \rightarrow (R_0)_{\mathfrak{m}}^{\wedge}$ be the composite  $W(k) \rightarrow W(k)\langle T_1^{\pm 1}, \ldots, T_d^{\pm 1}\rangle \rightarrow (R_0)_{\mathfrak{m}}^{\wedge}$, which is compatible with $\varphi$. Since $W(k) \rightarrow W(k_{\mathfrak{m}})$ is \'etale, $f$ factors uniquely through $W(k) \rightarrow W(k_{\mathfrak{m}}) \rightarrow (R_0)_{\mathfrak{m}}^{\wedge}$. By unicity, $W(k_{\mathfrak{m}}) \rightarrow (R_0)_{\mathfrak{m}}^{\wedge}$ is compatible with $\varphi$. Since $p \notin \mathfrak{m}^2$, $\{p\}$ can be extended to a minimal set generating $\mathfrak{m}$, and the map $W(k_{\mathfrak{m}}) \rightarrow (R_0)_{\mathfrak{m}}^{\wedge}$ extends to an isomorphism 
\[
W(k_{\mathfrak{m}})[\![t_1, \ldots, t_d]\!] \xrightarrow{\cong} (R_0)_{\mathfrak{m}}^{\wedge}.
\]
Furthermore, by Lemma~\ref{lem-recenter}, $t_1, \ldots, t_d$ can be chosen such that $\varphi(t_i)$ has zero constant term for each $i$ (where $\varphi$ on $W(k_{\mathfrak{m}})[\![t_1, \ldots, t_d]\!]$ is given by the above isomorphism).

Now, let $\mathfrak{P} \subset \mathfrak{S}[p^{-1}]$ be any maximal ideal. Then the prime ideal $\mathfrak{q} = \mathfrak{S} \cap \mathfrak{P}$ is maximal among the prime ideals of $\mathfrak{S}$ not containing $p$. Thus, $\mathfrak{n} \coloneqq \sqrt{\mathfrak{q}+p\mathfrak{S}}$ is a maximal ideal of $\mathfrak{S}$. Let $\mathfrak{S}_{\mathfrak{n}}^{\wedge}$ be the $\mathfrak{n}$-adic completion of the localization $\mathfrak{S}_{\mathfrak{n}}$. By the above discussion, $\mathfrak{S}_{\mathfrak{n}}^{\wedge} \cong W(k_{\mathfrak{n}})[\![t_1, \ldots, t_d]\!][\![u]\!]$ for some finite extension $k_{\mathfrak{n}}$ of $k$ and $t_1, \ldots, t_d$ such that $ W(k_{\mathfrak{n}}) \hookrightarrow \mathfrak{S}_{\mathfrak{n}}^{\wedge}$ is compatible with $\varphi$ and $\varphi(t_i)$ has zero constant term for each $i$. 

Let $\mathfrak{M}_{\mathfrak{n}} \coloneqq  \mathfrak{S}_{\mathfrak{n}}^{\wedge}\otimes_{\mathfrak{S}}\mathfrak{M}$ equipped with the induced tensor-product Frobenius. By Lemma~\ref{lem:rational-projectivity-power-series-case}, $\mathfrak{M}_{\mathfrak{n}}[p^{-1}]$ is projective over $\mathfrak{S}_{\mathfrak{n}}^{\wedge}[p^{-1}]$. Let $\mathfrak{P}_{\mathfrak{n}} \subset \mathfrak{S}_{\mathfrak{n}}^{\wedge}[p^{-1}]$ be a maximal ideal lying over $\mathfrak{P} \subset \mathfrak{S}[p^{-1}]$. Note that the natural map on localizations
\[
(\mathfrak{S}[p^{-1}])_{\mathfrak{P}} \rightarrow (\mathfrak{S}_{\mathfrak{n}}^{\wedge}[p^{-1}])_{\mathfrak{P}_{\mathfrak{n}}}
\]
is classically faithfully flat. Since $(\mathfrak{S}_{\mathfrak{n}}^{\wedge}[p^{-1}])_{\mathfrak{P}_{\mathfrak{n}}}\otimes_{\mathfrak{S}_{\mathfrak{n}}^{\wedge}[p^{-1}]}\mathfrak{M}_{\mathfrak{n}}[p^{-1}]$ is finite projective over $(\mathfrak{S}_{\mathfrak{n}}^{\wedge}[p^{-1}])_{\mathfrak{P}_{\mathfrak{n}}}$, we deduce that $(\mathfrak{S}[p^{-1}])_{\mathfrak{P}}\otimes_{\mathfrak{S}[p^{-1}]}\mathfrak{M}[p^{-1}]$ is projective over $(\mathfrak{S}[p^{-1}])_{\mathfrak{P}}$. This holds for any maximal ideal $\mathfrak{P} \subset \mathfrak{S}[p^{-1}]$, so $\mathfrak{M}[p^{-1}]$ is projective over $\mathfrak{S}[p^{-1}]$.
\end{proof}

\subsection{Crystalline representations and Breuil--Kisin modules in the CDVR case} \label{sec:cryst-rep-cdvr}

We follow Notation~\ref{notation:L}. In particular, recall that $\mathcal{O}_L$ denotes the $p$-adic completion of $R_{(\pi)}$. Then $L$ is a complete discrete valuation field whose residue field has a finite $p$-basis given by $\{T_1, \ldots, T_d\}$. We first consider crystalline representations of $\mathcal{G}_{\mathcal{O}_L}=\Gal(\overline{L}/L)$, and study certain properties of the associated Breuil--Kisin modules. By abuse of notation, we also write $G_L$ and $G_{\widetilde{L}_{\infty}}$ for the Galois groups $\mathcal{G}_{\mathcal{O}_L}$ and $\mathcal{G}_{\widetilde{\mathcal{O}}_{L, \infty}}$, respectively (see \eqref{eq:R-infinity} for the definition of $\calO_{\widetilde{L}_\infty}$). 
 
Fix a non-negative integer $r$. Let $V$ be a crystalline $\mathbf{Q}_p$-representation of $G_{L}$ with Hodge--Tate weights in $[0, r]$. 
By \cite[Prop.~4.17]{brinon-trihan}, there exists an $\mathfrak{S}_L$-module $\mathfrak{M}_L$ satisfying the following properties:
\begin{itemize}
\item $\mathfrak{M}_L$ is finite free over $\mathfrak{S}_L$;
\item $\mathfrak{M}_L$ is equipped with a $\varphi$-semi-linear endomorphism $\varphi_{\mathfrak{M}_L}\colon \mathfrak{M}_L \rightarrow \mathfrak{M}_L$ with $E$-height $\leq r$;
\item Set 
\[
M_L \coloneqq \mathcal{O}_{L_0}\otimes_{\varphi,\mathcal{O}_{L_0}}\mathfrak{M}_L/u\mathfrak{M}_L
\]
and equip it with the induced tensor-product Frobenius. We have a natural isomorphism of $L_0$-modules $M_L[p^{-1}] \cong D_{\mathrm{cris}}^{\vee}(V)$ compatible with Frobenii. 
Via this isomorphism, $M_L[p^{-1}]$ admits a topologically quasi-nilpotent connection $\nabla_{\fkM_L}$.
\end{itemize}
We call the triple $(\fkM_L,\varphi_{\fkM_L},\nabla_{\fkM_L})$ the \emph{Breuil--Kisin module associated with $V$}. Note that \cite{brinon-trihan} considers $\fkM_L/u\fkM_L$ instead of the Frobenius pullback $M_L$. However, we have a natural isomorphism of $L_0$-modules $M_L[p^{-1}] \xrightarrow{\cong}(\fkM_L/u\fkM_L)[p^{-1}]$ compatible with Frobenii. Following \cite{kim-groupscheme-relative}, we use $M_L$ since it is more suitable when we consider the filtration.

Let $S_L$ be the $p$-adically completed divided power envelope of $\mathfrak{S}_L$ with respect to $(E(u))$. The Frobenius on $\mathfrak{S}_L$ extends uniquely to $S_L$. For each integer $i \geq 0$, let $\mathrm{Fil}^i S_L$ be the PD-filtration of $S_L$ as before. Let $N_u\colon S_L \rightarrow S_L$ be the $\mathcal{O}_{L_0}$-linear derivation given by $N_u(u) = -u$. We also have a natural integrable connection $\nabla\colon S_L \rightarrow S_L\otimes_{\mathcal{O}_{L_0}} \widehat{\Omega}_{\mathcal{O}_{L_0}}$ given by the universal derivation on $\mathcal{O}_{L_0}$, which commutes with $N_u$. 

Set
\[
\mathscr{M}_L \coloneqq S_L\otimes_{\varphi, \mathfrak{S}_L} \mathfrak{M}_L
\]
equipped with the induced Frobenius. If
we let $q\colon S_L \twoheadrightarrow \mathcal{O}_{L_0}$ denote the $\varphi$-compatible projection given by $u \mapsto 0$, it induces the projection $q\colon \mathscr{M}_L \twoheadrightarrow M_L$.

We define two filtrations on $\mathscr{M}_L[p^{-1}]$ and study their compatibility. Let 
\[
\mathscr{D}_L \coloneqq S_L[p^{-1}]\otimes_{L_0} {D_{\mathrm{cris}}^{\vee}(V)}. 
\]
By the above isomorphism $M_L[p^{-1}] \cong D_{\mathrm{cris}}^{\vee}(V)$ and Lemma~\ref{lem:frobenius-compatible-section}, we have a $\varphi$-equivariant identification $\mathscr{M}_L[p^{-1}]=\mathscr{D}_L$. Let $N_u\colon \mathscr{D}_L \rightarrow \mathscr{D}_L$ be the $L_0$-linear derivation given by $N_{u, S_L}\otimes 1$, and let $\nabla\colon \mathscr{D}_L \rightarrow \mathscr{D}_L\otimes_{\mathcal{O}_{L_0}} \widehat{\Omega}_{\mathcal{O}_{L_0}}$ be the connection given by $\nabla_{S_L}\otimes 1+1\otimes \nabla_{D_{\mathrm{cris}}^{\vee}(V)}$. Define a decreasing filtration on $\mathscr{D}_L$ by $S_L[p^{-1}]$-submodules $\mathrm{Fil}^i \mathscr{D}_L$, inductively as follows: $\mathrm{Fil}^0 \mathscr{D}_L = \mathscr{D}_L$ and 
\[
\mathrm{Fil}^{i+1} \mathscr{D}_L = \{x \in \mathscr{D}_L \mid N_u(x) \in \mathrm{Fil}^i \mathscr{D}_L, ~q_{\pi}(x) \in \mathrm{Fil}^{i+1} (L\otimes_{L_0}D_{\mathrm{cris}}^{\vee}(V))\},
\]
where $q_{\pi}\colon \mathscr{D}_L \rightarrow L\otimes_{L_0}D_{\mathrm{cris}}^{\vee}(V)$ is the map induced by $S_L[p^{-1}] \rightarrow L, ~u \mapsto \pi$. The following is proved in \cite{moon-strly-div-latt-cryst-cohom-CDVF}.
 Note that \cite[\S 4.1]{moon-strly-div-latt-cryst-cohom-CDVF} assumes $p > 2$ and $r\leq p-2$, but the results we will cite in this subsection hold without these assumptions.

\begin{lem}[({\cite[Lem.~4.2]{moon-strly-div-latt-cryst-cohom-CDVF}})] \label{lem:CDVF-Griffiths-transversality}
The connection $\nabla$ on $\mathscr{D}_L$ satisfies the Griffiths transversality:
\[
\nabla(\mathrm{Fil}^{i+1} \mathscr{D}_L) \subset \mathrm{Fil}^{i} \mathscr{D}_L \otimes_{\mathcal{O}_{L_0}} \widehat{\Omega}_{\mathcal{O}_{L_0}}.
\]	
\end{lem}

For the second filtration, let
\[
\mathrm{F}^i \mathscr{M}_L[p^{-1}] \coloneqq \{x \in \mathscr{M}_L[p^{-1}] \mid (1\otimes \varphi_{\mathfrak{M}_L})(x) \in (\mathrm{Fil}^i S_L[p^{-1}])\otimes_{\mathfrak{S}_L} \mathfrak{M}_L\}.
\]
We will see that these two filtrations coincide under the identification $\mathscr{M}_L[p^{-1}]=\mathscr{D}_L$ and thus $\nabla_{\fkM_L}$ satisfies the $S_L$-Griffiths transversality. For this, consider the base change along $\mathcal{O}_{L_0} \rightarrow W(k_g)$ as in Notation~\ref{notation:L}. Note that $W(k_g)$ is a complete discrete valuation ring with \emph{perfect} residue field.
Let $S_g$ be the $p$-adically completed divided power envelope of $\mathfrak{S}_g \coloneqq W(k_g)[\![u]\!]$ with respect to $(E(u))$. It is equipped with $\varphi$, PD-filtration, and $N_u$ similarly as above. 
Let 
\[
\mathfrak{M}_g = \mathfrak{S}_g\otimes_{\mathfrak{S}_L}\mathfrak{M}_L, \quad
\mathscr{M}_g = S_g\otimes_{\varphi, \mathfrak{S}_g}\mathfrak{M}_g, \quad\text{and}\quad
\mathscr{D}_g = S_g[p^{-1}]\otimes_{W(k_g)}D_{\mathrm{cris}}^{\vee}(V|_{G_{K_g}}).
\]
We can identify $\mathscr{M}_g[p^{-1}] = \mathscr{D}_g$ compatibly with $\varphi$, and define two filtrations $\mathrm{Fil}^i \mathscr{D}_g$ and $\mathrm{F}^i \mathscr{M}_g[p^{-1}]$ similarly as above. By the proof of \cite[Cor.~3.2.3]{liu-semistable-lattice-breuil}, we have
\[
\mathrm{Fil}^i\mathscr{D}_g = \mathrm{F}^i \mathscr{M}_g[p^{-1}].
\]
Note also $K_{0,g}\otimes_{L_0}D_{\cris}^\vee(V)\xrightarrow{\cong}D_{\cris}^\vee(V|_{G_{K_g}})$ by \cite[4B]{Ohkubo}.

\begin{lem} \label{lem:CDVF-filtration-compatibility}
Under the $\varphi$-equivariant identification $\mathscr{M}_L[p^{-1}] = \mathscr{D}_L$, we have
\[
\mathrm{F}^i \mathscr{D}_L = \mathrm{Fil}^i \mathscr{D}_L.
\]	
In particular, the triple $(\fkM_L,\varphi_{\fkM_L},\nabla_{\fkM_L})$ is a quasi-Kisin module of $E$-height $\leq r$ over $\fkS_L$.
\end{lem}

\begin{proof}
We consider $\mathscr{D}_L$ as a $S_L[p^{-1}]$-submodule of $\mathscr{D}_{g}$ via $S_g\otimes_{S_L} \mathscr{D}_L = \mathscr{D}_{g}$. 
Recall that $\pi$ is a uniformizer of $\calO_L$, and let $e=[L:L_0]$.
Note that any $x \in S_L$ can be written as $x = \sum_{i\geq 0} \frac{E(u)^i}{i!} (\sum_{j=0}^{e-1} a_{ij}u^j)$ for some $a_{ij} \in \mathcal{O}_{L_0}$ (with $a_{ij} \rightarrow 0$ $p$-adically as $i \rightarrow \infty$). Furthermore, the $a_{ij}$'s can be seen to be uniquely determined by inductively setting $u = \pi$.  The analogous statement holds for the elements in $S_g$, and thus we have $S_L \cap \mathrm{Fil}^i S_{g} = \mathrm{Fil}^i S_L$.

Let $x \in \mathrm{Fil}^i \mathscr{D}_L$. Since $\mathrm{Fil}^i \mathscr{D}_L \subset \mathrm{Fil}^i \mathscr{D}_{g} = \mathrm{F}^i \mathscr{D}_{g}$, we have 
\[
(1\otimes \varphi)(x) \in (\mathrm{Fil}^i S_{g}[p^{-1}])\otimes_{\mathfrak{S}_L} \mathfrak{M}_L.
\] 	
Since $S_L[p^{-1}] \cap \mathrm{Fil}^i S_{g}[p^{-1}] = \mathrm{Fil}^i S_L[p^{-1}]$ and $\mathfrak{M}_L[p^{-1}]$ is projective over $\mathfrak{S}_L[p^{-1}]$, we deduce $(1\otimes \varphi)(x) \in (\mathrm{Fil}^i S_L[p^{-1}])\otimes_{\mathfrak{S}_L} \mathfrak{M}_L$ by Lemma~\ref{lem:intersection-modules-flat-base-change} (i). Thus, $\mathrm{Fil}^i \mathscr{D}_L \subset \mathrm{F}^i \mathscr{D}_L$.

Conversely, let $x \in \mathrm{F}^i \mathscr{D}_L$. Note that 
\[
(L\otimes_{L_0}D_{\mathrm{cris}}^{\vee}(V)) \cap \mathrm{Fil}^i (K_g\otimes_{W(k_g)}D_{\mathrm{cris}}^{\vee}(V|_{G_{K_g}})) = \mathrm{Fil}^i (L\otimes_{L_0}D_{\mathrm{cris}}^{\vee}(V)).
\]
Hence we deduce by induction on $i$ that $\mathscr{D}_L \cap \mathrm{Fil}^i \mathscr{D}_{g} = \mathrm{Fil}^i \mathscr{D}_L$. Since $\mathrm{F}^i \mathscr{D}_L \subset \mathrm{F}^i \mathscr{D}_{g} = \mathrm{Fil}^i \mathscr{D}_{g}$, we have $x \in \mathscr{D}_L \cap \mathrm{Fil}^i \mathscr{D}_{g} = \mathrm{Fil}^i \mathscr{D}_L$. Hence, $\mathrm{F}^i \mathscr{D}_L \subset \mathrm{Fil}^i \mathscr{D}_L$. The second assertion follows from the first and Lemma~\ref{lem:CDVF-Griffiths-transversality}: 
$N_u(\mathrm{Fil}^{i+1} \mathscr{D}_L) \subset \mathrm{Fil}^{i} \mathscr{D}_L$ by definition, and it is straightforward to check $\partial_u(\mathrm{Fil}^{i+1} \mathscr{D}_L) \subset \mathrm{Fil}^{i} \mathscr{D}_L$ by induction.
\end{proof}

Next we will explain how to recover $V$ from $\mathscr{D}_L$ as a representation of $G_L$.
Note that the embedding $\mathfrak{S}_L \rightarrow W(\mathcal{O}_{\overline{L}}^{\flat})$ given in \S~\ref{sec:etale phi-module} extends to $S_L \rightarrow \mathbf{A}_{\mathrm{cris}}(\mathcal{O}_{\overline{L}})$, which is compatible with $\varphi$, filtrations, and $G_{\widetilde{L}_{\infty}}$-actions. Consider an $\mathbf{A}_{\mathrm{cris}}(\mathcal{O}_{\overline{L}})$-semi-linear $G_{\widetilde{L}_{\infty}}$-action on $\mathbf{A}_{\mathrm{cris}}(\mathcal{O}_{\overline{L}})[p^{-1}]\otimes_{S_L} \mathscr{D}_L$ given by the $G_{\widetilde{L}_{\infty}}$-action on $\mathbf{A}_{\mathrm{cris}}(\mathcal{O}_{\overline{L}})[p^{-1}]$ and the trivial $G_{\widetilde{L}_{\infty}}$-action on $\mathscr{D}_L$. We will extend this action to a $G_L$-action as follows (see \cite[\S 4]{moon-strly-div-latt-cryst-cohom-CDVF}).
For each $i = 1, \ldots, d$, write $N_{T_i}\colon \mathscr{D}_L \rightarrow \mathscr{D}_L$ for the derivation given by $\nabla\colon \mathscr{D}_L \rightarrow \mathscr{D}_L\otimes_{\mathcal{O}_{L_0}} \widehat{\Omega}_{\mathcal{O}_{L_0}}\cong \bigoplus_{i = 1}^d \mathscr{D}_L \cdot d\log{T_i}$ composed with the projection to the $i$-th factor. Note $N_{T_i} = T_i\partial_{T_i}$ for the derivation $\partial_{T_i}\colon \mathscr{D}_L \rightarrow \mathscr{D}_L$ in \S~\ref{sec:quasi-kisin-mod-rational-descent-data}. For $\sigma \in G_L$, denote
\[
\underline{\varepsilon}(\sigma) \coloneqq \frac{\sigma([\pi^\flat])}{[\pi^\flat]} \quad\text{and }\quad\underline{\mu_i}(\sigma) \coloneqq \frac{\sigma([T_i^\flat])}{[T_i^\flat]} \quad(i = 1, \ldots, d).
\]
Note that $\log(\underline{\varepsilon}(\sigma))$ and $\log(\underline{\mu_i}(\sigma))$ lie in $\mathrm{Fil}^1 \mathbf{A}_{\mathrm{cris}}(\mathcal{O}_{\overline{L}})$. 
For any element $a\otimes x\in \mathbf{A}_{\mathrm{cris}}(\mathcal{O}_{\overline{L}})[p^{-1}]\otimes_{S_L} \mathscr{D}_L$, define
\begin{equation} \label{eq:CDVF-Galois-action}
\sigma(a\otimes x) = \sum \sigma(a) \gamma_{i_0}(-\log(\underline{\varepsilon}(\sigma)))\gamma_{i_1}(\log(\underline{\mu_1}(\sigma)))\cdots \gamma_{i_d}(\log(\underline{\mu_d}(\sigma))) \cdot N_u^{i_0}N_{T_1}^{i_1}\cdots N_{T_d}^{i_d}(x)
\end{equation}
where the sum ranges over the multi-index $(i_0, i_1, \ldots, i_d)$ of non-negative integers. This sum converges since $\nabla_{\mathscr{D}_L}$ is topologically quasi-nilpotent and since $\gamma_j(-\log(\underline{\varepsilon}(\sigma)))$, $\gamma_j(\log(\underline{\mu_i}(\sigma))) \rightarrow 0$ $p$-adically as $j \rightarrow \infty$. It follows from standard computations that this gives a well-defined $\mathbf{A}_{\mathrm{cris}}(\mathcal{O}_{\overline{L}})$-semi-linear $G_L$-action compatible with $\varphi$. This $G_L$-action preserves the filtration since $\log(\underline{\varepsilon}(\sigma))$, $\log(\underline{\mu_i}(\sigma))\in \mathrm{Fil}^1 \mathbf{A}_{\mathrm{cris}}(\mathcal{O}_{\overline{L}})$ and since $N_u$ and $\nabla$ satisfy the Griffiths transversality by definition and Lemma~\ref{lem:CDVF-Griffiths-transversality}.

Let 
\[
V(\mathscr{D}_L) \coloneqq \mathrm{Hom}_{S_L, \mathrm{Fil}, \varphi}(\mathscr{D}_L, \mathbf{A}_{\mathrm{cris}}(\mathcal{O}_{\overline{L}})[p^{-1}]). 
\]
Using the identification 
\[
\mathrm{Hom}_{S_L, \mathrm{Fil}, \varphi}(\mathscr{D}_L, \mathbf{A}_{\mathrm{cris}}(\mathcal{O}_{\overline{L}})[p^{-1}]) = \mathrm{Hom}_{S_L, \mathrm{Fil}, \varphi}(\mathbf{A}_{\mathrm{cris}}(\mathcal{O}_{\overline{L}})[p^{-1}]\otimes_{S_L}\mathscr{D}_L, \mathbf{A}_{\mathrm{cris}}(\mathcal{O}_{\overline{L}})[p^{-1}]),
\]
we define the $G_L$-action on $V(\mathscr{D}_L)$ by 
setting $\sigma(f)(x) = \sigma(f(\sigma^{-1}(x)))$ for any $x \in \mathbf{A}_{\mathrm{cris}}(\mathcal{O}_{\overline{L}})[p^{-1}]\otimes_{S_L}\mathscr{D}_L$.

\begin{prop}[(cf.~{\cite[\S 4]{moon-strly-div-latt-cryst-cohom-CDVF}})]\label{prop-VDL=V}
There is a natural $G_L$-equivariant isomorphism 
\[
V(\mathscr{D}_L) \cong V.
\]
\end{prop}

\begin{proof}
This is proved in \cite[\S 4]{moon-strly-div-latt-cryst-cohom-CDVF}, and we sketch the proof here. We first study how the above constructions are related to \'etale $\varphi$-modules. If we let $\mathcal{M}_L =  \mathcal{O}_{\mathcal{E}, L}\otimes_{\mathfrak{S}_L}\mathfrak{M}_L$ with the induced $\varphi$, then $\mathcal{M}_L$ is an \'etale $\varphi$-module over $\mathcal{O}_{\mathcal{E}, L}$. Consider the $G_{\widetilde{L}_{\infty}}$-equivariant map
\[
\mathrm{Hom}_{\mathfrak{S}_L, \varphi}(\mathfrak{M}_L, \widehat{\mathfrak{S}}^{\mathrm{ur}}_L) \rightarrow T^{\vee}(\mathcal{M}_L) = \mathrm{Hom}_{\mathcal{O}_{\mathcal{E}, L}, \varphi}(\mathcal{M}_L, \widehat{\mathcal{O}}_{\mathcal{E}, L}^{\mathrm{ur}})
\]
induced by the embedding $\widehat{\mathfrak{S}}^{\mathrm{ur}}_L \rightarrow \widehat{\mathcal{O}}_{\mathcal{E}, L}^{\mathrm{ur}}$. By \cite[Lem.~4.6]{moon-strly-div-latt-cryst-cohom-CDVF}, this map is an isomorphism.

The embedding $\varphi\colon \widehat{\mathfrak{S}}^{\mathrm{ur}}_L \rightarrow \mathbf{A}_{\mathrm{cris}}(\mathcal{O}_{\overline{L}})$ induces a natural $G_{\widetilde{L}_\infty}$-equivariant injective map
\[
\mathrm{Hom}_{\mathfrak{S}_L, \varphi}(\mathfrak{M}_L, \widehat{\mathfrak{S}}^{\mathrm{ur}}_L)[p^{-1}] \rightarrow V(\mathscr{D}_L)
\]
by Lemma~\ref{lem:frobenius-compatible-section}. On the other hand, any $f \in V(\mathscr{D}_L)$ induces a $\varphi$-equivariant map $f'\colon D_{\mathrm{cris}}^{\vee}(V) \rightarrow \mathbf{B}_{\mathrm{cris}}(\mathcal{O}_{\overline{L}})$ via the map $D_{\mathrm{cris}}^{\vee}(V) \rightarrow \mathscr{D}_L$. 
We see that $f'$ is also compatible with filtration, since $\mathbf{B}_{\mathrm{cris}}(\mathcal{O}_{\overline{L}}) \cong \mathbf{B}_{\mathrm{cris}}(\mathcal{O}_{\overline{K_g}})$ and the induced map $D_{\mathrm{cris}}^{\vee}(V|_{G_{K_g}}) \rightarrow \mathbf{B}_{\mathrm{cris}}(\mathcal{O}_{\overline{K_g}})$ is compatible with filtration by the proof of \cite[Lem.~8.1.2]{breuil-representations} and  \cite[\S 3.4]{liu-semistable-lattice-breuil}. So we obtain a natural injective map
\[
V(\mathscr{D}_L) \rightarrow \mathrm{Hom}_{\mathrm{Fil}, \varphi}(D_{\mathrm{cris}}^{\vee}(V), \mathbf{B}_{\mathrm{cris}}(\mathcal{O}_{\overline{L}})).
\]
Since $\mathrm{Hom}_{\mathfrak{S}_L, \varphi}(\mathfrak{M}_L, \widehat{\mathfrak{S}}^{\mathrm{ur}}_L)[p^{-1}]$ and $V$ are $\Q_p$-vector spaces of the same dimension, it suffices to show that $\mathrm{Hom}_{\mathrm{Fil}, \varphi}(D_{\mathrm{cris}}^{\vee}(V), \mathbf{B}_{\mathrm{cris}}(\mathcal{O}_{\overline{L}}))$ admits a $G_L$-action compatibly with $V(\mathscr{D}_L)$ and there exists a natural isomorphism $\mathrm{Hom}_{\mathrm{Fil}, \varphi}(D_{\mathrm{cris}}^{\vee}(V), \mathbf{B}_{\mathrm{cris}}(\mathcal{O}_{\overline{L}})) \cong V$ as $G_L$-representations. 

Write $D = D_{\mathrm{cris}}^{\vee}(V)$ for simplicity. 
Lemma~\ref{lem:cryst-period-ring} yields a $\mathbf{B}_{\mathrm{cris}}(\mathcal{O}_{\overline{L}})$-linear isomorphism
\[
\mathbf{B}_{\mathrm{cris}}(\mathcal{O}_{\overline{L}})\{X_1, \ldots, X_d\} \cong \mathbf{OB}_{\mathrm{cris}}(\mathcal{O}_{\overline{L}})
\]
sending $X_i$ to $T_i\otimes 1-1\otimes [T_i^\flat]$. The projection
\[
\mathrm{pr}\colon \mathbf{OB}_{\mathrm{cris}}(\mathcal{O}_{\overline{L}}) \rightarrow \mathbf{B}_{\mathrm{cris}}(\mathcal{O}_{\overline{L}})
\]
given by $X_i = T_i-[T_i^\flat] \mapsto 0$ induces the projection
\[
\mathrm{pr}\colon\mathbf{OB}_{\mathrm{cris}}(\mathcal{O}_{\overline{L}}) \otimes_{\iota_1,L_0}D  \rightarrow  \mathbf{B}_{\mathrm{cris}}(\mathcal{O}_{\overline{L}})\otimes_{\iota_2,L_0}D
\] 
compatible with Frobenii and filtrations (after tensoring with $L$ over $L_0$). Here, $\iota_1\colon L_0 \rightarrow  \mathbf{OB}_{\mathrm{cris}}(\mathcal{O}_{\overline{L}})$ is given by $T_i \mapsto T_i\otimes 1$, and $\iota_2\colon L_0 \rightarrow \mathbf{B}_{\mathrm{cris}}(\mathcal{O}_{\overline{L}})$ is given by $T_i \mapsto [T_i^\flat]$.

We define a $\mathbf{B}_{\mathrm{cris}}(\mathcal{O}_{\overline{L}})$-linear section $s$ to $\mathrm{pr}$ as follows. For $x \in D$, let
\[
s(x) = \sum (-1)^{i_1+\cdots+ i_d} \gamma_{i_1}(\log{(\frac{T_1}{[T_1^\flat]})})\cdots \gamma_{i_d}(\log{(\frac{T_d}{[T_d^\flat]})}) N_{T_{1}}^{i_1}\cdots N_{T_d}^{i_d} (x)
\]
where the sum ranges over the multi-index $(i_1, \ldots, i_d)$ of non-negative integers.The map $s$ is a well-defined section, and it induces an isomorphism
\[
s\colon  \mathbf{B}_{\mathrm{cris}}(\mathcal{O}_{\overline{L}})\otimes_{\iota_2,L_0}D \xrightarrow{\cong} ( \mathbf{OB}_{\mathrm{cris}}(\mathcal{O}_{\overline{L}})\otimes_{\iota_1, L_0}D)^{\nabla = 0}
\] 
of $\mathbf{B}_{\mathrm{cris}}(\mathcal{O}_{\overline{L}})$-modules, compatibly with filtrations and $\varphi$ (see \cite[\S 4.1]{moon-strly-div-latt-cryst-cohom-CDVF}).  Moreover, if we define $G_L$-action on $\mathbf{B}_{\mathrm{cris}}(\mathcal{O}_{\overline{L}})\otimes_{\iota_2,L_0}D $ by 
\[
\sigma(a\otimes x) = \sum \sigma(a)\gamma_{i_1}(\log(\underline{\mu_1}(\sigma)))\cdots \gamma_{i_d}(\log(\underline{\mu_d}(\sigma))) \cdot N_{T_1}^{i_1}\cdots N_{T_d}^{i_d}(x)
\]  
for $\sigma \in G_L$ and $a\otimes x \in \mathbf{B}_{\mathrm{cris}}(\mathcal{O}_{\overline{L}})\otimes_{\iota_2,L_0}D $, then by \textit{loc. cit.}, the map 
\[V(\mathscr{D}_L) \rightarrow \mathrm{Hom}_{\mathrm{Fil}, \varphi}(D, \mathbf{B}_{\mathrm{cris}}(\mathcal{O}_{\overline{L}})) = \mathrm{Hom}_{\mathrm{Fil}, \varphi}( \mathbf{B}_{\mathrm{cris}}(\mathcal{O}_{\overline{L}})\otimes_{\iota_2,L_0}D, \mathbf{B}_{\mathrm{cris}}(\mathcal{O}_{\overline{L}}))
\] 
is $G_L$-equivariant, and $s$ induces a $G_L$-equivariant isomorphism
\[
\mathrm{Hom}_{\mathrm{Fil}, \varphi}(D, \mathbf{B}_{\mathrm{cris}}(\mathcal{O}_{\overline{L}})) \cong \mathrm{Hom}_{\mathrm{Fil}, \varphi, \nabla}(D, \mathbf{OB}_{\mathrm{cris}}(\mathcal{O}_{\overline{L}})) = V.
\]
This shows that $V(\mathscr{D}_L) \cong V$ as representations of $G_L$.
\end{proof}

Since $\mathrm{Hom}_{\mathfrak{S}_L, \varphi}(\mathfrak{M}_L, \widehat{\mathfrak{S}}^{\mathrm{ur}}_L)[p^{-1}] \cong V$ as $G_{\widetilde{L}_\infty}$-representations, we have a natural map $\mathfrak{M}_L \rightarrow \widehat{\mathfrak{S}}^{\mathrm{ur}}_L\otimes_{\mathbf{Z}_p}V^{\vee}$.
Via the embedding $ \widehat{\mathfrak{S}}^{\mathrm{ur}} \xrightarrow{\varphi} \mathbf{B}_{\mathrm{cris}}(\mathcal{O}_{\overline{L}}) \hookrightarrow \mathbf{OB}_{\mathrm{cris}}(\mathcal{O}_{\overline{L}})$, this induces a map $\mathscr{M}_L[p^{-1}] \rightarrow V^{\vee}\otimes_{\mathbf{Q}_p} \mathbf{OB}_{\mathrm{cris}}(\mathcal{O}_{\overline{L}})$. Composing this with the section $M_L[p^{-1}] \rightarrow \mathscr{M}_L[p^{-1}]$ in Lemma~\ref{lem:frobenius-compatible-section}, we obtain a $\varphi$-compatible map 
\[
M_L[p^{-1}] \rightarrow  \mathbf{OB}_{\mathrm{cris}}(\mathcal{O}_{\overline{L}})\otimes_{\mathbf{Q}_p}V^{\vee}.
\]
Write $D = D_{\mathrm{cris}}^{\vee}(V)$ as before. If we compose the above map with
\[
 \mathbf{OB}_{\mathrm{cris}}(\mathcal{O}_{\overline{L}})\otimes_{\mathbf{Q}_p}V^{\vee} \xrightarrow{\alpha_{\mathrm{cris}}^{-1}}  \mathbf{OB}_{\mathrm{cris}}(\mathcal{O}_{\overline{L}})\otimes_{\iota_1,L_0}D \xrightarrow{\mathrm{pr}}  \mathbf{B}_{\mathrm{cris}}(\mathcal{O}_{\overline{L}})\otimes_{\iota_2,L_0}D,
\]
then we obtain a $\varphi$-compatible map
\[
M_L[p^{-1}] \rightarrow  \mathbf{B}_{\mathrm{cris}}(\mathcal{O}_{\overline{L}})\otimes_{\iota_2,L_0}D.
\]

We will use the following proposition in \S~\ref{sec:quasi-kisin-mod-connection}.

\begin{prop} \label{prop:CDVF-map-on-Dcris}
The image of the above map $M_L[p^{-1}] \rightarrow  \mathbf{B}_{\mathrm{cris}}(\mathcal{O}_{\overline{L}})\otimes_{\iota_2,L_0}D$ lies in
\[
D = L_0\otimes_{L_0}D \subset \mathbf{B}_{\mathrm{cris}}(\mathcal{O}_{\overline{L}})\otimes_{\iota_2,L_0}D .
\]
Furthermore, the induced map $M_L[p^{-1}] \rightarrow D$ is an isomorphism of $L_0$-modules.
\end{prop}

\begin{proof}
The construction of $G_L$-equivariant isomorphisms
\[
V(\mathscr{D}_L) \xrightarrow{\cong} \mathrm{Hom}_{\mathrm{Fil}, \varphi}(D, \mathbf{B}_{\mathrm{cris}}(\mathcal{O}_{\overline{L}})) \xrightarrow{\cong}V
\]
and diagram chasing implies that the above map injects into $D \subset  \mathbf{B}_{\mathrm{cris}}(\mathcal{O}_{\overline{L}})\otimes_{\iota_2,L_0}D$. Since $M_L[p^{-1}]$ and $D$ are $L_0$-vector spaces of the same dimension, the induced map $M_L[p^{-1}] \rightarrow D$ is an isomorphism. 
\end{proof}

\subsection{Construction of the quasi-Kisin module I: definition of \texorpdfstring{$\fkM$}{M}} \label{sec:quasi-kisin-mod-construction}

We now work over a general base ring: consider $R$ and $\fkS=R_0[\![u]\!]$ as in Set-up~\ref{set-up:base ring} and Notation~\ref{notation:Sigma}. 

Let $V$ be a crystalline $\Q_p$-representation of $\calG_R$ with Hodge--Tate weights in $[0, r]$, and let $T$ be a $\mathbf{Z}_p$-lattice of $V$ stable under $\calG_R$-action. Let 
\[
\mathcal{M} \coloneqq \mathcal{M}^{\vee}(T) = \mathrm{Hom}_{\mathcal{G}_{\tilde{R}_\infty}}(T, \widehat{\mathcal{O}}_{\mathcal{E}}^{\mathrm{ur}})
\]
be the associated \'etale $(\varphi, \mathcal{O}_{\mathcal{E}})$-module. In the following, we will construct a quasi-Kisin module over $\mathfrak{S}$ of $E$-height $\leq r$ associated with $T$.

Consider the base change along the map $R \rightarrow \mathcal{O}_{L}$ as in Notation~\ref{notation:L}. If we consider $T$ as a representation of $G_L\coloneqq \calG_{\calO_L}$ via $G_L \rightarrow \mathcal{G}_R$, then $T$ is a $\Z_p$-lattice in a crystalline $G_L$-representation with Hodge--Tate weights in $[0, r]$. Note that $\mathcal{M}_{L} \coloneqq \mathcal{O}_{\mathcal{E}, L}\otimes_{\mathcal{O}_{\mathcal{E}}}\mathcal{M} \cong \mathrm{Hom}_{G_{\widetilde{L}_\infty}}(T, \widehat{\mathcal{O}}_{\mathcal{E}, L}^{\mathrm{ur}})$ as \'etale $(\varphi, \mathcal{O}_{\mathcal{E}, L})$-modules. 

\begin{lem} \label{lem:kisin-mod-imperfect}
There exists a unique $\mathfrak{S}_L$-submodule $\mathfrak{M}_{L}$ of $\mathcal{M}_{L}$ stable under Frobenius such that the following properties hold.
\begin{itemize}
\item $\mathfrak{M}_{L}$ with the induced Frobenius is a quasi-Kisin module over $\mathfrak{S}_{L}$ of $E$-height $\leq r$. Furthermore, $\mathfrak{M}_{L}$ is free over $\mathfrak{S}_{L}$;
\item $\mathcal{O}_{\mathcal{E}, L}\otimes_{\mathfrak{S}_{L}} \mathfrak{M}_{L} = \mathcal{M}_{L}$;
\item if we let $M_L = \mathcal{O}_{L_0}\otimes_{\varphi,\mathcal{O}_{L_0}}\mathfrak{M}_L/u\mathfrak{M}_L $, then $M_L[p^{-1}] \cong D_{\mathrm{cris}}^{\vee}(V|_{G_L})$ compatibly with Frobenius and connection.	
\end{itemize}
\end{lem}

\begin{proof}
By \cite[Cor.~4.18]{brinon-trihan} and Lemma~\ref{lem:CDVF-filtration-compatibility}, there exists a quasi-Kisin module $\mathfrak{N}$ over $\mathfrak{S}_{L}$ of $E$-height $\leq r$ such that $\mathfrak{N}$ is free over $\mathfrak{S}_{L}$ and $\mathcal{O}_{L_0}[p^{-1}]\otimes_{\varphi,\mathcal{O}_{L_0}}\mathfrak{N}/u\mathfrak{N} \cong D_{\mathrm{cris}}^{\vee}(V|_{G_L})$ compatibly with Frobenii and connections. By \cite[Lem.~4.2.9]{gao-integral-padic-hodge-imperfect}, $\mathfrak{M}_{L} \coloneqq \mathfrak{N}[p^{-1}] \cap \mathcal{M}_{L}$ satisfies the required properties. The uniqueness also follows from the cited lemma.
\end{proof}

\begin{construction}\label{construction:quasi-Kisin module}
Let $T$ be a crystalline $\Z_p$-representation of $\calG_R$ with Hodge--Tate weights in $[0,r]$ and keep the notation as above. We set
\[
\mathfrak{M} \coloneqq 
\mathfrak{M}(T)\coloneqq\mathfrak{M}_{L} \cap \mathcal{M} \subset \mathcal{M}_{L}.
\]
This is an  $\fkS$-module. Moreover, since $\mathfrak{M}_{L}$ and $\mathcal{M}$ are $p$-adically complete and torsion free, so is $\mathfrak{M}$.
\end{construction}

We will show that $\mathfrak{M}$ is a quasi-Kisin module over $\mathfrak{S}$ of $E$-height $\leq r$ satisfying $\calO_{\calE}\otimes_{\fkS}\fkM\cong \calM$ and $\fkS_L\otimes_{\calS}\fkM\cong \fkM_L$.

\begin{prop} \label{prop:finitely-generated-intersection}
The $\fkS$-module $\mathfrak{M}$ is finitely generated. Furthermore, we have $\mathfrak{M}[u^{-1}] \cap \mathfrak{M}[p^{-1}] = \mathfrak{M}$.	
\end{prop}

\begin{proof}
Note that the cokernels of the maps $\mathfrak{S}_{L} \rightarrow \mathcal{O}_{\mathcal{E}, L}$ and $\mathcal{O}_{\mathcal{E}} \rightarrow \mathcal{O}_{\mathcal{E}, L}$ are $p$-torsion free. So the maps $\mathfrak{M}_{L}/p\mathfrak{M}_{L} \rightarrow \mathcal{M}_{L}/p\mathcal{M}_{L}$ and $\mathcal{M}/p\mathcal{M} \rightarrow \mathcal{M}_{L}/p\mathcal{M}_{L}$ are injective. By the proof of \cite[Lem.~4.1]{liu-moon-rel-crys-rep-p-div-gps-small-ramification}, the intersection $\mathfrak{M}_{L}/p\mathfrak{M}_{L} \cap \mathcal{M}/p\mathcal{M}$ inside $\mathcal{M}_{L}/p\mathcal{M}_{L}$ is finite over $\mathfrak{S}$. To show that $\mathfrak{M}$ is finite over $\mathfrak{S}$, it suffices to prove that the natural map $\mathfrak{M}/p\mathfrak{M} \rightarrow \mathfrak{M}_{L}/p\mathfrak{M}_{L} \cap \mathcal{M}/p\mathcal{M}$ is injective since $\fkS$ is noetherian and $\fkM$ is $p$-adically complete.

We have 
\[
p\mathfrak{M}_{L} \cap \mathfrak{M} = p\mathfrak{M}_{L}\cap \mathcal{M} \subset p\mathcal{M}_{L} \cap \mathcal{M}.
\]	
Since $p\mathcal{O}_{\mathcal{E}, L}\cap \mathcal{O}_{\mathcal{E}} = p\mathcal{O}_\mathcal{E}$ and $\mathcal{M}$ is classically flat over $\mathcal{O}_\mathcal{E}$, we have $p\mathcal{M}_{L} \cap \mathcal{M} \subset p\mathcal{M}$. Thus,
\[
p\mathfrak{M}_{L} \cap \mathfrak{M} \subset p\mathfrak{M}_{L} \cap p\mathcal{M} = p\mathfrak{M},
\]
where the last equality follows from the $p$-torsion freeness of $\mathcal{M}_{L}$. Thus the map $\mathfrak{M}/p\mathfrak{M} \rightarrow \mathfrak{M}_{L}/p\mathfrak{M}_{L}$ is injective. Since this map factors as $\mathfrak{M}/p\mathfrak{M} \rightarrow \mathfrak{M}_{L}/p\mathfrak{M}_{L} \cap \mathcal{M}/p\mathcal{M} \rightarrow\mathfrak{M}_{L}/p\mathfrak{M}_{L}$, we deduce the desired injectivity.

For the second part, since $\mathfrak{M}$ is torsion free, we have $\displaystyle \mathfrak{M} \subset \mathfrak{M}[u^{-1}] \cap \mathfrak{M}[p^{-1}]$. On the other hand,
\[
\mathfrak{M}[u^{-1}] \cap \mathfrak{M}[p^{-1}]	 \subset \mathfrak{M}_{L}[u^{-1}] \cap \mathfrak{M}_{L}[p^{-1}] = \mathfrak{M}_{L},
\]
and thus
\[
\mathfrak{M}[u^{-1}] \cap \mathfrak{M}[p^{-1}] \subset \mathfrak{M}[u^{-1}] \cap \mathfrak{M}_{L} \subset \mathcal{M} \cap \mathfrak{M}_{L} = \mathfrak{M}.
\]
\end{proof}

Since the Frobenius endomorphisms on $\mathcal{M}$ and $\mathfrak{M}_{L}$ are compatible with that on $\mathcal{M}_{\mathcal{O}_L}$, we have an induced Frobenius $\varphi_{\mathfrak{M}} \colon \mathfrak{M} \rightarrow \mathfrak{M}$.

\begin{prop} \label{prop:finite-E-height}
The $\fkS$-module $\mathfrak{M}$ with Frobenius has $E$-height $\leq r$.
\end{prop}

\begin{proof}
Since the composite of maps $\mathfrak{S}\otimes_{\varphi, \mathfrak{S}}\mathfrak{M} \rightarrow \mathcal{O}_{\mathcal{E}}\otimes_{\varphi, \mathcal{O}_\mathcal{E}}\mathcal{M} \xrightarrow{1\otimes\varphi} \mathcal{M}$ is injective, $1\otimes \varphi_{\mathfrak{M}} \colon \mathfrak{S}\otimes_{\varphi, \mathfrak{S}}\mathfrak{M} \rightarrow \mathfrak{M}$ is injective. Consider the natural map $\mathfrak{S}\otimes_{\varphi, \mathfrak{S}}\mathcal{O}_{\mathcal{E}} \rightarrow \mathcal{O}_{\mathcal{E}}\otimes_{\varphi, \mathcal{O}_{\mathcal{E}}} \mathcal{O}_{\mathcal{E}}$. Since $\mathfrak{S}\otimes_{\varphi, \mathfrak{S}}\mathcal{O}_{\mathcal{E}}$ and $\mathcal{O}_{\mathcal{E}}\otimes_{\varphi, \mathcal{O}_{\mathcal{E}}} \mathcal{O}_{\mathcal{E}}$ are $p$-adically complete and $p$-torsion free by Lemma~\ref{lem:Frobenius-faith-flat} and since the induced map $\mathfrak{S}/(p)\otimes_{\varphi, \mathfrak{S}}\mathcal{O}_{\mathcal{E}} \rightarrow \mathcal{O}_{\mathcal{E}}/(p)\otimes_{\varphi, \mathcal{O}_{\mathcal{E}}} \mathcal{O}_{\mathcal{E}}$ is an isomorphism, the map $\mathfrak{S}\otimes_{\varphi, \mathfrak{S}}\mathcal{O}_{\mathcal{E}} \rightarrow \mathcal{O}_{\mathcal{E}}\otimes_{\varphi, \mathcal{O}_{\mathcal{E}}} \mathcal{O}_{\mathcal{E}}$ is an isomorphism. Thus, the map
\[
\mathfrak{S}\otimes_{\varphi, \mathfrak{S}}\mathcal{M} \rightarrow \mathcal{O}_{\mathcal{E}}\otimes_{\varphi, \mathcal{O}_{\mathcal{E}}} \mathcal{M}
\]
is an isomorphism. On the other hand, since $R_0/pR_0$ has a finite $p$-basis which is also a $p$-basis of $\mathcal{O}_{L_0}/p\mathcal{O}_{L_0}$, the natural map $\mathfrak{S}/(p)\otimes_{\varphi, \mathfrak{S}}\mathfrak{S}_{L} \rightarrow \mathfrak{S}_{L}/(p)\otimes_{\varphi, \mathfrak{S}_{L}}\mathfrak{S}_{L}$ is an isomorphism. Since $\mathfrak{S}\otimes_{\varphi, \mathfrak{S}}\mathfrak{S}_{L}$ and $\mathfrak{S}_{L}\otimes_{\varphi, \mathfrak{S}_{L}}\mathfrak{S}_{L}$ are $p$-adically complete, the map $\mathfrak{S}\otimes_{\varphi, \mathfrak{S}}\mathfrak{S}_{L} \rightarrow \mathfrak{S}_{L}\otimes_{\varphi, \mathfrak{S}_{L}}\mathfrak{S}_{L}$ is an isomorphism. Hence
\[
\mathfrak{S}\otimes_{\varphi, \mathfrak{S}}\mathfrak{M}_{L} \rightarrow \mathfrak{S}_{L}\otimes_{\varphi, \mathfrak{S}_{L}}\mathfrak{M}_{L}
\]
is an isomorphism. 

Now, let $x \in \mathfrak{M}$. There exists a unique $y_1 \in \mathcal{O}_{\mathcal{E}}\otimes_{\varphi, \mathcal{O}_{\mathcal{E}}}\mathcal{M} = \mathfrak{S}\otimes_{\varphi, \mathfrak{S}} \mathcal{M}$ such that $(1\otimes\varphi)(y_1) = E(u)^rx$. On the other hand, there exists a unique $y_2 \in \mathfrak{S}_{L}\otimes_{\varphi, \mathfrak{S}_{L}}\mathfrak{M}_{L} = \mathfrak{S}\otimes_{\varphi, \mathfrak{S}}\mathfrak{M}_{L}$ such that $(1\otimes\varphi)(y_2) = E(u)^rx$. Hence, we have
\[
y_1 = y_2 \in (\mathfrak{S}\otimes_{\varphi, \mathfrak{S}}\mathcal{M}) \cap (\mathfrak{S}\otimes_{\varphi, \mathfrak{S}}\mathfrak{M}_{L}) = \mathfrak{S}\otimes_{\varphi, \mathfrak{S}} \mathfrak{M}
\]	 
by Lemma~\ref{lem:intersection-modules-flat-base-change} (i) since $\varphi\colon \mathfrak{S} \rightarrow \mathfrak{S}$ is classically flat.
\end{proof}

Next we will show that $\mathfrak{M}$ satisfies $\calO_{\calE}\otimes_{\fkS}\fkM\cong \calM$ and $\fkS_L\otimes_{\calS}\fkM\cong \fkM_L$. For this, we consider another description of $\mathfrak{M}$ as an inverse limit of $p$-power torsion $\mathfrak{S}$-modules as follows. Let 
\[
\mathfrak{M}_{L, n} \coloneqq \mathfrak{M}_{L}/ p^n\mathfrak{M}_{L}, \quad \mathcal{M}_n \coloneqq \mathcal{M} / p^n \mathcal{M},\quad\text{and}\quad \mathcal{M}_{L, n} \coloneqq \mathcal{M}_{L}/p^n\mathcal{M}_{L}. 
\]
Then $\mathfrak{M}_{L, n}$ and $\mathcal{M}_n$ are submodules of $\mathcal{M}_{L, n}$, and we set 
\[ 
\mathfrak{M}_{(n)} \coloneqq \mathfrak{M}_{L, n} \cap \mathcal{M}_n \subset \mathcal{M}_{L, n}.
\]
For any positive integers $i > j$, let $q_{i,j}$ denote the natural projection $\mathcal{M}_i \to \mathcal{M}_j$ given by modulo $p^j$, as well as its restriction $q_{i, j} \colon \mathfrak{M}_{(i)} \rightarrow \mathfrak{M} _{(j)}$. Note that $\mathcal{M}_{i-j}$ is naturally isomorphic to $p^j\mathcal{M}_i$. We have the commutative diagram 
\[
\xymatrix{ \ker(q_{i , j})\ar@{^{(}->}[r]\ar@{^{(}->}[d] & \mathfrak{M}_{(i)}\ar[r]^{q_{i, j}}\ar@{^{(}->}[d] & \fkM_{(j)}\ar@{^{(}->}[d] \\
\mathcal{M}_{i-j} \ar@{^{(}->}[r] & \mathcal{M}_i \ar@{->>}[r] ^{q_{i,j}} & \mathcal{M}_j.}
\]

\begin{lem} \label{lem:quasi-kisin-mod-inv-lim}
We have a natural isomorphism
\[
\mathfrak{M} \cong \varprojlim_n \mathfrak{M}_{(n)}.
\]	
\end{lem}

\begin{proof}
Recall that $\mathfrak{M}$ is $p$-adically complete. By a similar argument as in the proof of Proposition~\ref{prop:finitely-generated-intersection}, the natural map $\mathfrak{M}/p^n\mathfrak{M} \rightarrow \mathfrak{M}_{(n)}$ is injective for each $n \geq 1$. So the induced map
\[
f \colon \mathfrak{M} \xrightarrow{\cong} \varprojlim_n \mathfrak{M}/p^n\mathfrak{M} \rightarrow \varprojlim_n \mathfrak{M}_{(n)}
\]	
is injective. On the other hand, let $x = (x_n)_{n \geq 1} \in \varprojlim_n \mathfrak{M}_{(n)}$. Note that $x_n$ lies in both $\mathfrak{M}_{L, n}$ and $\mathcal{M}_n$ as an element in $\mathcal{M}_{L, n}$. Thus,
\[
x \in (\varprojlim_n \mathfrak{M}_{L, n}) \cap (\varprojlim_n \mathcal{M}_n) \subset  \varprojlim_n \mathcal{M}_{L, n},
\]
i.e., $x \in \mathfrak{M}_{L} \cap \mathcal{M} \subset \mathcal{M}_{L}$. This implies that $x$ lies in the image of $f$, and thus $f$ is surjective.
\end{proof}

\begin{rem} \label{rem:small-ramification-case}
Suppose $r = 1$ and $e < p-1$. By the above lemma and \cite[Prop.~4.3, 4.5]{liu-moon-rel-crys-rep-p-div-gps-small-ramification}, $\mathfrak{M}$ is projective over $\mathfrak{S}$. 
For general $r \geq 0$, as noted in \cite[Rem~4.6]{liu-moon-rel-crys-rep-p-div-gps-small-ramification}, the $\fkS$-module $\mathfrak{M}$ is projective when $er < p-1$. In particular, when $r = 0$, $\mathfrak{M}$ is projective for any $e$. 
\end{rem}

\begin{prop} \label{prop:torsion-mod-properties}
The following properties hold for $\fkM_{(n)}$:
\begin{enumerate}
\item $\mathfrak{M}_{(n)}$ is a finitely generated $\mathfrak{S}$-module;
\item $\mathfrak{M}_{(n)} [u^{-1}] \cong \mathcal{M}_n$ and $\mathfrak{S}_{L} \otimes_{\mathfrak{S}}\mathfrak{M}_{(n)} \cong \mathfrak{M}_{L, n}$;
\item $\mathfrak{M}_{(n)}$ has $E$-height $\leq r$.  
\end{enumerate}
\end{prop}

\begin{proof} Since the composite of maps $\mathfrak{S}\otimes_{\varphi, \mathfrak{S}}\mathfrak{M}_{(n)} \rightarrow \mathcal{O}_{\mathcal{E}}\otimes_{\varphi, \mathcal{O}_\mathcal{E}}\mathcal{M}_n \xrightarrow{1\otimes\varphi} \mathcal{M}_n$ is injective, $1\otimes\varphi \colon \mathfrak{S}\otimes_{\varphi, \mathfrak{S}}\mathfrak{M}_{(n)} \rightarrow \mathfrak{M}_{(n)}$ is injective. Thus, all statements follow from the same argument as in the proofs of \cite[Lem.~4.1, 4.2]{liu-moon-rel-crys-rep-p-div-gps-small-ramification} (where the case $r = 1$ is studied).
\end{proof}

Consider the set $\mathscr{A}_n$ consisting of $\mathfrak{S}$-submodules $\mathfrak{N}$ of $\fkM_{(n)}$ that are stable under $\varphi$, have $E$-height $\leq r$, and satisfy $\mathfrak{N}[u^{-1}] = \mathcal{M}_n$. Note that $\mathscr{A}_n$ is non-empty by the above proposition.
Let 
\[
\mathfrak{M}_{(n)}^\circ \coloneqq \bigcap_{\mathfrak{N} \in \mathscr{A}_n} \mathfrak{N} \subset \mathfrak{M}_{(n)}.
\]

\begin{lem} \label{lem:min-model} 
The following properties hold for $\fkM_{(n)}^\circ$:
\begin{enumerate}
    \item $\mathfrak{M}_{(n)}^\circ \in \mathscr{A}_n$;
    \item $ \mathfrak{M}^\circ_{(n)} \subset q_{n+1, n}( \mathfrak{M}^\circ_{(n+1)})$. 
\end{enumerate}	
\end{lem}

\begin{proof}
(i) Let $e \coloneqq [K : K_0]$ be the ramification index. We first show that for each fixed $n$, there exists an integer $s = s(n)$ such that $u^s\mathfrak{M}_{(n)} \subset \mathfrak{N} \subset \mathfrak{M}_{(n)}$ for all $\mathfrak{N} \in \mathscr{A}_n$. Choose an integer $a = a(n) \geq r$ such that $E(u)^{a}\equiv u ^{ea}\mod p^n$. Let $\mathfrak{N} \in \mathscr{A}_n$ and $\mathcal{L} \coloneqq \mathfrak{M}_{(n)} / \mathfrak{N}$. Without loss of generality, assume $\mathcal{L} \neq 0$. 
Note that $\mathfrak{M}_{(n)}$ and $\mathfrak{N}$ have $E$-height $\leq r$ and thus $E$-height $\leq a$. Hence we have unique $\mathfrak{S}$-linear maps $\psi_{\mathfrak{M}_{(n)}}\colon \mathfrak{M}_{(n)} \rightarrow \varphi^* \mathfrak{M}_{(n)}$ and $\psi_{\mathfrak{N}} \colon \mathfrak{N} \rightarrow \varphi^* \mathfrak{N}$ such that $\psi_{\mathfrak{N}} \circ (1 \otimes \varphi_{\mathfrak{N}}) = u^{ea} \text{Id}_{\varphi^* \mathfrak{N}}$ and $\psi_{\mathfrak{M}_{(n)}} \circ (1 \otimes \varphi_{\mathfrak{M}_{(n)}}) = u^{ea} \text{Id}_{\varphi^* \mathfrak{M}_{(n)}}$. 

The exact sequence $0 \rightarrow \mathfrak{N} \rightarrow \mathfrak{M}_{(n)} \rightarrow \mathfrak{L} \to 0$ induces the commutative diagram with exact rows:  
\[
\xymatrix{ 0 \ar[r] & \varphi^* \mathfrak{N}\ar[d]^{1 \otimes \varphi_{\mathfrak{N}}}\ar[r] & \varphi^*\mathfrak{M}_{(n)} \ar[r]\ar[d]^{1 \otimes \varphi_{\mathfrak{M}_{(n)}}} & \varphi^* \mathfrak{L}\ar[d]^{1 \otimes \varphi_{\mathfrak{L}}}\ar[r] & 0 \\ 
 0 \ar[r] &  \mathfrak{N}\ar[r]\ar[d]^{\psi_{\mathfrak{N}}} &  \mathfrak{M}_{(n)} \ar[r]\ar[d]^{\psi_{\mathfrak{M}_{(n)}}} &  \mathfrak{L}\ar[d]^{\psi_{\mathfrak{L}}}\ar[r] & 0 \\ 0 \ar[r] & \varphi^* \mathfrak{N}\ar[r] & \varphi^*\mathfrak{M}_{(n)} \ar[r] & \varphi^* \mathfrak{L}\ar[r] & 0.}
\]
Here, $1 \otimes \varphi_{\mathfrak{L}}$ and $\psi_{\mathfrak{L}}$ are the maps induced by $1\otimes \varphi_{\mathfrak{M}_{(n)}}$ and $\psi_{\mathfrak{M}_{(n)}}$ respectively. We have $\psi_{\mathfrak{L}} \circ (1 \otimes \varphi_{\mathfrak{L}}) = u^{ea} \text{Id}_{\varphi ^* \mathfrak{L}}$. 

We will show that $u^s\mathfrak{L}=0$ for $s=\lceil \frac{ea+p}{p -1} \rceil$.
Since $\mathfrak{N}[u^{-1}] = \mathcal{M}_n= \mathfrak{M}_{(n)}[u^{-1}]$, $\mathfrak{L}$ is killed by some $u$-power. Take an integer $l \geq 1$ such that $u^l \mathfrak{L} = 0$ and $u^{l-1}\mathfrak{L} \neq 0$. Pick $x \in \mathfrak{L}$ so that $u^{l-1} x \not = 0$. Set $y =1 \otimes x \in \varphi^* \mathfrak{L}$. Then $u^{pl} y = 1\otimes u^l x = 0$ but $u^{p(l-1)} y = 1 \otimes (u^{l-1} x ) \neq 0$, since $\varphi\colon \mathfrak{S} \rightarrow \mathfrak{S}$ is classically faithfully flat by Lemma~\ref{lem:Frobenius-faith-flat}. Let $z = (1 \otimes \varphi_{\mathfrak{L}})(y) \in \mathfrak{L}$. Since $u^l\mathfrak{L}=0$, we have $u^l z = 0$ and thus
\[
0 = \psi_{\mathfrak{L} } (u^l z) = u^l(\psi_{\mathfrak{L}} \circ  (1 \otimes\varphi_{\mathfrak{L}})) (y)  = u^{ea+l} y.
\]
So $ea+l > p(l-1)$, i.e., $l < \frac{ea +p}{p -1}$. Hence $u^s\mathfrak{L}=0$ for $s=\lceil \frac{ea+p}{p -1} \rceil$. This implies that $u^{s} \mathfrak{M}_{(n)} \subset \mathfrak{M}^\circ_{(n)} \subset \mathfrak{M}_{(n)}$, and $\mathfrak{M}^\circ_{(n)}[u^{-1}] = \mathcal{M}_n$. 
 
It remains to show that $\mathfrak{M}^\circ_{(n)}$ has $E$-height $\leq r$. Let $x \in \mathfrak{M}^\circ_{(n)}$. We need to show there exists $y \in \varphi^* \mathfrak{M}^\circ_{(n)}$ such that $(1 \otimes \varphi)(y) = E(u)^r x$. For each $\mathfrak{N} \in \mathscr{A}_n$, we have $x \in \mathfrak{N}$, and there exists $y \in \varphi^* \mathfrak{N}$, which is unique as an element of $\varphi^* \mathcal{M}_n$, such that $(1\otimes\varphi)(y) = E(u)^r x$. Since $\varphi\colon \mathfrak{S} \rightarrow \mathfrak{S}$ is finite free by Lemma~\ref{lem:Frobenius-faith-flat}, we deduce
\[
y \in \bigcap_{\mathfrak{N} \in \mathscr{A}_n} (\mathfrak{S} \otimes_{\varphi, \mathfrak{S}}\mathfrak{N}) = \mathfrak{S} \otimes_{\varphi, \mathfrak{S}} ( \bigcap_{\mathfrak{N} \in \mathscr{A}_n} \mathfrak{N}) = \varphi^* \mathfrak{M}^\circ_{(n)}.  
\]
 
(ii) Since $\mathfrak{M}^\circ_{(n+1)}[u^{-1}] = \mathcal{M}_{n+1}$ and $q_{{n+1}, n} (\mathcal{M}_{n +1}) = \mathcal{M}_n$, we have $q_{n+1, n} (\mathfrak{M}^\circ_{(n+1)})[u^{-1}] = \mathcal{M}_n$. So it suffices to show that $q_{n+1, n}(\mathfrak{M}^\circ_{(n+1)})$ has $E$-height $\leq r$. Let $\mathfrak{K} = \Ker (q_{n+1 , n})$. We have the commutative diagram with exact rows:
\[
\xymatrix{ 0 \ar[r] & \varphi^*\mathfrak{K}\ar[r]\ar[d] & \varphi^* \mathfrak{M}^\circ_{(n +1)}\ar[r]\ar[d]^{1 \otimes \varphi} & \varphi^*  q_{n+1, n}(\mathfrak{M}^\circ_{(n+1)})\ar[r]\ar[d] ^{1 \otimes \varphi} &  0 \\ 
 0 \ar[r] & \mathfrak{K}\ar[r] & \mathfrak{M}^\circ_{(n+1)}\ar[r]^-{q_{n+1, n}} &  q_{n+1, n}(\mathfrak{M}^\circ_{(n+1)})\ar[r] &  0 
 }
\]
Since $q_{n+1, n}(\mathfrak{M}^\circ_{(n+1)})\subset \mathcal{M}_n$, the rightmost vertical map is injective. From the first part, $1 \otimes \varphi$ in the middle column has cokernel killed by $E(u)^r$. Hence the cokernel of $1\otimes \varphi$ in the rightmost column is killed by $E(u)^r$.  
\end{proof}

\begin{prop} \label{prop:base-change-compatibility}
The natural maps 
\[
\mathcal{O}_{\mathcal{E}}\otimes_{\mathfrak{S}}\mathfrak{M} \rightarrow \mathcal{M}\quad\text{and}
\quad \mathfrak{S}_{L}\otimes_{\mathfrak{S}}\mathfrak{M} \rightarrow \mathfrak{M}_{L}.
\]	
are isomorphisms.
Moreover, $\mathfrak{M}$ is projective away from $(p, E)$, the Frobenius $\varphi_{\mathfrak{M}}$ has $E$-height $\leq r$, and $\mathrm{rank}_{\mathfrak{S}[p^{-1}]} \mathfrak{M}[p^{-1}] = \mathrm{rank}_{R_0[p^{-1}]} D_{\mathrm{cris}}^{\vee}(V)$. 
\end{prop}

\begin{proof}
We first prove $\mathcal{O}_{\mathcal{E}}\otimes_{\mathfrak{S}}\mathfrak{M} \cong \mathcal{M}$. Since $\mathfrak{M}\otimes_{\mathfrak{S}}\mathcal{O}_{\mathcal{E}}$ and $\mathcal{M}$ are $p$-adically complete and $\mathcal{M}$ is $p$-torsion free, it suffices to show that the induced map
\[
f \colon  \mathcal{O}_{\mathcal{E}}\otimes_{\mathfrak{S}}\mathfrak{M}/p\mathfrak{M} \cong (\mathfrak{M}/p\mathfrak{M})[u^{-1}] \rightarrow \mathcal{M}_1
\]	
is an isomorphism. It is shown in the proof of Proposition~\ref{prop:finitely-generated-intersection} that $\mathfrak{M}/p\mathfrak{M} \rightarrow \mathfrak{M}_{(1)}$ is injective. Hence $f$ is injective. By Lemmas~\ref{lem:quasi-kisin-mod-inv-lim} and \ref{lem:min-model} (ii), we have $\mathfrak{M}/p\mathfrak{M} \supset \mathfrak{M}^{\circ}_{(1)}$, and thus $f$ is surjective by Lemma~\ref{lem:min-model} (i).

For the second isomorphism, note that $\mathfrak{S}_{L}\otimes_{\mathfrak{S}}\mathfrak{M}/p\mathfrak{M} \rightarrow \mathfrak{S}_{L}\otimes_{\mathfrak{S}}\mathfrak{M}_{(1)} \cong \mathfrak{M}_{L, 1}$ is injective since $\mathfrak{S} \rightarrow \mathfrak{S}_{L}$ is classically flat. Since $\mathfrak{S}_{L}\otimes_{\mathfrak{S}}\mathfrak{M}$ is $p$-adically complete and $\mathfrak{M}_{L}$ is $p$-torsion free, $\mathfrak{S}_{L}\otimes_{\mathfrak{S}}\mathfrak{M} \rightarrow \mathfrak{M}_{L}$ is injective. In particular, $\mathfrak{S}_{L}\otimes_{\mathfrak{S}}\mathfrak{M}$ is a finite torsion free $\mathfrak{S}_{L}$-module. Furthermore, since $\mathfrak{S} \rightarrow \mathfrak{S}_{L}$ is classically flat, we have
\[
(\mathfrak{S}_{L}\otimes_{\mathfrak{S}}\mathfrak{M})[u^{-1}] \cap (\mathfrak{S}_{L}\otimes_{\mathfrak{S}}\mathfrak{M})[p^{-1}] = \mathfrak{M}\otimes_{\mathfrak{S}}\mathfrak{S}_{L}
\]
by Lemma~\ref{lem:intersection-modules-flat-base-change} (i) and Proposition~\ref{prop:finitely-generated-intersection}. Thus, $\mathfrak{S}_{L}\otimes_{\mathfrak{S}}\mathfrak{M}$ is a finite free $\mathfrak{S}_{L}$-module. 

Since $\mathfrak{M}$ has $E$-height $\leq r$ by Proposition~\ref{prop:finite-E-height}, $\mathfrak{S}_{L}\otimes_{\mathfrak{S}}\mathfrak{M}$ with the induced Frobenius has $E$-height $\leq r$. We have
\begin{align*}
 \mathcal{O}_{\mathcal{E},L}\otimes_{\mathfrak{S}_{L}}(\mathfrak{S}_{L}\otimes_{\mathfrak{S}}\mathfrak{M}) &\cong \mathcal{O}_{\mathcal{E}, L}\otimes_{\mathcal{O}_{\mathcal{E}}}(\mathcal{O}_{\mathcal{E}}\otimes_{\mathfrak{S}}\mathfrak{M})\\
	&\cong \mathcal{O}_{\mathcal{E}, L}\otimes_{\mathcal{O}_{\mathcal{E}}}\mathcal{M} \cong  \mathcal{O}_{\mathcal{E}, L}\otimes_{\mathfrak{S}_{L}}\mathfrak{M}_{L}.
\end{align*}
Hence, we deduce $\mathfrak{S}_{L}\otimes_{\mathfrak{S}}\mathfrak{M} \cong \mathfrak{M}_{L}$ by \cite[Prop.~4.2.5, 4.2.7]{gao-integral-padic-hodge-imperfect}. 

Finally, we deduce from the first isomorphism and Propositions~\ref{prop:rational-projectivity-etale-over-torus-case} and \ref{prop:finite-E-height} that $\mathfrak M$ is projective away from $(p, E)$, the Frobenius $\varphi_{\mathfrak{M}}$ has $E$-height $\leq r$, and 
\[
\mathrm{rank}_{\mathfrak{S}[p^{-1}]} \mathfrak{M}[p^{-1}] = \mathrm{rank}_{\mathcal{O}_{\mathcal{E}}} \mathcal{M} = \mathrm{rank}_{\mathbf{Q}_p} V =  \mathrm{rank}_{R_0[p^{-1}]} D_{\mathrm{cris}}^{\vee}(V).
\]
\end{proof}

\subsection{Construction of the quasi-Kisin module II: definition of \texorpdfstring{$\nabla$}{nabla}} \label{sec:quasi-kisin-mod-connection}

We further suppose that either $R$ is small over $\calO_K$ or $R=\calO_L$ (Assumption~\ref{assumption:base-ring-sec-3.4}).

Let $M = R_0\otimes_{\varphi, R_0}\mathfrak{M}/u\mathfrak{M}$. The $R_0$-module $M$ is equipped with the induced tensor-product Frobenius. We will construct a natural $\varphi$-equivariant isomorphism $M[p^{-1}] \cong D_{\mathrm{cris}}^{\vee}(V)$, via which we define $\nabla$ on $M[p^{-1}]$. Consider the $\varphi$-equivariant map $R_0 \rightarrow W(k_g)$ as in Notation~\ref{notation:L}, which naturally factors as $R_0 \rightarrow \mathcal{O}_{L_0} \rightarrow W(k_g)$.

\begin{lem} \label{lem:etale-galois-isom}
The natural $\mathcal{G}_{\tilde{R}_{\infty}}$-equivariant map
\[
\mathrm{Hom}_{\mathfrak{S}, \varphi}(\mathfrak{M}, \widehat{\mathfrak{S}}^{\mathrm{ur}}) \rightarrow T^{\vee}(\mathcal{M}) = \mathrm{Hom}_{\mathcal{O}_\mathcal{E}, \varphi}(\mathcal{M}, \widehat{\mathcal{O}}_{\mathcal{E}}^{\mathrm{ur}}) \cong T
\]	
is an isomorphism.
\end{lem}

\begin{proof}
For each $\mathfrak{p} \in \mathcal{P}$, consider the base change along $\overline{R}^{\wedge} \rightarrow (\overline{R}_{\mathfrak{p}})^\wedge\rightarrow (\mathcal{O}_{\overline{K_g}})^{\wedge}$ as in Notation \ref{notation:L}. The induced map 
\[
\mathrm{Hom}_{\mathfrak{S}, \varphi}(\mathfrak{M}, \widehat{\mathfrak{S}}_g^{\mathrm{ur}}) \rightarrow \mathrm{Hom}_{\mathcal{O}_\mathcal{E}, \varphi}(\mathcal{M}, \widehat{\mathcal{O}}_{\mathcal{E}, g}^{\mathrm{ur}})
\]
is an isomorphism by \cite[\S B Prop.~1.8.3]{fontaine-p-adic-rep-i}. Since $\widehat{\mathfrak{S}}^{\mathrm{ur}} = \widehat{\mathcal{O}}_{\mathcal{E}}^{\mathrm{ur}} \cap W(\overline{R}^\flat)$ by definition, we deduce the statement from Lemma~\ref{lem:intersection-witt-rings}.
\end{proof}

\begin{prop} \label{prop:comparing-quasi-kisin-mod-with-Dcris} 
Suppose that $R$ satisfies Assumption~\ref{assumption:base-ring-sec-3.4}. There exists a natural $\varphi$-compatible isomorphism $M[p^{-1}] \xrightarrow{\cong}D_{\mathrm{cris}}^{\vee}(V)$.	
\end{prop}

\begin{proof}
Let $\mathscr{M} \coloneqq S\otimes_{\varphi, \mathfrak{S}}\mathfrak{M}$ equipped with the induced Frobenius. Consider the $\varphi$-compatible projection $S \twoheadrightarrow R_0$ given by $u \mapsto 0$. This induces the projection $q\colon \mathscr{M} \twoheadrightarrow M$.
By Lemma~\ref{lem:frobenius-compatible-section}, Propositions~\ref{prop:rational-projectivity-etale-over-torus-case} and \ref{prop:finite-E-height}, the projection $q$ admits a unique $\varphi$-compatible section $s\colon M[p^{-1}] \rightarrow \mathscr{M}[p^{-1}]$, and $1\otimes s\colon S[p^{-1}]\otimes_{R_0[p^{-1}]}M[p^{-1}] \rightarrow \mathscr{M}[p^{-1}]$ is an isomorphism.

We first construct a $\varphi$-equivariant map $M[p^{-1}] \rightarrow D_{\mathrm{cris}}^{\vee}(V)$ similarly as in \S~\ref{sec:cryst-rep-cdvr}. By Lemma~\ref{lem:etale-galois-isom}, we have a natural map $\mathfrak{M} \rightarrow \widehat{\mathfrak{S}}^{\mathrm{ur}}\otimes_{\mathbf{Q}_p}V^{\vee}$. This extends to a map $\mathscr{M} \rightarrow \mathbf{OB}_{\mathrm{cris}}(\overline{R})\otimes_{\mathbf{Q}_p}V^{\vee}$ via the embedding $\widehat{\mathfrak{S}}^{\mathrm{ur}} \xrightarrow{\varphi} \mathbf{B}_{\mathrm{cris}}(\overline{R}) \hookrightarrow  \mathbf{OB}_{\mathrm{cris}}(\overline{R})$. So by composing with the section $s\colon M[p^{-1}] \rightarrow \mathscr{M}[p^{-1}]$ given in Lemma~\ref{lem:frobenius-compatible-section}, we get a $\varphi$-compatible map 
\begin{equation} \label{eq:map-from-M}
M[p^{-1}] \rightarrow \mathbf{OB}_{\mathrm{cris}}(\overline{R})\otimes_{\mathbf{Q}_p}V^{\vee}.
\end{equation}
By Lemma~\ref{lem:cryst-period-ring}, we have the projection $\mathrm{pr}\colon \mathbf{OB}_{\mathrm{cris}}(\overline{R}) \rightarrow \mathbf{B}_{\mathrm{cris}}(\overline{R})$ given by $T_i\otimes 1-1\otimes [T_i^\flat] \mapsto 0$. This induces the projection 
\[
\mathrm{pr}\colon \mathbf{OB}_{\mathrm{cris}}(\overline{R})\otimes_{\iota_1,R_0}D_{\mathrm{cris}}^{\vee}(V) \rightarrow   \mathbf{B}_{\mathrm{cris}}(\overline{R})\otimes_{\iota_2,R_0}D_{\mathrm{cris}}^{\vee}(V),
\]
where $\iota_1\colon R_0 \rightarrow \mathbf{OB}_{\mathrm{cris}}(\overline{R})$ is the natural map given by $T_i \mapsto T_i\otimes 1$ and $\iota_2\colon R_0 \rightarrow \mathbf{B}_{\mathrm{cris}}(\overline{R})$ is the embedding given by $T_i \mapsto [T_i^\flat]$.

If we compose the map \eqref{eq:map-from-M} with
\[
 \mathbf{OB}_{\mathrm{cris}}(\overline{R})\otimes_{\mathbf{Q}_p}V^{\vee} \xrightarrow{\alpha_{\mathrm{cris}}^{-1}} \mathbf{OB}_{\mathrm{cris}}(\overline{R})\otimes_{\iota_1,R_0}D_{\mathrm{cris}}^{\vee}(V) \xrightarrow{\mathrm{pr}} \mathbf{B}_{\mathrm{cris}}(\overline{R})\otimes_{\iota_2,R_0}D_{\mathrm{cris}}^{\vee}(V),
\]
then we obtain a $\varphi$-equivariant map
\[
f\colon M[p^{-1}] \rightarrow \mathbf{B}_{\mathrm{cris}}(\overline{R})\otimes_{\iota_2,R_0}D_{\mathrm{cris}}^{\vee}(V) .
\]

By Proposition~\ref{prop:CDVF-map-on-Dcris}, for each $\mathfrak{p} \in \mathcal{P}$ as in Notation \ref{notation:L}, we have
\[
f(M[p^{-1}]) \subset D_{\mathrm{cris}}^{\vee}(V|_{G_L}) = L_0\otimes_{R_0}D_{\mathrm{cris}}^{\vee}(V) \subset \mathbf{B}_{\mathrm{cris}}(\mathcal{O}_{\overline{L}})\otimes_{\iota_2,R_0}D_{\mathrm{cris}}^{\vee}(V).
\]
We claim that 
\[
\mathbf{B}_{\mathrm{cris}}(\overline{R}) \cap L_0 = R_0[p^{-1}]
\] 
as subrings of $\prod_{\mathfrak{p}\in \mathcal{P}} \mathbf{B}_{\mathrm{cris}}(\mathcal{O}_{\overline{L}})$. Since $\mathcal{O}_{L_0} \subset \mathbf{A}_{\mathrm{cris}}(\mathcal{O}_{\overline{L}})$ for each $\mathfrak{p} \in \mathcal{P}$, it suffices to show
\[
\mathbf{A}_{\mathrm{cris}}(\overline{R}) \cap \mathcal{O}_{L_0} = R_0.
\]
We clearly have $R_0 \subset \mathbf{A}_{\mathrm{cris}}(\overline{R}) \cap \mathcal{O}_{L_0}$. Let $x \in \mathbf{A}_{\mathrm{cris}}(\overline{R}) \cap \mathcal{O}_{L_0}$. Then 
\[
\theta(x) \in \overline{R}^{\wedge} \cap \mathcal{O}_{L_0} = R_0 \subset \prod_{\mathfrak{p} \in \mathcal{P}} \mathcal{O}_{\overline{L}}^{\wedge},
\]
which implies $x \in R_0$. This shows the claim.

Hence, we have by Lemma~\ref{lem:intersection-modules-flat-base-change} (i) that
\[
f(M[p^{-1}]) \subset (\mathbf{B}_{\mathrm{cris}}(\overline{R})\otimes_{R_0[p^{-1}]} D_{\mathrm{cris}}^{\vee}(V)) \cap (L_0\otimes_{R_0[p^{-1}]}D_{\mathrm{cris}}^{\vee}(V)) =  D_{\mathrm{cris}}^{\vee}(V),
\]
since $D_{\mathrm{cris}}^{\vee}(V)$ is projective over $R_0[p^{-1}]$. This gives a natural $\varphi$-equivariant map $f\colon M[p^{-1}] \rightarrow D_{\mathrm{cris}}^{\vee}(V)$.

It remains to show that $f$ is an isomorphism. Note that by Proposition~\ref{prop:CDVF-map-on-Dcris}, it suffices to consider the case where $R_0$ is the $p$-adic completion of an \'etale extension of $W(k)\langle T_1^{\pm 1}, \ldots, T_d^{\pm 1}\rangle$. By Proposition~\ref{prop:base-change-compatibility}, $M[p^{-1}]$ is projective over $R_0[p^{-1}]$ of rank equal to $\rank D_{\mathrm{cris}}^{\vee}(V)$. Moreover, by Proposition~\ref{prop:CDVF-map-on-Dcris}, the map
\[
L_0\otimes_{R_0}M[p^{-1}]  \rightarrow D_{\mathrm{cris}}^{\vee}(V|_{G_L}) = L_0\otimes_{R_0}D_{\mathrm{cris}}^{\vee}(V)
\]
induced by $f$ is an isomorphism. In particular, $f\colon M[p^{-1}] \rightarrow D_{\mathrm{cris}}^{\vee}(V)$ is injective. Let $I \subset R_0[p^{-1}]$ be the invertible ideal given by the determinant of $f$. Since $1\otimes\varphi\colon \varphi^* M[p^{-1}] \rightarrow M[p^{-1}]$ and $1\otimes\varphi\colon \varphi^*D_{\mathrm{cris}}^{\vee}(V) \rightarrow D_{\mathrm{cris}}^{\vee}(V)$ are isomorphisms, we have 
\[
\varphi(I)R_0[p^{-1}] = IR_0[p^{-1}].
\]
So $I = R_0[p^{-1}]$ by Proposition~\ref{prop:Frob-stable-invertible-ideal-etale-over-torus-case} (which is based on Lemma~\ref{lem:Frob-stable-invertible-ideal-power-series-case}) below, and $f$ is an isomorphism.
\end{proof}

\begin{lem} \label{lem:Frob-stable-invertible-ideal-power-series-case}
Let $k_1$ be a perfect field of characteristic $p$, and let $A = W(k_1)[\![t_1, \ldots, t_d]\!]$ be a power-series ring. Suppose that $A$ is equipped with a Frobenius endomorphism $\varphi$ extending the Witt vector Frobenius on $W(k_1)$. Let $I \subset A[p^{-1}]$ be an invertible ideal such that $IA[p^{-1}] = \varphi(I)A[p^{-1}]$. Then $I = A[p^{-1}]$.
\end{lem}

\begin{proof}
Suppose $I \neq A[p^{-1}]$. Since $A$ is a UFD, so is $A[p^{-1}]$. Hence $I$ is principal, say, generated by $x$. Since $p$ is a prime element of $A$, we may choose $x$ so that $x \in A \smallsetminus pA$ and $x$ is not a unit in $A$. Write $\varphi(x) = x^p+py$ for some $y \in A$. Since $p \nmid x$, we deduce from $\varphi(I)A[p^{-1}] = IA[p^{-1}]$ that $y = xz$ for some $z \in A$. Thus,
\[
\varphi(x) = x(x^{p-1}+pz).
\]
$x^{p-1}+pz$ is not a unit in $A$ since it lies in the ideal $(x, p)$ which is contained in the maximal ideal of $A$. Note that $p \nmid (x^{p-1}+pz)$ as elements in $A$. 
Thus, $x^{p-1}+pz$ is not a unit in $A[p^{-1}]$, which contradicts $\varphi(I)A[p^{-1}] = IA[p^{-1}]$. Hence, $I = A[p^{-1}]$.
\end{proof}

\begin{prop} \label{prop:Frob-stable-invertible-ideal-etale-over-torus-case}
Suppose $R_0$ is the $p$-adic completion of an \'etale extension of $W(k)\langle T_1^{\pm 1}, \ldots, T_d^{\pm 1}\rangle$. Let $I \subset R_0[p^{-1}]$ be an invertible ideal such that $\varphi(I)R_0[p^{-1}] = IR_0[p^{-1}]$. Then $I = R_0[p^{-1}]$.
\end{prop}

\begin{proof}
Let $\mathfrak{P} \subset R_0[p^{-1}]$ be any maximal ideal. Then the prime ideal $\mathfrak{q} = R_0 \cap \mathfrak{P}$ is maximal among the prime ideals of $R_0$ not containing $p$, and $\mathfrak{n} \coloneqq \sqrt{\mathfrak{q}+pR_0}$ is a maximal ideal of $R_0$. Let $(R_0)_{\mathfrak{n}}^{\wedge}$ be the $\mathfrak{n}$-adic completion of the localization $(R_0)_{\mathfrak{n}}$. As in the proof of Proposition~\ref{prop:rational-projectivity-etale-over-torus-case}, $(R_0)_{\mathfrak{n}}^{\wedge} \cong W(k_{\mathfrak{n}})[\![t_1, \ldots, t_d]\!]$ for some finite extension $k_{\mathfrak{n}}$ of $k$ and $t_1, \ldots, t_d$ such that $W(k_{\mathfrak{n}}) \hookrightarrow (R_0)_{\mathfrak{n}}^{\wedge}$ is compatible with $\varphi$.

Since the natural map $R_0[p^{-1}] \rightarrow (R_0)_{\mathfrak{n}}^{\wedge}[p^{-1}]$ is classically flat, we have $I(R_0)_{\mathfrak{n}}^{\wedge}[p^{-1}] = (R_0)_{\mathfrak{n}}^{\wedge}[p^{-1}]$ by Lemma~\ref{lem:Frob-stable-invertible-ideal-power-series-case}. Let $\mathfrak{P}_{\mathfrak{n}} \subset (R_0)_{\mathfrak{n}}^{\wedge}[p^{-1}]$ be a maximal ideal lying over $\mathfrak{P} \subset R_0[p^{-1}]$. Note that the natural map on localizations
\[
(R_0[p^{-1}])_{\mathfrak{P}} \rightarrow ((R_0)_{\mathfrak{n}}^{\wedge}[p^{-1}])_{\mathfrak{P}_{\mathfrak{n}}}
\]
is classically faithfully flat. Since $I((R_0)_{\mathfrak{n}}^{\wedge}[p^{-1}])_{\mathfrak{P}_{\mathfrak{n}}} = ((R_0)_{\mathfrak{n}}^{\wedge}[p^{-1}])_{\mathfrak{P}_{\mathfrak{n}}}$, we deduce that $I(R_0[p^{-1}])_{\mathfrak{P}} = (R_0[p^{-1}])_{\mathfrak{P}}$. This holds for any maximal ideal $\mathfrak{P} \subset R_0[p^{-1}]$, so $I = R_0[p^{-1}]$.
\end{proof}

By Proposition~\ref{prop:comparing-quasi-kisin-mod-with-Dcris}, the connection on $D_{\cris}^\vee(V)$ defines a connection
\[
\nabla_{\mathfrak{M}} \colon M[p^{-1}] \rightarrow M[p^{-1}]\otimes_{R_0} \widehat{\Omega}_{R_0}.
\]

Finally, we will show that $\nabla_\fkM$ satisfies the $S$-Griffiths transversality. For this, we study  the compatibility between two filtrations as in \S~\ref{sec:cryst-rep-cdvr}. By Lemma~\ref{lem:frobenius-compatible-section} and Proposition~\ref{prop:comparing-quasi-kisin-mod-with-Dcris}, we have a natural $\varphi$-equivariant isomorphism
\[
\mathscr{M}[p^{-1}] = S[p^{-1}]\otimes_{\varphi, \mathfrak{S}} \mathfrak{M} \cong S[p^{-1}]\otimes_{R_0} D_{\mathrm{cris}}^{\vee}(V).
\]
Let $\mathscr{D} \coloneqq S[p^{-1}]\otimes_{R_0} D_{\mathrm{cris}}^{\vee}(V)$, and identify $\mathscr{D} = \mathscr{M}[p^{-1}]$ via the above isomorphism. Let $N_u\colon \mathscr{D} \rightarrow \mathscr{D}$ be the $R_0$-linear derivation given by $N_{u} = N_{u, S}\otimes 1$, and let $\nabla\colon \mathscr{D}\rightarrow \mathscr{D}\otimes_{R_0} \widehat{\Omega}_{R_0}$ be the connection given by $\nabla_S\otimes 1+1\otimes\nabla_{D_{\mathrm{cris}}^{\vee}(V)}$. We consider two filtrations on $\mathscr{D} = \mathscr{M}[p^{-1}]$. For the first filtration, set $\mathrm{Fil}^0 \mathscr{D} = \mathscr{D}$, and inductively define for $i \geq 1$
\[
\mathrm{Fil}^i \mathscr{D} \coloneqq \{x \in \mathscr{D} ~|~ N_u(x) \in \mathrm{Fil}^{i-1} \mathscr{D}, ~~q_{\pi}(x) \in \mathrm{Fil}^i (R\otimes_{R_0}D_{\mathrm{cris}}^{\vee}(V))\},
\]
where $q_{\pi}\colon \mathscr{D} \rightarrow  R\otimes_{R_0}D_{\mathrm{cris}}^{\vee}(V)$ is the map given by $u \mapsto \pi$. For the second filtration, let
\[
\mathrm{F}^i \mathscr{M}[p^{-1}] \coloneqq \{x \in S[p^{-1}]\otimes_{\varphi, \mathfrak{S}} \mathfrak{M} ~|~ (1\otimes \varphi_{\mathfrak{M}})(x) \in (\mathrm{Fil}^i S[p^{-1}])\otimes_{\mathfrak{S}} \mathfrak{M}\}.
\]

\begin{lem} \label{lem:filtration-compatibility-Griffiths-transversality}
We have
\[
\mathrm{F}^i \mathscr{M}[p^{-1}] = \mathrm{Fil}^i \mathscr{D}.
\]	
Furthermore, 
\[
\nabla(\mathrm{Fil}^i \mathscr{D}) \subset \mathrm{Fil}^{i-1} \mathscr{D}\otimes_{R_0} \widehat{\Omega}_{R_0}.
\]
In particular, $\nabla$ satisfies the $S$-Griffiths transversality.
\end{lem}

\begin{proof}
By Proposition~\ref{prop:base-change-compatibility}, the first part follows from a similar argument as in the proof of Lemma~\ref{lem:CDVF-filtration-compatibility} using the base change along $R_0 \rightarrow W(k_g)$. The second part on the Griffiths transversality follows from a similar argument as in the proof of \cite[Lem.~4.2]{moon-strly-div-latt-cryst-cohom-CDVF}. Note that $N_u(\mathrm{Fil}^i \mathscr{D}) \subset \mathrm{Fil}^{i-1} \mathscr{D}$ by definition, and it is straightforward to check $\partial_u(\mathrm{Fil}^i \mathscr{D}) \subset \mathrm{Fil}^{i-1} \mathscr{D}$ by induction. 
\end{proof}

Combining this with Proposition~\ref{prop:base-change-compatibility}, we conclude the following.

\begin{prop} \label{prop:m-quasi-kisin-mod}
Suppose that $R$ satisfies Assumption~\ref{assumption:base-ring-sec-3.4}. With the above structures, $\mathfrak{M}$ is a quasi-Kisin module over $\mathfrak{S}$ of $E$-height $\leq r$.
\end{prop}

\subsection{Proof of the second part of Theorem~\ref{thm:main}} \label{sec:equivalence-categories}

Throughout this subsection, we suppose that $R$ satisfies Assumption~\ref{assumption:base-ring-sec-3.4}.

Let $V$ be a finite free $\mathbf{Q}_p$-representation of $\mathcal{G}_R$, which is crystalline with Hodge--Tate weights in $[0, r]$. Let $T \subset V$ be a $\mathbf{Z}_p$-lattice stable under the $\mathcal{G}_R$-action, and let $\mathfrak{M}$ be the quasi-Kisin module over $\mathfrak{S}$ of $E$-height $\leq r$ associated with $T$ as in Construction~\ref{construction:quasi-Kisin module} and Proposition~\ref{prop:m-quasi-kisin-mod}. 

By Propositions~\ref{prop:rational-descent-datum-over-A(2)} and \ref{prop:p=2-rational-descent-datum-over-A(2)}, we have a rational Kisin descent datum
\[
f\colon \mathfrak{S}^{(1)}[p^{-1}]\otimes_{p_1,\mathfrak{S}}\mathfrak{M}  \xrightarrow{\cong}  \mathfrak{S}^{(1)}[p^{-1}]\otimes_{p_2,\mathfrak{S}}\mathfrak{M}.
\]
 On the other hand, since $T$ is a finite free $\mathbf{Z}_p$-representation of $\mathcal{G}_R$,  \cite[Cor.~3.8]{bhatt-scholze-prismaticFcrystal} (see also \cite[Thm~3.2]{min-wang-rel-phi-gamma-prism-F-crys}) gives an isomorphism of $\mathfrak{S}^{(1)}[E^{-1}]^{\wedge}_p$-modules
\[
f_1\colon \mathfrak{S}^{(1)}[E^{-1}]^{\wedge}_p\otimes_{p_1,\mathcal{O}_{\mathcal{E}}}\mathcal{M} \xrightarrow{\cong} \mathfrak{S}^{(1)}[E^{-1}]^{\wedge}_p\otimes_{p_2,\mathcal{O}_{\mathcal{E}}}\mathcal{M}
\]
satisfying the cocycle condition over $\mathfrak{S}^{(2)}[E^{-1}]^{\wedge}_p$. Here, by Proposition~\ref{prop:base-change-compatibility}, $\mathcal{M} = \mathcal{O}_{\mathcal{E}}\otimes_{\mathfrak{S}}\mathfrak{M}$ is the \'etale $\varphi$-module finite projective over $\mathcal{O}_{\mathcal{E}}$ associated with $T$ (contravariantly). 
Note $\mathfrak{S}^{(1)}[E^{-1}]^{\wedge}_p\otimes_{p_j,\mathfrak{S}}\mathfrak{M} = \mathfrak{S}^{(1)}[E^{-1}]^{\wedge}_p\otimes_{p_j,\mathcal{O}_{\mathcal{E}}}\mathcal{M}$ since $\calO_\calE\otimes_{\fkS}\fkM=\calM$.

\begin{prop}\label{prop:compatibility of f and f_1}
Under the identification 
\[
\mathfrak{S}^{(1)}[E^{-1}]^{\wedge}_p[p^{-1}]\otimes_{p_j,\mathfrak{S}}\mathfrak{M}  = \mathfrak{S}^{(1)}[E^{-1}]^{\wedge}_p[p^{-1}]\otimes_{p_j,\mathcal{O}_{\mathcal{E}}}\mathcal{M} ,
\]
the maps $\mathrm{id}_{\mathfrak{S}^{(1)}[E^{-1}]^{\wedge}_p[p^{-1}]}\otimes_{\mathfrak{S}^{(1)}[p^{-1}]}f$ and $\mathrm{id}_{\mathfrak{S}^{(1)}[E^{-1}]^{\wedge}_p[p^{-1}]}\otimes_{\mathfrak{S}^{(1)}[E^{-1}]^{\wedge}_p}f_1$
coincide.
\end{prop}

To show Proposition~\ref{prop:compatibility of f and f_1}, let us consider the base change along $R_0 \rightarrow \mathcal{O}_{L_0}$ as before. Recall that $(\mathfrak{S}_L, (E)) = (\mathcal{O}_{L_0}[\![u]\!], (E))$ is a prism in $R_{\Prism}$ with the structure map $R \rightarrow \calO_L=\mathfrak{S}_L/(E)$. Let $(\mathfrak{S}_L^{(1)}, (E))$ be the self-product of $(\mathfrak{S}_L, (E))$ in $(\mathcal{O}_L)_{\Prism}$. Considering $(\mathfrak{S}_L^{(1)}, (E))$ as a prism in $R_{\Prism}$, the maps $f$ and $f_1$ induce the descent data
\[
f_L\colon \mathfrak{S}_L^{(1)}[p^{-1}]\otimes_{p_1,\mathfrak{S}_L}\mathfrak{M}_L \xrightarrow{\cong} \mathfrak{S}_L^{(1)}[p^{-1}]\otimes_{p_2,\mathfrak{S}_L}\mathfrak{M}_L
\]
and
\[
f_{1, L}\colon \mathfrak{S}_L^{(1)}[E^{-1}]^{\wedge}_p\otimes_{p_1,\mathfrak{S}_L}\mathfrak{M}_L  \xrightarrow{\cong} \mathfrak{S}_L^{(1)}[E^{-1}]^{\wedge}_p\otimes_{p_2,\mathfrak{S}_L}\mathfrak{M}_L,
\]
respectively. Here, $\mathfrak{M}_L = \mathfrak{S}_L\otimes_{\mathfrak{S}}\mathfrak{M}$ by Proposition~\ref{prop:base-change-compatibility}. Since the map $\mathfrak{S}^{(1)} \rightarrow \mathfrak{S}_L^{(1)}$ is injective, 
Proposition~\ref{prop:compatibility of f and f_1} follows if we show that $f_L$ and $f_{1, L}$ coincide over $\mathfrak{S}_L^{(1)}[E^{-1}]^{\wedge}_p[p^{-1}]$.
For this, we need the following proposition.

\begin{prop} \label{prop:Drinfeld}
There exists an $\mathfrak{S}_L$-submodule $\mathfrak{N}_L \subset \mathfrak{M}_L$ with $\mathfrak{N}_L[p^{-1}] = \mathfrak{M}_L[p^{-1}]$ such that $f_L$ induces an isomorphism of $\fkS_L^{(1)}$-modules
\[
f_L\colon \mathfrak{S}_L^{(1)}\otimes_{p_1,\mathfrak{S}_L}\mathfrak{N}_L \xrightarrow{\cong} \mathfrak{S}_L^{(1)}\otimes_{p_2,\mathfrak{S}_L}\mathfrak{N}_L.
\]
 Furthermore, $\mathfrak{N}_L$ can be chosen to be finite free over $\mathfrak{S}_L$ of $E$-height $\leq r$ and stable under $\varphi_{\mathfrak{M}_L}$.
\end{prop}

\begin{proof}
Since $p_i\colon \mathfrak{S}_L \rightarrow \mathfrak{S}_L^{(1)}$ is classically faithfully flat by Lemma~\ref{lem:AtoA2A3-faithful-flat}, the first part follows directly from the proof of \cite[Prop.~3.6.5]{drinfeld-stackycrystal}. We recall some points here. Note that for any $\fkS_L$-submodule $\mathfrak{N}_L \subset \mathfrak{M}_L$, the induced map $p_i^\ast \mathfrak{N}_L \rightarrow p_i^\ast \mathfrak{M}_L$ for $i = 1, 2$ is injective, where $p_i^\ast\fkN_L$ denotes $\fkS_L^{(1)}\otimes_{p_i,\fkS_L}\fkN_L$. 
Take an integer $n \geq 0$ such that $p^n f_L$ maps $p_1^*\mathfrak{M}_L$ into $p_2^* \mathfrak{M}_L$. It suffices to find an $\mathfrak{S}_L$-submodule $\mathfrak{N}_L \subset \mathfrak{M}_L$ such that $p^n\mathfrak{M}_L \subset \mathfrak{N}_L$ and $f_L$ maps $p_1^* \mathfrak{N}_L$ to $p_2^* \mathfrak{N}_L$; then it follows from the cocycle condition on $f_L$ for $\mathfrak{M}_L[p^{-1}]$ that the induced map $f_L\colon p_1^* \mathfrak{N}_L \rightarrow p_2^* \mathfrak{N}_L$ is an isomorphism (see also the proof of \cite[Thm.~1.9]{ogus-F-converg-isocryst-de-rham-cohom-II}). The map $f_L$ induces a map
\[
f_L\colon p_1^*(\mathfrak{M}_L/p^n\mathfrak{M}_L) \rightarrow p_2^*(p^{-n}\mathfrak{M}_L/\mathfrak{M}_L).
\]
This induces a morphism
\[
\beta\colon \mathfrak{M}_L/p^n\mathfrak{M}_L \rightarrow \Phi(p^{-n}\mathfrak{M}_L/\mathfrak{M}_L)
\]
where $\Phi \coloneqq (p_1)_*p_2^*$. Let $\mathfrak{N}_L$ be the kernel of the composite 
\[
\mathfrak{M}_L \twoheadrightarrow \mathfrak{M}_L/p^n\mathfrak{M}_L \xrightarrow{\beta} \Phi(p^{-n}\mathfrak{M}_L/\mathfrak{M}_L).
\]
Then $f_L(p_1^*\mathfrak{N}_L) \subset p_2^*\mathfrak{M}_L$, and by the proof of \cite[Prop.~3.6.5]{drinfeld-stackycrystal} (cf. also the proof of \cite[Thm.~1.9]{ogus-F-converg-isocryst-de-rham-cohom-II}), we further have $f_L(p_1^*\mathfrak{N}_L) \subset p_2^*\mathfrak{N}_L$.

Since $f_L$ is compatible with Frobenius, $\beta$ as above is compatible with $\varphi$. Thus, $\mathfrak{N}_L$ constructed as above is stable under $\varphi$. Consider the exact sequence
\[
0 \rightarrow \mathfrak{N}_L \rightarrow \mathfrak{M}_L \rightarrow \beta(\mathfrak{M}_L) \rightarrow 0
\]
where $\beta(\mathfrak{M}_L) \subset \Phi(p^{-n}\mathfrak{M}_L/\mathfrak{M}_L)$ denotes the image of $\mathfrak{M}_L$ under the above composite. Under the classically faithfully flat base change along $\mathfrak{S}_L \rightarrow \mathfrak{S}_g$, this induces an exact sequence
\[
0 \rightarrow \mathfrak{S}_g\otimes_{\mathfrak{S}_L}\mathfrak{N}_L  \rightarrow \mathfrak{S}_g\otimes_{\mathfrak{S}_L}\mathfrak{M}_L \rightarrow \mathfrak{S}_g\otimes_{\mathfrak{S}_L}\beta(\mathfrak{M}_L) \rightarrow 0.
\]
Note that $\mathfrak{S}_g\otimes_{\mathfrak{S}_L}\mathfrak{M}_L$ is a Kisin module of $E$-height $\leq r$ that is finite free over $\mathfrak{S}_g$. Furthermore, $\beta(\mathfrak{M}_L)$ is $u$-torsion free since $p^{-n}\mathfrak{M}_L/\mathfrak{M}_L$ is finite free over $\mathfrak{S}_L/p^n\mathfrak{S}_L$, and so $\mathfrak{S}_g\otimes_{\mathfrak{S}_L}\beta(\mathfrak{M}_L)$ is $u$-torsion free. Thus, $\mathfrak{S}_g\otimes_{\mathfrak{S}_L}\beta(\mathfrak{M}_L)$ is a torsion Kisin module over $\mathfrak{S}_g$ of $E$-height $\leq r$ by \cite[Prop.~2.3.2]{liu-fontaineconjecture}. Then $\mathfrak{S}_g\otimes_{\mathfrak{S}_L}\mathfrak{N}_L$ is a Kisin module of $E$-height $\leq r$ finite free over $\mathfrak{S}_g$ by \cite[Cor.~2.3.8]{liu-fontaineconjecture}. Since $\mathfrak{S}_L \rightarrow \mathfrak{S}_g$ is classically faithfully flat, $\mathfrak{N}_L$ is finite free over $\mathfrak{S}_L$ and has $E$-height $\leq r$. 
\end{proof}

\begin{proof}[Proof of Proposition~\ref{prop:compatibility of f and f_1}]
By the above proposition, $\mathcal{N}_L \coloneqq \mathcal{O}_{\mathcal{E}, L}\otimes_{\mathfrak{S}_L}\mathfrak{N}_L $ is an \'etale $\varphi$-module finite free over $\mathcal{O}_{\mathcal{E}, L}$, and $f_L$ induces an isomorphism
of $\mathfrak{S}_L^{(1)}[E^{-1}]^{\wedge}_p$-modules
\[
f_L\colon \mathfrak{S}_L^{(1)}[E^{-1}]^{\wedge}_p\otimes_{p_1,\mathcal{O}_{\mathcal{E}, L}}\mathcal{N}_L \xrightarrow{\cong} \mathfrak{S}_L^{(1)}[E^{-1}]^{\wedge}_p\otimes_{p_2,\mathcal{O}_{\mathcal{E}, L}}\mathcal{N}_L 
\]
satisfying the cocycle condition over $\mathfrak{S}_L^{(2)}[E^{-1}]^{\wedge}_p$. As in the proof of Proposition~\ref{prop:etale-realization}, this corresponds to a finite free $\mathbf{Z}_p$-representation $T'$ of $G_L$. Furthermore, by \cite[Cor.~3.7, Ex.~3.5]{bhatt-scholze-prismaticFcrystal}, $T'$ is determined by the $G_L$-action on $ W(\mathcal{O}_{\overline{L}}^{\flat}[(\pi^{\flat})^{-1}])\otimes_{\mathcal{O}_{\mathcal{E}, L}}\mathcal{N}_L$. On the other hand, note that $(\mathbf{A}_{\mathrm{cris}}(\mathcal{O}_{\overline{L}}), (p)) \in R_{\Prism}$ similarly as in Example \ref{eg:prism-OAcris}, and the composite $S \rightarrow S_L \rightarrow \mathbf{A}_{\mathrm{cris}}(\mathcal{O}_{\overline{L}})$ gives a map of prisms $(S, (p)) \rightarrow (\mathbf{A}_{\mathrm{cris}}(\mathcal{O}_{\overline{L}}), (p))$ over $R$. Thus, by the construction of the descent datum $f$ and definition of $f_S$ in Construction \ref{construction:rational S-descent datum}, the $G_L$-action on $\mathbf{A}_{\mathrm{cris}}(\mathcal{O}_{\overline{L}})[p^{-1}]\otimes_{\varphi,\mathfrak{S}_L} \mathfrak{N}_L \cong \mathbf{A}_{\mathrm{cris}}(\mathcal{O}_{\overline{L}})[p^{-1}] \otimes_{S_L}\mathscr{D}_L$ (with $\mathscr{D}_L = S_L[p^{-1}]\otimes_{\varphi, \mathfrak{S}_L} \mathfrak{M}_L$) is given by
\[
\sigma(a\otimes x) = \sum \sigma(a) \partial_u^{j_0}\partial_{T_1}^{j_1}\cdots \partial_{T_d}^{j_d} (x) \cdot \gamma_{j_0}(\sigma([\pi^\flat])-[\pi^\flat])\prod_{i=1}^d \gamma_{j_i}(\sigma([T_i^\flat])-[T_i^\flat]),
\]
for $\sigma \in G_L$ and $a\otimes x\in \mathbf{A}_{\mathrm{cris}}(\mathcal{O}_{\overline{L}})[p^{-1}]\otimes_{S_L}\mathscr{D}_L$ (where the sum goes over the multi-index $(j_0, \ldots, j_d)$ of non-negative integers). By \cite[\S 8.1]{Li-Liu-prismatic-cohom}, this is the same as the $G_L$-action given by Equation \eqref{eq:CDVF-Galois-action}, and it is proved in \S~\ref{sec:cryst-rep-cdvr} that this gives a $\mathbf{Q}_p$-representation of $G_L$ isomorphic to $T[p^{-1}] = V$. Thus, $T[p^{-1}] \cong T'[p^{-1}]$ as representations of $G_L$. This proves the claim that $f_L$ and $f_{1,L}$ coincide over $\mathfrak{S}_L^{(1)}[E^{-1}]^{\wedge}_p[p^{-1}]$ and thus Proposition~\ref{prop:compatibility of f and f_1}.
\end{proof}

By Proposition~\ref{prop:compatibility of f and f_1}, we see that the descent data $f$ and $f_1$ induce a map
\[
f\colon \mathfrak{S}^{(1)}\otimes_{p_1,\mathfrak{S}}\mathfrak{M} \rightarrow (\mathfrak{S}^{(1)}[p^{-1}]\otimes_{p_2,\mathfrak{S}}\mathfrak{M}) \cap ( \mathfrak{S}^{(1)}[E^{-1}]^{\wedge}_p\otimes_{p_2,\mathfrak{S}}\mathfrak{M}).
\]

\begin{proof}[End of the proof of Theorem~\ref{thm:main} (ii)]
By Lemma~\ref{lem:intersectin of two base chage over S^(2)}, we see that $f$ and $f_1$ induce a map
\[
f_{\mathrm{int}}\colon \mathfrak{S}^{(1)}\otimes_{p_1,\mathfrak{S}}\mathfrak{M} \rightarrow \mathfrak{S}^{(1)}\otimes_{p_2,\mathfrak{S}}\mathfrak{M}.
\]
By applying a similar argument to $f^{-1}$ and $f_1^{-1}$, we deduce that $f_{\mathrm{int}}$ is an isomorphism. Namely, $f_\mathrm{int}$ is a descent datum. Since $f$ is compatible with $\varphi$, so is $f_\mathrm{int}$. 
By Proposition~\ref{prop:equivalence-to-descent-datum}, the triple $(\fkM,\varphi,f_\mathrm{int})$ gives rise to a completed prismatic $F$-crystal $\calF$ on $R$ with $\fkM=\calF_{\fkS}$.
It is straightforward to see $T(\calF)=T$.
Hence the functor $T$ in Theorem~\ref{thm:main} is essentially surjective (when $R$ satisfies Assumption~\ref{assumption:base-ring-sec-3.4}).
\end{proof}

\begin{rem}\label{rem:how to recover Dcris-2}
We continue the discussion in Remark~\ref{rem:how to recover Dcris-1}. For $\calF\in\mathrm{CR}^{\wedge,\varphi}(R_\Prism)$, let $V = T(\calF)[p^{-1}] \in \mathrm{Rep}_{\mathbf{Z}_p, \geq 0}^{\mathrm{cris}}(\mathcal{G}_R)$. Consider the $\varphi$-equivalent $R_0[p^{-1}]$-linear isomorphism
\[
h\colon (R_0\otimes_{\varphi,R_0}\calF_\fkS/u\calF_\fkS)[p^{-1}]\cong D_\cris^\vee(V)
\]
in Remark~\ref{rem:how to recover Dcris-1}. Note that $\mathcal{F}_S[p^{-1}] \cong S[p^{-1}]\otimes_{\varphi, \mathfrak{S}}\mathcal{F}_{\mathfrak{S}}$ by Lemma~\ref{lem:completed-crystals-basic-properties} (iv) for the map of prisms $(\mathfrak{S}, E) \rightarrow (S, (p))$. So $h$ induces a $\varphi$-compatible isomorphism $\mathcal{F}_S[p^{-1}] \cong S[p^{-1}]\otimes_{R_0}D_\cris^\vee(V)$ by Lemma~\ref{lem:frobenius-compatible-section}. 

Since $\calF\in\mathrm{CR}^{\wedge,\varphi}(R_\Prism)$, we have an isomorphism of $\mathfrak{S}^{(1)}$-modules
\[
\mathfrak{S}^{(1)}\otimes_{p_1,\mathfrak{S}}\mathcal{F}_{\mathfrak{S}} \xrightarrow{\cong} \mathfrak{S}^{(1)}\otimes_{p_2,\mathfrak{S}}\mathcal{F}_{\mathfrak{S}}.
\]
Under the map $\varphi\colon \mathfrak{S}^{(1)} \rightarrow S^{(1)}$, this induces an isomorphism of $S^{(1)}[p^{-1}]$-modules
\begin{equation} \label{eq:S(2)-isom}
f_{S}\colon S^{(1)}[p^{-1}]\otimes_{p_1,S}\mathcal{F}_{S} \xrightarrow{\cong} S^{(1)}[p^{-1}]\otimes_{p_2,S}\mathcal{F}_{S}.
\end{equation}
Note that $f_S$ reduces to the identity after the base change along $S^{(1)} \rightarrow S$ and it satisfies the cocycle condition over $S^{(2)}$.
Let $\nu\colon R_0\otimes_W R_0 \rightarrow R_0$ be the multiplication, and let $R_0^{(1)}$ be the $p$-adically completed divided power envelope of $R_0\otimes_W R_0$ with respect to $\mathrm{Ker} (\nu)$. We also write $\nu\colon R_0^{(1)} \rightarrow R_0$ for the induced map. Consider the map $S^{(1)} \rightarrow R_0^{(1)}$ given by $u, y \mapsto 0$. From the isomorphism \eqref{eq:S(2)-isom}, we obtain an isomorphism of $R_0^{(1)}[p^{-1}]$-modules
\[
f_{R_0}\colon R_0^{(1)}[p^{-1}]\otimes_{p_1, R_0}D_\cris^\vee(V) \xrightarrow{\cong} R_0^{(1)}[p^{-1}]\otimes_{p_2, R_0}D_\cris^\vee(V) 
\]
such that $f_{R_0}$ reduces to the identity after the base change along $\nu$ and it satisfies a similar cocycle condition.
Since $\widehat{\Omega}_{R_0} \cong \mathrm{Ker} (\nu) / (\mathrm{Ker} (\nu))^{[2]}$ where $(\mathrm{Ker} (\nu))^{[2]}$ denotes the divided square of $\mathrm{Ker} (\nu)$, the isomorphism $f_{R_0}$ gives an integrable connection $\nabla\colon D_\cris^\vee(V) \rightarrow D_\cris^\vee(V)\otimes_{R_0} \widehat{\Omega}_{R_0}$. On the other hand, we have the natural integrable connection on $D_\cris^\vee(V)$ induced by that on $\OB_{\cris}(\overline{R})$ (cf. \S~\ref{sec:crystalline representations}). In the proof of Theorem~\ref{thm:main} (ii) on essential surjectivity, the isomorphism \eqref{eq:S(2)-isom} is obtained by Construction~\ref{construction:rational S-descent datum} using the natural connection on $D_\cris^\vee(V)$ as in \S~\ref{sec:quasi-kisin-mod-connection}. Thus, $\nabla$ given above agrees with the natural connection on $D_\cris^\vee(V)$.

Define the filtration on $\mathcal{F}_S[p^{-1}]$ by
\[
\mathrm{F}^i \mathcal{F}_S[p^{-1}] \coloneqq \{x \in S[p^{-1}]\otimes_{\varphi, \mathfrak{S}}\mathcal{F}_{\mathfrak{S}} ~|~ (1\otimes \varphi)(x) \in \mathrm{Fil}^i S[p^{-1}]\otimes_{\mathfrak{S}} \mathcal{F}_{\mathfrak{S}}\}.
\]
Let $\mathrm{Fil}^i(R\otimes_{R_0} D_\cris^\vee(V))$ be the quotient filtration given by $\mathrm{F}^i \mathcal{F}_S[p^{-1}]$ under the map $q_{\pi}\colon S \rightarrow R$, $u \mapsto \pi$. By the proof of Theorem~\ref{thm:main} (ii), Lemma~\ref{lem:filtration-compatibility-Griffiths-transversality} and a similar argument as in the proof of \cite[Prop. 6.2.2.3]{breuil-representations}, this quotient filtration agrees with the natural filtration on $R\otimes_{R_0} D_\cris^\vee(V)$ as in \S~\ref{sec:crystalline representations}. In this way, we can directly obtain the filtered $(\varphi,\nabla)$-module $(D_\cris^\vee(T(\calF)[p^{-1}]), \nabla, \Fil^i(R\otimes_{R_0}D_\cris^\vee(T(\calF)[p^{-1}])))$ from $\calF$.
\end{rem}

\begin{cor} \label{cor:CDVF-equivalence}
 The \'etale realization functor gives an equivalence of categories from $\mathrm{Vect}^{\varphi}_{[0, r]}((\mathcal{O}_L)_{\Prism})$ to $\mathrm{Rep}_{\mathbf{Z}_p, [0, r]}^{\mathrm{cris}}(G_L)$.
\end{cor}

\begin{proof}
By Remark~\ref{rem:vector-bundle-CDVF-case}, the category $\mathrm{CR}^{\wedge,\varphi}_{[0,r]}((\mathcal{O}_L)_{\Prism})$ is equal to $\mathrm{Vect}^{\varphi}_{[0, r]}((\mathcal{O}_L)_{\Prism})$. So the statement follows from Theorem~\ref{thm:main}.
\end{proof}

\begin{rem}
As a corollary, we can deduce that the construction of Brinon--Trihan in \cite{brinon-trihan} is independent of the choice of a uniformizer and the Kummer tower.
\end{rem}

\appendix
\section{Crystalline local systems}\label{sec:crystalline local systems}

Let $\fkX$ be a smooth $p$-adic formal scheme over $\calO_K$ and let $X$ denote its adic generic fiber.
In this appendix, we define the notion of crystalline local systems on $X$, which is used in \S~\ref{sec:globalization}. 
The definition of crystalline local systems goes back to the work \cite[V f)]{faltings} of Faltings. 
Tan and Tong \cite{Tan-Tong} also define crystalline local systems in the unramified case $\calO_K=W$ and prove that their definition agrees with the one given by Faltings. 
Since we also work on the ramified case, we give a minimal foundation that generalizes part of the work of Tan and Tong.

For our purpose, we work in two steps: when there exists a smooth $p$-adic formal scheme $\fkX_0$ over $W$ such that $\fkX\cong \fkX_0\otimes_W\calO_K$, we define the pro-\'etale sheaf $\calOB_\cris$ on $X$ and use it to define  crystalline local systems. Note that this assumption is satisfied Zariski locally, e.g., by considering a Zariski open covering consisting of small affines. In the general case, we define crystalline local systems via gluing.

Let $X_\proet$ denote the pro-\'etale site defined in \cite[\S 3]{scholze-p-adic-hodge} and \cite{Scholze-p-adicHodgeerrata}. It admits the morphism of site $\nu\colon X_\proet\rightarrow X_\et$.

\begin{defn}
We introduce sheaves on $X_\proet$.
\begin{enumerate}
    \item (\cite[Def.~4.1, 5.10]{scholze-p-adic-hodge}). Set
    \[
\calO_X^{+}\coloneqq \nu^{-1}\calO_{X_\et}^{+}, \quad
\widehat{\calO}_X^+\coloneqq \varprojlim \calO_X^+/p^n,\quad\widehat{\calO}_X\coloneqq \widehat{\calO}_X^+[p^{-1}],\quad\text{and}\quad
\widehat{\calO}_{X^\flat}^+\coloneqq \mspace{-6mu}\varprojlim_{\Phi\colon x\to x^p}\mspace{-6mu}\widehat{\calO}_X^+/p.
\]
    \item (\cite[Def.~6.1]{scholze-p-adic-hodge}. Set $\bA_{\mathrm{inf}}\coloneqq W(\widehat{\calO}_{X^\flat}^+)$ and $\mathbb{B}_{\mathrm{inf}}\coloneqq \bA_{\mathrm{inf}}[p^{-1}]$.
We have ring morphisms $\theta\colon \bA_{\mathrm{inf}}\rightarrow \widehat{\calO}_X^+$ and $\theta\colon \bB_{\mathrm{inf}}\rightarrow \widehat{\calO}_X$.
    \item (\cite[Def.~2.1]{Tan-Tong}).
    Let $\bA_\cris^0$ be the PD-envelope of $\bA_{\mathrm{inf}}$ with respect to the ideal sheaf $\Ker \theta$, and set $\bA_\cris\coloneqq \varprojlim \bA_\cris^0/p^n$.
    Note that the series $t\coloneqq \log [\varepsilon]$ converges and is a nonzero-divisor in $\bA_\cris|_{X_{\overline{K}}}$. See \cite[(2A.6), Cor.~2.24]{Tan-Tong}.
\end{enumerate}
\end{defn}

Now assume that $\fkX$ admits a \emph{$W$-model}, namely, there exists a smooth $p$-adic formal scheme $\fkX_0$ over $W$ such that $\fkX\cong \fkX_0\otimes_W\calO_K$. 
Let $X_0$ denote the adic generic fiber of $\fkX_0$. Hence we have a canonical identification $X\cong X_0\times_{\Spa(W[p^{-1}],W)}\Spa(K,\calO_K)$.
In \cite[\S~2B]{Tan-Tong}, Tan and Tong defined the structural crystalline period sheaves $\calOA_{\cris,X_0}$ and $\calOB_{\cris,X_0}$ on $(X_0)_\proet$. We define structural crystalline sheaves on $X_\proet$ in a similar way.

\begin{defn}[(cf.~{\cite[\S~2B]{Tan-Tong}})]
Consider the morphisms of sites
\[
w\colon X_\proet\rightarrow X_\et\rightarrow \fkX_\et\rightarrow (\fkX_0)_\et.
\]
Define sheaves $\calO_X^{\ur+}$ and $\calO_X^{\ur}$ on $X_\proet$ by
\[
\calO_X^{\ur+}\coloneqq\calO_X^{\ur/\fkX_0+}\coloneqq w^{-1}\calO_{(\fkX_0)_\et}\quad\text{and}\quad
\calO_X^{\ur}\coloneqq\calO_X^{\ur/\fkX_0}\coloneqq w^{-1}\calO_{(\fkX_0)_\et}[p^{-1}].
\]
Set $\calOA_{\mathrm{\inf}}\coloneqq \calO_X^{\ur+}\otimes_{\Z}\bA_{\mathrm{inf}}$.
By extending the scalars, we have an $\calO_X^{\ur+}$-algebra morphism $\theta_X\colon \calOA_{\mathrm{inf}}\rightarrow \widehat{\calO}_X^+$.
Define $\calOA_\cris$ to be the $p$-adic completion of the PD-envelope $\calOA_\cris^0$ of $\calOA_{\mathrm{inf}}$ with respect to the ideal sheaf $\Ker \theta_X$. Note that $\calOA_\cris$ is an $\bA_\cris$-algebra.
Set 
\[
\calOB_\cris^+\coloneqq \calOA_\cris[p^{-1}]\quad\text{and}\quad 
\calOB_\cris\coloneqq \calOB_\cris^+[t^{-1}].
\]
Here the sheaf $\calOB_\cris^+|_{X_{\overline{K}}}[t^{-1}]$ on $X_{\proet/X_{\overline{K}}}$ naturally descends to a sheaf on $X_\proet$ and $\calOB_\cris^+[t^{-1}]$ denotes the corresponding sheaf.
These sheaves are equipped with a decreasing filtration and a connection that satisfy the Griffiths transversality, which we omit to explain. See the remark below.
 \end{defn}

\begin{rem}
The definitions of our sheaves $\calOA_{\mathrm{inf}}$ and $\calOA_\cris^0$ are slightly different from the ones given by Tan and Tong: we use $\otimes_\Z$ instead of $\otimes_W$ to define $\calOA_{\mathrm{inf}}$. However, our $\calOA_{\cris}$ still coincides with theirs in the unramified case $\calO_K=W$ since $k$ is perfect. Hence it follows that $\calOA_{\cris}\cong \calOA_{\cris,X_0}|_{X_\proet}$ and $\calOB_{\cris}\cong \calOB_{\cris,X_0}|_{X_\proet}$. In particular, one can define the additional structures on $\calOA_\cris$ and $\calOB_\cris$ directly from \cite{Tan-Tong}.
\end{rem}

\begin{prop}[(cf.~{\cite[Cor.~2.19]{Tan-Tong}})]\label{prop:Tan-Tong Cor.2.19}
Let $\fkU_0=\Spf R_0 \in (\fkX_0)_\et$ be affine such that $R_0$ is connected and small over $W$.
With the notation as in \S~\ref{sec-basering}, set $R=R_0\otimes_W\calO_K$ and $U=\Spa(R[p^{-1}],R)$, and let $\overline{U}\in X_\proet$ denote the affinoid perfectoid corresponding to the pro-\'etale cover $(\overline{R}[p^{-1}],\overline{R})$ of $(R[p^{-1}],R)$. Then there is a natural isomorphism of $R_0\otimes_W \B_{\cris}(\overline{R})$-modules 
\[
\OB_\cris(\overline{R})\xrightarrow{\cong}\calOB_\cris(\overline{U})
\]
that is strictly compatible with filtrations.
Moreover, for every $i>0$ and $j\in\Z$, we have
\[
H^i(\overline{U},\calOB_\cris)=H^i(\overline{U},\Fil^j\calOB_\cris)=0.
\]
\end{prop}

\begin{proof}
Note that we have natural identifications $\overline{R}=\overline{R_0}$ and $\OB_\cris(\overline{R})=\OB_\cris(\overline{R_0})$.
Now the proposition is nothing but \cite[Cor.~2.19]{Tan-Tong} for $\fkU_0$ and $\overline{U}\in (X_0)_\proet$. Note that \emph{loc. cit.} only claims that the map $\OB_\cris(\overline{R})\rightarrow\calOB_\cris(\overline{U})$ is an isomorphism of $R_0\otimes_W\B_{\cris}(\overline{\calO_K})$-modules, but its proof together with \cite[Cor.~2.8]{Tan-Tong} shows that the map is indeed an isomorphism of $R_0\otimes_W \B_{\cris}(\overline{R})$-modules.
\end{proof}

\begin{rem}
 The modules $\OB_\cris(\overline{R})$ and $\calOB_\cris(\overline{U})$ admit an action of $\calG_R$ (or even $\calG_{R_0}$). The above isomorphism is compatible with the Galois actions. It is also compatible with the restriction along any \'etale morphism $\Spf R_0'\rightarrow \Spf R_0$.
\end{rem}

\begin{rem}\label{rem:Tan-Tong Lem.2.18}
 By the same argument, we also have a description of 
$\calOA_\cris(\overline{U})$ similar to \cite[Lem.~2.18]{Tan-Tong}. 
\end{rem}

We now explain the crystalline formalism.
Let $\mathrm{Loc}_{\Z_p}(X)$ (resp.~$\mathrm{ILoc}_{\Z_p}(X)$) denote the category of \'etale $\Z_p$-local systems (resp.~ \'etale isogeny $\Z_p$-local systems) on $X$. See \cite[\S~1.4, 8.4]{kedlaya-liu-relative-padichodge} for the precise formulation. By \cite[Prop.~8.2]{scholze-p-adic-hodge}, $\mathrm{Loc}_{\Z_p}(X)$ is equivalent to the category of $\widehat{\Z}_p$-local systems on $X_\proet$. We also note that if $\fkX=\Spf R$ is connected and affine with $X=\Spa(R[p^{-1}],R)$, then there are equivalences of categories 
\[
\mathrm{Loc}_{\Z_p}(X)\cong \mathrm{Rep}_{\mathbf{Z}_p}^{\mathrm{pr}}(\calG_R)\quad\text{and}\quad
\mathrm{ILoc}_{\Z_p}(X)\cong \mathrm{Rep}_{\mathbf{Q}_p}(\calG_R).
\]

\begin{defn}[(cf.~{\cite[Def.~3.12]{Tan-Tong}})]\label{defn:crystalline local systems when W-model exists}
Keep the assumption on the existence of $\fkX_0$.
 For an \'etale isogeny $\Z_p$-local system $\bL$ on $X$ with corresponding $\widehat{\Q}_p$-local system $\widehat{\bL}$ on $X_{\proet}$, we set
\[
 D_{\cris}(\bL)\coloneqq w_\ast(\calOB_{\cris}\otimes_{\widehat{\Q}_p}\widehat{\bL})
\quad\text{and}\quad
\Fil^i D_{\cris}(\bL)\coloneqq w_\ast(\Fil^i\calOB_{\cris}\otimes_{\widehat{\Q}_p}\widehat{\bL}).
\]
Note that these are sheaves of $\calO_{(\fkX_0)_\et}[p^{-1}]$-modules.

We say that $\bL$ is \emph{crystalline} (with respect to $\fkX_0$) if 
\begin{enumerate}
 \item the $\calO_{(\fkX_0)_\et}[p^{-1}]$-modules $D_{\cris}(\bL)$ and $\Fil^i D_{\cris}(\bL)$ ($i\in\Z$) are all coherent, and
 \item the adjunction morphism
\begin{equation}\label{eq:comparison map on proetale site}
\calOB_{\cris} \otimes_{\calO_X^{\ur}[p^{-1}]}w^{-1}D_{\cris}(\bL)
\rightarrow \calOB_{\cris}\otimes_{\widehat{\Q}_p}\widehat{\bL}    
\end{equation}
is an isomorphism of $\calOB_{\cris}$-modules.
\end{enumerate}
\end{defn}

\begin{rem}
In the unramified case, Tan and Tong \cite[Def.~3.10]{Tan-Tong} define crystalline local systems using the notion of association with a convergent filtered $F$-isocrystal, and they prove that their definition is equivalent to conditions (i) and (ii) above in \cite[Prop.~3.13]{Tan-Tong}.
\end{rem}

\begin{lem}[(cf.~{\cite[Lem.~3.14]{Tan-Tong}})]\label{lem:Tan-Tong Lem.3.14}
Assume that $\fkX$ admits a $W$-model $\fkX_0$ and let $\bL\in \mathrm{ILoc}_{\Z_p}(X)$.
For each small and connected affine formal scheme $\fkU_0=\Spf R_0$ that is \'etale over $\fkX_0$, set $R\coloneqq R_0\otimes_W\calO_K$ and $U\coloneqq \Spa(R[p^{-1}],R)$, and let $V_U$ denote the $\Q_p$-representation of $\calG_R$ corresponding to $\bL|_U$.
Then there exist natural isomorphisms of $R_0[p^{-1}]$-modules
\[
D_\cris(\bL)(\fkU_0)\xrightarrow{\cong}D_\cris(V_U)\quad\text{and}\quad 
(\Fil^i D_\cris(\bL))(\fkU_0)\xrightarrow{\cong}\Fil^i D_\cris(V_U)\quad (i\in \Z).
\]
Moreover, if we write $\overline{U}\in X_\proet$ for the affinoid perfectoid attached to $(\overline{R}[p^{-1}],\overline{R})$, then the evaluation of the adjunction morphism \eqref{eq:comparison map on proetale site} at $\overline{U}$ coincides with 
\[
\alpha_{\cris}(V_U)\colon \OB_{\cris}(\overline{R})\otimes_{R_0[p^{-1}]}D_{\cris}(V_U) \rightarrow \OB_{\cris}(\overline{R})\otimes_{\Q_p}V_U
\]
under the identification $D_\cris(\bL)(\fkU_0)\cong D_\cris(V_U)$.
\end{lem}

\begin{proof}
The proof in \cite[Lem.~3.14]{Tan-Tong} also works in the current setting if one uses Proposition~\ref{prop:Tan-Tong Cor.2.19} in place of \cite[Cor.~2.19]{Tan-Tong}. The second assertion follows from the construction.
\end{proof}

\begin{prop}\label{prop:def of crystalline local systems}
Assume that $\fkX$ admits a $W$-model $\fkX_0$ and let $\bL\in \mathrm{ILoc}_{\Z_p}(X)$. Then $\bL$ is crystalline with respect to $\fkX_0$ in the sense of Definition~\ref{defn:crystalline local systems when W-model exists} if and only if there exists an \'etale covering $\{\fkU_{\lambda,0}\rightarrow \fkX_0\}$ of small and connected affine $\fkU_{\lambda,0}=\Spf R_{\lambda,0}$ such that the $\Q_p$-representation $V_\lambda$ of $\calG_{R_\lambda}$ corresponding to $\bL|_{\Spa(R_\lambda[p^{-1}],R_\lambda)}$ is crystalline in the sense of Definition~\ref{defn:crystalline representations}, where $R_\lambda\coloneqq R_{\lambda,0}\otimes_W\calO_K$.
In particular, the notion of crystalline local systems on $X$ does not depend on the choice of a $W$-model of $\fkX$. 
\end{prop}

\begin{proof}
The necessity follows from Lemma~\ref{lem:Tan-Tong Lem.3.14}. For the sufficiency, observe that both of conditions (i) and (ii) in Definition~\ref{defn:crystalline local systems when W-model exists} can be verified locally on $(\fkX_0)_\et$. So we may assume $\fkX_0=\Spf R_{\lambda,0}$ for some $\lambda$. To simplify the notation, write $R_0$ for $R_{\lambda,0}$ and $V$ for $V_\lambda$.

First we verify condition (i). Since the proof is similar, we only show that $D_{\cris}(\bL)$ is a coherent $\calO_{(\fkX_0)_\et}[p^{-1}]$-module. Take any connected and affine $\fkU_0=\Spf R_0'\in (\fkX_0)_\et$. We need to show that the natural morphism 
\begin{equation}\label{eq:in the proof of prop:def of crystalline local systems}
R_0'[p^{-1}]\otimes_{R_0[p^{-1}]}D_\cris(\bL)(\fkX_0)\rightarrow D_{\cris}(\bL)(\fkU_0)
\end{equation}
is an isomorphism.
Set $R'\coloneqq R_0'\otimes_W\calO_K$. Then $R'$ is connected and small over $\calO_K$. By Lemma~\ref{lem:Tan-Tong Lem.3.14}, we have identifications $D_\cris(\bL)(\fkX_0)\xrightarrow{\cong}D_\cris(V)$ and $D_\cris(\bL)(\fkU_0)\xrightarrow{\cong}D_\cris(V|_{\calG_{R'}})$. Since $V$ is crystalline, the map \eqref{eq:in the proof of prop:def of crystalline local systems} is an isomorphism by Lemma~\ref{lem:base change map for Dcris}.
Now that we have verified condition (i), condition (ii) follows from the proof of \cite[Cor.~3.15 and 3.16]{Tan-Tong} with Remark~\ref{rem:Tan-Tong Lem.2.18}, Proposition~\ref{prop:Tan-Tong Cor.2.19}, and Lemma~\ref{lem:Tan-Tong Lem.3.14} in place of Lemma~2.18, Corollary~2.19, and Lemma~3.14 of \cite{Tan-Tong}.
This completes the proof of the sufficiency. The last assertion follows from \cite[Prop.~8.3.5]{brinon-relative}.
\end{proof}

With these preparations, we define the notion of crystalline local systems via gluing.

\begin{defn}
Let $\fkX$ be a smooth $p$-adic formal scheme over $\calO_K$ and let $X$ denote its adic generic fiber.
An \'etale isogeny $\Z_p$-local system $\bL$ on $X$ is said to be \emph{crystalline} if there exists an open covering $\fkX=\bigcup_\lambda \fkU_{\lambda}$ such that each $\fkU_\lambda$ admits a $W$-model and such that for each $\lambda$, $\bL|_{U_\lambda}$ is crystalline in the sense of Definition~\ref{defn:crystalline local systems when W-model exists} where $U_\lambda$ denotes the adic generic fiber of $\fkU_\lambda$.
By Proposition~\ref{prop:def of crystalline local systems}, this definition coincides with Definition~\ref{defn:crystalline local systems when W-model exists} when $\fkX$ itself admits a $W$-model.

An \'etale $\Z_p$-local system $\bL$ on $X$ is said to be \emph{crystalline} if the associated isogeny $\Z_p$-local system is crystalline.
\end{defn}

\begin{rem}
One could define crystalline local systems by introducing a period sheaf $\calOB_{\mathrm{max},K}$ on $X_{\proet}$ that generalizes the period ring $\mathrm{A}_{\max}(R)[p^{-1},t^{-1}]$ appearing in the proof of \cite[Prop.~8.3.5]{brinon-relative}. This period sheaf is defined without fixing a $W$-model of $\fkX$ and thus one could bypass the gluing approach.
\end{rem}

\bibliographystyle{amsplain}
\bibliography{library}
	
\end{document}